%% file: Dgeo.tex

\input amssym.def
\input amssym.tex

\def\item#1{\vskip1.3pt\hang\textindent {\rm #1}}

\def\itemitem#1{\vskip1.3pt\indent\hangindent2\parindent\textindent {\rm #1}}

\newskip\litemindent
\litemindent=0.7cm  
\def\Litem#1#2{\par\noindent\hangindent#1\litemindent
\hbox to #1\litemindent{\hfill\hbox to \litemindent
{\ninerm #2 \hfill}}\ignorespaces}

\tolerance=300
\pretolerance=200
\hfuzz=1pt
\vfuzz=1pt

\magnification=1000         
\hoffset=0in
\voffset=0.5in

\hsize=5.0 true in    
\vsize=8.1 true in

\parindent=25pt
\mathsurround=1pt
\parskip=1pt plus .25pt minus .25pt
\normallineskiplimit=.99pt

\countdef\revised=100
\mathchardef\emptyset="001F 
\chardef\ss="19
\def\3{\ss}
\def\anf{$\lower1.2ex\hbox{"}$}
\def\frac#1#2{{#1 \over #2}}
\def\>{>\!\!>}
\def\<{<\!\!<}

\def\ssssarr{\hbox to 15pt{\rightarrowfill}}
\def\sssarr{\hbox to 20pt{\rightarrowfill}}
\def\ssarr{\hbox to 30pt{\rightarrowfill}}
\def\sarr{\hbox to 40pt{\rightarrowfill}}
\def\arr{\hbox to 60pt{\rightarrowfill}}
\def\larr{\hbox to 60pt{\leftarrowfill}}
\def\Arr{\hbox to 80pt{\rightarrowfill}}

\def\ad{\mathop{\rm ad}\nolimits}

\def\Ad{\mathop{\rm Ad}\nolimits}

\def\Aut{\mathop{\rm Aut}\nolimits}

\def\Bil{\mathop{\rm Bil}\nolimits}

\def\Der{\mathop{\rm Der}\nolimits}
\def\deg{\mathop{\rm deg}\nolimits}

\def\Diff{\mathop{\rm Diff}\nolimits}

\def\Exp{\mathop{\rm Exp}\nolimits}
\def\ev{\mathop{\rm ev}\nolimits}

\def\End{\mathop{\rm End}\nolimits}

\def\Hom{\mathop{\rm Hom}\nolimits}%
\def\id{\mathop{\rm id}\nolimits} 
\def\im{\mathop{\rm im}\nolimits}


\def\OO{\mathop{\rm O{}}\nolimits}


\def\sgn{\mathop{\rm sgn}\nolimits}
\def\SL{\mathop{\rm SL}\nolimits}

\def\Sym{\mathop{\rm Sym}\nolimits}

\def\supp{\mathop{\rm supp}\nolimits}



\def\0{{\bf 0}}
\def\1{{\bf 1}}

\def\g{{\frak g}}
\def\gl{{\frak {gl}}}
\def\h{{\frak h}}

\def\m{{\frak m}}

\def\q{{\frak q}}

\def\C{{{\Bbb C}{\mskip+2mu}}} 
\def\K{{{\Bbb K}{\mskip+2mu}}} 

\def\R{{\Bbb R}} 
\def\Z{{\Bbb Z}} 
\def\N{{\Bbb N}}

\def\K{{\Bbb K}}

\def\:{\colon}  
\def\.{{\cdot}}
\def\|{\Vert}
\def\bsk{\bigskip}

\def\giantskip{\vskip2\bigskipamount}
\def\gsk{\giantskip}

\def\msk{\medskip}

\def\ssk{\smallskip}

\def\bbr{\bigbreak}
\def\giantbreak{\par \ifdim\lastskip<2\bigskipamount \removelastskip
         \penalty-400 \giantskip\fi}

\def\nin{\noindent}
\def\cen{\centerline}
\def\pagebreak{\vskip 0pt plus 0.0001fil\break}
\def\linebreak{\break}

\def\hat{\widehat}

\def\eps{\varepsilon}
\def\epsilon{\varepsilon}

\def\nin{\noindent}
\def\oline{\overline}
\def\uline{\underline}
\def\pder#1,#2,#3 { {\partial #1 \over \partial #2}(#3)}
\def\pde#1,#2 { {\partial #1 \over \partial #2}}
\def\phi{\varphi}


\def\tilde{\widetilde}

\font\ninerm=cmr9
\font\eightrm=cmr8

\font\eightbf=cmbx8


\font\smc=cmcsc10
\font\bfone=cmbx10 scaled\magstep1 
\font\bftwo=cmbx10 scaled\magstep2 

\def\qed{{\unskip\nobreak\hfil\penalty50\hskip .001pt \hbox{}\nobreak\hfil
          \vrule height 1.2ex width 1.1ex depth -.1ex
           \parfillskip=0pt\finalhyphendemerits=0\medbreak}\rm}

\def\qeddis{\eqno{\vrule height 1.2ex width 1.1ex depth -.1ex} $$
                   \medbreak\rm}

\def\Lemma #1. {\bigbreak\vskip-\parskip\noindent{\bf Lemma #1.}\quad\it}

\def\Sublemma #1. {\bigbreak\vskip-\parskip\noindent{\bf Sublemma #1.}\quad\it}

\def\Proposition #1. {\bigbreak\vskip-\parskip\noindent{\bf Proposition #1.}
\quad\it}

\def\Corollary #1. {\bigbreak\vskip-\parskip\nin{\bf Corollary #1.}
\quad\it}

\def\Theorem #1. {\bigbreak\vskip-\parskip\noindent{\bf Theorem #1.}
\quad\it}

\def\Definition #1. {\rm\bigbreak\vskip-\parskip\noindent{\bf Definition #1.}
\quad}

\def\Remark #1. {\rm\bigbreak\vskip-\parskip\noindent{\bf Remark #1.}\quad}

\def\Example #1. {\rm\bigbreak\vskip-\parskip\noindent{\bf Example #1.}\quad}
\def\Examples #1. {\rm\bigbreak\vskip-\parskip\noindent{\bf Examples #1.}\quad}

\def\Problems #1. {\bigbreak\vskip-\parskip\noindent{\bf Problems #1.}\quad}
\def\Problem #1. {\bigbreak\vskip-\parskip\noindent{\bf Problem #1.}\quad}
\def\Exercise #1. {\bigbreak\vskip-\parskip\noindent{\bf Exercise #1.}\quad}

\def\Conjecture #1. {\bigbreak\vskip-\parskip\noindent{\bf Conjecture #1.}\quad}

\def\Proof#1.{\rm\par\ifdim\lastskip<\bigskipamount\removelastskip\fi\smallskip
            \noindent {\bf Proof.}\quad}

\def\Axiom #1. {\bigbreak\vskip-\parskip\noindent{\bf Axiom #1.}\quad\it}

\def\Satz #1. {\bigbreak\vskip-\parskip\noindent{\bf Satz #1.}\quad\it}

\def\Korollar #1. {\bbr\vskip-\parskip\nin{\bf Korollar #1.} \quad\it}

\def\Folgerung #1. {\bbr\vskip-\parskip\nin{\bf Folgerung #1.} \quad\it}

\def\Folgerungen #1. {\bbr\vskip-\parskip\nin{\bf Folgerungen #1.} \quad\it}

\def\Bemerkung #1. {\rm\bigbreak\vskip-\parskip\noindent{\bf Bemerkung #1.}
\quad}

\def\Beispiel #1. {\rm\bigbreak\vskip-\parskip\noindent{\bf Beispiel #1.}\quad}
\def\Beispiele #1. {\rm\bigbreak\vskip-\parskip\noindent{\bf Beispiele #1.}\quad}
\def\Aufgabe #1. {\rm\bigbreak\vskip-\parskip\noindent{\bf Aufgabe #1.}\quad}
\def\Aufgaben #1. {\rm\bigbreak\vskip-\parskip\noindent{\bf Aufgabe #1.}\quad}

\def\Beweis#1. {\rm\par\ifdim\lastskip<\bigskipamount\removelastskip\fi
           \smallskip\noindent {\bf Beweis.}\quad}

\def\date{\ifcase\month\or January\or February \or March\or April\or May
\or June\or July\or August\or September\or October\or November
\or December\fi\space\number\day, \number\year}

\input Dgeo0000.tex 

\input Dgeo000.tex  

\nopagenumbers

\def\title{Title ??}
\def\author{Author ??}

\def\thanks#1{\footnote*{\eightrm#1}}

\def\rightheadline{\hfil{\eightrm\title}\hfil\tenbf\folio}
\def\leftheadline{\tenbf\folio\hfil{\eightrm\author}\hfil}
\headline={\vbox{\line{\ifodd\pageno\rightheadline\else\leftheadline\fi}}}

\def\firstheadline{}
\def\firstfootline{\cen{\rm\folio}}

\def\seite #1 {\pageno #1
               \headline={\ifnum\pageno=#1 \firstheadline
               \else\ifodd\pageno\rightheadline\else\leftheadline\fi\fi}
               \footline={\ifnum\pageno=#1 \firstfootline\else{}\fi}}

\newdimen\dimenone
 \def\checkleftspace#1#2#3#4{
 \dimenone=\pagetotal
 \advance\dimenone by -\pageshrink   
 \ifdim\dimenone>\pagegoal          
   \else\dimenone=\pagetotal
        \advance\dimenone by \pagestretch
        \ifdim\dimenone<\pagegoal
          \dimenone=\pagetotal
          \advance\dimenone by#1         
          \setbox0=\vbox{#2\parskip=0pt                
                     \hyphenpenalty=10000
                     \rightskip=0pt plus 5em
                     \noindent#3 \vskip#4}    
        \advance\dimenone by\ht0
        \advance\dimenone by 3\baselineskip   
        \ifdim\dimenone>\pagegoal\vfill\eject\fi
          \else\eject\fi\fi}


\def\subheadline #1{\nin\bigbreak\vskip-\lastskip
      \checkleftspace{0.9cm}{\bf}{#1}{\medskipamount}
          \indent\vskip0.7cm\centerline{\bf #1}\medskip}

\def\lsubheadline #1 #2{\bigbreak\vskip-\lastskip
      \checkleftspace{0.9cm}{\bf}{#1}{\bigskipamount}
         \vbox{\vskip0.7cm}\cen{\bf #1}\msk \cen{\bf #2}\bsk}

\def\sectionheadline #1{\bigbreak\vskip-\lastskip
      \checkleftspace{1.1cm}{\bf}{#1}{\bigskipamount}
         \vbox{\vskip1.1cm}\cen{\bfone #1}\bsk}

\def\lsectionheadline #1 #2{\bigbreak\vskip-\lastskip
      \checkleftspace{1.1cm}{\bf}{#1}{\bigskipamount}
         \vbox{\vskip1.1cm}\cen{\bfone #1}\msk \cen{\bfone #2}\bsk}

\def\lchapterheadline #1 #2{\bigbreak\vskip-\lastskip\indent\vskip3cm
                       \cen{\bftwo #1} \msk \cen{\bftwo #2} \gsk}
\def\llsectionheadline #1 #2 #3{\bigbreak\vskip-\lastskip\indent\vskip1.8cm
\cen{\bfone #1} \msk \cen{\bfone #2} \msk \cen{\bfone #3} \nobreak\bsk\nobreak}


\newtoks\literat
\def\[#1 #2\par{\literat={#2\unskip.}%
\hbox{\vtop{\hsize=.15\hsize\nin [#1]\hfill}
\vtop{\hsize=.82\hsize\nin\the\literat}}\par
\vskip.3\baselineskip}

\def\references{
\sectionheadline{\bf References}
\frenchspacing

\entries\par}

\mathchardef\emptyset="001F

\def\abstract #1{{\narrower\baselineskip=10pt{\noindent
\eightbf Abstract.\quad \eightrm #1 }
\bigskip}}

\def\Box #1 { \msk\par\nin 
\centerline{
\vbox{\offinterlineskip
\hrule
\hbox{\vrule\strut\hskip1ex\hfil{\smc#1}\hfill\hskip1ex}
\hrule}\vrule}\msk }

\def\adots{\mathinner{\mkern1mu\raise1pt\vbox{\kern7pt\hbox{.}}
                        \mkern2mu\raise4pt\hbox{.}
                        \mkern2mu\raise7pt\hbox{.}\mkern1mu}}


\input xy
\xyoption{all}

\def\Gl{\mathop{\rm Gl}\nolimits}
\def\dd{{\rm d}}

\def\X{{\frak X}}
\def\pr{\mathop{\rm pr}\nolimits}
\def\Diff{\mathop{\rm Diff}\nolimits}

\def\Gm{\mathop{\rm Gm}\nolimits}
\def\Fl{\mathop{\rm Fl}\nolimits}
\def\brs{\mathop{\rm brs}\nolimits}

\def\alt{\mathop{\rm alt}\nolimits}

\def\1{{\bf 1}}

\def\Gl{\mathop{\rm Gl}\nolimits}
\def\dd{{\rm d}}
\def\mod{{\rm mod}}

\def\X{{\frak X}}
\def\pr{\mathop{\rm pr}\nolimits}
\def\Diff{\mathop{\rm Diff}\nolimits}

\def\brs{\mathop{\rm brs}\nolimits}

\def\alt{\mathop{\rm alt}\nolimits}

\def\Gm{\mathop{\rm Gm}\nolimits}
\def\Sm{\mathop{\rm Sm}\nolimits}
\def\Vsm{\mathop{\rm Vsm}\nolimits}
\def\Shi{\mathop{\rm Shi}\nolimits}

\def\gm{{\frak {gm}}}

\def\Mult{\mathop{\rm Mult}\nolimits}
\def\Part{\mathop{\rm Part}\nolimits}

\def\eps{{\epsilon}}


\input Dgeo00.tex

\input Dgeo1.tex

\input Dgeo2.tex


\input Dgeo3.tex


\input Dgeo4.tex

\input Dgeo5.tex


\bye

%% file: Dgeo0000.tex

\pageno=-2006

$$
\phantom{xx}
$$

\vskip 50mm

\centerline
{\bf Wolfgang Bertram}

\medskip
\centerline{\bf  (Institut Elie Cartan, Nancy)}

\vskip 15mm

\centerline{\bfone 
Differential Geometry over General Base Fields and Rings} 

\vskip 105mm

\centerline{Version to appear in the Memoirs of the AMS.}

\bigskip

\centerline{Last updated: november}

\vfill \eject

$$
\phantom{xx}
$$

\vfill \eject

%% file: Dgeo000.tex

\pageno=-5
\def\title{Differential Geometry over General Base Fields and Rings} 
\def\author{Wolfgang Bertram}
\def\date{19.10.2006}

\def\leftheadline{\tenbf\folio\hfil{\eightrm\title}\hfil}

\nin {\bfone Contents}

\vskip 50mm
\nin Introduction
\dotfill{1}

\bigskip

\nin I. Basic notions

\msk
1. Differential calculus
\dotfill{14}

\ssk
2. Manifolds
\dotfill{20}

\ssk
3. Tangent bundle and general fiber bundles \dotfill{22}

\ssk
4. The Lie bracket of vector fields \dotfill{25}

\ssk
5. Lie groups and symmetric spaces: basic facts \dotfill{30}

\bigskip

\nin II. Interpretation of tangent objects via scalar extensions

\msk
6. Scalar extensions. I: Tangent functor and dual numbers \dotfill{36}
\ssk

7. Scalar extensions. II: Higher order tangent functors \dotfill{42}

\ssk
8. Scalar extensions. III: Jet functor and truncated polynomial rings
\dotfill{50}

\bigskip

\nin III. Second  order differential geometry

\msk
9. The structure of the tangent bundle of a vector bundle
\dotfill{57}

\ssk
10. Linear connections. I: Linear structures on bilinear bundles
\dotfill{61}

\ssk
11. Linear connections. II: Sprays
\dotfill{68}

\ssk
12. Linear connections. III: Covariant derivative 
\dotfill{71}

\ssk
13. Natural operations. I: Exterior derivative of a one-form
\dotfill{73}

\ssk
14. Natural operations. II: The Lie bracket revisited
\dotfill{75}

\bigskip

\nin IV. Third and higher order differential geometry

\msk
15. The structure of $T^k F$:  Multilinear bundles
\dotfill{79}

\ssk
16. The structure of $T^k F$: Multilinear connections
\dotfill{83}

\ssk
17. Construction of multilinear connections
\dotfill{87}

\ssk
18. Curvature \dotfill{91}

\ssk
19. Linear structures on jet bundles \dotfill{95}

\bigskip

\centerline
{\eightrm TABLE OF CONTENTS (continued)}

\bigskip
20. Shifts and symmetrization \dotfill{98}

\ssk
21. Remarks on differential operators and symbols \dotfill{102}

\ssk
22. The exterior derivative \dotfill{106} 

\bigskip


\nin V. Lie Theory 

\msk
23. The three canonical connections of a Lie group \dotfill{110}

\ssk
24. The structure of higher order tangent groups  \dotfill{116}

\ssk
25. Exponential map and Campbell-Hausdorff formula  \dotfill{124}

\ssk
26. The canonical connection of a symmetric space \dotfill{128}

\ssk
27. The higher order tangent structure of symmetric spaces \dotfill{134}

\bigskip
 \nin VI. Groups of diffeomorphisms and the exponential jet

\msk
28. Group structure on the space of sections of $T^k M$ 
\dotfill{139}

\ssk
29. The exponential jet for vector fields \dotfill{144}

\ssk
30. The exponential jet of a symmetric space \dotfill{148}

\ssk
31. Remarks on the exponential jet of a general connection \dotfill{151}

\ssk
32. From germs to jets and from jets to germs \dotfill{153}

\bigskip
\nin Appendix L. Limitations \dotfill{156}

\bigskip
\nin Appendix G. Generalizations \dotfill{159}

\bigskip
\nin
Appendix: Multilinear Geometry

\msk
BA. Bilinear algebra \dotfill{161}

\ssk
MA. Multilinear algebra \dotfill{168}

\ssk
SA. Symmetric and shift invariant multilinear algebra \dotfill{182}

\ssk
PG. Polynomial groups \dotfill{192}

\bigskip
\nin References \dotfill{199}

\vfill \eject

\nin {\bfone Preface}

\vskip 10mm

\rightline{\it Die Regeln des Endlichen gelten im Unendlichen weiter.}

\msk

\rightline{(G.\ W.\ Leibniz, 1702)}

\vskip 10mm
Classical, real finite-dimensional differential geometry and  
Lie theory can be generalized in several directions. For instance,
it is interesting to look for an {\it infinite dimensional} theory, 
or for a theory which works over other base fields, such as the $p$-adic
numbers or other topological fields. 
The present work aims at both directions at the same time.
Even more radically, we will develop the theory over certain base
{\it rings} so that, in general, there is no dimension at all.
Thus, instead of opposing ``finite" and ``infinite-dimensional"
or ``real" and ``$p$-adic" theories, we will simply speak of 
``general" Lie theory or differential geometry. 
Certainly, this is a daring venture, and therefore I should first
explain the origins and the aims of the project and say what can be
done and what can {\it not} be done. 

\ssk Preliminary versions of
the present text have  appeared as a series of five preprints
during the years 2003 -- 2005
at the Institut Elie Cartan, Nancy (still available on the web-server
{\tt http://www.iecn.u-nancy.fr/Preprint/}). The title of this work
is, of course, a reference to Helgason's by now classical monograph
 [Hel78] from which
I first learned the theory of differential geometry, Lie groups and
 symmetric spaces. Later I  became interested in the  more 
algebraic aspects of the theory,
related to the interplay between Lie- and Jordan-theory (see [Be00]).
However, rather soon I began to regret that differential geometric 
methods seemed to be limited to the case of the base field of real
or complex numbers and to the case of ``well-behaved'' topological
vector spaces as model spaces, such as finite-dimensional or
 Banach spaces (cf.\ the remark in [Be00, Ch.XIII] on
 ``other base fields''). Rather surprisingly (at least for me),
these limitations were entirely overcome  by  the
joint paper [BGN04] with Helge Gl\"ockner
and Karl-Hermann Neeb where differential calculus, manifolds and Lie groups
have been introduced in a very general setting, including the case of manifolds
modelled on
arbitrary topological vector spaces over any non-discrete topological
field, and even over topological modules over
certain topological base rings.
In the joint paper [BeNe05] with Karl-Hermann Neeb, a good deal of the
results from [Be00] could be generalized to this very general framework,
leading, for example, to a rich supply of infinite-dimensional symmetric
spaces.

\ssk
As to differential geometry on infinite dimensional manifolds, I used to have
the impression that its special flavor is due to its, sometimes rather
sophisticated, functional analytic methods. 
On the other hand, it seemed obvious that the ``purely differential''
aspects of differential geometry are algebraic in nature and thus should
be understandable without all the functional analytic flesh around it.
In other words, one should be able to formalize the fundamental 
differentiation process such that its general algebraic structure becomes 
visible. This is indeed possible, and one possibility is opened by the
general differential calculus mentioned above: for me, the most striking
result of this calculus is the one presented in Chapter 6 of this work
(Theorem 6.2), saying that the usual ``tangent functor'' $T$ of differential
geometry can, in this context, be interpreted as a functor of scalar extension
from the base field or ring $\K$ to the {\it tangent ring} $T\K = \K \oplus
\eps \K$ with relation $\eps^2 = 0$, sometimes also called the {\it dual 
numbers over $\K$}. Even though the result itself may not seem surprising
(for instance, in algebraic geometry it is currently used), it seems that
the proof in our context is the most natural one that one may imagine.
In some sense, all the rest
 of the text can be seen as a try to exploit the
consequences of this theorem. 

\ssk
I would like to add a remark on terminology.
Classical ``calculus" has two faces, namely {\it differential calculus}
and {\it integral calculus}. In spite of its name, 
classical ``differential geometry" also has these two faces.
However, in order to avoid misunderstandings, I propose  to
clearly distinguish between ``differential" and ``integral" geometry,
the former being the topic of the present work, the latter being 
restricted to special contexts where various kinds of integration
theories can be implemented (cf.\ Appendix L). In a nutshell, differential
geometry in this sense is the theory of $k$-th order
Taylor expansions, for any $k \in \N$,
of objects such as manifolds, vector bundles or Lie groups, without
requiring convergence of the Taylor series. In the French literature,
this approach
is named by the somewhat better adapted term
 ``d\'eveloppement limit\'e" (cf., e.g., [Ca71]),
which could be translated by ``limited expansion". 
Certainly, to a mathematician used to differential geometry in real life,
such a rigorous distinction between ``differential'' and 
``integral'' geometry may seem pedantic and purely formal -- probably, 
most readers would prefer to call  ``formal differential geometry" what
remains when integration is forbidden.  
However, the only ``formal" construction that we use is the
one of the tangent bundle, and hence 
the choice of terminology depends on wether
one considers the tangent bundle as a genuine object of differential 
geometry or not.
Only in the very last chapter (Chapter 32) we pass the frontier separating
true manifolds from formal objects by taking the projective limit
$k \to \infty$ of our $k$-th order Taylor expansions, thus making the link
with approaches by formal power series.

\msk The present text is completely self-contained, but by no 
means                                      
it is an exhaustive treatise on the topics mentioned in the title.
Rather, I wanted to present a first approximation of a new theory that,
in my opinion, is well-adapted to capture some very general aspects
of differential calculus and differential geometry and thus might be
interesting for mathematicians and physicists working on topics that
are related to these aspects.
				      
\msk
{\sl Acknowledgements.} 
I would not have dared to tackle this work if I had not had teachers
who, at a very early stage, attracted my interest to foundational questions of
calculus. In particular, although I do not follow the lines of thought
of non-standard analysis, I am indebted to Detlef Laugwitz from
whose reflections on the foundations of analysis I learned a lot.
(For the adherents of non-standard analysis, I add the remark that
the field ${ }^*\R$ of non-standard numbers is of course 
admitted as a possible base field of our theory.) 

I would also like to thank
Harald L\"owe for his careful proofreading of
the first part of the text, and all future readers for indicating to me 
errors that, possibly, have remained.

\vskip 10mm

Nancy, february 2005.

\vfill \eject

\nin {\bf Abstract.}
The aim of this work is to lay the foundations of differential
geometry and Lie theory
over the general class of topological base fields and -rings
for which a differential calculus has been developed in [BGN04],
without any restriction on the dimension or on the characteristic. 
 Two basic features distinguish our approach from the
classical real (finite or infinite dimensional) theory, namely
the interpretation of tangent- and jet functors as functors
 of {\it scalar extensions} and the introduction of 
{\it multilinear bundles} and {\it multilinear connections}
 which generalize the concept of vector bundles and linear connections.
\footnote{}{

\nin Received by the editors: February 9, 2005

\ssk
\nin {\bf AMS subject classification} (2000):

\ssk \nin
Primary:
22E65,
53B05,
53C35,
58A05,
58A20,
53B05, 
58A32 

\ssk
\nin Secondary:
22E20, 
14L05,
15A69,
51K10,
58B20

\bigskip
\nin {\bf Key words:}
Lie groups over topological fields,
infinite-dimensional Lie groups,
diffeomorphism groups,
formal group, general and special multilinear group,
symmetric space,
differentiable manifold, tangent functor, 
linear and multilinear bundle,
jet, 
linear connection, multilinear connection,
 curvature, 
exponential map,  exponential jet,
Campbell-Hausdorff formula,
multilinear algebra,
scalar extension,
dual numbers,
synthetic differential geometry,
Weil functor

\bigskip
\nin {\bf Author's address:}

\ssk
Wolfgang Bertram

Universit\'e Henri Poincar\'e (Nancy I)

Institut Elie Cartan 

B.P.\ 239

54506 Vand\oe uvre-l\`es-Nancy Cedex

France

e-mail: {\tt bertram@iecn.u-nancy.fr}

http: {\tt http://www.iecn.u-nancy.fr/ $\tilde{ }$ bertram}

}

\vfill \eject

$$ \phantom{xx} $$

\vfill \eject

%% file: Dgeo00.tex

\def\title{INTRODUCTION}
\def\author{Wolfgang Bertram}
\def\date{08.11.2006}

\def\rightheadline{\hfil{\eightrm \title}\hfil\tenrm\folio}
\def\leftheadline{\tenrm\folio\hfil{\eightrm
DIFFERENTIAL GEMETRY OVER GENERAL BASE FIELDS AND RINGS}\hfil}

\pageno=1

\sectionheadline
{0. Introduction}

The classical setting of differential geometry is the framework of
finite-dimen\-sio\-nal real manifolds. Later,  parts of the theory have
been generalized to various kinds of infinite dimensional manifolds
(cf.\ [La99]),
to analytic manifolds over ultrametric fields (finite-dimensional
or Banach, cf.\ [Bou67] or [Se65]) or to more formal
concepts going beyond the notion of manifolds such as
``synthetic differential geometry'' and its models, the 
so-called ``smooth toposes'' (cf.\ [Ko81] or [MR91]).
The purpose of the present work is to unify and generalize some
aspects of differential
geometry in the context of manifolds over very general base fields
and even -rings, featuring fundamental structures which presumably
will persist even in further generalizations beyond the manifold
context. We have divided the text into seven main parts of which we
now give a more detailed description.

\subheadline{I. Basic Notions}

\msk
\nin {\bf 0.1.} {\sl Differential calculus.}
Differential geometry may be seen as the ``invariant theory
of differential calculus",
and {\it differential calculus} deals, in one way or another,
with the ``infinitesimally small", either in the language
of limits or in the language of infinitesimals  -- in the words of G.W.
Leibniz: ``Die Regeln des Endlichen gelten im Unendlichen 
weiter"\footnote{$^1$}{``The rules of the finite continue to
hold in the infinite.''} (quoted from [Lau86, p.\ 88]). 
In order to formalize this idea, take a topological commutative unital ring
$\K$ (such as $\K=\R$ or any other model of the continuum that
you prefer) and let  ``das Endliche" (``the finite'') be represented
by the invertible scalars $t \in \K^\times$ and ``das unendlich (Kleine)''
(``the infinitely (small)'') by the non-invertible scalars. 
Then ``Leibniz principle"   may be interpreted by the assumption that
$\K^\times$ is dense in $\K$: the behaviour of a continuous mapping
at the infinitely small is determined by its behavior on the ``finite
universe''.
For instance, if $f: V \supset U \to W$ is a map defined on an open
subset of a topological $\K$-module, then for all invertible scalars
 $t$ in some neighborhood of $0$, the ``slope"
(at $x \in U$ in direction $v \in V$ with parameter $t$)
$$
f^{[1]}(x,v,t) := { f(x+tv) - f(x) \over t }
\eqno (0.1)
$$
 is defined.
We say that $f$ is {\it of class $C^1$ (on $U$)} if the slope
extends to a {\it continuous} map $f^{[1]}:U^{[1]} \to W$,
where $U^{[1]} = \{ (x,v,t) \vert \, x+tv \in U \}$.
By Leibniz' principle, this extension is unique if it exists.
Put another way, the fundamental relation
$$
f(x+tv)-f(x) = t \cdot   f^{[1]} (x,v,t)
\eqno (0.2)
$$
continues to hold for all $(x,v,t) \in U^{[1]}$.
The {\it derivative of $f$ at $x$ (in direction $v$)}
 is the continuation of the
slope to the most singular of all infinitesimal elements, namely
to $t=0$:
$$
df(x)v:=f^{[1]}(x,v,0).
\eqno (0.3)
$$
Now all the basic rules of differential calculus are
easily  proved,
from the Chain Rule via Schwarz' Lemma up to Taylor's formula
(see Chapter 1 and [BGN04]) -- the reader who compares the
proofs with the ``classical'' ones will find that our
approach at the same time simplifies and generalizes this
part of calculus.

\msk
\nin {\bf 0.2.} {\sl Differential geometry versus integral
geometry.}
The other half of Newton's and Leibniz' ``calculus" is 
{\it integral} calculus. For the time being, the simple
equation (0.2) has no integral analog -- in other words,
we cannot reverse the process of differentiation over
general topological fields or rings (not even over the real numbers
if we go too far beyond the Banach-space set-up). Thus
in general  we will
not even be able to determine the space of all solutions of
 a ``trivial" differential equation
such as $df=0$, and therefore no integration will appear
in our approach to differential geometry:
there will be no existence theorems for flows, for geodesics,
or for parallel transports, no exponential map, no Poincar\'e lemma, no
Frobenius theorem... 
what remains is really {\it differential} geometry, in
contrast to the mixture of integration and differentiation
that we are used to from the finite-dimensional real theory.
At a first glance, the remaining structure may look very poor;
but we will find that it is in fact still tremendously complicated
and rich. In fact, we will see that many classical notions and
results that usually are formulated on the global or  local level
or at least on the level of germs, continue to make sense on the level
of {\it jets (of arbitrary order)}, and although we have in general
no notion of convergence, we can take some sort of a limit of these
objects, namely a {\it projective limit}.  

\subheadline{II. Interpretation of tangent objects via scalar extensions}

\msk \nin
{\bf 0.3.}  {\sl Justification of infinitesimals.}
We define {\it manifolds} and {\it tangent bundles} in the
classical way using charts and atlasses (Ch.\ 2 and 3);
then, intuitively, we may think of the tangent space $T_p M$ of 
$M$ at $p$ as a ``(first order) 
infinitesimal neighborhood" of $p$. This idea may be formalized by
writing, with respect to some fixed chart of $M$, 
a tangent vector at the point $p$ in the form
$p + \eps v$, where $\epsilon$ is a ``very small" quantity  
(say, Planck's constant),
thus expressing that the tangent vector is ``infinitesimally
small" compared to elements of $M$ (in a chart).
The property of being ``very small" can mathematically be
expressed by requiring that $\epsilon^2$ is zero, compared
with all ``space quantities". 
This suggests that, if $M$ is a manifold modelled on the $\K$-module
$V$, then the tangent bundle should be modelled on the space
$V \oplus \eps V := V \times V$, which is considered as a module over
the ring $\K[\eps]=\K \oplus \eps \K$
of {\it dual numbers over $\K$} (it is constructed from 
$\K$ in a similar
way as the complex numbers are constructed from $\R$,
replacing the condition $i^2=-1$ by $\epsilon^2=0$). 
All this would be really meaningful if  $TM$ were a manifold
not only over $\K$, but also over the extended ring $\K[\epsilon]$.
In Chapter 6 we prove that this is indeed true:
the ``tangent functor" $T$ can be seen as a functor of scalar
extension by dual numbers (Theorem 6.2). Hence one may use the ``dual
number unit" $\epsilon$ when dealing with tangent bundles 
with the same right as the imaginary unit $i$ when dealing with
complex manifolds.  The proof of Theorem 6.2 is conceptual and
allows to understand why dual numbers naturally appear in this
context: the basic idea is that one can ``differentiate
the definition of differentiability" --  we write (0.2) as a
commutative diagram, apply the tangent functor to this diagram
and get a diagram of the same form, but taken over the tangent
ring $T\K$ of $\K$. The upshot is now that $T\K$ is nothing
but $\K[\epsilon]$.
				  
\ssk
Theorem 6.2 may be seen as a bridge from differential geometry
to algebraic geometry where the use of dual numbers for  
modelizing tangent objects is a standard technique (cf.\ [GD71], 
p.\ 11), 
and also to classical geometry (see article of Veldkamp in [Bue95]).
It is all the more surprising that in most textbooks
on differential geometry no trace of the dual number unit $\epsilon$
can be found -- although it clearly leaves a trace which is
well visible already in the usual real framework: in the same way 
as a complex structure
induces an almost complex structure on the underlying real manifold,
a dual number structure induces a tensor field of endomorphisms
having the property that its square is zero; let us call this an
``almost dual structure". Now, there is a 
{\it canonical} almost dual structure
on every tangent bundle (the almost 
trivial proof of this is given in Section 4.6).
 This canonical almost dual structure
is ``integrable" in the sense that its kernel is a distribution
of subspaces that admits integral submanifolds -- namely the
fibers of the bundle $TM$. Therefore Theorem 6.2 may be seen as
an analog of the well-known theorem of Newlander and Nirenberg: our integrable
almost dual structure induces on $TM$ the structure of
a manifold over the ring $\K[\epsilon]$.

\ssk Some readers may wish to avoid the use of rings which are not fields,
or even to stay in the context of the real base field. In principle,
all our results that do not directly involve dual numbers can be proved 
in a ``purely real" way, i.e., by interpreting $\eps$ just as a formal
(and very useful !) label which helps to distinguish two
copies of $V$ which play different r\^oles (fiber and base). But in the
end, just as the ``imaginary unit" $i$ got its well-deserved place in
mathematics, so will the ``infinitesimal unit" $\eps$.

\msk \nin                                      
{\bf 0.4.} {\sl Further scalar extensions.}
The suggestion to use dual numbers in differential geometry is not
new -- one of the earliest steps in this direction was by
A. Weil ([W53]); one of the most recent is [Gio03]. 
However, most of the proposed constructions  are so complicated
that one is discouraged  to iterate them. But this is
exactly what makes the dual number formalism so useful:
for instance,
if $TM$ is the scalar extension of $M$ by $\K[\epsilon_1]$,
then the double tangent bundle $T^2 M:=T(TM)$ is simply the scalar
extension of $M$ by the ring
$$
TTK:=T(T\K)=\K[\epsilon_1][\epsilon_2]  \cong \K \oplus \epsilon_1 \K \oplus
\epsilon_2 \K \oplus \epsilon_1 \epsilon_2 \K,
\eqno (0.4)
$$
and so on for all higher tangent bundles $T^k M$.
As a matter of fact,  most of the important notions of differential
geometry deal, in one way or another, with the second order
tangent bundle $T^2 M$ (e.g. Lie bracket, exterior derivative, connections)
or with $T^3 M$ (e.g. curvature) or even higher order tangent bundles
(e.g. covariant deriviate of curvature). Therefore 
second and third order differential geometry really is the 
central part of all differential geometry, and finding a good notation
concerning second and third order tangent bundles becomes a necessity.
Most textbooks, if at all $TTM$ is considered, use a component notation
in order to describe objects related to this bundle.
In this situation, the use of {\it different} symbols
$\epsilon_1,\epsilon_2,\ldots$
 for the infinitesimal units of the various scalar extensions
is a great notational progress, combining algebraic rigour and
transparency. It becomes clear that   
many structural features of $T^k M$ are simple consequences of
corresponding structural features of the rings $T^k \K = \K[\epsilon_1,
\ldots,\epsilon_k]$ (cf. Chapter 7). For instance, on both
objects there is a canonical action of the permutation group $\Sigma_k$;
the subring $J^k \K$ of $T^k \K$ fixed under this action is called a {\it jet
ring} (Chapter 8) -- in the case of characteristic zero this is
a truncated polynomial ring $\K[x]/(x^{k+1})$.
Therefore the subbundle $J^k M$
of $T^k M$ fixed under the action of $\Sigma_k$
can be seen as the scalar extension of $M$ by the jet ring $J^k \K$;
it can be interpreted as the {\it bundle of $k$-jets of (curves in) $M$}. 
The bundles $T^k M$ and $J^k M$ are the stage on which
higher order differential geometry is played.

\subheadline{III. Second  order differential geometry}

\msk \nin
{\bf 0.5.} {\sl Bilinear  bundles and connections.}
The theory of {\it linear connections} and their {\it curvature}
is a chief object of general differential geometry in our sense.
When studying higher tangent bundles $T^2M, T^3M,\ldots$ or
tangent bundles $TF, T^2 F,\ldots$ of linear (i.e. vector-)
 bundles $p:F \to M$,
one immediately encounters the problem that these bundles are
{\it not linear  bundles over $M$}.
One may ask, naively: ``if it is not a linear bundle, so what
is it then ?'' The answer can be found, implicitly, in all
textbooks where the structure of $TTM$ (or of $TF$) is discussed in
some detail
(e.g. [Bes78], [La99], [MR91]), but it seems that so far
no attempt has been made to separate clearly  linear algebra from 
differential calculus. Doing this, one is lead to a purely
algebraic concept which we call {\it bilinear space} (cf. Appendix
Multilinear Geometry). A {\it bilinear bundle} is simply a bundle
whose fibers carry a structure of bilinear space which is preserved
under change of bundle charts (Chapter 9). For instance, $TTM$ and $TF$ 
are bilinear bundles. Basically, a bilinear space $E$
is given by a whole collection of linear (i.e., $\K$-module-) 
structures, parametrized
by some other space $E'$, satisfying some purely algebraic axioms.
 We then say that these linear structures are {\it bilinearly 
related}. A {\it linear connection} is simply given by singling out,
in a fiberwise smooth way, in each fiber {\it one} linear structure among
all the bilinearly related linear structures.

\ssk
Of course, in the classical case our definition of a linear connection
coincides with the usual ones -- in fact, in the literature one finds
many different definitions of a connection, and, being aware of the
danger that we would be exposed to\footnote{$^1$}{[Sp79, 
p.\ 605]: ``I personally feel that the next
person to propose a new definition of a connection should be
summarily executed.''},
it was not our aim to add a new item to this long list.
Rather, we hope to have found the central item around which one
can organize the spider's web of related notions such as
{\it sprays, connectors, connection one-forms} and 
{\it covariant derivatives} (Chapters 10, 11, 12).
 We resisted the temptation to start with general  Ehresmann
connections on general fiber bundles. 
In our approach,  principal bundles play a much less important r\^ole 
than in the usual framework --
a simple reason for this is that the general (continuous) linear 
automorphism group
$\Gl_\K(V)$ of a  topological $\K$-module is
in general no longer a Lie group, and hence the prime example of a principal
bundle, the frame bundle, is not at our disposition.

\subheadline{ IV. Third and higher order differential geometry}

\msk \nin
{\bf 0.6.} {\sl Multilinear bundles and multilinear connections.} 
The step from second order to third order geometry is important and
conceptually not easy since it is at this stage that {\it non-commutativity}
enters the picture: just as $TTM$ is a bilinear bundle, the higher order
tangent bundles $T^k M$ carry the structure of a {\it multilinear bundle}
over $M$. The fibers are {\it multilinear spaces} (Chapter MA), which 
again are defined to be a set $E$ with a 
whole collection of linear structures
(called {\it multilinearly related}), parametrized by some space $E'$,
and a {\it multilinear connection} is given by fixing {\it one} of these
multilinearly related structures (Chapters 15 and 16).
In case $k=2$, $E'$ is an {\it abelian} (in fact, vector)
group, and therefore the space of linear connections
is an affine space over $\K$. For $k > 2$, the space $E'$ is a {\it
non-abelian}
group, called the {\it special multilinear group} $\Gm^{1,k}(E)$.

\msk \nin
{\bf 0.7.} {\sl Curvature.}
Our strategy of analyzing the higher order tangent bundle $T^k M$ is
by trying to define, in a most canonical way, a sequence of multilinear
connections on $TTM, T^3M,\ldots$, starting with a linear connection
in the sense explained above. In other words, we want to define
vector bundle structures over $M$ on $TTM, T^3 M,\ldots$, 
in such a way that they define, by restriction, also vector bundle 
structures over $M$ on
the jet bundles $J^1 M,J^2 M,\ldots$ which are compatible with taking 
the projective limit $J^\infty M$. 
Thus on the filtered, non-linear bundles $T^k M$ and $J^k M$ 
we wish to define a sequence of graded and linear structures.

\ssk
It turns out that there is a rather canonical procedure how to define
such a sequence of linear structures: we start with a linear connection
$L$ and define suitable derivatives $DL,D^2L,\ldots$ (Chapter 17).
Unfortunately, since second covariant derivatives do in general not commute,
the multilinear connection $D^k L$ on $T^{k+1} M$ is in general not
invariant under the canonical action of the permutation group
$\Sigma_{k+1}$. Thus, in the worst case, we get $(k+1)!$ different
multilinear connections in this way. Remarkably enough, they all
can be restricted to the jet bundle $J^{k+1} M$ (Theorem 19.2), but
in general define still $k !$ different linear structures there.
Thus, for $k=2$, we get $2$ linear structures on $TTM$: their difference
is the {\it torsion} of $L$. For $k=3$, we get $6$ linear structures on
$T^3 M$, but their number reduces to $2$ if $L$ was already chosen
to be torsionfree. In this case, their ``difference" is a tensor
field of type $(2,1)$, agreeing with the classical {\it curvature tensor}
of $L$ (Theorem 18.3). The induced linear structure on $J^3 M$ is then
unique, which is reflected by {\it Bianchi's identity} (Chapter 18 and
Section SA.10). For general $k \geq 3$, we get a whole bunch of
{\it curvature operators}, whose combinatorial structure seems to be
an interesting topic for further work. 

\ssk
If $L$ is flat (i.e, it has vanishing torsion and curvature tensors),
then the sequence $DL,D^2 L,\ldots$ is indeed invariant under permutations
and under the so-called {\it shift operators} (Chapter 20), which are
necessary conditions for integrating them to a ``canonical chart"
associated to a connection (Theorem 20.7). For general connections,
the problem of defining a permutation and shift-invariant sequence of
multilinear connections is far more difficult (see Section 0.16, below).

\msk \nin {\bf 0.8.} {\sl Differential operators, duality and 
``representation theory."} A complete theory of differential geometry
should not exclude the cotangent-``functor"
$T^*$, which in some sense plays a more important r\^ole in most modern
theories than the tangent functor $T$. In fact, the duality between
$T$ and $T^*$ corresponds to a duality between {\it curves} $\gamma:\K \to M$
and {\it  (scalar) functions} $f:M \to \K$.  
If the theory is rather based on curves,
then the tangent functor will play a dominant role since $\gamma'(t)$
is a tangent vector, and if the theory is rather based on functions,
then $T^*$ will play a dominant role since $df(x)$ is a cotangent vector.
Under the influence of contemporary algebraic geometry, most
modern theories adopt the second point of view. By choosing the first
one, we do not want to take a decision on ``correctness" of one
point of view or the other, but, being interested in the ``covariant theory",
 we just simplify our life and avoid technical problems such as
double dualization and topologies on dual spaces which inevitably arise if
one tries to develop a theory of the tangent functor in an approach based on
scalar functions. This allows
to drop almost all restrictions on the model space in infinite dimension
and opens the theory for the use of base fields of positive characteristic
and even of base rings (giving the whole strength to the point of view
of scalar extensions).

\ssk
However,
although we do not define, e.g., tangent vectors as differential operators
acting on functions, they may be {\it represented} by such operators.
Similarly, connections may be represented by covariant derivatives
of sections, and so on. In other words, one can develop a {\it representation
theory}, representing geometric objects from the ``covariant setting",
in a contravariant way,
by operators on functions and sections -- some sketchy remarks on
 the outline of such
a theory are given in Chapter 20, the theory itself remains to be worked out.
In Chapter 21, we define a particularly important differential operator,
namely the {\it exterior derivative} of forms.

\subheadline{V. Lie Theory }

\msk \nin
{\bf 0.9.} {\sl The Lie algebra of a Lie group.}
By definition, {\it Lie groups} are manifolds with a smooth group structure
(Chapter 5). Then, as in the classical real case,
the Lie algebra of a Lie group can be defined in
various ways: 

\ssk
\item{(a)} one may define the Lie bracket via the bracket of
(left or right) invariant vector fields (as is done in most textbooks), 
or 
\item{(b)} one may express it  by a ``Taylor expansion" of 
the group commutator $[g,h]=ghg^{-1}h^{-1}$ (as is done in [Bou72] and
[Se65]; the definition by
deriving the adjoint representation $\Ad:G \to \Gl( \g)$ 
is just another  version of this).

\ssk \nin
For convenience of the reader, we start with the first definition
(Chapter 5), and show (Theorem 23.2) that it agrees with the following 
version of the second definition:
as in the classical theory, the tangent bundle $TG$ of a Lie group $G$
carries again a canonical Lie group structure, and by induction, all
higher order tangent bundles $T^k G$ are then Lie groups. The projection
$T^k G \to G$ is a Lie group homomorphism, giving rise to an exact sequence
of Lie groups
$$
\matrix{1 & \to & (T^k G)_e & \to & T^k G & \to & G & \to 1 \cr}.
\eqno (0.5)
$$
For $k=1$, the group $\g := T_e G$ is simply the tangent space at the
origin with vector addition. For $k=2$, the group $(T^2 G)_e$ is no
longer abelian, but it is two-step nilpotent, looking like a Heisenberg
group: there are three abelian subgroups $\eps_1 \g, \eps_2 \g,\eps_1 \eps_2
\g$, all canonically isomorphic to $(\g,+)$, 
such that the commutator of elements from the first two of these subgroups
can be expressed by the Lie bracket in $\g$ via the relation
$$
\eps_1 X \cdot \eps_2 Y \cdot (\eps_1 X)^{-1} \cdot (\eps_2 Y)^{-1} =
\eps_1 \eps_2 [X,Y].
\eqno (0.6)
$$

\ssk \nin {\bf 0.10.} {\sl The Lie bracket of vector fields.}
Definition (a) of the Lie bracket in $\g$ requires a foregoing
definition of the Lie bracket of vector fields. In our general framework, 
this is more complicated than in the classical context because we do
not have a faithful representation of the space $\X(M)$
of vector fields on $M$ by 
differential operators acting on functions (cf.\ Section 0.8). 
It turns out (Chapter 14) that a conceptual
definition of the Lie bracket of vector fields is completely analogous
to Equation (0.6), where $\g$ is replaced by $\X(M)$ and $(T^2 G)_e$ 
by the space $\X^2(M)$ of smooth sections of the projection $TTM \to M$.
The space $\X^2(M)$ carries a canonical group structure, again looking
like a Heisenberg group (Theorem 14.3).
Strictly speaking, for our approach to Lie groups, Definition (a) of the
Lie bracket is not necessary. The main use of left- or right-invariant
vector fields in the classical theory is to be integrated to a flow
which gives 
the exponential map. In our approach, this is not possible, and we have
to use rather different strategies. 

\msk \nin
{\bf 0.11.} {\sl Left and right trivializations and analog of the 
Campbell-Hausdorff formula.} 
Following our general philosophy of multilinear connections, we try
to linearize $(T^k G)_e$ by introducing a suitable ``canonical chart''
 and to describe its group structure by a multilinear formula. 
This is possible in two canonical ways: via the zero section, the
exact sequences (0.5) split, and hence we may write $T^k G$ as
a semidirect product by letting act $G$ from the right or from the left.
By induction, one sees in this way that $(T^k G)_e$ can be written
via left- or right trivialization as an iterated 
semidirect product involving $2^k - 1$ copies of $\g$. Relation (0.6)
generalizes to a general ``commutation relation" for this group,
so that we may give an explicit formula for the group structure
on the space $(T^k \g)_0$ which is just a
certain direct sum of $2^k-1$ copies of $\g$. 
The result (Theorem 24.7) may be considered as a rather primitive
version of the Campbell-Hausdorff formula. It has the advantage
that its combinatorial structure is fairly transparent and that it
works in arbitrary characteristic, and it has the drawback that
inversion is more complicated than just taking the negative (as in
the Campbell-Hausdorff group chunk). 
It should be interesting to study the combinatorial aspects of this
formula in more detail, especially for the theory of Lie groups
in positive characteristic.

\msk \nin
{\bf 0.12.} {\sl The exponential map.}
As already mentioned, one cannot construct an
exponential mapping $\exp:\g \to G$ in our general context.
However, the groups $(T^k G)_e$ have a nilpotent Lie algebra $(T^k \g)_0$,
and thus one would expect that a polynomial exponential map
$\exp_k:(T^k \g)_0 \to (T^k G)_e$ (with polynomial inverse $\log_k$) 
should exist.
This is indeed true, provided that $\K$ is a field of characteristic zero
 (Theorem 25.2). Our proof of this fact is purely
algebraic (Appendix PG, Theorem PG.6)
and uses only the fact that the above mentioned 
analog of the Campbell-Hausdorff 
formula defines on $(T^k G)_e$ the structure of a 
{\it polynomial group} (see Section 0.18 below).
By restriction to the jet bundle $J^k G$, we get an exponential map $\exp_k$
for the nilpotent groups $(J^k G)_e$, which form of projective system so
that, as projective limit, we get an exponential map $\exp_\infty$ of
$(J^\infty G)_e$ (Chapter 32). If there is a notion of convergence
such that $\exp_\infty$ and $\log_\infty$ define convergent series, 
then they can indeed be used to define a ``canonical chart" and thus
a ``canonical atlas" of $G$. 

\msk \nin
{\bf 0.13.} {\sl Symmetric spaces.}
The definition and basic theory of {\it symmetric spaces} follows the one
given by O.\ Loos in [Lo69] (see also [BeNe05], where 
examples are constructed):  a symmetric space is a manifold
$M$ equipped with a family of symmetries $\sigma_x:M \to M$, $x \in M$,
such that the map $(x,y) \mapsto \mu(x,y):=\sigma_x(y)$ is smooth and
satisfies certain natural axioms  (Chapter 5). In [Lo69], the higher order
theory of symmetric spaces is developed by using the {\it linear} bundles
$T^{(k)} M$
whose sections are $k$-th order differential operators. 
Again following our general philosophy, we replace these bundles by the
iterated tangent bundles $T^k M$ which are not linear, but have the
advantage to be themselves
symmetric spaces in a natural way. In particular, we
construct the {\it canonical connection} of $M$ in terms of the
symmetric space $TTM$ (Chapter 26)
and express its curvature by the symmetric space
structure on $T^3M$ (Chapter 27). 
The construction of a family $\exp_k:(T^k (T_o M))_0 \to (T^k M)_o$
of exponential mappings (whose projective limit $\exp_\infty$ can
again be seen as a formal power series version of the missing exponential
map $\Exp:T_o M \to M$) uses some general results to be given in 
Part VI (below) and is achieved in Chapter 30.

\subheadline{VI. Groups of diffeomorphisms and the exponential jet}

\msk \nin
{\bf 0.14.} {\sl The group $\Diff(M)$.} If $M$ is a manifold, then the group
 $G:=\Diff(M)$  of all diffeomorphisms of  $M$ is in general not a Lie
group. However, quite often it is useful to think of $G$ 
as a Lie group with Lie algebra $\X(M)$, the Lie algebra of vector
fields on $M$, and exponential map associating to a vector field $X$
the flow $\Fl_t^X \in G$ at time $t=1$ (whenever possible).
Then, according to our approach to Lie theory, we would like 
to analyze the group $G$ by studying the higher order
tangent group $T^k G$ and its fiber $(T^k G)_e$ over 
the origin $e \in G$. 
It turns out that this approach makes perfectly sense for any
manifold $M$ over a general base field or -ring $\K$: 
we prove that the group $\Diff_{T^k \K}(T^k M)$ of all diffeomorphisms
of $T^k M$ which are smooth over the iterated tangent ring $T^k \K$,
behaves in all respects like the $k$-th order tangent group of $G$.
More precisely, there is an exact sequence of groups
$$
\matrix{1& \to& \X^k(M) & \to &\Diff_{T^k \K}(T^k M)& \to& \Diff(M)& \to & 1
\cr}
\eqno {(0.7)}
$$
which is the precise analog of (0.5). 
In particular, $\X^k(M)$, the space of smooth sections of 
the projection $T^k M \to M$, carries a canonical group structure which 
turns it into a polynomial group (Theorem 28.3).
For $\K=\R$ and $M$ finite-dimensional, it is
known that $\X^k(M)$ carries a natural group structure (this result is
due to Kainz and Michor, [KM87, Theorem 4.6], see also [KMS93, Theorem
37.7]), but the interpretation via the sequence (0.7) seems to be entirely
new.  A splitting of the exact sequence (0.7) is simply given by
$\Diff(M) \to \Diff_{T^k \K}(T^k M)$, $g \mapsto T^k g$, and thus 
$\Diff_{T^k \K}(T^k M)$ is
a semidirect product of $\Diff(M)$ and $\X^k(M)$. 
As for Lie groups, we can now define left- and right-trivializations
and have an analog of the Campbell-Hausdorff formula for the group
structure on $\X^k(M)$ (Theorem 29.2).  

\msk \nin {\bf 0.15.} {\sl Flow of a vector field.}
Even in the classical situation, not every vector field $X \in \X(M)$ can
be integrated to a global flow $\Fl_t^X:M \to M$, $t \in \R$.
In our general set-up, this cannot be done even locally.
But to some extent, the following construction can be seen as
an analog of the flow of a vector field:
as remarked above, $\X^k(M)$ carries the structure of a polynomial group,
and thus by Theorem PG.6, if $\K$ is a field of characteristic zero,
there is a bijective exponential map $\exp_k$ associated to this group.
Since $\X^k(M)$ plays the r\^ole of $(T^k G)_e$ for $G=\Diff(M)$,
the ``Lie algebra" of $\X^k(M)$ should be $(T^k \g)_0$, for $\g=\X(M)$.
As a $\K$-module, this is simply the space of section of the 
{\it axes-bundle} $A^k M$ which, by definition,
is  a certain direct sum over $M$ of $2^k-1$ copies of $TM$.
Thus $\exp_k$ will be a bijection from the space of 
sections of $A^k M$ to the space of sections of $T^k M$.
But using a canonical diagonal imbedding $TM \to A^k M$, every vector
field $X$ gives rise to a section $\delta X$ of $A^k M$, then via  
the exponential map to a section $\exp_k(\delta X)$ of $T^k M$ which finally
corresponds to a diffeomorphism $\Phi:T^k M \to T^k M$.
On the other hand,
if $X$ admits a flow $\Fl_t^X:M \to M$ 
in the classical sense (e.g., if $X$ is a complete vector field 
on a finite-dimensional
real manifold), then $T^k (\Fl_1^X):T^k M \to T^k M$ also defines a
diffeomorphism of $T^k M$. These two diffeomorphisms are closely
related to each other (see Chapter 29).
 Summing up, we have constructed an algebraic 
substitute of the ``missing exponential map" $\X(M) \to \Diff(M)$.

\msk \nin
{\bf 0.16.} {\sl The exponential jet of a connection.}
In this final part we come back to the problem left open in
Part IV, namely to the definition of a canonical sequence of
permutation- and shift-invariant multilinear connections
(Section 0.7).

\ssk In order to motivate our approach,
assume $M$ is a real finite-dimensional (or Banach) manifold 
equipped with an affine connection on the tangent bundle $TM$.
Then, for every $x \in M$, the exponential map $\Exp_x:U_x \to M$ of the
connection
is a local diffeomorphism from an open neighborhood of the origin in
$T_x M$ onto its image in $M$. Applying the $k$-th order tangent functor,
we get $T^k \Exp_x:T^k (U_x) \to T^k M$, inducing  bijections of fibers
over $0_x$ and $x$,
$$
(T^k \Exp_x)_{0_x}:\quad (A^k M)_x:= (T^k(T_x M))_{0_x} \to (T^k M)_x.
\eqno (0.8)
$$
As above, $(A^k M)_x$ is a certain direct sum of copies of $T_x M$,
the fiber of the  ($k$-th order)
axes-bundle of $M$. The bijections of fibers fit 
together to a smooth bundle isomorphism $T^k \Exp:A^k M \to T^k M$,
which we call the {\it exponential jet of the connection}.
In fact, this isomorphism is a multilinear connection in the sense
explained in Section 0.6, and has the desired properties of {\it shift
invariance} and {\it total
symmetry under the permutation group}. Therefore it is of basic
interest to construct the exponential jet $\Exp^k:A^k M \to T^k M$
 in a ``synthetic" way, 
starting with a connection on $TM$ and without ever using the map
$\Exp_x$ itself, which doesn't exist in our general framework.
In Chapters 30 and 31, we describe such a construction for Lie
groups and symmetric spaces, and then outline its definition for
a general (torsionfree) connection. 
Finally, in Chapter 32 we explain the problem of integrating jets to
{\it germs}: given an {\it infinity jet} (a projective limit of $k$-jets,
such as $\Exp^\infty:= \lim_{\leftarrow } \Exp^k$), one would like to
define a corresponding germ of a map, such as, for example,
the germ of a ``true'' exponential map. But this is of course a problem
belonging to integral geometry and not to differential geometry in our sense.

\subheadline{VII. Multilinear geometry}

\nin
The appendix on ``Multilinear geometry'' gathers all purely algebraic
results that are used in the main text. Two topics are treated:

\ssk
\item{(1)} We introduce the concept of a {\it multilinear
space} (over a commutative base ring $\K$); this is the algebraic model
for the fibers of higher order tangent bundles $T^k M$ over the base $M$
or, slightly more general, for the fibers of $T^{k-1} F$, where $F$ is a
vector bundle over $M$. 
\item{(2)} We introduce {\it polynomial groups}; these are essentially
formal groups (in the power-series approach, cf.\ [Haz78]) whose power
series are finite, i.e.\ polynomial.

\ssk
\nin The link between (1) and (2) comes from the fact that
the transition functions of bundles like $T^k F$ over $M$
are not linear but {\it multilinear} (Chapter 15), and they
belong to a group which can be defined in a purely algebraic way, called
the {\it general multilinear group}. Now, it turns out that this group
is a polynomial group, and hence we can apply  Lie theory in a
purely algebraic set-up. Conversely, much of the structure of a multilinear
space is encoded in the general multilinear group.
We now describe the concepts (1) and (2) in some more detail.

\msk \nin
{\bf 0.17.} {\sl Multilinear spaces.}
 The ingredients of the purely algebraic setting  are: a family
$(V_\alpha)_{\alpha \in I_k}$ of $\K$-modules (called a {\it cube
of $\K$-modules}), where $I_k$ is the
lattice of all subsets of the finite set $\N_k=\{ 1,\ldots,k \}$
($k \in \N$),
the {\it total space} $E:=\oplus_{\alpha \in I_k \setminus \emptyset} 
V_\alpha$ on
which the {\it special multilinear group} $\Gm^{1,k}(E)$ acts in a non-linear
way (essentially,  by multilinear maps in the components $V_\alpha$ of $E$).
We say that a structure on $E$ is {\it intrinsic} if it is invariant under
the action of the general multilinear group. There is no intrinsic linear
structure on $E$, but the collection of {\it all} linear structures obtained
by pushing around the initial linear structure on $\oplus_\alpha V_\alpha$
is an intrinsic object. Thus a multilinear space is not a $\K$-module but
rather a collection of many different $\K$-module structures on a given
set -- this point of view has already turned out to be very useful 
in relation with Jordan theory and geometry ([Be02] and [Be04], which 
somewhat was at the origin of the approach developed here).
In case $k=2$, the theory is particularly simple, and therefore we
first give an independent treatment of {\it bilinear spaces} (Chapter BA):
in this case, one can easily give explicit formulae for the various
$\K$-module structures (Section BA.1), thanks to the fact that in this case
(and only in this case) the special linear group is abelian.
The algebraically trained reader may start directly by reading
Chapter MA.  In the case of general $k \geq 2$,
 one should think of $\Gm^{1,k}(E)$ as a nilpotent Lie group (cf.\
Chapter PG).

\msk
We say that a multilinear space  is {\it of tangent type} if 
 all $\K$-modules $V_\alpha$ are canonically isomorphic to some fixed space
$V$ and hence are canonically isomorphic among each
other -- the terminology is motivated by the fact that the fibers of
higher order tangent bundle $T^k M$ have this property.
Then the symmetric group $\Sigma_k$  
acts on the total space $E$, where the action is
 induced by the canonical action of permutations on the set $\N_k$. 
This action is ``horizontal" in the sense that it preserves the
cardinality of subsets of $\N_k$.
 There is also another class of symmetry operators
acting rather ``diagonally" on the lattice of subsets of $\N_k$
 and called {\it shift operators}.
Invariance properties under these two kinds of operators and their
interplay are investigated in Chapter SA.

\ssk
Let us add some words on the relation of our concepts with
{\it tensor products}. By their universal property, tensor products
transform multilinear algebra into linear algebra.
In a way, our concept has the same purpose, but, whereas the tensor
product produces a new underlying set in order to reduce everything to
{\it one} linear structure, we work with {\it one} underlying set and produce a
new supply of linear structures.
Thus, in principle, it is possible to reduce our set-up into the linear
algebra of one linear structure
 by introducing some big tensor product space in a functorial
way such that the action of the general multilinear group becomes linear --
we made some effort to provide the necessary formulas by introducing a
suitable ``matrix calculus" for $\Gm^{1,k}(E)$. 
However, the definition of the
``big" tensor space is combinatorially rather complicated and destroys
the ``geometric flavor" of the non-linear picture; instead of simplifying
the theory, it causes unnecessary problems. This becomes clearly visible
if the reader compares our non-linear picture
 with the concept of an {\it osculating bundle of a vector bundle}
 introduced by F.W.\
Pohl in [P62]: with great technical effort, he constructs a {\it linear}
bundle which corresponds to $T^k F$ in the same way as the big tensor
space corresponds to a multilinear space, and finally admits
(loc.\ cit.\  p.\ 180): ``In spite of several attempts, we have not
been able to reduce even the [construction of the first osculating
bundle] to something familiar."

\ssk
Another advantage of our concept is that it behaves nicely with respect
to topology: since we do not change the underlying set $E$, we have
no problem in defining a suitable topology on $E$
 if $\K$ is a topological ring and the $V_\alpha$ are topological
$\K$-modules. For tensor products, the situation is much worse:
in general, topological tensor products of
topological $\K$-modules have  very bad properties which
make them practically useless (general topological tensor
products are not even associative, even over $\K=\R$, see [Gl04b]).

\msk \nin {\bf 0.18.} {\sl Polynomial groups.}
The general multilinear group $\Gm^{1,k}(E)$ is a special instance of a
{\it polynomial group}: it is at the same a group and a $\K$-module
such that the group multiplication is polynomial, and there is a common
bound on the degree of all iterated product maps.
One should think of a polynomial group as a nilpotent Lie group together
with some global chart in which the group operations are polynomial.
To such groups we may apply all results on {\it formal (power series) groups}
(see, e.g., [D73], [Haz78] or [Se65]), but since we are in a polynomial
setting, everything converges globally, and thus our results have
a direct  interpretation in terms of mappings.
 In particular, we can associate a Lie algebra to a polynomial
group and define (of $\K$ is a field of 
 characteristic zero) an exponential map and a logarithm
(Theorem PG.6).
Our proof of this result is much more elementary and closer to ``usual''
analysis on Lie groups than the proofs given for general formal groups
in [Bou72] or [Se65] -- we work with {\it one-parameter subgroups} 
essentially as in the differentiable case and 
do not need to invoke the Campbell-Hausdorff
group chunk nor the universal enveloping algebra.
Viewing formal power series groups as projective limits of polynomial
groups, our result also implies corresponding results for general formal
groups.

\subheadline{Limitations, Generalizations}

\msk
\nin {\bf 0.19.} {\sl Limitations.}
The purpose of the short Appendix L is to illustrate our distinction
between ``differential" and ``integral" geometry by a list of
results and notions that can {\it not} be treated in a purely differential
context. On the other hand, a good deal of these items can be generalized
from the finite-dimensional real set-up to surprisingly general
situations (i.e., for rather big classes of topological base fields
and of topological vector spaces) -- many of these very recent results
are due to H.\ Gl\"ockner (see references).

\ssk \nin {\bf 0.20.} {\sl Generalizations.}
Our methods apply to more general situations than just to manifolds over
topological fields and rings. Some of these possibilities of generalization
are discussed in Appendix G.

\subheadline{Related literature}

It is impossible to take account of all relevant literature on such
a vast area as Differential Geometry and Lie Theory, and I apologize 
for possibly omitting to give reference where it would have been 
due. 

\msk
{\sl Real, finite and infinite dimensional differential geometry.}
For geometry on mani\-folds modelled on Banach-spaces, [La99] is my basic
reference (the treatment of the second order tangent bundle and the
Dombrowski Splitting Theorem 10.5 are motivated by the corresponding
section of [La99]). From work of K.-H.\ Neeb (cf.\ [Ne02b] and
[Ne06]) I learned that
many notions of infinite-dimensional differential geometry from [La99] 
can be generalized to manifolds modelled on general locally convex vector
spaces which need not even be complete. Finally, although it treats only
the finite-dimensional case, the monograph ``Natural Operations in
Differential Geometry" by I.\ Kol\'a\v r, P.\ Michor and J.\ Slov\'ak
[KMS93] turned out to be most closely related to the spirit of the
present work. The reason for this is that an operation, if it is really
``natural", is robust and vital enough to survive in most
general contexts. In particular, Chapter VIII on ``Product Preserving
Functors" of [KMS93] had an important influence on my work. That chapter, 
in turn, goes back to the very influential paper ``Th\'eorie des points
proches sur les vari\'et\'es diff\'erentiables" by A.\ Weil [W53]. 

\ssk {\sl ``Synthetic differential geometry" and ``smooth toposes".}
These theories spring, in one way or another, from Weil's above mentioned 
paper. They propose a unification of differential geometry and algebraic 
geometry in very general categories which contain objects like
$\Diff(M)$ and permit to use dual numbers as a means to modelize 
``tangent" objects (see [Ko81], [Lav87] or [MR91]).
A similar approach is developed by the
russian school of A.\ M.\ Vinogradov (cf.\ [ALV91], [Nes03]).
It would be interesting to investigate the relations between those
approaches and the one presented here. 
It seems that these theories are ``dual" to ours in the sense explained
in Section 0.8: whereas we base our approach on the covariant tangent
functor $T$, those approaches are based on functions and thus are
contravariant in nature. In principle, it seems conceivable to
represent (in the sense of Section 0.8) our geometric objects in
``smooth toposes over $\K$" (cf.\ Appendix G). 

\ssk {\sl Formal groups and other ``formal" theories.}
Our version of Lie theory has much in common
with the theory of {\it formal (power series) groups} (see [Haz78]),
and our version of differential geometry is related to ``formal"
concepts such as {\it formal geometry} (see [CFT01]).
As above, this may also be considered from a ``contravariant"
point of view, and then corresponds to the approach to formal
groups via Hopf algebras (see [D73]). Again, it should be interesting
to interprete the second approach as a ``representation" of the first one.

\subheadline{Leitfaden}

The chapters on Lie theory can be read without having to go through
all the differential geometric parts, and conversely.
The reader who is mainly interested in Lie groups may 
read Chapters 5 -- 9 -- 23 -- 15 -- 24 -- PG -- 25,
and the  reader interested in symmetric spaces could read 
Chapters 5 -- 9 -- 26 -- 15 -- 27 -- 28 -- 29 -- 30. For such a reading,
the interpretation of the tangent functor as a functor of scalar
extension is not indispensable, provided one gives a proof of
 Theorem 7.5  not making use of dual numbers (which is possible). 

\subheadline{Notation}

Throughout the text, $\K$ denotes a commutative ring with unit $1$.
In certain contexts (e.g., symmetric spaces) we will assume that $2$
is invertible in $\K$ or that $\K$ is a field (Section PG) or even that $\K$
a field of characteristic zero (when talking about exponential mappings).
Some more specific notation concerning the following topics is explained in
the corresponding subsections:

\msk
\item{--} differential calculus: {\bf 1.3}
\item{--} dual numbers, iterated dual numbers, jet rings: {\bf 6.1}, {\bf 7.3}, {\bf 7.4}, {\bf 8.1}
\item{--} partitions (set theory): {\bf 7.4}, {\bf MA.1}, {\bf MA.3}, 
{\bf MA.4}
\item{--} bilinear maps, bilinear groups: {\bf BA.1}, {\bf BA.2}
\item{--} multilinear maps, multilinear groups: {\bf MA.5}

\vfill \eject 

%% file: Dgeo1.tex

\def\title{I. BASIC NOTIONS} 
\def\rightheadline{\hfil{\eightrm \title}\hfil\tenrm\folio}

\sectionheadline{I. Basic notions}

This part contains the ``conventional'' aspects of the theory:
having recalled the basic facts on differential calculus, we define
manifolds and bundles in the usual way by charts and atlasses and we 
define the Lie bracket of vector fields
in the old-fashioned way via its chart 
representation (the reader who prefers a more sophisticated definition
is referred to Chapter 14). This is used in order to define, in Chapter 5,
the {\it Lie functor} assigning a Lie algebra to a Lie group and a 
Lie triple system to a symmetric space.

\subheadline{1. Differential calculus}

\nin {\bf 1.1.} {\sl Topological rings and modules.}
By {\it topological ring} we mean a ring $\K$ with unit $1$ 
together with a topology (assumed to be Hausdorff)
such that the ring operations $+:\K \times \K \to \K$ and $\cdot:
\K \times \K \to \K$ are continuous and such that the set
$\K^\times$ of invertible elements is open in $\K$ and the inversion
map $i:\K^\times \to \K$ is continuous.
If $\K$ is an addition a field, it will be called a {\it topological
field}. In the sequel
 we will always  assume that $\K$ is a topological ring such
that the set $\K^\times$ of invertible elements
is dense in $\K$. (If $\K$ is a field, this means that $\K$
is not discrete.)

A module $V$ over a topological ring $\K$ is called a 
{\it topological $\K$-module} if the structure maps
$V \times V \to V$ and $\K \times V \to V$ are continuous.
We will assume that all our topological modules are Hausdorff.

\msk \nin {\bf 1.2.} {\sl The class $C^0$ and the ``Determination Principle''.}
For topological $\K$-modules $V$, $W$ and 
an open set $U \subset V$, we denote by $C^0(U,W):=C(U,W)$
the space of continuous maps $f:U \to W$.
Since $\K^\times$ is assumed to be dense in $\K$,
the value at zero of a continuous map $f:I \to W$, defined
on an open neighborhood $I$ of zero in $\K$, is already determined
by the values of $f$ at invertible elements, i.e., on $I \cap \K^\times$.
We call this the {\it determination  principle}.

\msk \nin
{\bf 1.3.} {\sl The class $C^1$ and the differential.}
We say that {\it $f:U \to W$ is $C^1(U,W)$} or just {\it of class $C^1$}
if there exists a $C^0$-map
$$
f^{[1]}:U \times V \times \K \supset  U^{[1]}:=
\{ (x,v,t) | \, x \in U, x+tv \in U \} \to W,
$$
such that
$$
f(x+tv)-f(x)=t \cdot  f^{[1]}(x,v,t) 
\eqno (1.1)
$$
whenever $(x,v,t) \in U^{[1]}$.
(It is an easy exercise in calculus that, for $\K=\R$,
$V=\R^n$ and $W=\R^m$ this coincides with the usual definition of the
class $C^1$; see [Be07] for a very elementary discussion of these topics,
and [BGN04] for more general results, including the infinite-dimensional
situation.)
The {\it differential of $f$ at $x$} is defined by
$$
df(x):V \to W, \quad v \mapsto
df(x)v:=f^{[1]}(x,v,0).
$$
Note that, by the Determination Principle,
 the map $f^{[1]}$ is uniquely determined by $f$ and
hence $df(x)$ is well-defined. The differential
$$
df:U \times V \to W, \quad (x,v) \mapsto df(x)v=f^{[1]}(x,v,0);
$$
is of class $C^0$ since so is $f^{[1]}$, and for all $v \in V$,
the {\it directional derivative in direction $v$} 
$$
\partial_v f: U \to W, \quad x \mapsto \partial_v f(x):= df(x)v
\eqno (1.2)
$$
is also of class $C^0$.
We define a $C^0$-map $Tf$, called the {\it tangent map}, by
$$
Tf: \, TU=U \times V \to W \times W, \quad
(x,v) \mapsto (f(x),df(x)v)
\eqno (1.3)
$$
and a $C^0$-map, called the {\it extended tangent map}, by
$$
\hat Tf: U^{[1]} \to W \times W \times \K, \quad
(x,v,t) \mapsto (f(x),f^{[1]}(x,v,t),t).
\eqno (1.4)
$$
The following first order differentiation rules are easily proved.

\Lemma 1.4. 
For all $x \in U$, $df(x):V \to W$ is a  $\K$-linear $C^0$-map.

\Proof. We have already remarked that $df$ and hence also $df(x)$ are $C^0$.
Let us prove additivity of $df(x)$. For $s\in \K$ sufficiently close to~$0$,
we have
$$
\eqalign{
s f^{[1]}(x,v+w,s) & = 
f(x+s(v+w))-f(x) \cr
&=  f(x+sv+sw)-f(x+sv)+f(x+sv)-f(x) \cr
&= sf^{[1]}(x+sv,w,s)+s f^{[1]}(x,v,s). \cr}
$$
For $s \in \K^\times$, we divide by $s$ and get
$$
f^{[1]}(x,v+w,s)=
f^{[1]}(x+sv,w,s)+f^{[1]}(x,v,s) .
$$
By the Determination Principle, this equality holds for all $s$ in a 
neighborhood of $0$, hence in particular for $s=0$, whence
$df(x)(v+w)=df(x)w + df(x)v$. Homogeneity of $df(x)$ is proved in the same way.
\qed

\Lemma 1.5.
If $f$ and $g$ are composable and of class $C^1$, then $g \circ f$ 
is $C^1$, and $T(g \circ f)=Tg \circ Tf$.

\Proof. 
Let $(x,y,t)\in U^{[1]}$. Then
$f(x+ty)=f(x)+tf^{[1]}(x,y,t)$,
hence 
$$
g(f(x+ty))-g(f(x)) =
g(f(x)+t f^{[1]}(x,y,t)) - g(f(x))
= t \cdot g^{[1]}(f(x),f^{[1]}(x,y,t),t),
$$
where $(x,y,t)\to g^{[1]}(f(x),f^{[1]}(x,y,t),t)$
is $C^0$ since it is a composition of $C^0$-maps.
Thus $g\circ f$ is $C^1$, with $(g\circ f)^{[1]}$ given by
$$
(g\circ f)^{[1]}(x,y,t)=g^{[1]}(f(x),f^{[1]}(x,y,t),t)\,.
$$
In particular, $d(g\circ f)(x)=dg(f(x))\circ df(x)$
for all $x\in U$, implying our claim.
\qed

\Lemma 1.6.
\item{(i)}
Multilinear maps of class $C^0$ are $C^1$ and are differentiated
as usual. In particular, if $f,g:U \to \K$ are $C^1$, then the
product $f \cdot g$ is $C^1$, and
$\partial_v(fg)=(\partial_v f)g+f\partial_vg$.
Polynomial maps $\K^n \to \K^m$ are always $C^1$ and are
differentiated as usual.
\item{(ii)}
Inversion $i:\K^\times \to \K$ is $C^1$, and
$di(x)v=-x^{-2}v$.
It follows that rational maps $\K^n \supset U \to \K^m$
are always $C^1$ and are differentiated as usual.
\item{(iii)}
The direct product of two $C^1$-maps is $C^1$, and
$$
T(f \times g) = Tf \times Tg,
$$
where we use the convention to identify, for a map 
$f:V \times W \to Z$, the tangent map
$Tf: V \times W \times V \times W \to Z \times Z$
with the corresponding map $V \times V \times W \times W \to
Z \times Z$ (and similarly for maps defined on open subsets of $V 
\times W$). Using this convention, the tangent map of
a diagonal map $\delta_V:V \to V \times V$, $x \mapsto (x,x)$ is 
the diagonal map $\delta_{V \times V}$.
\item{(iv)}
If $f:V_1 \times V_2 \supset U \to W$ is $C^1$, then
the {\it rule on partial derivatives} holds:
$$
df(x_1,x_2)(v_1,v_2)= d_1 f(x_1,x_2)v_1 + d_2 f(x_1,x_2)v_2.
$$

\Proof. One uses the same arguments as in usual differential calculus
(cf.\ [BGN04, Section 2]).
\qed

\nin The rule on partial derivatives may also be written, by introducing
for $(x_1,x_2) \in U$, the ``left'' and ``right translations'',
$$
l_{x_1}(x_2):= r_{x_2}(x_1):= f(x_1,x_2).
$$
Then 
$$
df(x_1,x_2)(v_1,v_2)=d(r_{x_2})(x_1)v_1 +  d(l_{x_1})(x_2)v_2 .
$$
The converse of the rule on partial derivatives holds if $\K$ is a field
(but we won't need it):
if $f$ is $C^1$ withe respect to the first and with respect to the second
variable (in a suitable sense), then  $f$ is of class $C^1$
(see [BGN04, Lemma 3.9]). 

\msk
\nin {\bf 1.7.} {\sl The classes $C^k$ and $C^\infty$.}
Let $f\!: V \supset U \to F$ be of class $C^1$.
We say that {\it $f$ is $C^2(U,F)$} or {\it of class $C^2$\/} 
if $f^{[1]}$ is $C^1$,
in which case we define $f^{[2]}:=(f^{[1]})^{[1]}\!:
U^{[2]}\to F$, where $U^{[2]}:=(U^{[1]})^{[1]}$.
Inductively, we say that $f$ is $C^{k+1}(U,F)$
or {\it of class $C^{k+1}$\/}
if $f$ is of class $C^k$
and $f^{[k]}\!: U^{[k]}\to F$
is of class $C^1$, in which case we
define $f^{[k+1]}:=(f^{[k]})^{[1]}\!: U^{[k+1]}\to F$
with $U^{[k+1]}:= (U^{[k]})^{[1]}$.
The map $f$ is called {\it  smooth\/} or {\it of class $C^\infty$\/}
if it is of class $C^k$ for each $k\in \N_0$.
Note that $U^{[k+1]}=(U^{[1]})^{[k]}$
for each $k\in \N_0$, and that
$f$ is of class $C^{k+1}$
if and only if $f$ is of class $C^1$
and $f^{[1]}$ is of class $C^k$;
in this case, $f^{[k+1]}=(f^{[1]})^{[k]}$.
Next we prove the basic higher order differentiation rules.

\Lemma 1.8.  If $f$ and $g$ are composable and of class $C^k$,
then $g \circ f$ is of class $C^k$, and  the generalized chain rule
$$
\hat T^k(g \circ f)=\hat T^k g \circ \hat T^k f
$$
holds. 

\Proof. For $k=1$, this is Lemma 1.5, and the general case is proved in
the same spirit by  induction on $k$ (cf.\ [BGN04, Prop.\ 4.5]).
\qed

\Lemma 1.9. 
If $f$ is of class $C^2$, then for all $x \in U$, $v,w \in V$,
$$
\partial_v \partial_w f(x) = \partial_w \partial_v f(x).
$$

\Proof. First of all, note that
$$
\eqalign{  
\partial_w \partial_v f(x) & = 
\lim_{t \to 0} {\partial_v f(x+tw)-\partial_v f(x) \over t} \cr
& = \lim_{t \to 0} \Big( \lim_{s \to 0}
{f(x+tw+sv)-f(x+tw)-f(x+sv)+ f(x) \over ts}  \Big)  \cr}
$$
where, for $t,s \in \K^\times$,
$$
\eqalign{
& {f(x+sv+tw)-f(x+tw)-f(x+sv)+ f(x) \over ts}    \cr
& \quad \quad ={1 \over t} \Big( f^{[1]}(x+tw,v,s) - f^{[1]}(x,v,s) \Big) 
= f^{[2]}\big( (x,v,s),(w,0,0),t\big) , \cr}
$$
whence, by continuity of $f^{[2]}$,
$\partial_w \partial_v f(x) = f^{[2]}\big( (x,v,0),(w,0,0),0\big)$.
Of course, we have also, by the same calculation as above, 
for $t,s \in \K^\times$,
$$
\eqalign{
& {f(x+sv+tw)-f(x+sv)-f(x+tw)+ f(x) \over ts}  \cr 
& \quad \quad = {1 \over s} \Big( f^{[1]}(x+sv,w,t) - f^{[1]}(x,w,t) \Big) 
= f^{[2]}\big( (x,w,t),(v,0,0),s\big) , \cr}
$$
and hence
$$
 f^{[2]}\big( (x,v,s),(w,0,0),t\big) = f^{[2]}\big( (x,w,t),(v,0,0),s\big).
$$
By continuity of $f^{[2]}$, this equality holds also for
 $t=s=0$, and this means that
 $\partial_w \partial_v f(x) =\partial_v \partial_w f(x)$.
\qed

\Corollary 1.10.
If $f$ is of class $C^k$ and $x \in U$, then the map
$$
d^k f(x):V^k \to W, \quad
(v_1,\ldots,v_k) \mapsto \partial_{v_1} \ldots \partial_{v_k} f(x)
$$
is a symmetric multilinear $C^0$-map.

\Proof. The maps $d^k$ are partial maps of $f^{[k]}$ and hence are $C^0$;
for instance,  we have seen above that
$d^2 f(x)(v,w)=f^{[2]}((x,w,0),(v,0,0),0)$.
Symmetry follows  from  Lemma 1.9, and
multilinearity now follows from Lemma 1.4.
\qed

\Theorem 1.11. {\rm
(The second order Taylor expansion.)} 
Assume $f:U \to W$ is of class $C^2$. Then for all $(x,h,t) \in U^{[1]}$,
$f$ has an expansion
$$
f(x+th)=f(x)+ t \, df(x)h + t^2 \, a_2(x,h) + t^2 \, R_2(x,h,t),
\eqno (1.5)
$$
where the remainder term $R_{2}(x,h,t)$ is $O(t)$, i.e., it is
 of class $C^0$  and
takes the value $0$ for $t=0$. The coefficient $a_2(x,h)$
is of class $C^0$ jointly in both variables and, in the
variable $h$, is a vector-valued homogeneous form of degree $2$, i.e.,
 $a_2(x, \cdot)$ is homogeneous of degree $2$, and the map
$(v,w) \mapsto a_2(x,v+w)-a_2(x,v)-a_2(x,w)$ is bilinear. More precisely,
we have 
$$
a_2(x,v+w)-a_2(x,v)-a_2(x,w)=d^2f(x)(v,w).
\eqno (1.6)
$$
In particular $2 \, a_2(x,h)= d^2 f(x)(h,h)$. 
Therefore, if $2$ is invertible in $\K$, the second order Taylor expansion
{\rm (1.5)} may also be written
$$
f(x+th)=f(x)+ t \, df(x)h + {t^2 \over 2} \, d^2 f(x)(h,h) + t^2 \, R_2(x,h,t).
\eqno (1.7)
$$

\Proof. Since $f$  is of class $C^2$,
applying  the fundamental relation (1.1) first to $f$ and then
to $f^{[1]}$, we get
$$
\eqalign{
f(x+th)&=f(x) + t f^{[1]}(x,h,t) \cr
&= f(x) + t df(x)h + t ( f^{[1]}(x,h,t)- f^{[1]}(x,h,0)) \cr
&= f(x) + t df(x)h + t^2 f^{[2]}((x,h,0),(0,0,1),t) \cr
&= f(x)+ t df(x)h + t^2  f^{[2]}((x,h,0),(0,0,1),0) \cr
&  \quad \quad 
+ t^2  \left(f^{[2]}((x,h,0),(0,0,1),t)-
 f^{[2]}((x,h,0),(0,0,1),0)\right). \cr}
$$
We let $a_2(x,h):=f^{[2]}((x,h,0),(0,0,1),0)$ and 
$$
R_2(x,h,t):=f^{[2]}((x,h,0),(0,0,1),t)-
 f^{[2]}((x,h,0),(0,0,1),0).
$$
These maps are all $C^0$ since so is $f^{[2]}$, and
$R_2(x,h,t)$  takes the value $0$ at $t=0$. 
One shows, as in usual differential calculus, that an expansion with these
properties (``d\'eveloppement limit\'e'' in French)
is unique (cf.\ [BGN04, Lemma 5.2]), 
and then deduces that $a_2(x,\cdot)$ is homogeneous of degree $2$.
Let us prove (1.6). From (1.5), together with linearity
of $df(x)$, we get
$$
\eqalign{
f(x+t(h_1+h_2))-f(x+th_1)-f(x+th_2)+f(x) & \cr
= \quad  t^2(a_2(x,h_1+h_2)-a_2(x,h_1)-a_2(x,h_2)) & + t^2O(t).
\cr}
$$
Thus, for $t \in \K^\times$,
$$
\eqalign{
{f(x+t(h_1+h_2))
-f(x+th_1)-f(x+th_2)+f(x) \over  t^2} & \cr
= \quad  a_2(x,h_1+h_2)-a_2(x,h_1)-a_2(x,h_2) & + O(t). \cr}
$$
For $t=0$, the $C^0$-extension of the left hand side equals
$d^2f(x)(h_1,h_2)$ (cf.\ proof of Lemma 1.9),
and this  implies (1.6).
Finally, since we already know that $a_2(x,h)$ is homo\-geneous
quadratic in $h$, we get for $h_1=h_2=h$ the relation
$2 \, a_2(x,h)=d^2 f(x)(h,h)$, which implies (1.7).
\qed

\Corollary 1.12. If $U$ is an open neighborhood of the origin in $V$ and
 $f:U \to W$ is of class $C^2$ and
 homogeneous of degree $2$ (i.e., $f(rx)=r^2 \, f(x)$
whenever this makes sense), then $f$ is a $W$-valued quadratic form.

\Proof. Clearly, $f(0)=0$.
Using the expansion (1.5) at $x=0$, we get
$$
t^2 f(h)=f(th) = t df(0)h + t^2 (a_2(0,h) + R_2(0,h,t)).
$$
Dividing, for $t \in \K^\times$, by $t$, we get an identity whose
both sides are $C^0$-functions of $t$. Letting $t=0$, we deduce that
$df(0)h=0$. Dividing again, for $t \in \K^\times$, by $t^2$, we
get $f(h)=a_2(0,h) + R_2(x,h,t)$. Both sides are $C^0$-functions of $t$,
and hence for $t=0$ we get $f(h)=a_2(0,h)$. According to (1.6),
$f$ is then a vector-valued quadratic form.
\qed

\nin
It is fairly obvious that the arguments given in the proof of 
Theorem 1.11 can be iterated
in  order to prove the general expansions, for $f$ of class $C^k$,
$$
\eqalign{
f(x+th) & =f(x)+\sum_{j=1}^k t^j a_j(x,h) + t^k R_{k}(x,h,t), \cr
& =f(x)+ \sum_{j=1}^k  {t^j \over j!} d^jf(x)(h,\ldots,h) +  t^k R_k(x,h,t)\cr}
$$
the latter provided the integers are invertible in $\K$ -- see [BGN04,
Chapter 5] for the details. In this work we will use only the second order
Taylor expansion. 
 Thus the present summary of basic results on differential calculus is 
sufficient for our purposes. For more  information, we refer to [BGN04];
see also Appendix L, Section L.3, for examples of some important categories
of topological fields and vector spaces and their main properties, and
Appendix G, Section G.1,
 for the generalization of the preceding set-up to general
``$C^0$-concepts'', as introduced in [BGN04].

\vfill \eject

\subheadline{2. Manifolds}

\nin {\bf 2.1.} {\sl Manifolds and manifolds with atlas.}
We fix a topological $\K$-module $V$ as
 ``model space'' of our manifold.
A {\it $C^k$-manifold with atlas} is a   
topological space $M$ together with a {\it $V$-atlas}
${\cal A} = (\phi_i,U_i)_{i \in I}$. 
This means that $U_i$, $i \in I$ is a covering of $M$ by
open sets, and
 $\phi_i:M \supset U_i \to \phi_i(U_i)
\subset V$ is a {\it chart}, i.e.
a homeomorphisms of the open set $U_i \subset M$ onto an open
set $\phi_i(U_i) \subset V$, and any
 two charts $(\phi_i,U_i), (\phi_j,U_j)$ are
{\it $C^k$-compatible} in the sense that
$$
\phi_{ij}:=
\phi_i \circ \phi_j^{-1}\vert_{\phi_j(U_i \cap U_j)}:
\phi_j(U_i \cap U_j) \to \phi_i(U_i \cap U_j)
$$
and its inverse $\phi_{ji}$ are of class $C^k$.
If $M$ is $C^k$ for all $k$, we say that $M$ is a 
smooth or $C^\infty$-manifold.
We see no reason to assume that the topology of  $M$ is
Hausdorff (compare [La99, p.\ 23] for this issue).
Also, we see no reason to assume that the atlas $\cal A$ is
{\it maximal}. Of course, any atlas can be completed to a
maximal one, if one wishes to do so. (See, however, Section 2.4 below,
for a word of warning concerning maximal atlasses.)
A {\it manifold over $\K$} is a manifold together with a
maximal atlas.

Sometimes it may be useful to consider also the disjoint union
of two manifolds (possibly modelled on non-ismorphic topological
modules) as a manifold; then one should call  {\it $V$-manifolds}
the class of manifolds defined before.

Let $M,N$ be $C^k$-manifolds with atlas.
 A map $f:M \to N$ is {\it of class $C^k$} if, 
 for all choices of charts
 $(\phi,U)$ of $M$ and $(\psi,W)$ of $N$,
$$
\psi \circ f \circ \phi^{-1}: 
\phi(U \cap f^{-1}(W)) \to \psi( W)
$$
is of class $C^k$. 
From the chain rule (Lemma 1.8) it follows that
  $C^k$-manifolds with atlas form a category.
We denote by
$\Diff(M)$ or $\Diff_\K(M)$ the automorphism group of $M$ in this
category, also called the group of {\it diffeomorphisms}.

In particular, all topological $\K$-modules 
are $C^\infty$-manifolds 
and we can define {\it smooth functions on $M$} to be smooth maps 
$f:M \to \K$.  The space $C^\infty(M)$ of smooth functions on $M$
 may be reduced to the  constants, and it may also happen that
$C^\infty(U_i)$ is reduced to the constants for all chart domains $U_i$
(e.g. case of topological vector spaces that admit no non-zero 
continuous linear forms, cf.\ [BGN04, Example 8.2]).
 Therefore it is no longer possible
to define differential geometric objects via their action  on 
smooth functions.

\msk
\nin {\bf 2.2.} {\sl Direct products.}
In the category of smooth manifolds over $\K$
 one can form direct products: given 
two $C^k$-manifolds with atlas $(M,{\cal A}), (N,{\cal B})$,
endow $M \times N$ with
the product topology, 
and the charts are given by the maps $\phi_i \times \psi_j$.
These  charts are again $C^k$-compatible and define an atlas
${\cal A} \times {\cal B}$.

\msk \nin
{\bf 2.3.} {\sl Submanifolds.}
For our purposes, the following strong definition of submanifolds
will be convenient (of course, one may consider also weaker
conditions):
if $E$ is a submodule of a topological $\K$-module $V$, 
then we say that $E$ is {\it admissible} if there exists
a complementary submodule $F$ such that the 
bijection
$$
E \times F \to V, \quad (e,f) \mapsto e+f
$$
is a homeomorphism. Thus $E$ and $F$ are closed and complemented submodules
of $V$.
Now a {\it submanifold} is defined as usual to be
a subset $N \subset M$ such that there exists a subatlas
of $\cal A$ having the property that $U_i \cap N$
is either empty or corresponds to an admissible
submodule of $V$.

\msk \nin {\bf 2.4.} {\sl Locality.}
For classical real manifolds, any topological refinement of the open chart
domain cover $(U_i)_{i \in I}$ of $M$ defines a new and equivalent
atlas of $M$. In particular, a maximal atlas of $M$ contains, for any point
$x \in M$, a system of open neighborhoods of $x$ as possible chart domains.  
This fact remains true  if $\K$ is a field, because of the
following Locality Lemma ([BGN04, Lemma 4.9]): {\it If  $\K$ is a field and
if $f$ is $C^k$ on each open set $U_j$ of an open cover
$(U_j)_{j \in J}$ of the open set
$U \subset V$, $f$ is also of class $C^k$ on $U$.}
For general base rings, this may become false. In this case, one should
rather think of a maximal atlas as something 
like a ``maximal bundle atlas'' of a fiber bundle, with fibers corresponding
to ideals of non-invertible elements of $\K$. 
As long as one minds this remark, the theory of manifolds over rings
is not really different from the theory of manifolds over fields.

\vfill \eject

\subheadline
{3. Tangent bundle and general fiber bundles}

\msk
\nin {\bf 3.1.} {\sl An equivalence relation describing $M$.}
Assume $M$ is a manifold with atlas modelled on $V$.
A point $p \in M$ is described in the form $p=\phi_i^{-1}(x)$
with $x \in \phi_i(U_i)$ and $i \in I$. In a different chart
it is given by $p=\phi_j^{-1}(y)$. In other words, 
$M$ is the set of equivalence classes $S/\sim$, where
$$
S:=\{ (i,x)  | \, x \in \phi_i(U_i) \} \subset I \times V,
$$
and $(i,x) \sim (j,y)$ iff $\phi_i^{-1}(x)=\phi_j^{-1}(y)$
iff $\phi_{ji}(x)=y$.
We write $p=[i,x] \in M = S/\sim$.

\msk \nin {\bf 3.2.} {\sl An equivalence relation describing $TM$.}
Next we define an equivalence relation on the set
$$
TS:= S \times V\subset I \times V \times V
$$
by:
$$
\eqalign{
(i,x,v) \sim (j,y,w) \quad  :& \Leftrightarrow  \quad
\phi_j \circ \phi_i^{-1}(x)=y, \,
d(\phi_j \circ \phi_i^{-1})(x)v=w  \cr
&  \Leftrightarrow  \quad
\phi_{ji}(x)=y, \, d\phi_{ji}(x)v=w. \cr}
$$
Since $\phi_{ii}=\id$, $\phi_{ij}\phi_{jk}=\phi_{ik}$,
the chain rule implies that this is an equivalence relation.
Again we denote equivalence classes by $[i,x,v]$, and we let
$$
TM:=TS/\sim.
$$
If $[i,x,v]=[j,y,w]$, then $[i,x]=[j,y]$, and hence the map
$$
\pi: TM \to M, \quad [i,x,v] \mapsto [i,x]
$$
is well-defined. For $p=[i,x] \in M$, we let
$$
T_p M := \pi^{-1}(p) =\{ [i,x,v] \in TM  \vert \, v \in V \}.
$$
The map
$$
T_x \phi_i^{-1}: V \to T_p M, \quad v \mapsto [i,x,v]
$$
is a bijection (it is surjective by definition and injective
 since the differentials are
bijections), and we can use it to define the structure of 
a topological $\K$-module 
on $T_pM$ which actually does not depend on $(i,x)$
because $T_y \phi_j^{-1} = T_x \phi_i^{-1} \circ d\phi_{ij}(y)$.
The $\K$-module $T_pM$ is called the {\it tangent space of $M$
at $x$}.

We define an atlas $T {\cal A}:=(T\phi_i)_{i \in I}$ on $T M$ by:
$$
TU_i:=\pi^{-1}(U_i), \quad  \quad T \phi_i:
TU_i \to V \times V, \,\, [i,x,v] \mapsto (x,v)
$$
Change of charts is now given by 
$$
T\phi_{ij}:  (x,v) \mapsto
(\phi_{ij}(x),d\phi_{ij}(x)v)
$$  
which is $C^{k-1}$ if $\phi_{ij}$ is 
$C^k$. Thus $(TM,T{\cal A})$ is a manifold with atlas. 

If $f:M \to N$ is $C^k$ we define its {\it tangent map} by
$$
Tf: TM \to TN, \quad [i,x,v] \mapsto [j,f_{ij}(x),d f_{ij}(x)v]
$$
where $f_{ij}=\psi_j \circ f \circ \phi_i^{-1}$ (supposed to be
defined on a non-empty open set). In other words,
$$
Tf=(T\psi_j)^{-1} \circ (f_{ij},df_{ij}) \circ T\phi_i
=(T\psi_j)^{-1} \circ Tf_{ij} \circ T\phi_i 
$$
with $Tf_{ij}$ as defined in (1.3).
This is well-defined, linear in fibers and $C^{k-1}$.
Clearly the functorial rules hold, i.e. we have defined
a covariant functor $T$ from the category of smooth manifolds over
$\K$ into itself. 

Using the product atlas from  2.2 on the direct product $M \times N$,
we see that there is a natural isomorphism
$$
T(M \times N) \cong TM \times TN
$$
which is directly constructed by using the equivalence relation 
defining $TM$.

If $f:M \to \K$ is a smooth function, then $T_x f:T_x M \to T_{f(x)} \K = \K$
gives rise to a function $TM \to \K$, linear in fibers, 
which we denote by $df$.
The product rule (Lemma 1.6 (i)) implies that $d(fg)=f\,dg+g\,df$.

\msk\nin 
{\bf 3.3.} {\sl General fiber bundles.}  
General fiber bundles over $M$ are defined following the same pattern as
above for the tangent bundle: assume $M$ is modelled on $V$ and let
$N$ some manifold modelled on a topological $\K$-module
$W$. Assume that, 
for all triples $(i,j,x) \in I \times I \times V$ such that
$\phi_i^{-1}(x) \in  \phi_j^{-1}(U_j)$, an element
$$
g_{ji}(x) \in {\rm Diff}(N)
$$
of the diffeomorphism group of $N$ 
is given such that the cocycle relations
$$
g_{ik}(x) = g_{ij}(\phi_{jk}(x)) \circ g_{jk}(x)  , \quad
g_{ii}(x)=\id_N
$$
are satisfied and such that 
$$
(x,w) \mapsto g_{ij}(x,w):= g_{ij}(x)w
$$
 is smooth wherever
defined. Then we define an equivalence relation on 
$S \times N$ by
$$
(i,x,v) \sim (j,y,w) \quad : \Leftrightarrow  \quad
\phi_{ji}(x)=y, \, \, g_{ji}(x)v=w.
$$
By the cocycle relations, this is indeed an equivalence relation,
and by the smoothness assumption, $F:=S \times N/\sim$ can be
turned into a manifold modelled on $V \times W$ and locally
isomorphic to $U_i \times N$ and such that
the projection $p:F \to M$, $[i,x,w] \mapsto [i,x]$ is a
well-defined smooth map whose fibers are all diffeomorphic  to $N$.
When we speak of a ``chart'' of a bundle, we always mean a bundle
chart, and we never consider maximal atlasses on bundles. 

{\it Homomorphisms} or {\it bundle maps} are pairs
$\tilde f:F \to F'$, $f:M \to M'$ of smooth maps such that
$p' \circ \tilde f = f \circ p$. Using charts, it is seen that
fibers of bundles are submanifolds, and
bundle maps induce by restriction smooth maps on fibers.

\msk \nin {\bf 3.4.} {\sl Vector bundles.}
If $N$ carries an additional structure ($\K$-module, affine space,
projective space...) and the $g_{ij}(x)$ respect this structure,
then each fiber also carries this structure.
In particular, if $N$ is a $\K$-module and the transition maps
respect this structure, then $F$ is called a
{\it vector bundle}. In this case, the bundle atlas is of the
form $\tilde A =(p^{-1}(U_i), \tilde \phi_i)_{i \in I}$ with
$$
\tilde \phi_i: p^{-1}(U_i) \to V \times W, \quad
[i,x,w] \mapsto (x,w).
$$
Change of charts is given by the transition functions
$$
\tilde g_{ij}(x,w)=(\phi_{ij}(x),g_{ij}(x)w).
$$
Homomorphisms are required to respect the extra structure on
the fibers. Thus homomorphisms of vector bundles are such
that restriction to fibers induces $\K$-linear maps.

If $F$ is a vector bundle over $M$, then there is a well defined map
$$
z: M \to F, \quad [i,x] \mapsto [i,x,0],
$$
called the {\it zero section}, and if $(\tilde f,f)$ is a homomorphism
of vector bundles, then $f$ can be recovered from $\tilde f$ via
$f = p \circ \tilde f \circ z$.

\msk
\nin
{\bf 3.5.} {\sl Direct sum of vector bundles.}
If $F$ and $E$ are fiber bundles over $M$, then, by definition,
$F \times_M E$ is the fiber bundle over $M$ whose fiber over
$x$ is $F_x \times E_x$. The transition functions are given by
$g_{ij}(x) \times f_{ij}(x)$. If $E$ and $F$ are
in addition vector bundles, then $E \times_M F$ is
also a vector bundle, denoted by $E \oplus_M F$,
with transition functions, in matrix form,
$$
g_{ij}(x) \oplus f_{ij}(x) = 
\pmatrix{g_{ij}(x)&0 \cr 0 & f_{ij}(x) \cr}.
$$
We do not define tensor products, dual or hom-bundles
of vector bundles
in our general context -- see Appendix L, Sections L1 and
L2 for explanatory comments.

\vfill \eject

\subheadline
{4. The Lie bracket of vector fields}

\msk \nin
{\bf 4.1.} {\sl Vector fields and derivations.}
A {\it section} of a vector bundle $F$ over $M$ is a smooth map
$X:M \to F$ such that $p \circ X = \id_M$. 
If $F$ is a vector bundle, then the sections 
of $F$ form a module, denoted by $\Gamma^\infty(F)$ or simply by
$\Gamma(F)$, over the ring $C^\infty(M)$. 
Sections of $TM$ are also called {\it vector fields},
 and we also use the classical notation $\X(M)$ for
$\Gamma(TM)$. 
In a chart $(U_i,\phi_i)$, vector fields can be identified with
smooth maps $X_i:V \supset \phi_i(U_i) \to V$, given by 
$$
X_i:=\pr_2 \circ T\phi_i \circ X \circ \phi_i^{-1}:
\phi_i^{-1}(U_i) \to U_i \to TU_i \cong U_i \times V \to V.
$$
Similarly, sections of an arbitrary vector bundle are locally
represented by smooth maps 
$$
X_i: V \supset \phi_i^{-1}(U_i) \to W.
$$
If the chart $(U_i,\phi_i)$ is fixed, for brevity of notation
 we will often suppress
the index $i$ and write the chart representation of $X$ in the form
$$
U \to U \times W, \quad x \mapsto x + X(x) \quad {\rm or} \quad
(x,X(x)).
$$

For a vector field $X:M \to TM$ and a smooth function $f:M \to \K$,
 recall that the differential
$df:TM \to \K$ is smooth and hence we can define a smooth
function $L_X f$ by $L_X f:=df \circ X$. Then we have the Leibniz rule:
$$
L_X(fg)=d(fg) \circ X=(fdg+gdf) \circ X=g L_Xf+f L_Xg.
$$
Thus the map $X \mapsto L_X$ is a $\K$-linear map into the space of
derivations of $C^\infty(M)$.
However, it will in general neither be injective nor surjective, not
even when restricted to suitable open subsets
 (cf.\ [BGN04, Examples 8.2 and 8.3]).
 Therefore it cannot be used to define the Lie algebra
structure on $\X(M)$.
In the following theorem, we define the Lie bracket by using
its expression in a chart;
an intrinsic definition of the Lie bracket needs
a closer look at the second tangent bundle $TTM$ and
is postponed to Chapter 14.

\Theorem 4.2.
There is a unique structure of a Lie algebra over $\K$ on $\X(M)$
such that for all $X,Y \in \X(M)$ and $(i,x) \in S$,
$$
[X,Y]_i(x) = dX_i(x) Y_i(x) - dY_i(x) X_i(x).
$$

\Proof.
Uniqueness is clear. Let us show that, on the intersection of two
chart domains, the bracket $[X,Y]$ is independent of the choice of
chart: assume $(i,x) \sim (j,y)$, i.e. $y=\phi_j \phi_i^{-1} x =
\phi_{ji}(x)$; then
$X_j(y)=d \phi_{ji}(x) X_i(x)$ or
$$
X_j \circ \phi_{ji} = T \phi_{ji} \circ X_i.
\eqno (4.1)
$$
We have to show that $[X,Y]$ has
the same transformation property under changes of charts. 
For simplicity, we let $\phi:=\phi_{ji}$ and calculate:
$$
\eqalign{
d Y_j(\phi(x)) X_j(\phi(x)) & =
d(T \phi \circ Y_i)(x) T\phi \circ X_i(x) \cr
&=d^2\phi(x)(X_i(x),Y_i(x)) + d\phi(x) dY_i(x) X_i(x).
\cr}
$$
We exchange the r\^oles of $X$ and $Y$ and take the difference of
the two equations thus obtained: we get, using Schwarz' lemma (Lemma 1.9),
$$
d Y_j(\phi(x)) X_j(\phi(x)) -d X_j(\phi(x)) Y_j(\phi(x))  =
d\phi(x) (d Y_i(x) X_i(x) - dX_i(x) Y_i(x))
$$
which had to be shown.
Summing up, the bracket operation $\X(M) \times \X(M) \to \X(M)$
is well-defined, and it clearly is $\K$-bilinear and antisymmetric
in $X$ and $Y$. For the case that $2$ is not invertible in $\K$,
note that we clearly also have the identity $[X,X]=0$.

All that remains to be proved is the Jacobi identity. 
This is done by a direct computation which involves only 
the composition rule and Schwarz' lemma: define a (chart dependent)  
``product" of $X_i$ and $Y_i$ by
$$
(X_i \cdot Y_i)(x):=dY_i(x) X_i(x).
$$
(Later this will be interpreted
 as $\nabla_{X_i}Y_i$, the canonical flat connection
of the chart $U_i$, see Chapter 10.)
Then, by a direct calculation, one shows that this product
is a {\it left symmetric} or {\it Vinberg algebra} (cf., e.g., [Koe69]), i.e.
it satisfies the identity
$$
X_i \cdot (Y_i \cdot Z_i)-(X_i \cdot Y_i)\cdot Z_i =
Y_i \cdot (X_i \cdot Z_i)-(Y_i \cdot X_i) \cdot Z_i.
$$
But it is immediately checked that for every left symmetric algebra,
the commutator $[X_i,Y_i]=X_i \cdot Y_i - Y_i \cdot X_i$
satisfies the Jacobi identity.
\qed

The Lie bracket is natural in the following sense:
assume $\phi:M \to N$ is a smooth map and $X \in \X(M)$, $Y \in \X(N)$.
We say that the pair $(X,Y)$ is {\it $\phi$-related} if
$$
Y \circ \phi = T\phi \circ X.
$$

\Lemma 4.3. If $(X,Y)$ and $(X',Y')$ are $\phi$-related, then 
so is $([X,X'],[Y,Y'])$. In particular, the diffeomorphism group
of $M$ acts by automorphisms of the Lie algebra $\X(M)$.

\Proof. This is the same calculation as the one given in the preceding
proof after Eqn. (4.1). \qed

Moreover, from the definitions it follows easily that
$$
\X(M) \to \Der(C^\infty(M)), \quad X \mapsto - L_X
$$
is a homomorphism of Lie algebras. (The sign in the definition of
the Lie bracket is a matter of convention.)

\msk
\nin {\bf 4.4.} {\sl Infinitesimal automorphisms.}
If $X$ is a vector field on $M$, then an {\it integral curve}
is a smooth map $\phi:\K \supset I \to M$ such that
$T_t \phi \cdot \1 = X(\phi(t))$.
In a chart, this is equivalent to the usual differential
equation $\phi'(t)=X(\phi(t))$.
In our general set-up, we have no existence or uniqueness
statements for solutions of ordinary differential equations,
and hence  we will never work with integral curves
or flows.

However, a very rough infinitesimal version of the flow of a vector field
can be defined (for more refined versions see Chapter 29):
we say that a diffeomorphism
$f:TM \to TM$ is an {\it infinitesimal automorphism} if
$f$ preserves fibers,  and  in each fiber $f$ acts by
translations. Clearly, this defines a group, denoted by
$$
{\rm InfAut}(M):=\bigl\{ f \in \Diff(TM) |   \,
  \forall p \in M:
\exists v_p \in T_p M : \forall v \in T_p M: \, f(v)=v+v_p   \bigr\}.
$$
Clearly,  $X(p):=v_p=f(0_p)$ defines a
map $X:M \to TM$ that can also be written $X=f \circ z$, where $z:M \to TM$
is the zero section, whence is a smooth section of $TM$, i.e., it
 is a vector field on $M$. Conversely,
if $X:M \to TM$ is a vector field on $M$, we define a map
$$
\tilde X:TM \to TM, \quad v \mapsto v + X(\pi(v))).
$$
In each fiber, $\tilde X$ acts by translations, and using
charts, we see that $\tilde X$ is smooth.
The inverse of $\tilde X$ is $\tilde{-X}$; more generally,
we have $\tilde{X+Y}=\tilde X \circ \tilde Y$, and
hence $\tilde X$ is an infinitesimal automorphism. Summing up, we have
defined a group isomorphism $(\X(M),+) \cong {\rm InfAut} (M)$.
Note that the same notions could be defined by replacing
$TM$ by an arbitrary vector bundle $F$ over $M$.
In Chapter 28, ``higher order versions'' of the group ${\rm InfAut}(M)$
will be defined.

\msk \nin {\bf 4.5.} {\sl Tensor fields, forms.}
Let $E$ and $F$ be vector bundles over $M$.
An {\it $E$-valued $k$-multilinear form on $F$}
is a smooth map $A: F \times_M \ldots \times_M F \to E$ ($k$ factors)
 such
that $p_E \circ  A = p_{\times^k_M F}$, i.e., $A$ 
maps fibers over $x$ to fibers over $x$, and 
$A_x:F_x \times \ldots \times F_x \to E_x$ is $\K$-linear in all $k$ variables.
If $F=TM$, then we speak also of {\it $E$-valued
forms}; if $F=TM$ and $E= M \times \K$ is the trivial
line bundle, then $A$ is called a {\it tensor field
of type $(k,0)$}, and if $F=E=TM$, then $A$ is
called {\it tensor field of type $(k,1)$}.
{\it Differential forms of degree $k$} are skew-symmetric tensor
fields of type $(k,0)$.
(We don't define tensor fields of type $(r,s)$ with $s>1$ since  
we have not defined the tensor product bundles; cf. Appendix L2.)
Note that we do not interprete forms as sections
of a vector bundle over $M$; for instance, one-forms are defined as maps
$\omega:TM \to \K$,
but $T^*M$ is not defined as a bundle (cf. Appendix L1).

\msk \nin
{\bf 4.6.} {\sl Almost (para-) complex and dual structures.}
A tensor field $J$ of type $(1,1)$, i.e., a smooth map $J:TM \to TM$
acting linearly on all fibers, such that
$J^2=-\1$ will be called an {\it almost complex structure}
(even if the base field is not real), and if $J^2 = \1$,
$J$ will be called a {\it polarization}
(or a {\it para-complex structure} in case $2$ is invertible
in $\K$ and both eigenspaces of $J_x$ are isomorphic as
submodules of $T_xM$ for all $x \in M$).
If $J^2 = 0$ and $\ker(J_x)=\im(J_x)$ for all $x \in M$,
then $J$ will be called an {\it almost dual structure}.

\msk \nin
{\bf 4.7.} {\sl The canonical almost dual structure on the
tangent bundle.}
The tangent bundle $TM$ of a manifold $M$
carries a natural almost dual structure.  This can be seen as follows:
recall that the transition maps of the bundle atlas of
$TM$ are given by
$T\phi_{ij}(x;v)=(\phi_{ij}(x),d\phi_{ij}(x)v)$.
Therefore, using the rule on partial derivatives, 
the transition maps of the bundle atlas of $TTM=T(TM)$
are given by
$$
\eqalign{
TT\phi_{ij}(x,v;x',v') & =
(T\phi_{ij}(x,v),d(T \phi_{ij})(x,v)(x',v')) \cr
& =(\phi_{ij}(x),d\phi_{ij}(x)v,d\phi_{ij}(x)x',
d\phi_{ij}(x)v' + d^2 \phi_{ij}(x)(v,x')) \cr
& =  (\phi_{ij}(x),d\phi_{ij}(x)v,
\pmatrix{d \phi_{ij}(x) & 0 \cr d^2\phi_{ij}(x)(v,\cdot) & d \phi_{ij}(x) \cr}
\pmatrix{x' \cr v' \cr} ). \cr}
$$
The matrix in the last line describes the differential
$d(T \phi_{ij})(x,v)$. It is clear that this matrix commutes with
the matrix
$$
Q := \pmatrix{0 & 0 \cr \1_V & 0 \cr}.
$$
Therefore, if $p=[i,x] \in M$ and $u =[i,x,v] \in T_pM$,
and we identify $T_u(TM)$ with $V \times V$ via the chart
$TT\phi_i$, then the endomorphism 
$\epsilon_u$ of $T_u(TM)$ defined by
$Q$ does not depend on the chart. Hence
$$
\epsilon :=(\epsilon_u)_{u \in TM}
$$
defines a tensor field of type $(1,1)$ on $TM$ such that
$\epsilon^2=0$ and the kernel of $\epsilon_u$ is exactly
$T_p M \cong T_u(T_pM) \subset T_u(TM)$.
This tensor field is natural in the sense that it is invariant under
{\it all} diffeomorphisms of the kind $Tf$ where $f$ is a diffeomorphism 
of $M$: this is proved by the above calculation, with $\phi_{ij}$
replaced by $f$.
(In [Bes78, p.\ 20/21] this tensor field, in a slightly different
presentation, is called the ``vertical endomorphism''.)

The tensor field $\epsilon$ is {\it integrable} in the following
sense: the distribution of subspaces given by $\ker \epsilon$
admits {\it integral submanifolds} which are simply the fibers
of the bundle $TM$. Comparing with the situation of an almost 
complex structure (replace the condition $\epsilon^2=0$ by
$\epsilon^2 = - \1$), one might conjecture that in analogy with
the theorem of Newlander and Nirenberg from the complex case,
also in this case it may be true that indeed $TM$ is already
a manifold defined over the ring $\K \oplus \epsilon \K$
with relation $\epsilon^2 = 0$, i.e.\ over the dual
numbers. In Chapter 6 we will show that this is indeed true.

\msk \nin {\bf 4.8.} {\sl Remarks on differential operators of degree zero.}
In the real finite-dimensional case, tensor fields can be
interpreted as  (multi-)differential operators of degree zero.
This is in general no longer possible in our set-up.
By definition, a {\it (multi-) differential operator of degree
zero} between vector bundles $E,F$ over $M$ is given by a
collection of maps from sections to sections,
for all open sets $U \subset M$, or at least for all chart domains $U$,
$$
\tilde A_U: \Gamma(F|_U) \times \ldots \times \Gamma(F|_U)
  \to \Gamma(F|_U), \quad
(\xi_1,\ldots,\xi_k) \mapsto \tilde A_U(\xi_1,\ldots,\xi_k)
\eqno (4.2)
$$
($k$ factors) such that

\ssk
\item{(1)} $\tilde A_U$ is ${\cal C}^\infty(U)$-multilinear,
\item{(2)} for $x \in U$, the value $\tilde A_U(\xi_1,\ldots,\xi_k)(x)$
depends only on the values of the sections $\xi_1,\ldots,\xi_k$ at the
point $x$.
\ssk

\nin
If $A: F \times_M \ldots \times_M F \to E$
 is a tensor field as in Section 4.5, then
clearly we obtain a multidifferential operator of degree zero by letting
$$
(\tilde A_U(\xi_1,\ldots,\xi_k))(x):=
A_x(\xi_1(x),\ldots,\xi_k(x)).
$$
Conversely, given a 
collection of mappings $\tilde A_U$
having Properties (1) and (2),  we can define a map
 $A: F \times_M \ldots \times_M F \to E$, multilinear in fibers and depending
smoothly on the base $M$,
in the following way: let $x \in M$, choose
a chart $U$ around $x$, extend  elements 
$v_1,\ldots,v_k\in F_x$ to constant sections
$\tilde v_1,\ldots,\tilde v_k$ of $F$ over $U$
(i.e. $(\tilde v_i)(x)=v_i$) and let
$$
A_x(v_1,\ldots,v_k):=(\tilde A_U(\tilde v_1,\ldots,\tilde v_k))(x).
\eqno (4.3)
$$
(This is well-defined: independence from 
$U$ and from the extension of elements in the fiber to local
sections follows from (2).)
However, we cannot guarantee 
smooth dependence on $v_1,\ldots,v_k$; this has to be
checked case by case in application situations.
For instance,  a one-form $A:TM \to \K$  corresponds to 
$\omega:= \tilde A_M:\X(M) \to C^\infty(M)$, and for vector fields 
$X,Y$, we let
$$
(L_X \omega)(Y):=L_X (\omega(Y)) -\omega([X,Y]).
\eqno (4.4)
$$
Then $L_X \omega$ satisfies Properties (1) and (2);
smooth dependence on $v=Y(x)$ is easily checked in a chart,
and hence $L_X \omega$ defines, via (4.3), a tensor field which we denote
by $L_X A$.
Similarly we could define the Lie derivative of an arbitrary
$(k,0)$- or $(k,1)$-form, and using this, 
the exterior derivative of a differential form
could be defined in the usual way. We will, however,  give
a more intrinsic definition of the exterior derivative later on (Chapters
13 and 22). 
Finally, let us remark that in some situations, Condition (2) is
 automatically satisfied:
it is automatic in case  $\K$ is a field and the fibers of $F$ are
finite-dimensional (essentially, the proof of [La99, Lemma VIII.2.3] applies
to this situation).
Moreover, if cutoff functions
in the sense of  [La99] exist, then it suffices to work with
$\tilde A_M$ instead of the whole collection of the $\tilde A_U$'s.
See Chapter 21 for more general remarks on differential operators.

\vfill \eject 

\subheadline
{5. Lie groups and symmetric spaces: basic facts}

\nin {\bf 5.1.} {\sl Manifolds with multiplication.}
A {\it product} or {\it multiplication map} 
on a manifold $M$ 
is a smooth binary map $m:M \times M \to M$, and 
{\it homomorphisms of manifolds with multiplication}
are smooth maps that are compatible with the respective 
multiplication maps.
{\it Left and right multiplication operators}, defined by
$l_x(y)=m(x,y)=r_y(x)$, are partial maps of $m$ and
hence smooth self maps of $M$. 
Applying the tangent functor to this situation, we see that
$(TM,Tm)$ is again a manifold with multiplication, and tangent
maps of homomorphisms are homomorphisms of the respective tangent
bundles.
The tangent map $Tm$ is given by the  formula
$$
T_{(x,y)}m(v,w) = T_{(x,y)}m((v,0_y)+(0_x,w))= T_x (r_y)v + T_y (l_x)w.
\eqno (5.1)
$$
Formula (5.1) is nothing but the rule on partial derivatives (Lemma 1.5 (iv)),
 written in the language of manifolds.
In particular, (5.1) shows that the canonical projection and the zero
section,
$$
\pi:TM \to M, \quad v \to \pi(v), \quad \quad
z:M \to TM, \quad p \mapsto 0_p
\eqno (5.2)
$$
are homomorphisms of manifolds with multiplication. We will always 
identify $M$ with the subspace $z(M)$ of $TM$. 
Then (5.1) implies that the operator of left multiplication by $p=0_p$
in $TM$ is nothing but $T(l_p):TM \to TM$, and similarly for
right multiplications.

\msk
\nin
{\bf 5.2.} {\sl Lie groups.}
A {\it Lie group over $\K$}
is a smooth $\K$-manifold $G$ carrying a group structure such that
 multiplication  $m:G \times G \to G$ and  inversion
$i:G \to G$ are smooth. Homomorphisms of Lie groups are smooth 
group homomorphisms. 
Clearly, Lie groups and their homomorphisms form a category in which
direct products exist.

Applying the tangent functor to the defining identities of the group
structure $(G,m,i,e)$, it is immediately seen that then $(TG,Tm,Ti,
0_{T_eG})$ is
again a Lie group such that $\pi:TG \to G$ becomes a homomorphism
of Lie groups and such that the zero section $z:G \to TG$ also
is a homomorphism of Lie groups.

\msk \nin
{\bf 5.3.} {\sl The Lie algebra of a Lie group.}
A vector field $X \in \X(G)$ is called {\it left invariant} if, for
all $g \in G$,
$X \circ l_g = T l_g \circ X$, and similarly we define {\it right invariant
vector fields}. If $X$ is right invariant, then
 $X(g)=X(r_g(e))=T_e r_g \cdot X(e)$; thus $X$ is uniquely
determined by the value $X(e)$, and thus the map
$$
\X(G)^{r_G} \to T_e G, \quad X \mapsto X(e)
\eqno (5.3)
$$
from the space of right invariant vector fields into $T_e G$ is injective.
It is also surjective: if $v \in T_e G$, then
left multiplication with $v$ in $TG$, $l_v:T G \to TG$, preserves
fibers (because $\pi$ is a homomorphism) and hence defines a vector field 
$$
v^R:=l_v \circ z: G \to TG, \quad 
g \mapsto l_v (0_g)= v \cdot 0_g = Tm(v,0_g)=T_e r_g(v)
$$ 
which is right  invariant since right
multiplications commute with left multiplications. 
Now, the space $\X(G)^{r_G}$ is a Lie subalgebra of $\X(M)$; 
this follows immediately from Lemma 4.3
because $X$ is right invariant iff the pair $(X,X)$ is $r_g$-related
for all $g \in G$. The space $\g:=T_eG$ with the Lie bracket defined by
$$
[v,w]:=[v^R, w^R]_e
$$
is called the {\it  Lie algebra of $G$}.

\Theorem 5.4. 
\item{(i)} The Lie bracket $\g \times \g \to \g$ is of class $C^0$.
\item{(ii)} For every homomorphism $f:G \to H$, the tangent map
$\dot f:=T_e f:\g \to \h$ is a homomorphism of Lie algebras.

\Proof. (i)
Pick a chart $\phi:U \to V$ of $G$ such that $\phi(e)=0$.
Since
$w^R(x)=Tm(w,x)$ depends smoothly on $(x,w)$, it is represented in
the chart by a smooth map (which again will be denoted by $w^R(x)$).
But this implies that
$[v^R,w^R](x)=d(v^R)(x) w^R(x)-d(w^R)(x)v^R(x)$
depends smoothly on $v,w$ and $x$ and hence
$[v,w]$ depends smoothly on $v,w$.

(ii)
First one checks that the pair of vector fields 
$(v^R,(\dot f v)^R)$ is $f$-related, and then
one applies Lemma 4.3 in order to conclude that
$\dot f[v,w]=[\dot f v, \dot f w]$.
\qed

The functor from Lie groups over $\K$ to $C^0$-Lie algebras
over $\K$ will be called the {\it Lie functor (for $\K$-Lie groups
and $\K$-Lie algebras)}.
We say that a $C^0$-Lie algebra over $\K$ is {\it integrable}
if there is a Lie group over $\K$ such that $\g={\rm Lie}(G)$.
The following is  the
{\it integrability problem for $\K$-Lie algebras}:

\msk
\item{IP1.} Which $\K$-Lie algebras are integrable ?
\item{IP2.} If a $\K$-Lie algebra is integrable, how can we describe
all equivalence classes of Lie groups belonging to this Lie algebra ?

\msk
\nin These are  difficult
problems which fall out of the scope of the present work.
Let us just mention that
all full matrix groups $\Gl(n,\K)$ are Lie groups (with atlas given
by the natural chart), as well as all
orthogonal groups $\OO(b,\K^n)$ corresponding to non-degenerate
forms $b:\K^n \times \K^n \to \K$ -- in these cases the manifold
structure can be defined  by ``Cayley's rational
chart'', cf.\ [Be00] and [BeNe05]. For $\SL(n,\K)$ one can define
an atlas via the Bruhat decomposition; thus all ``classical groups
over $\K$'' are Lie groups, and their Lie algebras are calculated in
the usual way, namely, they are subalgebras of some matrix algebra
$\gl(V)$ with the usual Lie bracket $[X,Y]=XY - YX$, and $V \cong \K^m$.
In fact, for $X \in \gl(V)$, the value of the right invariant vector field
$X^R$ at $g \in \Gl(V)$ is $X^R(g)=X g$ (product in $\End(V)$), and hence
the Lie bracket can be calculated in the natural chart $\End(V)$:
$$
[X,Y]= [X^R,Y^R]_e = d(X^R)(e) \cdot Y^R(e) - d(Y^R)(e) \cdot X^R(e)=
XY-YX.
$$
(This explains our choice of defining the Lie algebra via {\it right}
invariant vector fields; the other choice would have given  the opposite
Lie algebra structure on $\gl(V)$. Since we do not invoke any contravariant
constructions, we can avoid the well-known ``sign dilemma'', cf.\ 
[Be00, p.\ 32/33], [Hel78, p.\ 122].)
If $V$ is not isomorphic to some $\K^m$, then in general there is no
Lie group structure on $\Gl(V)$ -- but a good substitute of general
linear groups is given by unit groups of ``continuous inverse algebras''
(cf.\ Section 25.6). See [Gl04e] for these and  more general 
constructions  of Lie groups over topological fields.

\msk
\nin {\bf 5.5.} {\sl Symmetric spaces.}
A {\it reflection space (over $\K$)}
is a smooth manifold with a multiplication map $m:M \times M \to M$
such that, for all $x,y,z \in M$, and writing $\sigma_x(y)=m(x,y)$
instead of $l_x(y)$,

\ssk
\item{(M1)} $m(x,x)=x$
\item{(M2)} $m(x,m(x,y))=y$, i.e. $\sigma_x^2 = \id_M$, 
\item{(M3)} $m(x,m(y,z))=m(m(x,y),m(x,z))$, i.e. $\sigma_x \in \Aut(M,m)$.

\ssk
\nin
(Reflection spaces -- ``Spiegelungsr\"aume'' in German -- have been 
introduced by O. Loos in [Lo67] in the finite dimensional real case.)
The left multiplication operator $\sigma_x$ 
is, by (M1)--(M3), an automorphism of 
order two fixing $x$; it is called the {\it symmetry around $x$}.
The ``trivial reflection space''
 $\sigma_x = \id_M$ for all $x$ is not excluded by the axioms
(M1)--(M3). We say that $(M,m)$ is a {\it symmetric space (over $\K$)}
if $(M,m)$ is a reflection space such that $2$ is invertible in $\K$
and the property

\msk
\item{(M4)} for all $x \in M$, $T_x(\sigma_x)=-\id_{T_xM}$
\msk

\nin holds. The assumption that $2$ is invertible in $\K$
guarantess that $0$ is the only fixed point of the differential
$T_x (\sigma_x):T_x M \to T_x M$, and without this assumption (M4)
would be useless. 
In the finite dimensional case over $\K=\R$, or, more generally, in
any context where we have an implicit function theorem at out disposition,
 (M4) implies that $x$ is an {\it isolated fixed point}
of $\sigma_x$ and hence our definition contains the one from [Lo69]
as a special case (see [Ne02a, Lemma 3.2] for the case of Banach symmetric
spaces). 
The group $G(M)$ generated by all $\sigma_x \sigma_y$
is a (normal) 
subgroup of $\Aut(M,m)$, called the {\it group of displacements}.
A distinguished point $o\in M$ is called a {\it base point}.
With respect to a base point, one defines the {\it quadratic representation}
$$
Q:=Q_o:M \to G(M), \quad x \mapsto Q(x):=\sigma_x \sigma_o
\eqno (5.4)
$$
and the {\it powers} (for $n\in \Z$)
$$
x^{-1}:=\sigma_o(x), \quad x^{2n}:=Q(x)^n.o, \quad x^{2n+1}:=Q(x)^n.x.
\eqno (5.5)
$$
By a straightforward calculation, one proves then the {\it fundamental
formula}
$$
Q(Q(x)y)=Q(x)Q(y)Q(x)
\eqno (5.6)
$$
and the {\it power associativity rules} (cf.\ [Lo69] or [Be00, Lemma I.5.6])
$$
m(x^n,x^m)=x^{2n-m},\quad (x^m)^n = x^{mn}.
\eqno (5.7)
$$
A symmetric space is called {\it abelian} or {\it commutative} if the group
$G(M)$ is commutative. For instance, any topological $\K$-module with the
product
$$
m(u,v)=u-v+u=2u-v
\eqno (5.8)
$$
is a commutative symmetric space. See Section 26.1 for some more basic
facts on commutative symmetric spaces.

\Proposition 5.6. Assume $(M,m)$ is a symmetric space over $\K$.
\item{(i)}
The tangent bundle $(TM,Tm)$ of a reflection
space is again a reflection space.
\item{(ii)}
The tangent bundle $(TM,Tm)$ of a symmetric
space is again a symmetric space.

\Proof. (i)
We express the identities (M1)--(M3) by commutative diagrams
to which we apply the tangent functor $T$.
Since $T$ commutes with direct products, we get the same diagrams
and hence the laws (M1)--(M3) for $Tm$ (cf.\ [Lo69] for the explicit
form of the diagrams).

(ii) We have to prove that (M4) holds for $(TM,Tm)$.
First of all, note that the fibers of $\pi:TM \to M$ (i.e. the
tangent spaces) are stable under $Tm$ because $\pi$ is a homomorphism.
We claim that for $v,w \in T_pM$ the explicit formula
$Tm(v,w)=2v-w$
holds (i.e. the structure induced on tangent spaces is the canonical
``flat'' symmetric structure (5.8) of an affine space). In fact, from (M1)
for $Tm$ we get $v=Tm(v,v)=T_p(\sigma_p)v+T_p(r_p)v=-v+T_p(r_p)v$, 
whence $T_p(r_p)v=2v$ and
$$
Tm(v,w)=T_p(\sigma_p)w+T_p(r_p)v = 2v-w.
$$

Now fix $o \in M$ and $v \in T_o M$. We choose $0_o$ as base point
in $TM$.
Then $Q(v)=\sigma_v \sigma_{0_o}$ is, by (M3), an automorphism of  
$(TM,Tm)$ such that $Q(v) 0_o = \sigma_v(0_o)=2v$.
But multiplication by $1 \over 2$,
$$
{1 \over 2}:TM \to TM, \quad u \mapsto {1 \over 2} u
$$
also is an automorphism of $(TM,Tm)$, as shows Formula (5.1).
Therefore the automorphism group of $TM$ acts transitively on fibers,
and after conjugation of $\sigma_v$ with $({1 \over 2} Q(v))^{-1}$ we
may assume that $v=0_o$.
But in this case the proof of our claim is easy:
we have $\sigma_{0_o}=T \sigma_o$, and 
since $T_o\sigma_o = - \id_{T_oM}$, 
the canonical identification $T_{0_o}(TM) \cong T_o M \oplus T_o M$
yields
$T_{0_o}(\sigma_{0_o}) = (- \id_{T_o M}) \times (- \id_{T_oM})
= - \id_{T_{0_o} TM}$, whence (M4).
\qed

The alert reader may have noticed that the preceding proof of (M4) already
contains the construction of a canonical connection on $TM$:
in fact, the argument of the proof shows that the ``vertical
space" $V_v:=T_v (T_p M) \subset T_v(TM)$ has a canonical complement
$H_v = {1 \over 2} Q(v) H_0$, where $H_0$ is one of the factors of
the canonical decomposition $T_{0_p}(TM) \cong T_p M \oplus T_p M$
-- see Chapter 26 for the further theory concerning this.

\msk \nin {\bf 5.7.} {\sl The algebra of derivations of $M$.}
A vector field $X:M \to TM$ on a symmetric space $M$ is called a
{\it derivation} if $X$ is also a homomorphism of symmetric spaces.
This can be rephrased by saying that $(X \times X,X)$ is
$m$-related. Lemma 4.3 therefore implies that the space
$\g \subset \X(M)$ of derivations is stable under the Lie bracket.
It is also easily checked that it is a $\K$-submodule of $\X(M)$,
and hence $\g \subset \X(M)$ is a Lie subalgebra.

\msk \nin {\bf 5.8.} {\sl The Lie triple system of a symmetric space with
base point.}
We fix a base point $o \in M$. 
The map $X \mapsto T \sigma_o \circ X \circ \sigma_o$ is a Lie algebra 
automorphism of $\X(M)$  of order 2 which stabilizes 
$\g$. We let
$$
\g=\g^+ \oplus \g^-, \quad
\g^{\pm} = \{ X \in \g \vert \, T \sigma_o \circ X \circ \sigma_o =
\pm X \}
$$
be its associated eigenspace decomposition (recall our assumption on
 $\K$ !).  The space $\g^+$ is a Lie subalgebra of $\X(M)$,
whereas $\g^-$ is only closed under the triple bracket
$$
(X,Y,Z) \mapsto [X,Y,Z]:=[[X,Y],Z].
$$

\Proposition 5.9.
\item{(i)}
The space $\g^+$ is the kernel of the evaluation map
${\rm ev}_o:\g \to T_oM$, $X \mapsto X(o)$.
\item{(ii)}
Restriction of ${\rm ev}_o$ yields a bijection
$\g^- \to T_oM$, $X \mapsto X(o)$.

\Proof.
(i) Assume $X \in \g^+$.
Then $T_o \sigma X(o)=X(\sigma_o(o))=X(o)$ implies $-X(o)=X(o)$
and hence $X(o)=0$. 
On the other hand, if $X(o)=0$, then
$X(\sigma_o(p))=X(m(o,p))=Tm(X(o),X(p))=Tm(0_o,X(p))=
T\sigma_o X(p)$,
whence $X \in \g^+$.

(ii)
By (i), $\g^- \cap \ker({\rm ev}_o) = \g^- \cap \g^+ = 0$, and hence
${\rm ev}_o:\g^- \to T_o M$ is injective. It is also surjective:
let $v \in T_o M$. Consider the map
$$
\tilde v = {1 \over 2} Q(v) \circ z:
M \to TM, \quad p \mapsto {1 \over 2} Q(v) 0_p = {1 \over 2}
Tm(v,Tm(0_o,0_p)).
\eqno (5.9)
$$
It is a composition of homomorphisms and hence is itself a
homomorphism from $M$ into $TM$. Moreover, as seen in the proof
of Proposition 5.6, $\tilde v(o)=v$. Thus we will be done if can show 
that $\tilde v \in \g^-$.
First of all, $\tilde v$ is a vector field since,  for all $p \in M$ and
$w \in T_p M$,
$Q(v)w \in T_{m(o,m(o,p))}M=T_p M$.
Finally,
$$
\eqalign{
T\sigma_o \circ \tilde v \circ \sigma_o & =
{1 \over 2} T\sigma_o \circ Q(v) \circ z \circ \sigma_o \cr 
& = {1 \over 2} Q(T \sigma_o v) \circ z =
{1 \over 2} Q(-v) \circ z = - \tilde v. \cr}
$$
(The notation $\tilde v$ is taken from [Lo69]. In later chapters and in
[Be00] we use also the notation $l(v):=\tilde v$ and consider
$l(v)$ as a ``vector field extension of $v$''.)
\qed

\nin The space $\m:=T_o M$ with triple bracket given by
$$
[u,v,w]:=[[\tilde u,\tilde v], \tilde w](o)
$$
is called the {\it Lie triple system (Lts) associated to $(M,o)$}.
It satisfies the identities of an abstract Lie triple system over
$\K$ (cf.\ [Lo69, p.\ 78/79] or [Be00]): for all $u,v,w \in \m$,

\msk
\item{(LT1)} $[u,u,v]=0$,
\item{(LT2)} $[u,v,w]+[v,w,u]+[w,u,v]=0$,
\item{(LT3)} the operator $R(u,v):=[u,v,\cdot]$ is a derivation of the
trilinear product $[\cdot,\cdot,\cdot]$.

\msk
\nin In fact, these identities are easy consequences of the defining
identities of a Lie algebra: they hold in $\g$ and hence also in
the subspace $\g^-$ of $\g$ which is stable under taking triple Lie
brackets.

\Theorem 5.10. Let $M$ be a symmetric space over $\K$ with base point $o$.
\item{(i)} The triple Lie bracket of the Lts $\m$ associated to $(M,o)$
is of class $C^0$.
\item{(ii)} If $\phi:M \to M'$ is a homomorphism of symmetric spaces 
such that $\phi(o)=o'$,
then $\dot \phi:=T_o \phi:\m \to \m'$ is a homomorphism of associated
Lts.

\Proof.
One uses the same arguments as in the proof of Theorem 5.4.
\qed

\msk
\nin {\bf 5.11.} {\sl Lie functor, group case, and the classical spaces.}
The functor described by the preceding theorem will be called
the {\it Lie functor for symmetric spaces and Lie triple systems
(defined over $\K$)}, and the same integrability problem as for
Lie groups arises. It generalizes the integrability problem for
Lie groups because every Lie group can be turned into a symmetric
space by letting
$$
\mu:G \times G \to G, \quad (x,y) \mapsto xy^{-1} x;
$$
then $(G,\mu)$ is a symmetric space, called  {\it of group type}
 (Properties (M1)--(M3) are
immediate; for (M4) one proves as usual in Lie group theory
 that $T_e j = - \id_{T_e G}$ for the inversion map $j:G \to G$);
the Lie triple system of $(G,\mu,e)$ is then $\g$ with the
triple bracket
$$
[X,Y,Z] = {1 \over 4} [[X,Y],Z]
$$
-- for proving this, one can use the same arguments as in [Lo69, p.\ 81].
Note that, in the group case, $\Aut(G,\mu)$ acts transitively on the
space since all left and right translations are automorphisms;
however, the action of the transvection group $G(G,\mu)$ is in general
not transitive -- 
e.g., take the Lie group $\Gl(n,\K)$;
we have $\sigma_x \sigma_y = l_x r_x l_{y^{-1}} r_{y^{-1}}$ in
terms of left and right translations, and hence the determinant
of $\sigma_x \sigma_y(z)$ is congruent to the determinant of $z$
modulo a square in $\K$. This shows that the orbit structure of
$G=\Gl(n,\K)$ under the action of $G(G,\mu)$ is at least as
complicated as the structure of $\K^\times/(\K^\times)^2$.

The ``classical spaces'' can all be constructed in the following
way: assume $G$ is a Lie group and $\sigma$ an involution of $G$;
we assume moreover that the {\it space of symmetric elements}
$$
M:= \{ g \in G | \, \sigma(g) = g^{-1} \}
\eqno (5.11)
$$
is a submanifold of $G$ in the sense of 2.3.
(The last assumption is not automatically satisfied since we do not
have an exponential map at our disposition; however, in all 
``classical cases'' mentioned below, it is easily checked either directly or
by remarking that in ``Jordan coordinates'' it is indeed automatic,
cf.\ [Be00].)
 Then $M$ is a subsymmetric
space of $(G,\mu)$ in the obvious sense. If $\g = \h \oplus \q$
is the decomposition of $\g={\rm Lie}(G)$ into $+1$- and $-1$-eigenspaces
of $T_e \sigma$, then, as in [Lo69, p.\ 82], it is seen that
 the Lie triple system of $(M,\mu,e)$
is $\q$ with $[X,Y,Z] = [[X,Y],Z]$. For an arbitrary base point
$g \in M$, the Lie triple system structure on $T_gM$
is described in the same way as for $g=e$,
but replacing $\sigma$ by the involution $I_g \circ \sigma$ where
$I_g(x)=gxg^{-1}$ is conjugation by $g$. 
In contrast to the group case, here the action of the automorphism
group $\Aut(M)$ is in general no longer transitive on $M$.
For instance, if $G=\Gl(n,\K)$  and $\sigma=\id$, then the orbits
in $M=\{ g \in G | \, g^2 = \1 \}$ are described by
the matrices $I_{p,q}=\bigl( {\1_p \atop 0}{0\atop -\1_q} \bigr)$,
$p+q=n$. 
Similarly, taking for $G$ a classical matrix group (in the sense explained
after Theorem 5.4) and for $\sigma$ a ``classical involution''
(see [Be00, Chapter I.6] for a fairly exhaustive list in the real case)
we get analogs of all matrix series from Berger's classification
of irreducible real symmetric spaces [B57], but where the number
of non-isomorphic types may be considerably bigger in case of 
base-fields or rings different from $\R$ or $\C$. 
All these symmetric spaces are ``Jordan symmetric spaces'' in the
sense of [BeNe05] where a very general construction of symmetric
spaces (of classical or non-classical type) is described.

\vfill \eject 

\def\title{II. TANGENT OBJECTS AND SCALAR EXTENSIONS} 
\def\rightheadline{\hfil{\eightrm \title}\hfil\tenrm\folio}

\sectionheadline
{II. Interpretation of tangent objects via scalar extensions}

\subheadline
{6. Scalar extensions. I: Tangent functor and dual numbers}

\msk \nin
{\bf 6.1.} {\sl The  ring of dual numbers over $\K$.}
The {\it ring of dual numbers over $\K$}
is the truncated polynomial ring $\K[x]/(x^2)$; it can also be
defined in a similar  way as complex numbers,
replacing the condition $i^2=-1$ by $\epsilon^2 = 0$: 
$$
\K[\epsilon]:= \K \oplus \epsilon \K, 
\quad (a + \epsilon b)(a' + \epsilon b')=aa' + \epsilon(ab'+ba').
\eqno (6.1)
$$
In analogy with the terminology for complex numbers, let us call
$z=x+\epsilon y$ a {\it dual number}, 
$x$ its {\it spacial part},
$y$ its {\it infinitesimal part}
and $\epsilon$ the {\it infinitesimal} or {\it dual number unit}.
A dual number is invertible iff its spacial part is invertible, and
 inversion is given by
$$
 (x+\epsilon y)^{-1} = {x-\epsilon y \over x^2} .
\eqno (6.2)
$$
It follows that $\K[\epsilon]$ is again a topological ring with
dense group of invertible elements.
Taking the spacial part is a ring homomorphism $\K[\epsilon] \to \K$,
and its kernel is the ideal $\epsilon \K$ of ``purely infinitesimal'' 
numbers.
If $\K$ is a field, then this is the unique maximal ideal of $\K[\epsilon]$. 
A matrix realization of $\K[\epsilon]$
in the algebra of $2 \times 2$-matrices over
$\K$ is given by
$$
\K[\epsilon] \cong \K \, \pmatrix{1 & 0 \cr 0 & 1 \cr} \oplus \K 
\pmatrix{0 & 0 \cr 1 & 0 \cr} .
\eqno (6.3)
$$
For any topological $\K$-module $V$, the scalar extension
$$
V \otimes_\K \K[\eps] = V \otimes_\K (\K \oplus \eps \K) =V  \oplus \eps V
$$
(with $\K[\eps]$-action given by $(a+\eps b)(v+\eps w)=ab+\eps(aw+bv)$)
is a topological $\K[\eps]$-module. In particular, it is a smooth manifold
over $\K[\eps]$. More generally, for every open set $U \subset V$, the set
$TU:=U \oplus \eps V$ is open in $V\oplus \eps V$ and hence is a manifold over
$\K[\eps]$.

\Theorem 6.2.
If $M$ is a manifold of class $C^{k+1}$ over $\K$, modelled on $V$,
then
$TM$ is, in a canonical way, a manifold of class $C^k$ over
$\K[\epsilon]$, modelled on the scalar extended $\K[\epsilon]$-module
$V \oplus \epsilon V$.
If $f:M \to N$ is of class $C^{k+1}$ over $\K$, then
$Tf:TM \to TN$ is of class $C^k$ over $\K[\epsilon]$.
Thus the tangent functor can be characterized as the ``functor 
of scalar extension" from
manifolds over $\K$ into manifolds over $\K[\epsilon]$ agreeing on
open submanifolds of $\K$-modules with the algebraic scalar extension
functor.

\Proof.
Change of charts on $TM$ is described by the transition maps
$T\phi_{ij}(x,v)=$
$(\phi_{ij}(x),$
$ d\phi_{ij}(x)v)$.
Therefore the theorem follows from the following claim:

\Theorem 6.3.
Assume $V$ and $W$ are topological $\K$-modules, $U$ open in $V$ and
 $f:V \supset U \to W$ is of class $C^{k+1}$ over $\K$. 
Then $Tf: V \oplus \eps V
 \supset U \times V \to W \times W = W \oplus \eps W$ is of class $C^k$
over $\K[\epsilon]$. More precisely, the ``higher order
difference quotient maps''
defined in {\rm Section 1.7} are related by
$$
(Tf)^{[k]} = T(f^{[k]}),
$$
where the difference quotient map on the
left hand side is taken with respect to $\K[\epsilon]$.

\Proof.
It suffices to prove the claim for $k=1$, i.e. that
$$
(Tf)^{[1]} = T(f^{[1]}),
$$
the general case
then follows by a straightforward induction.
Assume $f:V \supset U \to W$ is a map of class $C^2$. 
We are going to apply the tangent functor $T$ to the  defining
identity
$$
t \cdot f^{[1]}(x,v,t)= f(x+tv) - f(x)
\eqno (!)
$$
of $f^{[1]}$. Let us 
write $(!)$ in the form of a commutative diagram:

\msk
$$
\matrix{
U^{[1]} & {\buildrel \beta \over \rightarrow} & U^{[1]} \times U &
{\buildrel \gamma \over \rightarrow} &
  U \times U 
\cr
\alpha \downarrow \phantom{\alpha} & & & & 
\phantom{\id_\K \times f} \downarrow f \times f  \cr
\K \times U^{[1]} & & & & W \times W
\cr
\id_\K \times f^{[1]} \downarrow \phantom{ \id_\K \times f}
& & &  & \phantom{\id_\K } \downarrow \phi
\cr
\K \times W & {\buildrel m_W \over \rightarrow} & W & = & W
 \cr}
$$

\msk \nin 
where $\alpha(x,v,t)=(t,x,v,t)$,
$\beta(x,v,t)=(x,v,t,x)$,
$\gamma(x,v,t,y)=(x+tv,y)$ and $\phi(u,v)=u-v$, i.e.

\msk
$$
\matrix{
(x,v,t) & \mapsto & (x,v,t,x) & \mapsto & (x+tv,x)
\cr
 \downarrow  & & & & \downarrow \cr
(t,x,v,t) & & & &  (f(x+tv),f(x)) 
\cr
\downarrow 
& & &  &  \downarrow 
\cr
(t,f^{[1]}(x,v,t)) 
 & \mapsto  & t f^{[1]}(x,v,t)   & 
=  &f(x+tv)-f(x)
 \cr}
$$

\ssk \nin 
Applying the tangent functor $T$ to this diagram,
we get a diagram of the same type, where
all maps are replaced by their tangent maps.
(Note that this is the same argument as in the proof that the tangent
group $TG$ of $G$ is a Lie group, resp.\ the analog for symmetric spaces,
cf.\ Chapter 5; it only uses that $T$ is a
covariant functor commuting with direct products and
diagonal imbeddings.)
Let us describe the tangent map 
 $Tm:T\K \times T\K \to T\K$ of $m:\K \times \K \to \K$:
for $t \in \K^\times$,
$$
m^{[1]}((x_1,x_2),(v_1,v_2),t)={1 \over t}
(m(x_1 + t v_1,x_2 + t v_2)-m(x_1,x_2)) = 
v_1 x_2 + x_1 v_2 + t v_1 v_2,
$$
and letting $t=0$, we get
$Tm((x_1,x_2),(v_1,v_2)) = (x_1 x_2, v_1 x_2 + x_1 v_2)$.
Comparing with (6.1), we see that $(T\K,Tm)$ is isomorphic
to the ring $\K[\epsilon]$ of dual numbers over $\K$.
For this reason, $\K[\epsilon]$ will also be called
the {\it tangent ring of $\K$}.
In a similar way, we calculate the tangent maps of the
 structural maps $m_V:\K \times V \to V$ and $a:V \times V \to V$
of the topological $\K$-module $V$ are $C^1$ and find that
$$
\eqalign{
Tm_V((r,s),(x,v)) & = (rx,rv+sx), \quad \cr
Ta((x,v),(x',v')) & =(x+x',v+v') \cr}
$$
which is also obtained by writing $r+\epsilon s$ for $(r,s)$
and $x+\epsilon v$ for $(x,v)$ and calculating in
the scalar extended module $V \oplus \eps V$.
Hence $TV \cong V \oplus \eps V$ as modules over the ring
$T\K \cong \K[\epsilon]$, and we call $V \oplus \eps V$ also 
the {\it tangent module of $V$}.
Summing up, when applying the tangent functor to (!),
all structural maps for $V,W$ are replaced by the corresponding
structural maps for the corresponding scalar extended modules,
 and $f$ and $f^{[1]}$ are
replaced by $Tf$ and $T(f^{[1]})$. 
This in turn means that the identity
$$
t \cdot T(f^{[1]})(x,v,t)= Tf(x+tv) - Tf(x)
\eqno (!!)
$$
holds for all $(x,v,t) \in TU \times TV \times T\K$ with
$x+tv \in TU$, and this means that $Tf$ is of class $C^1$ over
the ring $T\K$, and that
$(Tf)^{[1]} = T(f^{[1]})$.
\qed

\msk \nin
{\bf 6.4.} {\sl Application to Lie groups and symmetric spaces.}
For any $\K$-Lie algebra $\g$, the scalar extended
algebra $\g \otimes_\K \K[\epsilon]$ can be described as
$$
\g \oplus \epsilon \g, \quad [X+\epsilon Y,X' + \epsilon Y'] =
[X,X'] + \epsilon([X,Y'] + [X',Y]), 
$$
i.e. one simply takes the $\epsilon$-bilinear bracket on
$\g \oplus \epsilon \g$ and takes account of the relation $\epsilon^2=0$.
Similarly, the scalar extended Lie triple system $\m \otimes_\K \K[\epsilon]
\cong \m \oplus \epsilon \m$
of a $\K$-Lie triple system $\m$ can be described.

\Theorem 6.5.
If $G$ is a Lie group over $\K$, then
$TG$ is a Lie group over $\K[\epsilon]$ with Lie algebra
$\g \oplus \epsilon \g$, and if 
 $M$ is a symmetric space over $\K$ with base point $o$ and
 associated Lie triple system $\m$, 
then $TM$ is a symmetric space over $\K[\epsilon]$ 
with base point $0_o$ and associated Lie triple system $\m \oplus \epsilon \m$.

\Proof.
If $M$ is a manifold with product map $m:M \times M \to M$, then, by
Theorem 6.2,
$Tm:TM \times TM \to TM$ is smooth over $\K[\epsilon]$.
Thus, if $G$ is a Lie group over $\K$, then $TG$ is a Lie group,
not only over $\K$, but also over $\K[\epsilon]$.
Moreover, all constructions used to define the Lie bracket of $G$
naturally commute with the functor of scalar extension by $\K[\epsilon]$,
and hence the Lie algebra of $TG$ is the scalar extended Lie
algebra $\g \otimes_\K \K[\epsilon]$ of $\g$. 
Similar arguments apply to symmetric spaces. 
\qed 

\msk \nin
{\bf 6.6.} {\sl Comparison of the structures of $T\K$ and of $TM$.}
The fiber over $0$ in the ring $T\K$ is the ideal
$\epsilon \K$; it is the kernel
of the canonical projection $T\K \to \K$, which is a ring
homomorphism. In other words, we have an exact sequence of non-unital rings
$$
0 \to \epsilon \K \to T\K \to \K \to 0.
\eqno (6.4)
$$
On the level of manifolds, this corresponds to
the inclusion of the tangent space $T_o M$ at a base point
$o \in M$, followed by the canonical projection:
$$
T_o M \to TM \to M.
$$
By a direct computation, one checks that for any $\lambda \in \K$, the map
$$
l_\lambda :\K[\epsilon] \to \K[\epsilon], \quad
r + \epsilon s \mapsto r + \epsilon \lambda s
$$
is a ring homomorphism and that
$$
l:\K \times T\K \to T\K, \quad (\lambda, z) \mapsto l_\lambda(z)
$$
is an action of $\K$ on $T\K$.
If $\lambda = -1$, we may call $l_{-1}$ the  ``spatial conjugation";
if $\lambda = 0$, $l_0$ corresponds  to the projection onto $\K$.
On the level of manifolds, we have corresponding smooth maps
$$
l_\lambda: TM \to TM, \, \, v \mapsto \lambda v, \quad \quad
l:\K \times TM \to TM
$$
where the product is taken in the $\K$-module $T_{\pi(v)} M$.
Finally, there is also an inclusion map
$\K \to \K[\epsilon]$, which corresponds to the zero section
$$
z_{TM}: M \to TM, \quad x \mapsto 0_x.
$$

\msk \nin {\bf 6.7.} {\sl The trivial scalar extension functor.}
Not every manifold over $\K[\epsilon]$ is (isomorphic to)
a tangent bundle, not even locally. In fact,
since $\epsilon \K$ is an ideal in $T\K$,
every $\K$-module is also a $\K[\epsilon]$-module, just by
letting act $\epsilon$ trivially
(composition of the action of $\K$ and the ring homomorphism
$\K[\epsilon] \to \K$, $r+\epsilon s \mapsto r$).
%
%
In a similar way, every manifold over $\K$ is also a manifold
over $\K[\epsilon]$, just by letting act $\epsilon$ trivially
everywhere; we call this {\it the trivial scalar extension functor}.
Therefore, in order to characterize the $\K[\eps]$-manifold structure
in Theorem 6.2, it was necessary to state explicitly that, on chart
domains, it coincides with the algebraic scalar extension.
Occasionally, it may be useful to consider $\K$-manifolds $M$
as $\K[\eps]$-manifolds in the trivial way. For instance, the projection
$\pi:TM \to M$ is then smooth over $\K[\epsilon]$ -- indeed, in a chart
the projection $V \oplus \eps V \to V$ is then  $\K[\eps]$-linear,
whence smooth over $\K[\eps]$.

\msk
Theorem 6.2 implies that $\Diff_\K(M) \to \Diff_{T\K}(TM)$,
$g \mapsto Tg$ is a well-defined imbedding. In fact, it is not too 
difficult to determine now the full group of $T\K$-diffeomorphisms of
$TM$:

\Theorem 6.8. Assume $M$ and $N$ are manifolds over $\K$.
\item{(1)}
 For any vector field $X:M \to TM$,
the infinitesimal automorphism $\tilde X:TM \to TM$,
$v \mapsto v + X(\pi(v))$ (cf.\ {\rm Section 4.4})
 is smooth over $\K[\epsilon]$.
\item{(2)}
Assume $F:TM \to TN$ is smooth over $T\K$. Then $F$ maps fibers to fibers
(i.e., there is a smooth $f: M \to N$ with $\pi_N \circ F = f \circ \pi_M$)
and $F$ acts $\K$-affinely from fibers to fibers.
\item{(3)}
The following is an exact and splitting sequence of groups:
$$ \matrix{
0 & \to & \X(M) & {\buildrel X \mapsto \tilde X \over \to} & \Diff_{T\K}(TM)&
{\buildrel F \mapsto f \over \to} &  \Diff_\K(M) & \to 1. \cr}
$$
In particular, any $T\K$-smooth diffeomorphism $F:TM \to TM$ may be written
in a unique way as $F=\tilde X \circ Tf$ for some $X \in \X(M)$ and
$f \in \Diff(M)$.

\Proof. (1) We use the following simple lemma:

\Lemma 6.9. Assume $F:V \oplus \eps V \supset U \oplus \eps V
 \to W \oplus \eps W$ is of the form
$F(a+\eps b)=\eps g(a)$ with $g:U \to W$ of class $C^\infty$ over $\K$.
 Then $F$ is of class $C^\infty$ over $T\K$.

\Proof. The claim can be checked by a direct computation. Here is a
computation-free argument:
scalar multiplication by $\eps$,
$$
l_\eps: W \oplus \eps W \to W \oplus \eps W, \quad x + \eps v \mapsto \eps x,
$$
is $\K[\eps]$-linear and continuous, whence of class $C^\infty$ over 
$\K[\eps]$. 
Now, we may write $F(a + \eps b)=\eps(g(a)+\eps dg(a)b) = l_\eps \circ Tg
(a + \eps b)$, and thus $F$ is seen to be a composition of $\K[\eps]$-smooth
maps. 
\qed

\nin Coming back to the proof of the theorem,
 in a bundle chart we write $\tilde X(x + \epsilon v)=x +
\epsilon v +\epsilon X(x) = \id(x + \eps v) + \epsilon X(x)$,
and hence $\tilde X$ is represented by the sum of the identity map and
a map having the form given in the preceding lemma. It follows that
$\tilde X$ is smooth over $\K[\eps]$.

\ssk
(2)
We represent $F$ with respect to charts and write again $F:TU \to TU'$ for the 
chart representation. Since $F$ is smooth over $T\K$, we have the following
second order Taylor expansion of $F$ at the point $x=x+\eps 0 \in TU$
(Theorem 1.11):
$$
F(x+\eps v)=F(x)+ \eps dF(x) v + \eps^2 a_2(x,v) +
\eps^2 R_2(x,v,\eps)  =F(x)+ \eps dF(x) v.
\eqno (6.5)
$$
Since $\eps dF(x)v \in \eps W$,
this shows that $F$ maps fibers to fibers. Next, we decompose
$F(x+\eps v)=F_1(x+\eps v)+\eps F_2(x + \eps v)$ into its ``spacial'' and
``infinitesimal part'' $F_i: TU \to W$. 
We let also $f:= F_1|_U:U \to W$, which is the chart representation of the map
 $f:=\pi_N \circ F \circ z_M: M \to N$. 
Over $\K$, we may write $dF=d(F_1 + \eps F_2)=dF_1 + \eps dF_2$, 
and hence we get from (6.5)
$$
\eqalign{ F(x+\eps v) & = F(x) + \eps dF(x) v 
= F_1(x) + \eps (dF_1(x) v + F_2(x)) \cr
& = f(x) + \eps df(x) v + \eps F_2(x) = Tf(x)v + \eps F_2(x). \cr}
$$  
Thus $F$ acts affinely on fibers (the linear part is $df(x)$ and the 
translation part is translation by $F_2(x)$).

\ssk (3)
The first map is well-defined by (1); clearly, it is an injective group
homomorphism. The second map is surjective since it admits the
section $g \mapsto Tg$. Clearly $\X(M)$ belongs to the kernel of the second
map. Conversely, assume $F \in \Diff_{T\K}(TM)$ is such that $f=\id_M$.
In a chart representation as in the proof of Part (2), this means that 
$F_1(x)=x$. But as seen above, then the linear part of $F|_{T_x M}$
is the identity, hence $F|_{T_x M}$ acts as a translation for all $x$,
and thus $F$ is an infinitesimal automorphism.  
This proves exactness of the given sequence. 

Finally, as for any split exact sequence, the middle group $\Diff_{T\K}(TM)$
can now be written as a semidirect product of $\X(M)$ and $\Diff(M)$,
implying the last statement.
\qed

\nin
In fact, the proof of Part (2) of the theorem shows a bit more than anounced,
namely that any $T\K$-smooth map $F:TM \to TN$ can be written in the
form $F(u)=Tf(u) + Y(\pi_M(u))$ where $Y:M \to TN$ is a ``vector field over
$f$'' (represented in a chart by $F_2|_U:U \to W$),
 i.e., $\pi_N \circ Y = f$. --
In Chapter 28 these facts will be generalized, and we will show that the
sequence from the preceding theorem is the precise analog of the 
sequence $0 \to \g \to TG \to G \to 1$ of a Lie group $G$.

\msk
\nin {\bf 6.10.} {\sl Remark on the interpretation of
the tangent functor as a contraction of the direct product functor.} 
Intuitively, the tangent bundle may be seen in the following way:
take the direct product $M \times M$, realize $M$ as the diagonal
in $M \times M$, take the vertical fiber over the point $x=(x,x)$
and enlarge it more and more (with fixed center $(x,x)$)
until you only see the ``infinitesimal neighborhood of $x$", and finally forget
the rest. What you get is the tangent bundle, seen as a ``contraction
of the direct product". This intuitive picture
can be put on a completely rigorous footing:
look at the family of rings $\K[x]/(x^2 - tx)$
with $t \in \K$; taking the classes of $1$ and $x$ as a $\K$-basis,
this ring is described by the product
$$
(a+b x)(c + dx) = ac + (ad+bc+t bd)x.
\eqno (6.6)
$$
A matrix realization of this ring
 in the algebra of $2 \times 2$-matrices over $\K$ is given by
$$
\K[x]/(x^2-tx) \cong \K \pmatrix{1 & 0 \cr 0 &1 \cr} \oplus \K 
\pmatrix{t & 0 \cr 1 & 0 \cr} .
$$
If $t$ is invertible in $\K$ (in particular, for $t=1$),
this ring is isomorphic to the direct product of rings $\K \times \K$:
take, in the ring $\K \times \K$ with
pointwise product $(x,y)(x',y')=(xx',yy')$, the new basis
$e=(1,1), f=(0,t)$; $e$ is the unit element in $\K \times \K$,
and $f^2 = tf$; hence the product in $\K \times \K$ is described by
$(x_1 e + v_1 f)(x_2 e + v_2 f)=
x_1 x_2 e +(x_1 v_2 + v_1 x_2 + t v_1v_2)f$, which is the same as
(6.6). If $t=0$ we get the dual numbers
over $\K$; they appear here as a ``contraction
for $t \to 0$" of the rings $\K \times \K \cong \K[x]/(x^2 - tx)$.
In a similar way, differential calculus can be seen as
a contraction of difference calculus. Namely,
recall from Equation (1.4) the definition of the
extended tangent map $\hat T f(x,v,t)=(f(x),f^{[1]}(x,v,t),t)$.
Since the functorial rule $\hat T(f \circ g)=\hat Tf \circ 
\hat T g$ holds, all arguments used in the proof of Theorem
6.3 can be adapted to the case of  the functor $\hat T( -- ,t)$.
(Note that, for invertible $t$, this functor is, by a change of
variables, isomorphic to the direct product functor.)
In particular, applying this functor to the ring $\K$ itself, we get
$$
\hat T(m)((x_1,x_2),(v_1,v_2),t))=
(x_1 x_2, v_1 x_2 + x_1 v_2 + t v_1 v_2, t).
$$
Fixing and suppressing the last variable, this leads 
to a product 
$(\K \times \K)\times (\K \times \K) \to \K \times \K$
defined by
$$
(x_1,v_1) \cdot (x_2,v_2):= \hat Tm ((x_1,x_2),(v_1,v_2),t)=
(x_1 x_2, v_1 x_2 + x_1 v_2 + t v_1 v_2).
$$
As in the case $t=0$, this defines on $\K \times \K$ the
structure of a topological ring; 
comparing with (6.5), we see that this is nothing but
the ring $\K[x]/(x^2 - tx)$.
Applying our functor $\hat T(--,t)$ to (!), we thus see 
 that for all $t \in \K$, $\hat T f(--,t)$ is of class $C^k$ 
 over the ring $\K/(x^2-tx)$.

For invertible $t$, and in particular for $t=1$, the preceding claims
can all be checked by  direct and elementary calculations.
Taking the ``limit case $t \to 0$", we may again deduce Theorem 6.3,
thus giving  another proof of it. However, proving the limit case Theorem 6.3
 by checking the definition of differentiability over $\K[\eps]$
directly, leads to fairly long and involved calculations.
(The reader may try this as an exercise, to get accustomed to the
definitions.)

\vfill \eject

\subheadline
{7. Scalar extensions. II: Higher order tangent functors}

\msk \nin
{\bf 7.1.} {\sl Higher order tangent bundles.}
For any manifold $M$, the {\it second order tangent bundle}
is $T^2M:=T(TM)$, and the {\it $k$-th order tangent 
bundle} is inductively defined by $T^k M:=T(T^{k-1}M)$.
A repeated application of Theorem 6.2 shows:

\Theorem 7.2.
The $k$-th order tangent bundle $T^kM$ is, in a canonical
way, a manifold over the ring of iterated dual numbers
$T^k \K = \K[\epsilon_1] \ldots [\epsilon_k]$,
and if $f:M \to N$ is smooth over $\K$, then $T^k f:T^kM \to T^k N$
is smooth over $T^k \K$.
\qed

\nin
Thus the higher order tangent functor $T^k$ can be seen as the
functor of scalar extension by the ring $T^k \K$.
We will add the index $\epsilon_i$ to the
tangent functor symbol $T$ if we want to indicate that we mean
scalar extension with ``new infinitesimal unit" $\epsilon_i$.
Psychologically, we associate $\eps_1$ to the ``first'' scalar extension
and $\eps_k$ to the ``last'' scalar extension; but we will see below that
all orders are equivalent, i.e., there is an action of the permutation
group $\Sigma_k$ by automorphisms of the whole structure. -- 
When working with chart representations, we may use, for  bundle chart
domains of $TM$ over $U \subset M$, resp.\ of $TTM$ over $U$,
 the following two equivalent notations
$$
\eqalign{
TU & = \{ (x, v) | \, x \in U, v\in V \} =
\{ x + \epsilon v | \, x \in U, v\in V \},  \cr
TTU & = \{ (x,v_1,v_2,v_{12}) | x \in U, v_1,v_2,v_{12} \in V \} \cr
& = \{ x + \epsilon_1 v_1 + \epsilon_2 v_2+ \epsilon_1 \epsilon_2
v_{12} | x \in U, v_1,v_2,v_{12} \in V \} \cr}
$$
but the more classical component notation will be gradually abandoned in
favor of the ``algebraic notation'' which turns out to be more suitable
(the component notation corresponds to the one used in [La99] and
in most other references where $TTM$ is treated).

\msk \nin {\bf 7.3.} {\sl Structure of $TT\K$ and of $TTM$.}
First of all, one should note that $TTM$ is a fiber bundle over
$M$, but it is {\it not} a vector bundle -- the transition maps,
calculated in Section 4.7, are not linear in fibers.
The ``algebraic structure'' of $TTM$
 will be discussed in Chapter 9. Here we focus on properties
that are directly related to properties of the ring
$$
TT\K  \cong (\K[\epsilon_1]) [\epsilon_2] =
\K \oplus \epsilon_1 \K \oplus \epsilon_2 \K \oplus \epsilon_1
\epsilon_2 \K,
$$
(relations: $\epsilon_1^2 = 0 = \epsilon_2^2$ and
$\epsilon_1 \epsilon_2 = \epsilon_2 \epsilon_1$).
We are going to discuss the following features of $TT\K$:
projections; injections; ideal structure; canonical automorphisms and
endomorphisms.

\ssk {\sl (A) Projections.}
The projections onto the ``spacial part'' corresponding to the
various extensions fit together to a commutative diagram
$$
\matrix{  & & \K[\epsilon_1][\epsilon_2] & & \cr
& \swarrow & & \searrow & \cr
\K[\epsilon_1] & & & & \K[\epsilon_2] \cr
&\searrow & & \swarrow & \cr
& & \K & & \cr} .
\eqno (7.1)
$$
For every $\K$-module $V$, we get by scalar extension a $TT\K$-module
$TTV$ and a diagram of $\K$-modules corresponding to (7.1),
which, on the level of manifolds, corresponds to the diagram
of vector bundle projections 
$$
\matrix{  & & TTM & & \cr
& \swarrow & & \searrow & \cr
T_{\epsilon_1} M & & & & T_{\epsilon_2}M \cr
&\searrow & & \swarrow & \cr
& & M & & \cr} .
\eqno (7.2)
$$
In a bundle chart,  (7.2) is represented by
$$
\matrix{  & & x+\epsilon_1 v_1 + \epsilon_2 v_2 + \epsilon_1 \epsilon_2
v_{12} & & \cr
& \swarrow & & \searrow & \cr
x + \epsilon_1 v_1 & & & & x + \epsilon_2 v_2 \cr
&\searrow & & \swarrow & \cr
& & x & & \cr} .
\eqno (7.3)
$$

\ssk {\sl (B) Injections.}
There is  a diagram of inclusions of rings
$$
\matrix{  & & \K[\epsilon_1][\epsilon_2]  & & \cr
& \nearrow & & \nwarrow & \cr
\K[\epsilon_1] & & & & \K[\epsilon_2] \cr
&\nwarrow & & \nearrow & \cr & & \K & & \cr} 
\eqno (7.4)
$$
which corresponds to a diagram of zero sections of vector bundles
$$
\matrix{  & & TTM & & \cr
& \nearrow & & \nwarrow & \cr
T_{\epsilon_1} M & & & & T_{\epsilon_2}M \cr
&\nwarrow & & \nearrow & \cr & & M & & \cr} .
\eqno (7.5)
$$

\ssk {\sl (C) Ideal structure.}  In the following, we abbreviate $R:= TT\K$.
The following are inclusions of ideals:
$$
\matrix{  &  & \epsilon_1 R = \epsilon_1 \K \oplus \epsilon_1 \epsilon_2 \K
= \eps_1 \K[\eps_2]
&  &  \cr
	  &\nearrow &       &\searrow & \cr
      \epsilon_1 \epsilon_2 R = \epsilon_1 \epsilon_2 \K & & & & R \cr
	  &\searrow &       &\nearrow & \cr
 &  & \epsilon_2 R =\epsilon_2 \K \oplus \epsilon_1 \epsilon_2 \K =
\eps_2 \K[\eps_1] & & \cr} .
\eqno (7.6)
$$
Here, also the composed inclusion $\epsilon_1 \epsilon_2 \K \subset R$
is an ideal. The three ideals from (7.6)
are kernels of three projections appearing in the following three
exact sequences of (non-unital) rings.
$$
\matrix{
0 & \to & \epsilon_2 R & \to & R &
{\buildrel p_1 \over \longrightarrow} & \K[\epsilon_1] & \to & 0, \cr
0 & \to & \epsilon_1 R &  \to & R &
{\buildrel p_2 \over \longrightarrow} & \K[\epsilon_2] & \to & 0, \cr
0 & \to & \epsilon_1 \epsilon_2 R & \to & R &
{\buildrel p_1 \times p_2\over \longrightarrow} &
\K[\epsilon_1] \times \K[\epsilon_2] & \to &  0. \cr}
\eqno (7.7)
$$ 
The first two sequences correspond to the two versions of the
vector bundle projection $TTM \to TM$ that appear in (7.2),
 where the kernel describes
just the fibers. The third sequence corresponds to a sequence of
 smooth maps
$$
TM \to TTM \to TM \times_M TM.
\eqno (7.8)
$$
The first map in (7.8) is the injection given in a chart by
$x+\epsilon v \mapsto x + \epsilon_1 \epsilon_2 v$ and the
second map in (7.8) is the surjection given
by 
$$
x +\epsilon_1 v_1 +\epsilon_2 v_2 + \epsilon_1
\epsilon_2 v_{12} \mapsto (x + \epsilon_1 v_1,x +\epsilon_2 v_2);
$$
 it is just the vector bundle direct sum of
the two projections $TTM \to TM$ from (7.7).
The injection from (7.8)
may intrinsically be defined as $\epsilon_2 \circ T_{\epsilon_1} z=
\epsilon_1 \circ T_{\epsilon_2} z$, where $z:M \to TM$
is the zero section and $\epsilon_i: TTM \to TTM$
is the almost dual structure corresponding to $\epsilon_i$
(cf. Section 4.7).
The image of this map is called the {\it vertical bundle},
$$
VM:=\epsilon_1 \epsilon_2 TM:= \epsilon_1 T_{\epsilon_2} z(TM) =
\epsilon_2 T_{\epsilon_1} z(TM);
\eqno (7.9)
$$
it is a subbundle of $TTM$ that carries the structure
of a vector bundle over $M$ isomorphic to $TM$. 

\ssk {\sl (D) Automorphisms, endomorphisms.}
An important  feature of $TT\K$
 is the ``exchange'' or ``flip automorphism" $\kappa$
 exchanging $\epsilon_1$ and $\epsilon_2$,
$$
\eqalign{
\kappa : TT\K  & \to TT\K, \cr
x_0 + x_1 \epsilon_1 + x_2 \epsilon_2 + x_{12} \epsilon_1 \epsilon_2  &\mapsto
x_0 + x_2 \epsilon_1 + x_1 \epsilon_2 + x_{12} \epsilon_1 \epsilon_2. \cr}
$$
It induces a $\K$-linear automorphism of $TTV$ and a canonical
diffeomorphism (over $\K$) of $TTM$, given in a chart by
$$
\eqalign{
\kappa:TTU & \to TTU, \cr
x+\epsilon_1 v_1+\epsilon_2 v_2 + \epsilon_1 \epsilon_2 v_{12} & \mapsto
x+\epsilon_1 v_2 +\epsilon_2 v_1 + \epsilon_1 \epsilon_2 v_{12}. \cr}
\eqno (7.10)
$$
This map is chart independent: the transition functions of $TTM$ are
calculated in Section 4.7; they are given by $T^2 \phi_{ij}$,
and because of the symmetry of $d^2 \phi_{ij}(x)$ (Schwarz' lemma), 
they commute with $\kappa
$. Thus $\kappa:TTM \to TTM$ does not depend on the chart.

Recall from Section 6.6 that the ring $\K$ acts by endomorphisms on $T\K$:
$$
l:\K \times T\K \to T\K, \quad (\lambda,a+\eps b) \mapsto a + \eps \lambda b
$$ 
which corresponds to the action $l:\K \times TM \to TM$ by scalar 
multiplication in tangent spaces.
By taking the tangent map $Tl$, we get two actions 
$T_{\eps_i} l: T\K \times TT\K \to TT\K$, $i=1,2$,
corresponding to two actions $T\K \times TTM \to TTM$ coming from the
two $T\K$-vector bundle structures of $TTM$ over $TM$. 
Both actions are interchanged by the flip $\kappa$. The combined action
$$
l^{(2)}:
T\K \times TTM \to TTM, \quad (\lambda,u) \mapsto T_{\eps_1} l (\lambda,
T_{\eps_2} l (\lambda,u)) \eqno (7.11)
$$
commutes with $\kappa$. In a chart, for $r,s \in \K$  and 
$u= x+\epsilon_1 v_1 + \epsilon_2 v_2 + \epsilon_1 \epsilon_2 v_{12}$,
$$
\eqalign{
T_{\eps_2}l(r+\eps_2 s,u) & =
x+\epsilon_1 v_1 + \epsilon_2 r v_2 + \epsilon_1
\epsilon_2 (r v_{12} + s v_1), \cr
l^{(2)}(r + \eps_2 , u) & = x + r \epsilon_1 v_1 + r \epsilon_2 v_2 +
 \epsilon_1 \epsilon_2 (r^2 v_{12} + rs v_1 + rs v_2).  \cr}
$$

\msk \nin
{\bf 7.4.} {\sl Structure of $T^k \K$ and of $T^k M$.}
The structure of the ring
$T^k \K = \K[\epsilon_1,\ldots,\epsilon_k]$ can be analyzed
in a similar way as we did above for $k=2$. In the following, we focus
on the case $k=3$, 
$$
T^3 \K =  \K \oplus \epsilon_1 \K \oplus
\epsilon_2 \K \oplus \epsilon_3 \K
\oplus \epsilon_1 \epsilon_2 \K \oplus \epsilon_1 \epsilon_3 \K
\oplus \epsilon_2 \epsilon_3 \K \oplus \epsilon_1 \epsilon_2 \epsilon_3 \K,
$$
with the obvious relations.

\ssk
{\sl (A) Projections.}
There is a commutative cube of ring projections:
$$
\xymatrix@ur@!0@C+10pt{
&\K[\epsilon_1,\epsilon_2] \ar@{-}[rr]\ar@{-}'[d][dd]
 & & \K[\epsilon_1,\epsilon_2,\epsilon_3] \ar@{-}[dd]
\\
\K[\epsilon_1] \ar@{-}[ur]\ar@{-}[rr]\ar@{-}[dd]
& & \K[\epsilon_1,\epsilon_3] \ar@{-}[ur]\ar@{-}[dd]
\\
& \K[\epsilon_2] \ar@{-}'[r][rr]
& & \K[\epsilon_2,\epsilon_3]
\\
\K \ar@{-}[rr]\ar@{-}[ur]
& & \K[\epsilon_3] \ar@{-}[ur]
}
\eqno (7.12)
$$
 On the level of manifolds, the
 corresponding diagram of vector bundle projections will be written 
$$
\xymatrix@dl@!0@C+10pt{
& T_{\eps_0}M \ar@{-}[rr]\ar@{-}'[d][dd] & & M \ar@{-}[dd]
\\
 T^2_{\eps_0 \epsilon_2}M \ar@{-}[ur]\ar@{-}[rr]\ar@{-}[dd]
& &  T_{\eps_2}M \ar@{-}[ur]\ar@{-}[dd]
\\
&  T^2_{\eps_0 \eps_1}M \ar@{-}'[r][rr] & & T_{\eps_1}M
\\
 T^3_{\eps_0 \eps_1 \eps_2} M \ar@{-}[rr]\ar@{-}[ur] & & 
T^2_{\eps_1 \eps_2} M \ar@{-}[ur]
}
$$

{\sl (B) Injections.}
The diagrams written above can also be read, from the base to the top,
as inclusions of subrings; on the level of manifolds, this corresponds to a
diagram of zero sections of vector bundles.
For general $k$, we have similar diagrams of projections and injections
 which could be represented by $k$-dimensional cubes.

\ssk {\sl (C) Ideal structure.}
There is also a cube of inclusions of ideals of the ring $R=T^3 \K$:
$$
\xymatrix@ur@!0@C+10pt{
&\epsilon_1 R \ar@{-}[rr]\ar@{-}'[d][dd]
 & & R \ar@{-}[dd]
\\
\epsilon_1 \epsilon_2 R \ar@{-}[ur]\ar@{-}[rr]\ar@{-}[dd]
& & \epsilon_2 R \ar@{-}[ur]\ar@{-}[dd]
\\
& \epsilon_1 \epsilon_3 R \ar@{-}'[r][rr]
& & \epsilon_3 R
\\
\epsilon_1 \epsilon_2 \epsilon_3 R \ar@{-}[rr]\ar@{-}[ur]
& & \epsilon_2 \epsilon_3 R \ar@{-}[ur]
}
\eqno (7.13)
$$
which on the level of manifolds corresponds to a cube of inclusions
of various ``vertical bundles''; but only the last and smallest
of these vertical bundles, $\epsilon_1 \epsilon_2 \epsilon_3 TM$,
is a vector bundle over $M$. 
For general $k$,
the first two exact sequences of (7.7) generalize, for $j=1,\ldots,k$, to
$$ \matrix{
0 & \to &  \epsilon_j R & \to &  R & {\buildrel p_j \over\to} &
\K[\epsilon_1,\ldots,\hat \epsilon_j,\ldots,\epsilon_k] &  \to &  0 \cr}
\eqno (7.14)
$$
where a ``hat'' means ``omit this coefficient''.
The general form of the third sequence from (7.7) is, for any
choice $1 \leq j_1 < \ldots < j_\ell \leq k$,
$$ \matrix{
0 & \to &
\bigcap_{i=1}^\ell \epsilon_{j_i} R  & \to &
R &  {\buildrel \prod p_{j_i} \over \longrightarrow}
\prod_{i=1}^\ell \K[\epsilon_1,\ldots,\hat \epsilon_{j_i},\ldots,\epsilon_k] 
& \to &  0. \cr}
\eqno (7.15)
$$
In particular, for the ``maximal choice'' $j_i=i$, $i=1,2,\ldots,k$, we get
$$ \matrix{
0& \to & \epsilon_{1} \ldots \epsilon_{k} R &\to& R & \to& 
\prod_{i=1}^k \K[\epsilon_1,\ldots,\hat \epsilon_j,\ldots,\epsilon_k]& \to& 0.
\cr}
\eqno (7.16)
$$
On the level of manifolds, (7.16) corresponds to the ``most vertical'' bundle
$$ \matrix{
\epsilon_1 \ldots \epsilon_k
TM &  \to & T^k M &  \to & \prod_{i=1}^k T^{k-1} M . \cr}
$$
There are also various projections $T^k \K \to T^\ell \K$ for all
$\ell=0,\ldots,k-1$; their kernels are certain sums of the ideals
considered so far. On the level of manifolds, this corresponds
to an iterated fibration associated to the projections $T^k M \to T^\ell M$.
In particular, for $\ell=0$, we obtain the ``augmentation ideal'' of $T^k \K$
containing
 all linear combinations of the $\epsilon_i$ ``without constant term''.
It corresponds to the fibers of $T^k M$ over $M$.

\ssk {\sl (D) Automorphisms, endomorphisms.} 
According to Diagram (7.12), $T^3 M$ can be seen in three different ways
as the second order tangent bundle of $TM$.
Thus we have three different versions of the canonical flip $\kappa$ (Section
7.3 (D)), corresponding to the canonical action of the three transpositions
$(12)$, $(13)$ and $(23)$ on $T^3M$.
The action of $(12)$ may also be interpeted as $T_{\eps_3} \kappa$ with
$\kappa :TTM \to TTM$ as in 7.3 (D).
 Altogether this gives us a canonical action of the permutation
group $\Sigma_3$ by $\K$-diffeomorphisms of $T^3 M$. In a chart, this action
is given by permuting the symbols $\eps_1,\eps_2,\eps_3$. Similarly,
for general $k$ we have an action of the permutation group $\Sigma_k$ on
$T^k M$. 

The action $l:\K \times TM \to TM$ gives rise to an action
$T^k l:T^k \K \times T^{k+1} M \to T^{k+1} M$ which is not invariant
under the $\Sigma_k$-action, but as for $k=2$ one can put these actions
together to define some invariant action -- see Section 20.4 for more
details.

\msk \nin
{\bf 7.4.} {\sl Transition functions and 
chart representation of higher order tangent maps.}
For a good understanding of the structure of $T^k M$
 we will need to know the structure of the transition functions
$T^k \phi_{ij}$ of $T^kM$, and more generally, the structure
of higher order tangent maps $T^k f$.
For $k=1$, recall that, for a smooth map $f:V \supset U \to W$ (which 
represents a smooth map $M \to N$ with respect to charts), we have
$$
Tf(x+\eps v)=f(x) + \eps df(x) v.
\eqno (7.17)
$$
Using the 
rule on partial derivatives (Lemma 1.6 (iv)), we may calculate
$d(Tf)$ and hence $TTf$. One gets
$$
\eqalign{
T^2 f(x+\epsilon_1 v_1+\epsilon_2 v_2+\epsilon_1 \epsilon_2 v_{12}) &=
f(x)+ \epsilon_1 df(x) v_1+ \epsilon_2 df(x)v_2 + \cr
& \quad \quad
\epsilon_1 \epsilon_2 \bigl( df(x) v_{12} + d^2f(x) (v_1,v_2)\bigr). \cr}
\eqno (7.18)
$$
Iterating this procedure once more, we get
$$
\eqalign{
 & T^3 f(x+\epsilon_1 v_1+\epsilon_2 v_2+\epsilon_3 v_3+
 \epsilon_1 \epsilon_2 v_{12} +
\epsilon_3 \epsilon_1 v_{13} + \epsilon_3\epsilon_2 v_{23} +
\epsilon_1 \epsilon_2 \epsilon_3 v_{123}) 
 \cr
& = 
f(x)+ \sum_{j=1,2,3} 
\epsilon_j df(x) v_j +
\sum_{1 \leq i < j \leq 3}
\epsilon_{i} \epsilon_{j} ( df(x) v_{i j} + d^2f(x)(v_{i},v_{j})) + \cr
& \quad \quad \quad \quad
\epsilon_1 \epsilon_2 \epsilon_3 \cdot  \bigl( df(x) v_{123} +
d^2f(x)(v_1,v_{23})+ d^2f(x)( v_2,v_{13})+  \cr
& \quad \quad \quad \quad \quad \quad \quad \quad
d^2f(x)( v_3,v_{12})+ d^3f(x)(v_1,v_2,v_3) \bigr).    
\cr}
$$
In order to state the general result, we introduce a multi-index notation
(cf.\ also Appendix MA, Sections MA.1 and MA.3). 
For  multi-indices $\alpha, \beta \in I_k:=\{ 0, 1 \}^k$ we define
$$
\epsilon^\alpha := \epsilon_1^{\alpha_1} \cdots \epsilon_k^{\alpha_k},
\quad \quad | \alpha | :=\sum_i \alpha_i.
$$
Moreover, we fix a total ordering on our index set $I_k$; for the moment,
any ordering would do the job, but for later purposes, we agree to take
the usual lexicographic ordering on $I_k=\{ 0, 1 \}^k$.
Then a general element in the fiber over $x$ can be written in the
form $x + \sum_{\alpha > 0} \epsilon^\alpha v_\alpha$
 with $v_\alpha \in V$, i.e., our chart domain is written in the form
$$
T^k U = \{ x + \sum_{\alpha \in I_k \atop \alpha > 0} \epsilon^\alpha v_\alpha|
\, x \in U, \, \forall \alpha: \, v_\alpha \in V \}.
$$
A {\it partition} of $\alpha \in I_k$ is an $\ell$-tuple
$\Lambda= (\Lambda^1,\ldots,\Lambda^\ell)$ with $\Lambda^i \in I_k$ such that
$0<\Lambda^1<\ldots<\Lambda^\ell$ and $\sum_i \Lambda^i = \alpha$; the
integer $\ell=\ell(\Lambda)$ is called the {\it length} of the partition. 
The set of
all partitions (of length $\ell$) of $\alpha$ is denoted by ${\cal P}(\alpha)$
(resp.\ by ${\cal P}_\ell(\alpha)$).

\Theorem 7.5. The higher order tangent maps $T^k f:T^k U \to T^k W$
of a $C^k$-map $f:U \to W$ are expressed in terms
of the usual higher differentials $d^j f:U \times V^j \to W$ 
via the following formula:
$$
\eqalign{
T^k f(x+\sum_{\alpha >0} \epsilon^\alpha v_\alpha) & =f(x)+
\sum_{\alpha \in I_k \atop \alpha > 0} \epsilon^\alpha
\bigl(\sum_{\ell =1}^{|\alpha|} \sum_{\Lambda \in {\cal P}_\ell(\alpha) }
d^\ell f(x)(v_{\Lambda^1},\ldots,v_{\Lambda^\ell}) \bigr)
\cr
& =f(x)+
\sum_{j=1}^k \sum_{\alpha \in I \atop |\alpha|=j} \epsilon^\alpha
\bigl(\sum_{\ell=1}^j \sum_{\Lambda^1 + \ldots + \Lambda^\ell = \alpha
\atop \Lambda^1 < \ldots < \Lambda^\ell }
d^\ell f(x)(v_{\Lambda^1},\ldots,v_{\Lambda^\ell}) \bigr).
\cr}
$$

\Proof. The claim is proved by induction on $k$.
 For $k=1$, this is simply formula (7.17).
Assume now that the formula holds for $k \in \N$.
We know that the iterated tangent map
$\tilde f:= T^k f:T^k U \to T^k W$ is smooth over $T^k \K$
and that $\tilde f|_U$ coincides, via the zero sections, with $f$, i.e.,
$\tilde f(x)=f(x) \in W$ for all $x \in U$.
In the following proof  we will use only these properties of $\tilde f$.
As in the proof of Theorem 6.8 (2), a second order Taylor expansion gives,
for $\eps$ being one the $\eps_i$,
$$
\tilde f(x+\eps v)=f(x)+ \eps d\tilde f(x)v + \eps^2 a_2(x,v) + \eps^2
R_2(x,v,\eps) = f(x) + \eps \partial_v f(x),
$$
which, of course, is again (7.17). Note that we used $d\tilde f(x)v=
\partial_v \tilde f(x) = \partial_v f(x) = df(x) v$, which holds
since $x,v \in V$ and $\tilde f|_U = f$.
 Now we are going to repeat this argument,
using that, if $\tilde f$ is $C^3$ over $TT\K$, then all maps $\partial_u f$,
$u \in TTU$, are $C^2$ over $TT\K$ and hence also over $T\K$:
$$
\eqalign{ \tilde f(x+\eps_1 v_1 + \eps_2 v_2 & \, + \eps_1 \eps_2 v_{12})  =
\tilde  f(x+\eps_1 v_1 + \eps_2 (v_2 + \eps_1 v_{12})) \cr
&= \tilde f(x + \eps_1 v_1) + \eps_2 (\partial_{v_2 + \eps_1 v_{12}} \tilde
f)(x+\eps_1 v_1) \cr
& =f(x) + \eps_1 \partial_{v_1} f(x) +  \eps_2
\partial_{v_2 + \eps_1 v_{12}} \tilde
f(x) + \eps_1 \eps_2 \partial_{v_1} \partial_{v_2 + \eps_1 
v_{12}} \tilde f(x) \cr
& = f(x) + \eps_1 \partial_{v_1} f(x) +  \eps_2 \partial_{v_2} 
f(x) + \eps_1 \eps_2 (\partial_{v_{12}} + \partial_{v_1} \partial_{v_2})f(x).
 \cr}
$$
The proves again (7.18), and again the resulting expression depends only
on $f|_U$. Now, for the general induction step, we let
$u=x+v = x + \sum_{\alpha \in I_{k+1}} \eps^\alpha v_\alpha  \in T^{k+1} U =
T(T^k U)$ and decompose $u=x+ v' + \eps_{k+1} v''$ correspondingly. 
We let also $\omega=(0,\ldots,0,1) \in I_{k+1}$ and identify $I_k$ with the
subset of $\alpha \in I_{k+1}$ such that $\alpha_{k+1}=0$. Then $\eps^\omega =
\eps_{k+1}$, and $I_{k+1}$ is the disjoint union of $I_k$ and $I_k + \omega$.
By induction, 
$$
\eqalign{\tilde f (& x+v)  =
\tilde f\big(x+\eps^\omega v_\omega + \sum_{\alpha \in I_k \atop \alpha >0}
\eps^\alpha (v_\alpha + \eps^\omega v_{\alpha + \omega}) \big) 
\cr
& =\tilde f(x+\eps^\omega v_\omega)+
\sum_{\alpha \in I_k \atop \alpha > 0} \epsilon^\alpha
\sum_{\ell=1}^{|\alpha|} \sum_{\Lambda \in {\cal P}_\ell(\alpha) } 
(\partial_{v_{\Lambda^1} + \eps^\omega v_{\Lambda^1 + \omega}} \cdots
\partial_{v_{\Lambda^\ell} + \eps^\omega v_{\Lambda^\ell + \omega}} \tilde f)
(x+\eps^\omega
v_{\omega}) 
\cr
& =
f(x) +\eps^\omega \partial_{v_\omega} f(x) +
\sum_{\alpha \in I_k \atop \alpha > 0} \epsilon^\alpha
\sum_{\ell=1}^{|\alpha|} \sum_{\Lambda \in {\cal P}_\ell(\alpha) }
(\partial_{v_{\Lambda^1} + \eps^\omega v_{\Lambda^1 + \omega}} \cdots
\partial_{v_{\Lambda^\ell} + \eps^\omega v_{\Lambda^\ell
 + \omega}} \tilde f)(x) 
\cr
& \quad \quad \quad \quad \quad   + 
\eps^\omega \bigl( \partial_{v_\omega} 
\sum_{\alpha \in I_k \atop \alpha > 0} \epsilon^\alpha
\sum_{\ell=1}^{|\alpha|} \sum_{\Lambda \in {\cal P}_\ell(\alpha) } 
(\partial_{v_{\Lambda^1} + \eps^\omega v_{\Lambda^1 + \omega}} \cdots
\partial_{v_{\Lambda^\ell} + \eps^\omega v_{\Lambda^\ell + \omega}} 
\tilde f) \bigr)(x)  .
\cr}
$$ 
Next, we use linearity of the map $y \mapsto \partial_y f$ in order to
expand the whole expression. Fortunately, due to the relation $(\eps^\omega)^2=
\eps_{k+1}^2=0$, many terms vanish, and the remaining sum is a sum over
partitions $\Omega$ of elements $\beta \in I_{k+1}$ of
the following three types:

\ssk
\item{(a)} 
$\Omega$ is a partition $\Lambda$ of some $\alpha \in I_k$,
\item{(b)}
$\Omega$ is obtained from a partition $\Lambda$ of $\alpha \in I_k$
by adjoining $\omega$ as last component (then $\ell(\Omega)=\ell(\Lambda)+1$),
\item{(c)}
$\Omega$ is obtained from a partition $\Lambda$ of $\alpha \in I_k$
by adding $\omega$ to some component of $\Lambda$ 
(then $\ell(\Omega)=\ell(\Lambda)$).

\ssk
\nin 
But it is clear that every partition $\Omega$ of any element 
$\beta \in I_{k+1}$ is exactly of one of the preceding three types, and
therefore $T^{k+1} f(x+v)$ is of the form given in the claim.
\qed

\nin As a consequence of the theorem, one can get another proof of the
existence of a canonical action of the permutation group $\Sigma_k$
on $T^k M$: in a chart, $\Sigma_k$ acts via its canonical action on
$\{0,1\}^k$, and the chart formula shows that $T^k$ commutes with this
action because of the symmetry of the higher differentials. --
Another consequence of the preceding proof is the following uniqueness
result on the ``scalar extension of $f$ by $T^k f$'':

\Theorem 7.6.
Assume $f:M \to N$ is smooth over $\K$ and 
$\tilde f:T^k M \to T^k N$ is a $T^k \K$-smooth extension of $f$ in
the sense that  $\tilde f \circ z_{T^kM} =
z_{T^k N} \circ f$:
$$
\matrix{T^k M & {\buildrel \tilde f \over \to} & T^k N \cr
\uparrow & & \uparrow \cr
M &{\buildrel f \over \to} & N \cr}
$$
Then $\tilde f = T^k f$.

\Proof. The preceding proof has shown that $\tilde f$ agrees with $T^k f$
in any chart representation, whence $\tilde f = T^k f$.
\qed

\vfill \eject 

\subheadline{8. Scalar extensions. III: Jet functor and
truncated polynomial rings}

\msk \nin {\bf 8.1.} {\sl Jet rings and truncated polynomial rings.}
Recall the canonical action of the permutation group
$\Sigma_k$ by automorphisms of the higher order tangent ring
$T^k \K = \K[\epsilon_1,\ldots,\epsilon_k]$.
The subring fixed under this action will be called the {\it ($k$-th
order) jet ring of $\K$}, denoted by $J^k \K$.
For $k=0,1$ we have $J^0 \K = \K$, $J^1 \K = T\K = \K[\epsilon]
= \K[x]/(x^2)$, and for $k=2,3$,
$$
\eqalign{
J^2 \K & = \{ a_0 + \epsilon_1 a_1 + \epsilon_2 a_1 +
\epsilon_1 \epsilon_2 a_2 \vert \, a_0,a_1,a_2 \in \K \} =
\K \oplus \K (\epsilon_1 + \epsilon_2)\oplus \K \epsilon_1 \epsilon_2, \cr
J^3 \K & = \K \oplus \K(\epsilon_1 + \epsilon_2+\epsilon_3)
\oplus \K (\epsilon_1 \epsilon_2 + \epsilon_2 \epsilon_3 + \epsilon_1
\epsilon_3) \oplus \K \epsilon_1 \epsilon_2 \epsilon_3. \cr}
$$
For general $k$, we introduce the notation
$$
\delta:=\delta_k:= \sum_{i=1}^k \epsilon_i, \quad \quad
\delta^{(2)}:= \sum_{1\leq i<j \leq k} \epsilon_i \epsilon_j,
\quad \quad 
\delta^{(j)}:=\delta_k^{(j)} := 
\sum_{1 \leq i_1 < \ldots < i_j \leq k} \epsilon_{i_1} \cdots \epsilon_{i_j}
\eqno (8.1)
$$
(the $j$-th elementary symmetric polynomial in the variables
 $\epsilon_1,\ldots,
\epsilon_k$) which defines a basis of the $\K$-module $J^k \K$:
$$
J^k \K = \K \oplus \K \delta \oplus
\K \delta^{(2)} \oplus 
 \ldots \oplus \K \delta^{(k)}.
\eqno (8.2)
$$
An element of $T^k \K$ and of $J^k \K$ is invertible iff its lowest
coefficient is invertible; therefore $J^k \K$ has a dense unit
group if $\K$ has a dense unit group.
The powers of $\delta$ are
$$
\eqalign{
\delta^2 & =(\sum_i \epsilon_i)^2 = \sum_{i,j} \epsilon_i \epsilon_j =
2 \sum_{i< j} \epsilon_i \epsilon_j = 2 \delta^{(2)}, \ldots \cr
\delta^j & = j ! \delta^{(j)}, \quad \ldots \quad,
\delta^k  = k ! \epsilon_1 \cdots  \epsilon_k, \quad 
\quad \delta^{k+1} = 0.\cr}
$$
More generally, we have the multiplication rule
$$
\delta^{(i)} \cdot \delta^{(j)} = \pmatrix{i+j \cr i \cr} \delta^{(i+j)}.
$$
If the integers $2,\ldots,k$ are invertible in $\K$, 
we can take $1,\delta,\delta^2,\ldots,\delta^k$ as a new $\K$-basis
in $J^k \K$:
$$
J^k \K = \K \oplus \K \delta \oplus \K \delta^2 \oplus 
\K \delta^3 \ldots \oplus \K \delta^k.
$$
For any $\K$, the ring
$$
\K[\delta_k] := \K[x]/(x^{k+1})
$$
is called a {\it truncated polynomial ring (over $\K$)}.
A $\K$-basis of $\K[x]/(x^{k+1})$ is given by the classes of the
polynomials $1,x,x^2,\ldots,x^k$. We denote these classes
by $\delta^0,\delta,\delta^2,\ldots,\delta^k$;
then we have the relation $\delta^{k+1}=0$, and elements are
multiplied by the rule
$$
\sum_{i=0}^k a_i \delta^i \cdot 
\sum_{j=0}^k b_j \delta^j =
\sum_{\ell=0}^k (\sum_{i+j=\ell} a_i b_j) \delta^\ell.
$$
These relations show that, if $2,\ldots,k$ are invertible in $\K$,
then the jet ring $J^k \K$ and the truncated polynomial ring
$\K[\delta_k]$ are isomorphic.
For general $\K$, the unit element of $\K[\delta_k]$ 
is $1 =\delta^0$, and, given $a=\sum_i a_i \delta^i$,
the condition $a b=1$ for $b=\sum_j b_j \delta^j$ is
a triangular system of linear equations whose matrix has a 
diagonal of coefficients all equal to $a_0$; this system can
be solved if and only if $a_0$ is invertible. In other words,
the unit group of the ring $J^k \K$ is the set of truncated
polynomials having an invertible lowest coefficient.
Clearly, this group is dense in $J^k \K$ if the unit group
of $\K$ is dense in $\K$.

\Theorem 8.2. If $M$ is a manifold over $\K$, we denote by
$$
J^k M := (T^k M)^{\Sigma_k}
$$  
the fixed point set of the canonical action of the permutation group
$\Sigma_k$ on $T^k M$. Then $J^k M$ is a subbundle of the bundle $T^k M$ 
over $M$, and it is a smooth manifold over the 
$k$-th order jet ring $J^k \K$. If $f:M \to N$ is smooth
over $\K$, then the restriction
$J^k f:J^k M \to J^k N$ is well-defined and is smooth over $J^k \K$.
Summing up, $J^k$ can be seen as the functor of scalar extension from
$\K$ to $J^k \K$.
If the integers $2,\ldots,k$ are invertible in $\K$, then
the corresponding statements are true with respect to the truncated
polynomial ring $\K[x]/(x^{k+1})$ instead of $J^k \K$.

\Proof. 
First of all we prove that $J^k M$ is a subbundle of $T^kM$. In a bundle chart 
over $U \subset M$, the fiber over $x$  is represented by
$$
(J^k U)_x = \{ x + \sum_{i=1}^k\epsilon_i v_1 + \sum_{i<j}
\epsilon_i \epsilon_j  v_2 +
\ldots + \epsilon_1 \cdots \epsilon_k v_k  \,  \vert \,
 v_1,\ldots,v_k \in V \}.
\eqno (8.3)
$$
This is a $\K$-submodule and corresponds to
the chart representation of a submanifold provided
we can find a complementary submodule.
In case $k=2$, this amounts to find a complementary submodule
of the diagonal in $V \times V$: for instance, we may take
$0 \times V$ (if $2$ in invertible in $\K$, we may also take
the antidiagonal); similarly for general $k$.
Thus $J^k M$ is a submanifold of $T^k M$ in the sense of 2.3;
from (8.3) it is then seen to be a subbundle. 
If $f:M \to M$ is a smooth map,
then $T^k f:T^k M \to T^k N$ is, due to the symmetry of higher
differentials, compatible with the action of
the permutation group $\Sigma_k$, and hence gives rise to a well-defined
and $\K$-smooth map
$$
J^k f: J^k M \to J^k N, \quad u \mapsto T^k f(u)
$$
such that we have the functorial rules $J^k(f \circ g)=J^k f \circ
J^k g$, $J^k \id_M = \id_{J^k M}$.
Smoothness of $J^k f$ over $J^k \K$ follows simply by restriction:
we know that
$T^k f:T^k U \to T^kW$ is smooth over $T^k \K$ (Theorem 7.2);
we write Condition (1.1) for differentiability over $T^k \K$ and
restrict to $\Sigma_k$-invariants and thus get Condition (1.1) for
differentiability ofer $J^k \K$. By induction, we see that $J^k f$
is of class $C^\ell$ over $J^k \K$ for all $\ell$.

Finally,
the last statement follows from the isomorphism of topological rings
$J^k \K \cong \K[x]/(x^{k+1})$ explained in Section 8.1.
\qed

\nin
The bundle $J^k M$ over $M$ is called the {\it $k$-th order jet bundle
of $M$} -- see Section 8.8  below for a discussion with more conventional
definitions of jets.

\Theorem 8.3. 
If $G$ is a Lie group over $\K$, then, for $k \in \N$,
$T^k G$, resp. $J^k G$ are Lie groups  over the rings $T^k \K$,
resp. over $J^k \K$,
and the corresponding Lie algebra is the scalar extension of
the Lie algebra $\g$ of $G$ by the ring $T^k \K$, resp. by $J^k \K$. 
Similarly,
if $M$ is a symmetric space over $\K$ with base point $o$, 
then, for $k\in \N$, $T^k M$, resp. $J^k M$,
are symmetric spaces over $T^k \K$, resp. over $J^k \K$, with
base point being the origin in the fiber over $o$.
The associated Lie triple system is the scalar extension
of the Lie triple system associated to $(M,o)$
by the ring $T^k \K$, resp. by $J^k \K$.
The various natural projections and inclusions are homomorphisms
of Lie groups, resp. of symmetric spaces.

\Proof. This is the generalization of Theorem 6.5, and its version for
$T^k \K$ follows by induction from that theorem. The version for $J^k \K$
then follows by restriction to $\Sigma_k$-invariants. \qed

\msk
\nin {\bf 8.4.} {\sl The structure of $J^k \K$ and of $J^k M$.}
Not all properties of the ring $T^k \K$ are invariant under the permutation
group, and hence the structure of $J^k \K$ is somewhat poorer than the
structure of $T^k \K$.

\ssk
{\sl (A) Projections.}
The various exact sequences $0\to \epsilon T^{k-1}\K \to
T^k \K \to T^{k-1} \K \to 0$ which are all of the
form $0\to \epsilon R \to R[\epsilon] \to R \to 0$
(where $\epsilon$ is one of the $\epsilon_i$) induce,
by restriction to invariants under permutations, an exact sequence
$$
0 \to \delta^{(k)} \K \to J^k \K {\buildrel p \over \to} J^{k-1} \K \to 0,
\quad \quad 
p(\sum_{i=0}^k \delta^{(i)} v_i) =  \sum_{i=0}^{k-1} \delta^{(i)} v_i.
\eqno (8.4)
$$
More generally, for all $\ell \leq k$ we have projections
$$
\pi_{\ell,k}: J^k \K \to J^{\ell} \K, \quad
\sum_{i=0}^k \delta^{(i)} v_i \mapsto  \sum_{i=0}^{\ell} \delta^{(i)} v_i
\eqno (8.5)
$$
which all come from the various projections $T^k\K \to T^\ell \K$.
On the space level, $\pi_{k,k+1}$
 corresponds to the {\it jet bundle projection}
$$
\matrix{
 J^{k+1} M & \subset & T^{k+1} M \cr
\downarrow & & \downarrow \cr
J^k M & \subset & T^k M \cr}
\eqno (8.6)
$$
and the projections $\pi_{\ell,k}:J^k M \to J^\ell M$ are given in a chart by
$$
J^{k} M \to J^\ell M, \quad
x+ \sum_{j=1}^{k} \delta^{(j)} v_j \mapsto
x+ \sum_{j=1}^{\ell} \delta^{(j)} v_j.
\eqno (8.7)
$$

{\sl (B) Non-existence of injections.} 
The ``zero-sections'' $z_j:\epsilon_j T^k M \to T^{k+1} M$  
of the vector bundles $p_j:T^{k+1} M \to T^k M$ are, for
$k \geq 1$, not compatible with the action of the symmetric
group: for instance, in a chart we have
$$
z_{k+1}(x + \epsilon_1 v_1 + \ldots + \epsilon_{k} v_{k}) =
x + \epsilon_1 v_1 + \ldots + \epsilon_{k} v_{k} + 0,
$$
hence the image of $z_{k+1}$ does not belong to $J^{k+1} M$.
Therefore, when restricting the vector bundle 
$T^{k+1}M \to T^k M$ to invariants under the symmetric group,
we obtain a ``vector bundle without zero section'', that is,
an {\it affine bundle} $J^{k+1} M \to J^k M$. Thus 
the fiber over $u \in J^k M$,
$$
(V^{k+1,1} M)_u = (J^{k+1} M)_{x+ \sum_{j=1}^{k} \delta^{(j)} v_j} 
= \{ x+ \sum_{j=1}^{k+1} \delta^{(j)} v_j \vert \, v_{k+1} \in V  \} 
\eqno (8.8)
$$
carries a canonical structure of an affine space over $\K$,
but not of a $\K$-module.

\ssk {\sl (C) Ideal structure.} The kernel of the projection
$\pr_{k,\ell}:J^k \K \to J^\ell \K$,
$$
J_{k,\ell}=\{ \sum_{i=\ell+1}^k \delta^{(i)} r_i | \forall i: \, r_i \in \K \},
\eqno (8.9)
$$
is an ideal of $J^k \K$, and instead of a cube of ideals as for $T^k \K$,
now we have a chain of ideals. For $\ell=0$ we get the ``augmentation ideal'',
i.e., the set of elements having lowest coefficient equal to zero;
if $\K$ is a field, this is the unique maximal ideal of $J^k \K$.
For $\ell=k-1$ we get the ideal $\K \delta^{(k)}$ which, on the space level,
corresponds to fibers of the affine bundle $J^k M \to J^{k-1} M$.

\ssk {\sl (D) Automorphisms, endomorphisms.} The only automorphisms or
endomorphisms from $T^k \K$ that suvive in a non-trivial way come from
the action (which has been defined in Section 7.3 (D))
 $l^{(k)}:T^k \K \times T^{k+1} M \to T^{k+1} M$.
In particular, there is a well-defined action $\K \times J^k M \to J^k M$.
In a chart, it is given by
$$
r .(x + \sum_{i=1}^k \delta^{(i)} v_i)= x+\sum_{i=1}^k r^i \delta^{(i)} v_i.
\eqno (8.10)
$$

\msk \nin {\bf 8.5.} {\sl Chart representation of $J^k f$.}
Let $f:M \to N$ be a smooth map. Using charts, $J^k M$ and 
$J^k N$ are represented in the form (8.3), and we may ask
for the chart expression of $J^k f$. For $k=2$, letting
 $v_1=v_2=:w$,  we get from (7.18):
$$
J^2 f(x+ \delta v +\delta^{(2)} w)=
f(x)+ \delta df(x)v +  \delta^{(2)} \bigl( df(x)w + d^2f(x) (v,v) \bigr),
\eqno (8.11)
$$
which implies, with
$\delta^2 = 2 \delta^{(2)}$, if $2$ is invertible in $\K$,
$$
J^2 f(x+\delta v_1+\delta^2 v_2)=
f(x)+ \delta df(x)v_1 + 
\delta^2 \bigl( df(x) v_2 + {1 \over 2} d^2f(x) (v_1,v_1) \bigr).
$$
For $k=3$ we get by restriction from the formula for $T^3 f$
$$
\eqalign{
J^3 f(x + \delta v_1 + & \delta^{(2)} v_2 + \delta^{(3)} v_3)  = \cr
& f(x) + \delta df(x)v_1 +  
\delta^{(2)}\bigl( df(x)v_2 + d^2 f(x)(v_1,v_1) \bigr) +
\cr
& \quad \quad \quad \quad \delta^{(3)} \bigl( 
df(x)v_3+3 d^2 f(x)(v_1,v_2) + d^3 f(x)(v_1,v_1,v_1) \bigr).
\cr}
$$
In particular, if $3$ is invertible in $\K$, we have
$$
\eqalign{
J^3 f(x + \delta v) & =
f(x) + \delta df(x)v  +
{1 \over 2} \delta^2 \, df(x)v  + {1 \over 3!} \delta^3 \, d^3 f(x)(v,v,v).
\cr}
$$

\Theorem 8.6.
If $f:U \to W$ is smooth over $\K$, then $J^k f:J^k U \to J^k W$ is
described in terms of the usual higher order differentials
$d^\ell f:U \times V^\ell \to W$ by
$$
J^k f(x + \sum_{j=1}^k \delta^{(j)} v_j ) =
f(x) +
\sum_{j=1}^k \delta^{(j)} \bigl( \sum_{\ell=1}^j 
\sum_{i_1 + \ldots + i_\ell = j \atop i_1 \leq \ldots \leq i_\ell}
C^j_{i_1,\ldots,i_\ell} d^\ell f(x)(v_{i_1},\ldots,v_{i_\ell}) \bigr),
$$
where $C^j_{i_1,\ldots,i_\ell}$ is the number of decompositions
of a set of $j$ elements into $\ell$
 disjoint subsets containing $i_1,i_2,\ldots,
i_\ell$ elements.  If the integers are invertible
in $\K$, then this can also be written
$$
J^k f(x + \sum_{j=1}^k \delta^{(j)} v_j ) =
f(x) +
\sum_{j=1}^k \delta^{(j)} (\sum_{\ell=1}^j 
\sum_{i_1 + \ldots + i_\ell = j}
{j! \over i_1! \ldots i_\ell ! \ell !} d^\ell f(x)(v_{i_1},\ldots,v_{i_\ell})).
$$
In particular,
$$
\eqalign{
J^k f(x + \delta v)& =f(x) + \sum_{j=1}^k \delta^{(j)} df(x)(v,\ldots,v)\cr
&=f(x) + \sum_{j=1}^k {1 \over j!}\delta^j  df(x)(v,\ldots,v) , \cr}
$$
the latter provided the integers are invertible in $\K$.

\Proof.
We restrict the formula from Theorem 7.5 to $\Sigma_k$-invariants:
 two multi-indices
$\alpha$, $\beta \in \{ 0,1 \}^k$ are conjugate under $\Sigma_k$ if
and only if $|\alpha| = |\beta|$; therefore we  have to take
the special case of Theorem 7.5 where,  for all $\alpha$ with $|\alpha| = j$,
the arguments $v_{\alpha}$ are equal to a given element $v_j \in V$.
Recall also that, by definition,
$\delta^{(j)} = \sum_{|\alpha|=j} \epsilon^\alpha$.
Then the claim follows by counting the number of terms that are equal.

Note that, by definition, the numbers $C^j_{i_1,\ldots,i_\ell}$ are integers, 
and hence the claim
is valid in arbitrary characteristic. See [Mac79, p.\ 22] for another
formula and some remarks on these combinatorial coefficients.
They are partition numbers, and there is no
``explicit formula'' for these constants.
But the index set of the sum in the first formula of the theorem
 also is a set of
 partitions, and these two combinatorial difficulties cancel out if we
take the  sum over {\it all} decompositions $i_1 + \ldots + i_\ell =j$
(with repetitions!) -- instead of counting  partitions of
$\{1,\ldots,j \}$ we count $\ell$-tuples of subsets, the $m$-th set having
$i_m$ elements (so the number of such $\ell$-tuples is
${j! \over i_1 ! \cdots i_\ell !}$),  and then divide
by $\ell!$ (provided the integers are invertible in $\K$)
which gives the second formula of the claim.
%
%
\qed 

\nin One should note that for practical computations the first formula of the
theorem is much more efficient since it has less terms.

\msk \nin
{\bf 8.7.} {\sl The composition rule for higher differentials.}
Assume $f$ and $g$ are smooth mappings of $\K$-modules such that
the composition $g \circ f$ is defined on a non-empty open set.
As an interesting consequence of the explicit formulae from Theorem
8.6 for
the $k$-jets of $f$ and $g$ we get the following explicit formula 
for the higher differentials of $g \circ f$ in terms of those of
$f$ and $g$:
$$
\eqalign{
d^j(g \circ f)(x)& ( v,\ldots ,v)   \cr
= & \sum_{\ell=1}^j 
\sum_{i_1 + \ldots + i_\ell = j \atop i_1 \leq \ldots \leq i_\ell}
C^j_{i_1,\ldots,i_\ell} d^\ell g(f(x))
\Bigl( d^{i_1} f(x)(v_{i_1}),\ldots,
d^{i_\ell} f(x)(v_{i_\ell}) \Bigr) \cr
= & \sum_{\ell=1}^j 
\sum_{i_1 + \ldots + i_\ell = j}
{j! \over i_1! \ldots i_\ell ! \ell !} d^\ell  g(f(x))
\Bigl(d^{i_1} f(x)(v_{i_1}),\ldots,
d^{i_\ell} f(x)(v_{i_\ell})\Bigr), \cr}
$$
where for simplicity we wrote $d^{i_m} f(x)(v_{i_m})$ instead of
 $d^{i_m} f(x)(v, \ldots, v)$.
In order to prove this formula, we write  the explicit formulae
for both sides of the equality $J^k(g\circ f)(x+\delta v) = 
J^kg(J^k f(x+\delta v))$ and compare terms having coefficient
$\delta^{(j)}$. 
The second version of the composition
rule (8.19) (in the finite-dimensional real case) can be found in
[Chap03, p.\ 172, ``formule de Faa de Bruno''] and
[Ho01, Satz 60.16]; the proofs given there are different, using Taylor
expansions in the classical sense.
See also [Mac79, p.\ 23, Ex.\ 12] for the one-dimensional version of this
formula.  

\msk \nin {\bf 8.8.}
{\sl Relation with classical definitions of jets.}
Usually, jets are defined via an equivalence relation (cf., e.g.,
[KMS93, Ch.\ IV] or [ALV91, p.\ 132]):
one says that two maps $f,g:M \to N$ 
{\it have the same $k$-jet at $x \in M$} if, in a chart, all derivatives
$d^j f(x)$ and $d^j g(x)$, $j=0,1,\ldots,k$, coincide.
Then we write $f \sim_x^k g$. We have the following relation with our
concepts of $T^k f(x)$ and $J^k f(x)$:
$$
f \sim_x^k g \quad \Leftrightarrow \quad (T^k f)_x = (T^k g)_x \quad
\Rightarrow \quad (J^k f)_x = (J^k g)_x,
$$
and the last implication becomes an equivalence if the integers are invertible
in $\K$. In fact, it is clear that the first two properties are equivalent
since Theorem 7.5 permits to calculate $T^k f(x)$ from the higher order
derivatives and vice versa.  If $J^k f(x)=J^k g(x)$,
it follows from Theorem 8.5
that we can recover the values of $d^j f(x)(v,\ldots,v)$ 
(for $j=0,1,\ldots,k$), and thus,
if $2,\ldots,k$ are invertible in $\K$, it follows by polarisation that
the differentials up to order $k$ at $x$ agree, i.e.
$T^k f(x)=T^k g(x)$ and  $f \sim_x^k g$. 
In particular, in the real  case our definition
coincides with the classical one. 
(For the case $k=2$, cf. also [Bes78, p.\ 20].) 
Another possible interpretation of the bundles $J^kM$ is 
via $k$-jets of {\it curves}: 
if $\alpha:I \to M$ is a curve defined on an open neighborhood $I$
of $0$ in $\K$ such that $\alpha(0)=x$, then
$(J^k \alpha)_0(\delta 1)$ is an element of $(J^k M)_x$, and if the integers
are invertible in $\K$, then every
element $u \in (J^k M)_x$ is obtained in this way: in fact,
writing, in a chart,
 $u=u_0 + \delta u_1 + \ldots + \delta^{(k)} u_k \in (J^k M)_x$, we let
$\alpha(t) = \sum_{\ell=0}^k {1 \over \ell !} t^\ell u_\ell$.
According to Theorem 8.6, we then have $J^k \alpha(0 + \delta 1)= u$. 
Summing up, in case of characteristic zero,
 both functors $T^k$ and $J^k$ may then be seen as ``faithful
representations of the abstract jet functor'', and hence are equivalent.
However, in practice there are important differences:

\ssk
\item{(a)} At a first glance, the functor $J^k$ seems preferable, since we get
rid of the big amount of redundant  information contained in $T^k$:
the number of variables is reduced  from $2^k$ to $k$.
\item{(b)} For a fixed pair $k \geq \ell$,
there is just {\it one} jet projection $\pi_{\ell,k}:J^kM \to J^\ell M$,
and thus $J,J^2,J^3,\ldots$ forms a projective system so that we can pass
to the projective limit $J^\infty$ (see Chapter 32). The ring $J^\infty \K$
is essentially the ring of {\it formal power series over $\K$}.
For $T^k$, there is no nice projective limit because of the more complicated
cube structure of the projections.
\item{(c)} On the other hand, 
the functor $T^k$ behaves well with respect to induction procedures because
$T^{k+1} = T \circ T^k$, whereas $J^{k+1}$ and $J^1 \circ J^k$ yield
quite different results. 
\item{(d)} The non-existence of zero-sections makes the functors $J^k$
conceptually more difficult than the functors $T^k$.
\item{(e)} Finally, in positive characteristic, $T^k$ may contain strictly
more information than $J^k$. 

\ssk \nin For all of these reasons, in the sequel we will work, as long
as we can, with $T^k$, and restrict only at the very last stage to
$\Sigma_k$-invariants, before we take the projective limit --
see Chapter 32.

\vfill \eject

\def\title{III. SECOND ORDER DIFFERENTIAL GEOMETRY} 
\def\rightheadline{\hfil{\eightrm \title}\hfil\tenrm\folio}

\sectionheadline{III. Second order differential geometry}

This part, although strongly motivated by Section 7.3 on the second order
tangent bundle $TTM$, is independent of Part II. In particular, it can
be read without having to bother about manifolds over rings, for readers
who prefer fields to rings. Then notation like $V \oplus \eps V$ or
$V \oplus \delta V$ should simply be interpreted as another way of
writing $V \oplus V$, with a formal label added in order to distinguish
both copies of $V$. 

\subheadline
{9. The structure of the tangent bundle of a
vector bundle}

\msk \nin {\bf 9.1.} {\sl Transition functions of the tangent bundle
  of a vector bundle.}
Recall from Section 4.7 the transition maps for $TTM$.
Let us make a similar calculation for
$TM$ replaced by a general vector bundle $p:F \to M$.
In order to keep notation simple, we use similar
conventions as in [La99]: a  bundle chart domain of $F$ is denoted
by $U \times W \subset V \times W$ and a bundle chart domain
of $TF$ by $T(U \times W)=U \times W \times V \times W
\subset V \times W \times V \times W$ with typical element
denoted by $(x,w;x',w')$ or $x+w+\epsilon(x'+w')$.
The transition maps of $F$ are, with the notation of Section 3.3,
$$
\tilde g_{ij}(x,w)=(\phi_{ij}(x),g_{ij}(x)w) = (\phi_{ij}(x),g_{ij}(x,w)).
$$
Using the rule on partial derivatives (Lemma 1.6 (iv)), 
it follows that the transition maps of $TF$ are
$$
\eqalign{
T\tilde g_{ij}(x,w;x',w') & =
(\tilde g_{ij}(x,w),d(\tilde g_{ij})(x,w)(x',w')) \cr
& =(\phi_{ij}(x),g_{ij}(x)w,d\phi_{ij}(x)x',
g_{ij}(x)w' + d_1( g_{ij}) (x,w)x') \cr}
\eqno (9.1)
$$
which may also be written
$$
\eqalign{
T\tilde g_{ij}(x+w+\epsilon(x'+w')) 
& =\phi_{ij}(x)+g_{ij}(x)w +  \cr
& \quad \quad   \epsilon(d\phi_{ij}(x)x' +
g_{ij}(x)w' +d_1(g_{ij}) (x,w)x'). \cr}
$$

\msk \nin
{\bf 9.2.} {\sl $TF$ as a fiber bundle over $M$.}
Let $p:F \to M$ be as above.
Denote by $Tp:TF \to TM$ its tangent map and by
$\pi:TF \to F$ the canonical projection. In bundle charts,
these projections are given by
$$
Tp(x+w+\epsilon(x'+w'))=x+\epsilon x', \quad \quad 
\pi(x+w+\epsilon(x'+w'))=x+w.
$$
The projections $p \circ \pi$ and $\pi_{TM} \circ Tp:TF \to M$
agree: we have the following commutative diagram of
vector bundle projections:
$$
\matrix{ & & TF &  & \cr
& \swarrow & & \searrow & \cr
F & & & & TM \cr
& \searrow & & \swarrow & \cr
& & M & & \cr}
$$

\ssk \nin {\bf 9.3.} {\sl $TF$ as a bilinear bundle.}
The fiber $(TF)_x$ over $x \in M$ does not carry a canonical
$\K$-module structure: in a chart, $(TF)_x$ is represented by
$E:= W \times V \times W$,
and the transition maps, calculated above, are {\it non-linear} in
$(w,x',w')$ (the last term
in (9.1) depends bilinearly on $(w,x')$ and not linearly).
Hence $TF \to M$  {\it is not a vector bundle} -- so the question 
arises: what is it then? Formula (9.1) contains the answer:
for fixed $x$, the transition functions can be written in the form
$$
f(w,x',w')=(h_1(w),h_2(x'),h_3(w' + b(w,x'))),
\eqno (9.2)
$$
with (continuous) invertible linear maps
$$
h_1 = h_3 = g_{ij}(x) \in \Gl(W), \quad
h_2 = d\phi_{ij}(x) \in \Gl(V),
$$
and a (continuous) bilinear map
$$
b:=b_x: W \times V \to W, \quad (w,x') \mapsto
g_{ij}(x)^{-1} \, (d_1(g_{ij})(x,w)x').
\eqno (9.3)
$$
In the terminology of Appendix BA, Section BA.2, this means that the
transition functions for fixed $x \in M$ belong to the {\it general bilinear
group} $\Gm^{0,2}(E)$ of $E=V_1 \times V_2 \times V_3$ (where $V_1=V_3=W$,
$V_2=V$). 
Every pair of bundle charts (say, corresponding to indices
$i,j \in I$) defined at $x$ identifies the fiber
$(TF)_x$ with $E$, in two different ways, such that transition is given by
$f$:
$$
\matrix{  &  & (TF)_x &  &  \cr
          & \phantom{i} {\buildrel (i) \over \swarrow} & & 
{\buildrel (j) \over \searrow} \phantom{ j} &  \cr
E & & {\buildrel f \over \longrightarrow } & & E \cr}
$$
The $\K$-module structure on $(TF)_x$ induced from  the chart indexed by
$(i)$ is not an intrinsic feature since it is different from the one
induced from the chart $(j)$. But the {\it class} of all $\K$-module
structures induced from all possible bundle charts is an intrinsic feature
of the fiber $(TF)_x$. As we have seen, the transition map $f$ belongs
to the group $\Gm^{0,2}(E)$, and hence the orbit under this group of
the $\K$-module structure induced by chart $(i)$ contains the one induced
by chart $(j)$. Therefore this orbit is an intrinsic structure of the fiber
$(TF)_x$. In the terminology of Appendix BA, this can be rephrased by saying
that the {\it set of linear structures on $(TF)_x$ that are bilinearly 
related to linear structures induced from charts} is an intrinsic feature.
These structures can all be described in a very explicit way: 
the linear part $h=(h_1,h_2,h_3)$ of the map $f$ 
preserves the linear structure of $E$, and therefore we may forget this 
linear part. Only the bilinear part (9.3) really is responsible for
non-linearity of the action. In other words, all linear structures $L^b$ that
are bilinearly related to the linear structure $L^0$
coming from the ``original chart'' are given by push-forward of $L^0$
via a map of the form
$$
f_b:W \times V \times W \to W \times V \times W, \quad
(w,x',w') \mapsto (w,x',w'+b(w,x'))
\eqno (9.4)
$$
with bilinear $b=b_x$. Addition and multiplication by scalars in the
$\K$-module $(E,L^b)$ are then given by the explicit Formulae (BA.1).
Using the terminology of Appendix BA, we may summarize:

\Proposition 9.4.
If $F$ is a vector bundle over $M$, then any two $\K$-module
structures on the fibers of $TF$ over $M$ induced by
bundle charts are bilinearly related to each other
(via continuous bilinear maps).
If $F=TM$, then they are related to each other via continuous
{\rm symmetric} bilinear maps.

\Proof. Only the claim on $TM$ remains to be proved.  In this case, (9.3) gives
$$
b_x(w,x')=d\phi_{ij}(x)^{-1} \, d^2\phi_{ij}(x) (w,x'),
\eqno (9.5)
$$
which is a {\it symmetric} bilinear map. 
\qed

\nin
{\bf 9.5.} {\sl ``Dictionary" between bilinear geometry
and differential geometry.}
It follows from the preceding proposition that any intrinsic object
of bilinear geometry (in the sense of Section BA.2) 
on the ``standard fiber'' $W \times V \times W$ of $TF$ over $M$
has an intrinsic interpretation on the fiber bundle $TF$ over
$M$. In the following, we present a table of such structures
together with their differential geometric,
 chart-independent definition.
Here, $z_F:M \to F$ denotes the zero section of $F$
and $z_{TF}:F \to TF$ the zero section of $TF$,
and 
$$
P:\pi \oplus Tp: TF \to F \oplus TM
$$
the vector bundle direct sum of the two projections $\pi:TF \to F$ and
$Tp:TF \to TM$.
A typical element $u \in F$ will be represented in
a chart by $(x,w)$, $w \in W$, and $0_x$ denotes the
zero vector in $F_x$. Notation concerning bilinear algebra
is explained in Appendix BA. 


\bigskip
\settabs 2 \columns
\+
$\quad \quad${\sl bilinear algebra} & {\sl differential geometry} \cr
\msk
\+
$\quad \quad W\times V \times W$, bilinear space  & 
$(TF)_x$, fiber over $x \in M$ \cr
\+
axes: & \cr
\+ 
$\quad \quad W \times 0 \times 0$ (first axis) &
$z_{TF}(F_x) \subset (TF)_x$
\cr
\+
$\quad \quad 0 \times V \times 0$ (second axis) &
$T(z_f)(T_xM) \subset (TF)_x$
\cr
\+
$\quad \quad 0 \times 0 \times W$ (third axis) &
$T_{0_x}(F_x) \subset (TF)_x$ 
\cr
\+
projections: & \cr
\+
$\quad \quad \pr_1:W \times V \times W \to W$ &
$\pi_x:(TF)_x \to F_x$
\cr
\+
$\quad \quad \pr_2:W \times V \times W \to V$ &
$(Tp)_x:(TF)_x\to T_x M$ \cr
\+
$\quad \quad \pr_{12}: W \times V \times W \to W \times V$ &
$P_x=\pi_x \oplus (Tp)_x:(TF)_x \to F_x \oplus T_xM$ \cr
\+
fibrations:& \cr
\+
$\quad \quad \{ w \} \times V \times W
= \pr_1^{-1} (w)$   &
$\pi_x^{-1}(u)$
\cr
\+
$\quad \quad W \times \{ v \} \times W
= \pr_2^{-1} (v)$   &
$(Tp)^{-1}(v)$
\cr
\+
$\quad \quad \{ w \} \times \{v \} \times W
= \pr_{12}^{-1} (w,v)$   &
$P^{-1}_x(u,v)$ 
\cr
\+
special endomorphisms: & \cr
\+
$\quad \quad \epsilon: W \times V \times W \to W \times V \times W$,
 &
$\K[\epsilon]$-module structure on
$T(F_x)$ (fiber of $Tp$)
\cr
\+
$\quad \quad \quad \quad (w,v,z) \mapsto (0,v,w)$ & 
\cr
\+ 
$\quad \quad$if $V=W$:
& case of tangent bundle, $F=TM$:
\cr
\+
$\quad \quad \quad \quad \kappa: (w,v,z) \mapsto (v,w,z)$ &
the ``canonical flip'' $\kappa:TTM \to TTM$ \cr

\bigskip
\nin
Later (Section 10.7) this table will be continued to a table
of {\it non-intrinsic} structures, i.e. structures that depend
on the choice of one linear structure among all the bilinearly
related structures (linear connection).

\msk \nin {\bf 9.6.} {\sl Bilinear bundles.}
A {\it bilinear bundle} is a fiber bundle $F$ over $M$
such that each fiber is a {\it bilinear space} in the sense of Appendix BA.
This means that vector bundles $F_1,F_2,F_3$ over $M$ together with
bundle injections $F_i \to F$ ($i=1,2,3$) are fixed
(called {\it axes of $F$}) such that the axes are {\it linear} subbundles
of $F$ (i.e., subbundles on which transition functions act linearly);
in all chart representations of $F$, the fiber $F_x$ over $x \in M$ is assumed
to be a
direct sum $E=(F_1)_x \oplus (F_2)_x \oplus (F_3)_x$, and
transition maps in the fiber over $x$ are assumed to
belong to the general bilinear group
 $\Gm^{0,2}(E)$. 
Proposition 9.4 says that the tangent bundle of a vector bundle
   is a bilinear bundle.  One may construct also
bilinear bundles that are not of this type.
Many definitions and statements of the subsequent sections generalize to
this framework; however, we will rather stay in the framework of tangent
bundles of vector bundles. -- 
Let us just add the remark that one can
 characterize bilinear bundles intrinsically, without referring to chart
conditions. We leave this is an exercise to the interested reader.


\vfill \eject 

\subheadline
{10. Linear connections. I: Linear structures on bilinear bundles}

\msk \nin
{\bf 10.1.} {\sl Linear structures on fiber bundles.}
A {\it linear structure $L=(a,m)$ on a fiber bundle
 $p:E \to M$} is given by  smooth bundle maps
$$
a: E \times_M E \to E, \quad m: \K \times  E \to E
$$
such that, for all $x \in M$, the fiber $E_x$ is turned
into a $\K$-module with linear structure $L_x$ given by
addition $a_x:E_x \times E_x \to E_x$
and multiplication by scalars $m_x: \K \times E_x \to E_x$.
Morphisms of fiber bundles with linear structure are
morphisms of fiber bundles that are fiberwise linear
(w.r.t.  the given linear structures $L,L'$).

\msk \nin
{\bf 10.2.} {\sl  Definition of linear connections and their morphisms.}
As in Section 9.1, we assume that $F$ is a vector bundle over
$M$ with projection $p:F \to M$.
A {\it (linear) connection}  on $F$ is a linear structure
$L$ on the fiber bundle $TF$ over $M$ such that, for all $x \in M$
and 
for every bundle chart (say, with index $j)$ around $x$, the structure $L_x$
is  bilinearly related to the linear structure $L_x^{(j)}$
induced on $(TF)_x$ by the  bundle chart of $TF$ with index $j$.
(This is well-defined: 
if $L_x$ is bilinearly related to {\it one} structure induced
from a bundle chart, then Proposition 9.4 implies that
it is bilinear related to {\it all}
structures induced from bundle charts because
``bilinearly related'' is an equivalence relation.)
A morphism of vector bundles $F,F'$ with linear connections $L,L'$
is a  smooth bundle map $\tilde f: F \to F'$ such that
$T\tilde f:TF \to TF'$ is fiberwise linear with respect to $L$ and $L'$;
following the usual convention, we will call them also
{\it affine (with respect to the given connections)}
(although the term {\it linear} would be more
appropriate).

\msk \nin
{\bf 10.3.} {\sl Connections on the tangent bundle.}
In case $F=TM$, we say that a connection
$L$ on $TM$ is {\it torsionfree} or {\it symmetric}
 if $L_x$ is bilinearly related via
a {\it symmetric} bilinear map to  $L_{\phi_j}$.
(As above, using the last assertion of Proposition 9.4,
it is seen that this condition is well-defined.) This means
that the action of the canonical flip $\kappa:TTM \to TTM$
is linear with respect to $L$ (cf. Section BA.5).
A smooth map $\phi:M \to M'$ is called
{\it affine (with respect to linear connections on $TM$ and on $TM'$)} 
if $T\phi:TM \to TM'$ is $L-L'$-affine;
i.e. if
$TT\phi:TTM \to TTM'$ is fiberwise linear with respect to $L$ and $L'$.
In case $\phi$ is an affine curve, 
i.e. $M=I \subset \K$ is open and equipped with the
connection induced by this chart, $\phi:I \to M'$ is also
called a {\it geodesic}.

\msk \nin
{\bf 10.4.} {\sl Chart representations and existence questions.}
Recall that the space of all bilinearly related structures on
a product of three $\K$-modules is an affine space over $\K$
(Section  BA.2). Applying this fiberwisely, we see that
the space of all affine connections on $F$, denoted by
${\rm Conn}(F)$, also is an affine space over $\K$.
More precisely, assume $L,L' \in {\rm Conn}(F)$.
Then, for all $x \in M$, $L_x$ and $L_x'$ are bilinearily
related to each other, which means that their difference
is a well-defined bilinear map
$$
b_x:V_1 \times V_2 \to V_3, \quad (u,v) \mapsto b_x(u,v)=\pr_3^{L}(u +_{L'} v)
$$
(Section BA.1), where the three axes are canonically given by
$V_1 = F_x$, $V_2=T_xM$, $V_3=\epsilon F_x$ (cf. Table 9.4).
In other words, 
$b=L-L'$ is a tensor field whose value at $x$ is a bilinear 
 map $F_x \times T_x M \to F_x$.
Conversely, adding such a tensor field $b$ to a linear connection
$L$ we get another linear connection $L'$. 

Thus ${\rm Conn}(F)$
 is an affine space over $\K$ (and over the ring ${\cal C}^\infty(M,\K)$).
Unfortunately, this affine space may be empty. However,
the space of connections of the bundle over a chart
domain is never empty because, by the very
definition of a connection, every chart induces a connection over
the chart domain. More precisely, fix a bundle chart $\tilde g_i=(\phi_i,g_i)$;
this chart induces a linear structure $L^i_x$ on the fiber
$(TF)_x$, given by the bundle chart $(TF)_x \cong W \times V \times 
\epsilon W$. Then, if $L$ is any other connection on $F$ defined over $U_i$,
the difference $L - L_i$ corresponds to a tensor field
 $b:=b^i:U \times W \times V \to W$, $(x,w,v) \mapsto b_x(w,v)$, 
which we call
 the {\it Christoffel tensor (of the connection in the given chart)}.
The linear structure $L_x$ is then explicitly given by
$+_{L_x} = +_{b_x}$ with $+_{b_x}$ given by Equation (BA.1), and 
similarly for multiplication by scalars:
$$
\eqalign{
(x,u,v,w)+_{L_x} (x,u',v',w') & =
(x,u+u',v+v',w+w'+b_x(u,v')+b_x(u',v)), \cr
r_{L_x} (x,v,w) & = (x,ru,rv,rw+(r^2 -r) b_x(u,v)). \cr}
\eqno (10.1)
$$
The calculation from Sections 9.1 and 9.2 gives us the relation
between Christoffel symbols $b^i$, $b^j$ belonging to different
bundle charts:
$$
b^j_{\phi_{ij}(x)}(g_{ij}(x)u,d\phi_{ij}(x)v) =
g_{ij}(x) b^i_x (u,v) + \partial_v  g_{ij}(x,u).
\eqno (10.2)
$$
In the sequel we will mainly work with connections over
a chart domain. Global connections, if they exist, can
be constructed in various ways:
a family of tensor fields $(b^i)_{i \in I}$ 
satisfying (10.2) defines a linear connection $L$; or, 
if our manifold $M$ has a partition of unity subordinate to
a cover by charts, 
then  the usual arguments show that
the chart connections can be patched together to define a global
connection on $M$. Finally, connections on Lie groups and symmetric
spaces will be constructed in an intrinsic way without
using chart arguments (see Chapters 23 and 26).
-- The following result generalizes a result due to Dombrowki, cf.\
 [La99, Th.\ X.4.3]:

\Theorem 10.5. {\rm (The Dombrowski Splitting Theorem.)}
A linear connection $L$ on a vector bundle $F$ induces on $TF$
a unique  structure of a vector bundle over $M$ such that
$$
\Phi: F \oplus_M TM \oplus_M F \to TF, \quad
(w_x, \delta_x,w_x') \mapsto z(w_x) + z(\delta_x) + \epsilon(z(w_x'))
$$
becomes an isomorphism of vector bundles over $M$,
where the sum is taken with respect to the linear structure $L_x$ on the
fiber $(TF)_x$.
Conversely, given an isomorphism of bilinear bundles 
$\Phi: F \oplus_M TM \oplus_M F \to TF$ whose restriction to the axes
$F$, $TM$  is just inclusion of axes, one recovers a
linear connection $L$ simply by push-forward of the canonical linear
structure $L_0$ on $F \oplus_M TM \oplus_M F$ by $\Phi$.

\Proof.
Choose a bundle chart $g_j$
and let $(b_x)_{x \in U}$ be the corresponding ``Christoffel symbols''.
Then $\Phi$ is represented in the chart by
$$
\Phi(x,w,x',w')= (x,w,0,0)+_{b_x} (x,0,x',0) +_{b_x} (x,0,0,w') =
(x,w,x',w'+b_x(w,x')),
$$
that is, in each fiber $\Phi$ is given by $f_{b_x}$
(with $f_b$ as in Section BA.1).
Hence $\Phi$ is invertible (the inverse being given in each fiber by
$f_{-b_x}$) and smooth in both directions.
It is fiberwise linear (by construction), hence an isomorphism of linear
bundles.

Conversely, if $\Phi$ is an isomorphism of bilinear bundles (in the sense
of Section 9.6), then it follows directly from the definitions that
the push-forward $\Phi_*(L_0)$ is a linear connection on $F$.
\qed

\nin
If a connection is fixed, then the isomorphism $\Phi$ from the
Dombrowski splitting will often be considered as an identification,
also by denoted by
$$
TF  \,\,  {\buildrel L \over \cong} \, \,  F \oplus TM \oplus \epsilon F
\eqno (10.3)
$$
(where $\epsilon F$ is just a copy of $F$,  the
prefix $\epsilon$ being added in order to distinguish the first and the
third factor). In the following two sections we explain the relation
with several classically well-known concepts that are related to connections.

\msk \nin
{\bf 10.6.} {\sl  Structures defined by a connection:
horizontal subspace and horizontal lift.}
Assume $L$ is a linear connection on a vector bundle $p:F \to M$.
Fix $x \in M$ and $u \in F_x$; in a chart we write $u=(x,w)$.
Recall from the table in Section 9.5 the following
intrinsic subspaces:  the {\it tangent space}
$T_u F = \{ w \} \times V \times W$, the
axis $T_x M = z(T_x M) = 0 \times V \times 0$, and
the {\it vertical subspace}  $V_u = \{ w \} \times 0 \times W$.
The parallel to the axis $0 \times V \times 0$ in $(TF)_x$
through $u=(w,0,0)$ depends on the linear structure $L_x$;
it is called the {\it horizontal subspace $H_u$} and is given in
a chart representation by (cf.\ Equation (BA.15))
$$
H_u = \{ \bigl( w,v,b_x(w,v) \bigr) \vert \, v \in V \}= \pr_{13}^{-1}(w,0).
\eqno (10.4)
$$
It is a subspace of the tangent space $T_u F$ such that
$$
T_u F = V_u \oplus H_u,
\eqno (10.5)
$$
and $Tp:TF \to TM$ induces a linear isomorphism
$H_u \to T_xM$ with inverse called the {\it horizontal lift} and
denoted by
$$
h_u: T_xM \to H_u \subset T_u F .
\eqno (10.6)
$$

\msk
\nin {\bf 10.7.} {\sl Structures defined by a connection:
connector and  connection one-form.}
The isomorphism (10.3) from the
Dombrowski splitting 
 gives rise to three bundle projectors over $M$:
$$
\pr_F:TF \to F, \quad \quad
 \pr_{TM}:TF \to TM, \quad \quad
 \pr_{\epsilon F}: TF\to \epsilon F.
$$
The first one is just $\pi_{TF}$, the second one is $Tp$, and the
third one really depends on the connection; it is called the
{\it connector} and in the literature is often denoted by 
$K: TF \to F$.  We prefer the notation
$$
K= \pr_{\epsilon F}: TF \to \epsilon F
$$
which permits to distinguish $K$ from the first projection.
 In a bundle chart, the third projection is given by
$$
K(x,u,v,w)=(x,w-b_x(u,v)).
$$
The connector $K$ has the usual properties
 known from real finite-dimensional differential
geometry (cf. e.g. [KMS93, p.\ 110], [Bes78, p.\ 38]): for all $x \in M$,

\ssk
\item{(1)}
$K_x:(TF)_x \to \epsilon F_x$ is {\it (intrinsically) bilinear} in the sense
explained in Prop.\ BA.7,
\item{(2)}
$K_x(z)=z$ for all $z$ belonging to the third axis $\epsilon F_x$ (vertical
space) of $(TF)_x$.

\ssk
\nin One can show that these properties characterize
the third projection. Instead with the connector $K$,
we could also work  with the projector
$\pr_{13}=\pr_{1} \times \pr_3$, i.e. with the bundle projector over $M$
$$
\Psi:=\pr_F \oplus \pr_{\epsilon F}: TF \to F \oplus \epsilon F
$$
which is called the {\it connection one-form}, given in a chart by
$$
\Psi(x,u,v,w)=(x,u,w-b_x(u,v)).
$$
Its kernel is the horizontal bundle, and its image can be seen
as the vertical bundle.
The map $\Psi$ can be seen as a fiber bundle map over $M$ or as
a vector bundle map over $F$ which is linear in fibers and
hence can be interpreted as a one-form over $F$.
 Moreover, composing with the canonical injection
$F \oplus \epsilon F \to TF$ (vertical bundle),
we get $\Psi$ in the version of a projector $\tilde \Psi:TF \to TF$
having as image the vertical bundle
(version used in [KMS93]).
Using these maps, the (inverse of the) isomorphism from the Dombrowski
splitting theorem can be written
$$
\Phi^{-1} = \pi \oplus Tp \oplus K = Tp \oplus \Psi :
TF \to F \oplus TM \oplus \epsilon F.
\eqno (10.7)
$$
The ``dictionary" between bilinear algebra and differential
geometry from Section 9.5 can now be continued -- this part concerns
``non-intrinsic structures":

\bigskip
\settabs 2 \columns
\+
$\quad${\sl bilinear algebra} & {\sl differential geometry} \cr
\bigskip
\+
fixing $L^b$, $b \in \Bil(W \times V, W)$ & fixing a connection $L$
\cr
\+
$f_b:W \times V \times W \to W \oplus_b V \oplus_b W$ &
$\Phi: F \oplus TM \oplus \epsilon F \to TF$ (Dombrowski)
\cr
\+
projectors: & \cr
\+
$\quad \pr_{13}:W \oplus_b V \oplus_b W \to W \times W$ &
$\Psi:TF \to F \oplus \epsilon F$ (connection one-form) \cr
\+
$\quad \pr_3: W \oplus_b V \oplus_b W \to W$ &
$K:TF \to \epsilon F$ (connector) \cr
\+
$\quad \pr_{23}:W \oplus_b V \oplus_b W \to V \times W$ &
$\pr_{23}:TF \to TM \oplus \epsilon F$ (no name)
\cr
\+
spaces, decompositions: & \cr
\+
$\quad H_w = \{ (w,v,b(w,v)) | \, v \in V \}$
& horizontal space $H_u = \Psi^{-1}((u,0))$
\cr
\+
$\quad V_w = \{ (w,0,z) | \, z \in W \}$ (intrinsic)
&
vertical space $V_u$ (intrinsic)
\cr
\+
$\quad \{w\} \times V \times W = V_w \oplus_b H_w$
&
$T_u F = V_u \oplus H_u$
\cr
\+
$\quad V \to H_w$, $v \mapsto (w,v,b(w,v))$
&
horizontal lift $h_u:T_x M \to H_u \subset T_u F$
\cr
\+
non-intrinsic exact sequences:
& 
\cr
\+
$\quad W \oplus_b V \oplus_b 0 \to W \times V \times W {\buildrel \pr_3 \over
\to} W$
&
$F \oplus TM \to TF {\buildrel K \over \longrightarrow} \epsilon F$
\cr
\+
$\quad V  \to W\times V \times W {\buildrel \pr_{13} \over \to} W 
\times W$ &
$TM \to TF {\buildrel \Psi \over \longrightarrow}
 F \oplus \epsilon F$ 
\cr
\+
Case $V = W$: & Case $F=TM$:
\cr
\+
$\quad$ torsion of $L^b$ (= $L^b - L^{\kappa.b}$) & torsion tensor
of the connection
\cr
\+
$\quad$ diagonal map $V \to (V \times V) \times V$
&
$S:TM \to TTM$
(spray, cf.\ Section 11.3.)
\cr

\bigskip \nin
{\bf 10.8.} {\sl Connections on bilinear bundles.}
The preceding definitions and facts can be generalized in a
straightforward way for arbitrary bilinear bundles (cf. Section 9.6):
a {\it linear connection on a bilinear bundle $F$} is
 given by a linear structure $L$
on $F$ over $M$ which is bilinearily related to structures induced by
bundle charts. 
The Dombrowski Spitting Theorem and the 
notions of connector, connection one-form, and so on, immediately
carry over to this context;  details are left to the reader.

\msk \nin
{\bf 10.9.} {\sl Direct sum of connections.}
From the classical linear algebra constructions of connections
(tensor products, hom-bundles, dual connection...), only the
{\it direct sum of connections} survives in our general context. In order to
give an intrinsic definition of the direct sum, we first describe some
canonical isomorphisms of bundles. For this, assume given three vector bundles
$F_1,F_2,F_3$ over $M$, with fibers modelled on $\K$-modules $W_1,W_2,W_3$.

 First of all, there is a canonical isomorphism
$$
T(F_1 \times_M F_2) \cong TF_1 \times_{TM} TF_2.
\eqno (10.8)
$$
In a bundle chart with base domain $U$,
 this amounts to the canonical identification
$$
T(U \times W_1 \times W_2) \cong TU \times TW_1 \times TW_2
$$
which is a special case of the general isomorphism $T(M \times N)
\cong TM \times TN$ for manifolds $M,N$. 

Next, $F_1 \times_M F_2$ is a vector bundle over $M$, but it may also
be considered as a vector bundle over $F_2$: the fiber over $g \in (F_2)_x$
with $x \in M$ is $\{ (f,g) | \, f \in (F_1)_x \} \subset 
(F_1)_x \times (F_2)_x$. Then 
$$
(F_1 \times_M F_2) \times_{F_2} (F_3 \times_M F_2)
$$
is a vector bundle over $F_2$, whereas over $M$, {\it a priori}, 
it is just a fiber bundle. But
$$
F_1 \times_M F_2 \times_M F_3 \to
(F_1 \times_M F_2) \times_{F_2} (F_3 \times_M F_2), \quad
(x;f,g,h) \mapsto (x,g;f,h)
\eqno (10.9)
$$
is a bijection of bundles over $M$
-- in a chart, this amounts to the canonical identification
$$
U \times (W_1 \times W_2 \times W_3) \cong (U \times W_2) \times
W_1 \times W_3.
$$
Now, the left hand side in (10.9) is a vector bundle {\it over $M$},
and hence, via the isomorphism (10.9), also the right hand side
carries a canonical vector bundle structure over $M$.

Next,
 assume that $F,G$  are vector bundles over $M$ with connections
$L_i$, $i=1,2$. Then, using the canonical isomorphisms described above
and the Dombrowski splittings (10.3) corresponding to $L_1$ and $L_2$, the 
following isomorphism defines a linear structure on $T(F \oplus_M G)$:
$$
\eqalign{
T(F \oplus_M G) & \cong TF \oplus_{TM} TG  \cr
& \cong (F \oplus_M TM \oplus_M \eps TF) \oplus_{TM} (
G \oplus_M TM \oplus_M \eps TG) \cr
& \cong F \oplus_M \eps F \oplus_M TM \oplus_M G \oplus_M \eps G. \cr}
\eqno (10.10)
$$
In a chart, this isomorphism is described, if $b_x$ and $a_x$ are
the Christoffel tensors on $F$, resp.\ $G$, 
$$
(x;\eps v,u,w, \eps u',\eps w') \mapsto
(x;\eps v,u,w,\eps (u'+b_x(v,u)), \eps(w'+a_x(v,w)),
$$
i.e.\ by the Christoffel tensor $c_x(v,(u,w))=b_x(v,w)+a_x(v,w)$.

\msk
\nin {\bf 10.10.} {\sl Connections and the tensor bundle.}
Recall from Section BA.10 that to a bilinear space
$E=V_1 \times V_2 \times V_3$ we can associate in an intrinsic
way the {\it linear} space
$Z:=V_1 \otimes V_2 \oplus V_3$. Doing this construction pointwise,
we can associate to every bilinear bundle $F$ an  ``algebraic linear bundle"
$$
Z:= Z(F): = F_1 \otimes F_2 \oplus F_3
\eqno (10.11)
$$
i.e., $Z$ is defined as a set, but without topology and hence
without manifold structure (cf. Appendix L.2). Nevertheless we may speak about
``bundle charts", which are just transition functions without
specified smoothness properties. As with tensor products in ordinary finite
dimensional real differential geometry, it is seen that there is
a canonical projection $Z(F) \to M$ such that the
fibers carry a well-defined $\K$-module structure, and hence $Z$
is a well-defined ``algebraic linear bundle''.
The map
$$
TF \to Z(F), \quad (x,w,v,w') \mapsto (x,w \otimes v + w')
$$
does not depend on the chart and defines a
bundle map from a bilinear bundle to a linear bundle;
this map is fiberwise intrinsically bilinear in the sense of Prop.\ BA.7.
 The sequence $F_3 \to F \to F_1 \times_M F_2$ gives rise to an exact
sequence of algebraic linear bundles over $M$
$$
0 \to F_3 \to Z(F)  \to F_1 \otimes F_2 \to 0.
\eqno (10.12)
$$
Linear connections now give rise to linear splittings of the 
sequence (10.12):
$$
L: F_1 \otimes F_2 \to Z(F)
\eqno (10.13)
$$
If $F=TM$, then (10.12) reads $0 \to TM \to Z(TM) \to TM \otimes TM \to 0$,
 and restricting to $\Sigma_2$-invariants, we get an exact sequence 
$0 \to TM \to Z(TM)^\kappa \to
S^2(TM) \to 0$. Torsionfree connections give rise to cross-sections
$\Gamma: S^2(TM) \to Z(TM)^\kappa$ of this sequence. In the finite-dimensional
real case, there is a one-to-one correspondence between such cross-sections
and connections --
this is the point of view on connections used in [Lo69].
See also [P62] and the remarks on differential operators and symbols in
Chapter 21.

\msk
\nin {\bf 10.11.} {\sl Comments on definitions of a linear connection.}
As mentioned in the introduction, 
one finds many different definitions of connections in the literature,
and it was not our aim to add a new item to this long list,
but rather to propose a way how to organize it:
both the axiomatics based on the connector (see, e.g.,
[Bes78],  [KMS93]) and the axiomatics based on the connection one-form
(see, e.g., [BGV92], [KMS93], [KoNo69]) focus on specific aspects
(projectors)
of the linear structure $L$ and are equivalent to ours.
Similarly, the definition in [La99, p.\ 104] takes some features of the
injection $F \times_M TM \to TF$ as axiomatics and again is
equivalent to ours. Using our concept of 
bilinear algebra, it is fairly easy to prove  that all the aforementioned
 definitions 
of a connection are equivalent, or even to give yet other axiomatic
characterizations. For instance,  a connection defines an isomorphism of
 the fibers of $(Tp)_x:(TF)_x \to T_xM$ with
$F_x \oplus \epsilon F_x \cong F_x \oplus F_x$, and hence we have
a whole $2\times 2$-matrix algebra $M(2,2;\K)$ acting
by bundle maps of $TF$ over $TM$,
this action depending only on the connection $L$:
$$
\epsilon_u =\pmatrix{0 & 0 \cr \1 & 0 \cr}, \quad
J_u = \pmatrix{0_W & \cr  & \1_W  \cr}, \quad
H_u=\pmatrix{\1_W &  \cr  & - \1_W \cr}, \quad
Q_u=\pmatrix{0 & \1 \cr 0 & 0 \cr}
$$
where $W \cong F_x$ and $u=x+\epsilon v \in TM$.
We propose to call this algebra the {\it connector algebra}.
Here, the first operator is independent of $L$ (the intrinsic
$\epsilon$-tensor), whereas the other operators do depend on
$L$ (except on the fiber of $TF$ over the zero vector $0_x \in T_xM$).
The operator $Q$ seems to be particularly interesting --
it is in a sense opposed to $\epsilon$:
 when $\epsilon$ is seen as sort of ``annihilator''
(making an object ``infinitesimally small''),
$Q$ is a sort of ``creator'' (magnifying an ``infinitesimally
small'' object by an ``infinitely large'' factor).
-- 
Finally, the relation with the definition of connections via
sprays and covariant derivatives is discussed in the following
two chapters.  

\vfill \eject  

\subheadline
{11. Linear connections. II: Sprays}

\nin {\bf 11.1.} {\sl Second order vector fields.} 
This chapter concerns linear connections on the tangent bundle: let $F=TM$.
Recall that the second order jet bundle $J^2 M$ is
the subbundle of $TF=TTM$ fixed under the action of the
permutation group $\Sigma_2$, and that the two projections
$p_j:TTM \to TM$ give rise to one jet projection $\pi_{1,2}:J^2M \to J^1M=TM$:
$$
\matrix{
J^2 M & \subset & TTM \cr
\phantom{\pi_{1,2}} \downarrow \pi_{1,2} &  & \phantom{p_j} \downarrow p_j \cr
J^1 M & = & TM .\cr}
\eqno (11.1)
$$
As we have seen in 8.4 (B), the jet bundle $J^2 M \to TM$ 
 does not have a canonical ``zero-section'' since
the zero-sections of $p_j$ are not compatible with the canonical
flip. Thus $\pi_{1,2}:J^2 M \to J^1 M$ is an affine bundle and not a 
vector bundle. We define 
a {\it second order vector field on $M$} to be a cross-section 
$X:J^1 M \to J^2 M$ of the jet projection $J^2 M \to J^1 M$. 
This defines a vector bundle structure on $J^2 M$ over $J^1 M$, simply by
considering $X$ as zero-section.
One may also say that $X:TM \to J^2 M \subset TTM$ is a vector field on $TM$,
taking values in the space of $\Sigma_2$-invariants. Thus a vector field
 $X:TM \to TTM$ is a second order vector field if and
only if in a bundle chart it is given by
$$
X(x+\epsilon v) = x + (\epsilon_1 + \epsilon_2) v + \epsilon_1
\epsilon_2 \xi(x,v)
\eqno (11.2)
$$
with some smooth map $\xi:U \times V \to V$.

\msk \nin
{\bf 11.2.} {\sl The second order vector field defined by a
linear connection.}
If $L$ is a linear connection on $TM$ and
$$
\epsilon_1 TM \times_M \epsilon_2 TM \times_M \epsilon_1
\epsilon_2 TM \to TTM, \quad
(x;\epsilon_1 u,\epsilon_2 v,\epsilon_1 \epsilon_2v) \mapsto
\epsilon_1 u +_L \epsilon_2 v +_L \epsilon_1 \epsilon_2v
$$
its corresponding Dombrowski splitting, then the diagonal imbedding
$$
S: TM \to  TM \times_M  TM \times_M TM \, \, {\buildrel L \over \cong} \, \,
TTM, \quad v \mapsto (v,v,0) = \epsilon_1 v +_L \epsilon_2 v
$$
is a second-order vector field. In a chart, it is given by
$$
S(x+\epsilon v) = x + (\epsilon_1 + \epsilon_2) v + \epsilon_1
\epsilon_2 b_x(v,v)
\eqno (11.3)
$$
where $b$ is the Christoffel tensor of $L$.
(Note that $S$ depends only on the quadratic map determined by $b$,
and hence, if $2$ is invertible in $\K$, we may w.l.o.g. $b$ assume
to be symmetric, i.e. $L$ to be torsionfree.)
For instance, if $M$ is open in the topological $\K$-module
 $V$ and $L$ is the canonical chart
connection ($b=0$), then $S(x,v)=(x,v,v,0)$ is the corresponding
second order vector field.

\msk \nin
{\bf 11.3.} {\sl Sprays.} 
Not every second order vector field comes from a linear connection
on $TM$ because the smooth map $\xi$ from (11.2) can be chosen 
arbitrarily. Maps of the special form $\xi(x,v)=b_x(v,v)$ can
be characterized intrinsically by the following {\it spray
condition}: let $r \in \K$ and recall from Section 7.3 (D) the intrinsically
defined action $l^{(2)}:\K \times TTM \to TTM$ which is given in a chart by
$r^2(x,u,v,w)=(x,ru,rv,r^2w)$ (and which can also be defined in the purely
algebraic context of bilinear spaces, see remark after Prop. BA.3).
Then, if $S$ comes from a linear connection $L$, we have, in a chart
$$
\eqalign{
S  (x+\eps r v) & = 
 x+(\epsilon_1 + \epsilon_2)rv +
\epsilon_1 \epsilon_2 b_x(rv,rv) \cr
& =x+r(\epsilon_1 + \epsilon_2)v +
r^2\epsilon_1 \epsilon_2 b_x(v,v)= l^{(2)}_r \circ S(x+ \eps v). \cr}
$$
By definition, a {\it spray} is a second order vector field satisfying the
``homogenity condition''
$X \circ r_{TM} = l^{(2)}_r \circ X$, for all $r \in \K$.

\Theorem 11.4. Assume $2$ is invertible in $\K$. Then
there is a canonical bijection between torsionfree connections
$L$ on $TM$ and sprays $S:TM \to J^2 M$. 

\Proof.
We have just seen that the second order vector field $S$ associated to
a connection is a spray. 
Conversely, if $S$ is a spray, then, writing
$S(x+\epsilon v)=x + \epsilon_1 v + \epsilon_2 v + \epsilon_1
\epsilon_2 q_x(v)$, we deduce that $q_x$ is smooth and 
homogeneous quadratic with respect to scalars
 and hence is, by Cor.\ 1.12, a vector valued quadratic form. 
If $2$ is invertible in $\K$, then the symmetric bilinear map $b_x$
can be recovered from $q_x$  by polarisation. Clearly, both constructions are
inverse to each other.             \qed

\nin The preceding result  generalizes a theorem of 
 Ambrose, Palais and Singer, cf.\ [MR91, 7.3].
In the case of general characteristic, the spray is equivalent
to what one might call the ``connection on the {\it quadratic} 
bundle $TTM$'' (cf. Section BA.9).

\msk \nin
{\bf 11.5.} {\sl Second order differential equations and geodesics.} 
If $X:TM \to J^2 M$ is a second order vector field, then, by definition,
a {\it geodesic with respect to $X$} is a curve $\alpha:\K \to M$ such that
 $\alpha':\K \to TM$ is an integral
curve of $S$, i.e. $\alpha''(t) = S(\alpha'(t))$.
Here, $\alpha'(t)=T_t \alpha \cdot 1=T\alpha (t + \eps 1)$
 and $\alpha''(t)=(\alpha')'(t)=TT\alpha(t + \delta 1 + \delta^{(2)} 1)$.
In a chart representation (11.2), the condition $\alpha''(t) = S(\alpha'(t))$
amounts to the second order differential equation
$\alpha''(t) = \xi(\alpha(t),\alpha'(t))$ (cf.\ [La99, p.\ 99]).

\ssk
On the other hand, if $X=S$ is the spray of a connection $L$, the notion
of  geodesic has already been defined in Section 10.3: it
is an affine map $\alpha:\K \to M$,
i.e., $J^2 \alpha : J^2 \K \to J^2 M$ is fiberwise linear.
This is equivalent to saying that the following diagram commutes:
$$
\matrix{J^2 \K & {\buildrel J^2 \alpha \over \to} & J^2 M \cr
S_\K \uparrow \phantom{X_\K} & & \phantom{X} \uparrow S \cr
J^1 \K & {\buildrel J^1 \alpha \over \to} & J^1 M, \cr}
$$
where $S_\K:T\K \to J^2\K$ is the spray corresponding to the canonical chart
connection of $\K$. 
In a bundle chart, we have,  with $\xi$ as in (11.2), using the spray 
property,
$$
\eqalign{
S \circ J^1 \alpha (x + \eps v) & = \alpha(t) + \delta \alpha'(t) v +
\delta^{(2)} \xi(\alpha(t),\alpha'(t)v) \cr
&  = \alpha(t) + \delta \alpha'(t) v +
\delta^{(2)} v^2 \xi(\alpha(t),\alpha'(t)), \cr
J^2 \alpha \circ S_\K(t + \eps v) & = \alpha(t) + \delta \alpha'(t) v +
\delta^{(2)} \alpha''(t) v^2, \cr}
$$
and hence the condition $J^2 \alpha \circ S_\K = 
S \circ J^1 \alpha$ reads in a chart
$\alpha''(t) = \xi(\alpha(t),\alpha'(t))$.
This is the same condition as above, and hence both notions of geodesic
agree. 
For $v \in T_x M$, we say that $S(v) \in (J^2 M)_x$
is the {\it geodesic $2$-jet of $v$}. 
Note that in general we have no existence or uniqueness statements on
geodesics; but the preceding arguments show that, if a geodesic
$\alpha_v$ with $\alpha_v(0)=x$ and $\alpha_v'(0)=v$ exists,
then $J^2 \alpha(0)=S(v) $.

\vfill \eject  

\subheadline
{12. Linear connections. III: Covariant derivative}
			
\msk \nin
{\bf 12.1.} {\sl The covariant derivative associated to a connection.}
We return to the case of
 a general vector bundle $F$ with linear connection $L$.
Assume $X:M \to TM$ is a vector field on $M$ and $Y:M \to F$
a section of $F$. If $K:TF \to \epsilon F$ is the connector
associated to a connection $L$ on $F$, we define
$$
\eqalign{
\nabla Y := K \circ TY: & \quad  TM \to \epsilon F, \cr
\nabla_X Y:= K \circ TY \circ X: & \quad  M \to \epsilon F. \cr}
$$

\Lemma 12.2.
Let $(b_x)_{x \in U}$ be the Christoffel tensor of the
connection $L$ in the chart $U$. Then
$$
(\nabla_X Y)(x) = dY(x) X(x) + b_x(Y(x),X(x)).
$$

\Proof. In a bundle chart, by some abuse of notation,
$X(x)=x + \eps X(x)$,
$Y(x)=(x,Y(x))$,
$TY(x+\eps v)=(x , Y(x)) + \eps ( v , dY(x)v))$,
$K((x,w) + \eps (x',w'))=x + \eps (w'+ b_x(w,x'))$, and thus
$$
\eqalign{
(\nabla_X Y)(x) & =x + \eps K \bigl( (x,Y(x)) + \eps(X(x),dY(x)X(x))\bigr) \cr
&  = x + \eps (dY(x)X(x) + b_x(Y(x),X(x))). \cr}
\qeddis

\Corollary 12.3.
If $L$ is a torsionfree connection on $TM$, 
then for all vector fields $X,Y$,
$$
[X,Y](x) = (\nabla_Y X - \nabla_X Y)(x).
$$

\Proof.
In a chart, the preceding lemma gives for the right hand side,
due to the symmetry of $b_x$, the expression
$dX(x) Y(x) - dY(x)X(x)$ which
 is precisely the definition of the Lie bracket.
\qed

\msk \nin {\bf 12.4.} {\sl General covariant derivatives.}
A {\it (general) covariant derivative} is a $\K$-bilinear map
from sections to sections,
$$
\nabla: \Gamma(TM) \times \Gamma(F) \to \Gamma(F)
$$
such that, for all smooth functions $f$ and all sections $X,Y$,
$$
\nabla_{fX} Y = f \nabla_X Y, \quad \quad
\nabla_X(fX)=Xf \cdot Y + f \nabla_X Y.
$$
It is easily proved that the covariant derivative associated to
a connection $L$ has these properties. The converse is not true:
not every general covariant derivative is associated to a linear
connection
-- see  [La99, p.\ 202 ff and p.\ 279/80] for a discussion (which can, {\it 
mutatis mutandis}, be applied to our situation) of the
infinite dimensional real case. In particular, if one wants to
define the curvature tensor in terms of covariant derivatives, 
it is not sufficient to use just the properties of a general covariant
derivative (see [La99, p.\ 231 ff]). Therefore, in Chapter 18, we will define 
the curvature tensor in a different way.

\msk \nin
{\bf 12.5.} {\sl Covariant derivative of a tensor field.}
For a tensor field of type $(k,0)$, say  $\omega: TM \times_M
\ldots \times_M TM \to \K$,  the covariant derivative 
$\nabla \omega$ can be defined along the lines of Remark 4.8.
E.g., for a one-form $\omega:TM \to \K$ we define
$$
(\nabla \omega)(X,Y):= L_X (\omega(Y)) - \omega(\nabla_X Y) ;
\eqno (12.1)
$$
we then have to check that the value at $x$ depends only on
the value of $X$ and $Y$ at $x$ and that the result depends
smoothly on $x$, $X(x)$ and $Y(x)$.
Then, by Remark 4.8, $\nabla \omega$ can be identified with 
a well-defined tensor field $TM \times_M TM \to \K$.
In the sequel we will not use this construction, but see
Section 13.4 for an analog of Corollary 12.3.

\vfill \eject

\subheadline
{13. Natural operations. I: Exterior derivative of a one-form}

\msk \nin
{\bf 13.1.} {\sl Bilinear maps defined on $TF$.}
If $p:F \to M$ is a vector bundle over $M$, a smooth bundle map
$\tilde f: TF \to E$, $f:M \to N$ into a vector bundle $E$ over $N$
will be called {\it (intrinsically) bilinear} if it is linear in fibers
both over $F$ and over $TM$. This is equivalent to saying that, for
all $x \in M$, $\tilde f_x:(TF)_x \to E_{f(x)}$ is intrinsically
bilinear in the sense of Proposition BA.7.
 From that proposition  it follows that $\tilde f$ is bilinear if,
 and only if, $\tilde f$ has a chart representation of the form
$$
\tilde f(x,w,v,w') = (f(x), g_x(w,v) + \lambda_x(w'))
\eqno (13.1)
$$
with bilinear $g_x$ and linear $\lambda_x$.
A bilinear map $\tilde f$ 
 is called {\it homogeneous} if it is constant on the vertical bundle
(fibers of $\pi \times_M Tp$),
and this is the case iff $\lambda = 0$ (see Cor.\ BA.8).
Thus the homogeneous bilinear maps are in one-to-one correspondence
with the fiberwise bilinear maps
$g:F \times_M TM \to E$
via $\tilde f = g \circ (Tp \times \pi)$:
$$
\matrix{TF & {\buildrel \tilde f \over \longrightarrow} \,\,  & E \cr
\downarrow   & {\buildrel g \over \nearrow} & \cr
TM \times_M F & & \cr}
$$

\ssk \nin
{\bf 13.2.}
{\sl The exterior derivative of a one-form.}
Recall that a {\it tensor field of type $(k,0)$}  is a smooth 
map $\omega: \times^k_M TM \to \K$, multilinear in fibers over $M$.
If $\omega$ is alternating in fibers, it is called a {\it $k$-form}.
A {\it one-form on a vector bundle $E$} is a smooth map $\omega: E \to \K$
such that $\omega_x:E_x \to \K$ is linear for all $x \in M$. Then
$$
\pr_2 T \omega: TE \to T\K \to \epsilon \K
$$
is (intrinsic) bilinear in the sense of 13.1: in fact, $T\omega$
is linear on fibers over $E$ since it is a tangent map, and $\pr_2 T \omega$
is linear in fibers over $TM$ since already $\omega$ was linear
in fibers over $M$. Equivalently, the bilinearity of $\pr_2 T\omega$ follows 
from the chart representation    
$$
 T\omega:TF \to T\K, \quad  (x,\epsilon v,w,\epsilon w') \mapsto
\omega(x,v) + \epsilon ({\partial_v }\omega(x,w) + \omega(x,w')).
$$
Comparing with (13.1) we see that $\pr_2 T\omega$ is bilinear, but
non-homogeneous (the term $\lambda_x(w')$ corresponds to $\omega(x,w')$).
For a general vector bundle $F$ we cannot extract a homogeneous
bilinear map from this. However, if $F=TM$, then
 we have the alternation map from intrinsic bilinear
maps $TTM \to \K$ to homogeneous ones:
$$
\eqalign{
\alt & (\pr_2 T \omega) (x+\epsilon_1 v_1  + \epsilon_2 v_2 +
\epsilon_1\epsilon_2v_{12}) \cr
& = (\pr_2 T\omega)(x+\epsilon_1 v_1 + \epsilon_2 v_2 +
\epsilon_1\epsilon_2 v_{12}) -
(\pr_2 T\omega)(x+\epsilon_1 v_2 + \epsilon_2 v_1 +
\epsilon_1\epsilon_2 v_{12}) \cr
& =
\partial_{v_1} \omega(x,v_2)-\partial_{v_2} \omega(x,v_1) \cr}
$$
is homogeneous (the value does not depend on $v_{12}$). As explained in 13.1,
homogeneous bilinear maps $\tilde f$
correspond to tensor fields $g$ of type $(2,0)$, and if $\tilde f$ is
skew-symmetric, then so is $g$. Hence we end up with
 $2$-form denoted by
$$
\dd \omega: TM \times_M TM \to \K, \quad
(x;v_1,v_2) \mapsto \alt(\pr_2 T\omega)(x+\epsilon_1 v_1 + \epsilon_2 v_2).
\eqno (13.2)
$$
In other terms, the diagram
$$
\matrix{TTM & {\buildrel \pr_2 T\omega \over \longrightarrow} & \K \cr
\downarrow & & \cr
TM \times_M TM & & \cr}
$$
can in general not be completed to a commutative triangle, but
$$
\matrix{TTM & {\buildrel {\rm alt} \pr_2 T\omega \over \longrightarrow}
 & \K \cr \downarrow & & \cr TM \times_M TM & & \cr}
$$
can, namely by $\dd \omega$.
The chart representation
$$
\dd \omega (x;v_1,v_2) = 
{\partial_{v_1} }\omega(x,v_2)  -
{\partial_{v_2} }\omega(x,v_1) 
\eqno (13.3)
$$
shows that $\dd \omega$ coincides with the usual exterior derivative
of a one-form. It is clear that in the preceding arguments the target
space $\K$ can be replaced by any other fixed topological $\K$-module
$W$, i.e. we can define $\dd \omega$ for any $1$-form $\omega:TM \to W$.

\msk \nin {\bf 13.3.} 
If $\omega=df$, i.e. $\omega(x,v)=T_xf \cdot v$,  then from (13.3) we get
$\dd \omega = \dd \dd f =0$ because the ordinary second 
differential $d^2 f(x)$ is symmetric.

\msk \nin {\bf 13.4.} {\sl Exterior derivative and connections.}
The chart formula (13.3) has the following generalization:
if $\nabla$ is the covariant derivative
of  an arbitrary torsionfree connection on $M$ and
$\omega$ a one-form on $M$, then for all $x \in M$
and $u,v \in T_x M$, we have
$$
\dd \omega(x;u,v)= (\nabla \omega)(x;u,v) - (\nabla \omega)(x;v,u).
\eqno (13.4)
$$
In fact,
by definition of the covariant derivative of a one-form (Section 12.5), 
$$
\eqalign{
\nabla \omega(X,Y) - \nabla \omega(Y,X) & =
X \omega(Y)-Y \omega(X) - \omega(\nabla_X Y)
+ \omega(\nabla_Y X) \cr
& =X \omega(Y)-Y \omega(X) - \omega([X,Y]) \cr}
$$
since $\nabla$ is torsionfree (Cor.\ 12.3). Now
fix $x \in M$ and $v,w \in T_x M$;   extend,
in a chart, $v$ and $w$ to constant vector fields $X,Y$ on the chart domain
 $U \subset M$; 
then $[X,Y]=0$ and hence from (13.3) we get
$$
\nabla \omega(x;v,w)- \nabla \omega(x;w,v) = \partial_v \omega(x,w) -
\partial_w \omega(x,v)= \dd \omega(x;v,w).
$$

\msk \nin {\bf 13.5.} {\sl Exterior derivative of a vector bundle
valued one-form.}
If we replace $\omega$ by a bundle-valued one-form
$\omega:TM \to E$, then in presence of a connection on $E$
with connector $K_E:TE \to \epsilon E$
we can define $\dd \omega := {\rm alt}(K_E \circ T\omega)$
which is skew-symmetric $TM \times TM \to \epsilon E$.
Note that the preceding definition in the scalar or vector valued case
is a special case of this, where $\pr_2$ is just the connector
of the canonical flat connection on $\K$, resp. on $X$.
An exterior derivative of bundle valued forms, in presence of
a connection, has also been defined in [Lav87, p.\ 165 ff].

\msk \nin {\bf 13.6.} {\sl Remark on the infinitesimal Stoke's theorem.}
In the context of synthetic differential geometry, one can prove
a purely infinitesimal version of Stoke's theorem, cf. [Lav87, p.\ 114 ff]
or [MR91, p.\ 134 ff]. It should be possible to state and prove such a version
also in the present context -- this will be taken up elsewhere.

\vfill \eject 

\subheadline{14. Natural operations. II: The Lie bracket revisited}

\msk \nin {\bf 14.1.}
{\sl Subgroups of the group of diffeomorphisms of $TTM$.}
The group $\Gamma:=\Diff(TTM)$ of diffeomorphisms
(over $\K$) of $TTM$ contains several subgroups which
are defined via the behaviour with respect to the 
projections $p_i:TTM \to TM$, $i=1,2$ and $p:TTM \to M$:

\ssk
\item{(1)} Let $\Gamma^+$ be the subgroup preserving fibers of
the projection $p:TTM \to M$ and
permuting the fibers of $p_1$ and those of $p_2$.
In a chart, $f \in \Gamma^+$ is represented by
$$
f(x+\eps_1 v_1 + \eps_2 v_2 + \eps_1 \eps_2 v_{12})=
 x + \eps_1 f_1(x,v_1) + \eps_2 f_2(x,v_2)+\eps_1 \eps_2 f_{12}
 (x, v_1, v_2 , v_{12})
$$
\item{(2)} Let $\Gamma^1$ be the subgroup of $\Gamma^+$ 
preserving all fibers of the projection
$p_1$. In a chart:
$$
f(x+\eps_1 v_1 + \eps_2 v_2 + \eps_1 \eps_2 v_{12})=
x+\eps_1 v_1 + \eps_2 f_2(x,v_2) + \eps_1 \eps_2 
f_{12}(x,v_1,v_2,v_{12})
$$
\item{(3)} Let $\Gamma^2$ be the subgroup of $\Gamma^+$
preserving fibers of $p_2$. In a chart,
$$
f(x+\eps_1 v_1 + \eps_2 v_2 + \eps_1 \eps_2 v_{12})=
x+\eps_1 f_1(x,v_1) + \eps_2 v_2 + \eps_1 \eps_2 
f_{12}(x, v_1, v_2, v_{12})
$$
\item{(4)} Let $\Gamma^{12}=\Gamma^1 \cap \Gamma^2$ be the
subgroup preserving all vertical spaces; in a chart,
$$
f(x+\eps_1 v_1 + \eps_2 v_2 + \eps_1 \eps_2 v_{12})=
x+\eps_1 v_1 + \eps_2 v_2 + \eps_1 \eps_2 f_{12}
(x, v_1, v_2, v_{12}).
$$

\msk
\nin
It is clear that $\Gamma^1$ and $\Gamma^2$ are normal subgroups
in $\Gamma^+$. Therefore, if $g \in \Gamma^1$ and $h \in \Gamma^2$,
it follows that the group commutator
$[g,h] = (ghg^{-1})h^{-1} = g (h g^{-1} h^{-1})$
belongs both to $\Gamma^1$ and to $\Gamma^2$ and hence belongs
to $\Gamma^{12}$. In Theorem 14.4 we will express the Lie bracket
of vector fields in terms of this commutator. 

\msk \nin
{\bf 14.2.} {\sl Sections of the second order tangent bundle.}
We denote by $\X^2(M):=\Gamma(M,T^2M)$ the space of smooth
sections of the second tangent bundle $T^2 M \to M$.
In a chart with chart domain $U \subset M$, we write
$T^2 U \cong U \times
\eps_1 V \times \eps_2 V \times \eps_1 \eps_2 V$, and a section
$X:U \to T^2 U$ is written in the form
$$
X(p) = p + \eps_1 X_1(p) + \eps_2 X_2(p) + \eps_1 \eps_2 X_{12}(p)
\eqno (14.1)
$$
with (chart-dependent) vector fields $X_\alpha: U \to V$.
In this notation, $X$ is a second-order vector field (cf.\ Section 11.1)
 if, and only if, in any
chart representation we have $X_1(p)=X_2(p)$ for all $p \in U$.
There are three canonical injections of the space of vector fields
 $\X(M)$ into $\X^2(M)$ which simply correspond to the
three canonical inclusions of axes $\iota_\alpha:TM \to TTM$,
$\alpha=01,10,11$, by letting $\eps^\alpha Y:= \iota_\alpha \circ Y:
M \to TTM$ for a vector field $Y:M \to TM$. 
In the notation (14.1), these are the sections 
$X$ having just one non-trivial component, i.e.,
with $X_\alpha = Y$ and $X_\beta =0$ for $\beta \not= \alpha$.

\Theorem 14.3.
\item{(1)}
There is a natural group structure on the space $\X^2(M)$, given in
a chart representation as above, by the formula
$$
\eqalign{
(X \cdot Y)(x) & =
x + \eps_1 (X_1(x)+Y_1(x)) + \eps_2 (X_2(x)+Y_2(x)) + \cr
& \quad \quad  \eps_1 \eps_2 
(X_{12}(x)+Y_{12}(x) + dX_1(x) Y_2(x) + dX_2(x) Y_1(x)). \cr}
$$
The space $\Gamma(M,J^2 M)$ of sections of $J^2 M$ is a subgroup of $\X^2(M)$,
and the three canonical injections $\X(M) \to \X^2(M)$ are group homomorphisms
(where $\X(M)$ is equipped with addition of vector fields).
\item{(2)}
Every $X \in \X^2(M)$ gives rise, in a natural way, to a diffeomorphism
$\tilde X:T^2 M \to T^2 M$, which, in a chart representation as
above, is described by
$$
\eqalign{
\tilde X(x+\eps_1v_1 + \eps_2 v_2+\eps_1 \eps_2 v_{12}) & =
x+\eps_1(v_1+X_1(x)) +\eps_2 (v_2+X_2(x)) + \cr
& \quad \eps_1 \eps_2 (v_{12} +
X_{12}(x) + dX_1(x) v_2 + dX_2(x) v_1). \cr}
$$
\item{(3)}
The map $\X^2(M) \to \Diff(TTM)$, $X \mapsto \tilde X$ is an injective
group homomorphism.

\Proof. 
(2) We fix a point $p \in M$ and define $\tilde X$ on a chart neighborhood
$T^2 U$ of $p$ by the formula given in the claim.
We then have to show that this definition is independent of the chosen chart.
This, in turn, amounts to proving that the map $X \mapsto \tilde X$ is
{\it natural} in the following sense: for all local diffeomorphisms $f$
defined on a neighborhood of $x$, we have $f_* \tilde X = \tilde{(f_* X)}$,
where $f_*$ is the natural action of $f$ on $\Diff(TTM)$, resp.\ on
$\X^2(M)$ by $f_* h = T^2 f \circ h \circ (T^2 f)^{-1}$, resp.\ by
$f_* X =T^2 f \circ X \circ f^{-1}$.
In order to get rid of the inverse, we show, more generally, that,
if $X,Y \in \X^2(M)$ are {\it $f$-related} (i.e.\
$T^2 f \circ X = Y \circ f$) for a smooth map $f:U \to M$, then
$T^2 f \circ \tilde X = \tilde Y \circ T^2 f$.
Now, 
$$
\eqalign{
T^2 f  (\tilde X & \,(x+\eps_1v_1 + \eps_2 v_2+\eps_1 \eps_2 v_{12}))  = \cr
& f(x)+\eps_1(df(x)(v_1+X_1(x))) + \eps_2 (df(x)(v_2+X_2(x)))  
+  \cr
& \quad \eps_1 \eps_2 \cdot \Bigl( df(x)(v_{12} +
X_{12}(x) + dX_1(x) v_2 + dX_2(x) v_1) + \cr
& \quad \quad \quad \quad 
d^2 f(x)(v_1 + X_1(x),v_2+X_2(x)) \Bigr), \cr
\tilde Y  ( T^2&  f  (x+\eps_1v_1 + \eps_2 v_2+\eps_1 \eps_2 v_{12}))  = \cr
& f(x) + \eps_1 (df(x)v_1+Y_1(f(x))) + \eps_2(df(x)v_2 + Y_2(f(x))) 
+  \cr
& \eps_1 \eps_2 \cdot \Bigl( d^2 f(x)(v_1,v_2) + df(x) v_{12} +
dY_1(f(x)) df(x) v_2 + dY_2(f(x)) df(x) v_1 \Bigr), \cr}
$$
and both expressions are equal if $X$ and $Y$ are $f$-related. 
Thus the map $\tilde X:TTM \to TTM$ is well-defined, and the construction
is natural in the
sense explained above. Let us prove that $\tilde X$ is a diffeomorphism.
In fact, a direct check shows that
its inverse is $\tilde{X^{-1}}$ with a section $X^{-1} :M \to TTM$ defined by
$$
\eqalign{
X^{-1}(x) & = x - \eps_1 X_1(x) - \eps_2 X_2(x) - \cr
& \quad \eps_1 \eps_2 \Bigl(X_{12}(x)
- (dX_1(x) X_2(x) + dX_2(x) X_1(x))\Bigr). \cr}
\eqno (14.2)
$$

\ssk (1)
By comparing the formulas from parts (1) and (2) of the claim, we see
that the product is defined by $(X \cdot Y)(x)=\tilde X(Y(x))$, and hence
is chart-independent.
By a straightforward computation, one shows that the product is associative:
$((X\cdot Y)\cdot Z)(x)=(X \cdot (Y \cdot Z))(x)$, and that
$X^{-1}$ as in (14.2) is an inverse of $X$ as in (14.1).
Moreover, it is clear from the explicit formula
 that, if $X,Y$ are sections of $J^2 M$, i.e.,
$X_1=X_2$ and $Y_1=Y_2$, then $X \cdot Y$
is again a section of $J^2 M$ over $U$,
 and similarly for the inverse, hence the
sections of $J^2 M$ form a subgroup of $\X^2 (M)$.
Finally, the explicit formula also shows that, if $X$ and $Y$        
are vector fields, seen as elements of $\X^2(M)$ in one of the three
canonical ways, then $X \cdot Y$ simply corresponds to the sum $X+Y$
of vector fields.

\ssk (3)
Choosing $Z$ such that $Z(x)=v$, we get
$$
\tilde{XY}(v)=((XY)Z)(x)=(X(YZ))(x)=\tilde X \circ \tilde Y(v).
$$
Clearly, $\tilde 0 = \id_{TTM}$, and hence we have a group action.
The action is faithful since, if  $\tilde X(x)=x$ for all $x\in U$,
then $X_1(x)=0=X_2(x)$ and thus also $X_{12}(x)=0$ for all $x \in U$.
\qed

\msk \nin For a conceptual, chart-independent definition of the
group structure on the space $\Gamma(M,T^k M)$ of sections
of $T^k M$ over $M$, see Theorem 28.2.

\Theorem 14.4. The group commutator in the group $\X^2(M)$ and the
Lie bracket in $\X(M)$ are related via
$$
[\eps_1 X, \eps_2 Y]_{\X^2(M)} = \eps_1 \eps_2 [X,Y]
$$
for all $X,Y \in \X(M)$. Equivalently, the group commutator in the group
$\Diff(TTM)$ and the Lie bracket in $\X(M)$ are related via
$$
[\tilde{\eps_1 X}, \tilde{\eps_2 Y}] = \tilde {\eps_1 \eps_2 [X,Y]}.
$$
  
\Proof. We fix $x \in M$.
 Let $X:M \to TM$ be a vector field. It gives rise, via the three
imbeddings $\X(M) \to \X^2(M)$, to three diffeomorphisms of $TTM$.
Specialising the chart formula from Theorem 14.3 (2), these three
diffeomorphisms are described by their values at the point
$u = x+\epsilon_1 v_1 + \epsilon_2 v_2 + \epsilon_1 \epsilon_2
v_{12}$:
$$
\eqalign{
 \tilde{\eps_1 X} (u)  & =
x+ \epsilon_1 (v_1+X(x)) + \epsilon_2 v_2  + \epsilon_1 \epsilon_2 
(v_{12} + dX(x) v_2) \cr 
\tilde{\eps_2 X} (u)  & =
x+ \epsilon_1 v_1+ \epsilon_2(v_2 + X(x))  + \epsilon_1 \epsilon_2 (v_{12} 
+ dX(x) v_1) \cr
\tilde{\eps_1 \eps_2 X} 
(u)  &  =
x+ \epsilon_1 v_1 + \epsilon_2 v_2  + \epsilon_1 \epsilon_2 
(v_{12} + X(x)).\cr}
\eqno (14.3)
$$
We observe that these diffeomorphisms belong to the groups
$\Gamma^2,\Gamma^1$ resp.\ $\Gamma^{12}$ defined in Section 14.1.
(Without using Theorem 14.3, these diffeomorphisms may be defined as follows:
the first two 
maps may also be seen as the tangent maps $T_{\eps_i} \tilde X$
of the infinitesimal automorphism $\tilde X:TM \to TM$,
and the third map is defined by translation in direction of the
vertical bundle $VM \subset TTM$, and such translations are canonical
in the bilinear space $(TTM)_x$, cf.\ Section BA.2.)
Thus, if $X$ and $Y$ are vector fields,
both diffeomorphisms $\tilde{\eps_1 X}$ and $\tilde{\eps_2 Y}$
preserve the fiber $(TTM)_x$ and are represented in a chart
$V \times V \times V$ of this fiber by bijections 
$g$ and $h$ of the fiber given by
$$
g(v,x',v')=(v+a,x',v'+\alpha(x')), \quad
h(v,x',v')=(v,x'+b,v'+\beta(v))
$$
with $a=X(x)$, $b=Y(x)$, $\alpha=dX(x)$, $\beta=dY(x)$.
Then we obtain for the commutator
$$
\eqalign{
[g,h](v,x',v')&=
gh(v-a,x'-b,v'-\alpha(x'-b)-\beta(v)) \cr
& = g(v-a,x',v'-\alpha(x'-b)\beta(v)+\beta(v-a)) \cr
&= g(v,x',v'-\alpha(x')+\alpha(b)-\beta(a)+\alpha(x')) \cr
&= (v,x',v'+\alpha(b)-\beta(a)). \cr}
$$
The value at $0_x$ is therefore 
$$
\eqalign{
[g,h](0,0,0) & =(0,0,\alpha(b)-\beta(a)) 
=(0,0,dX(x)Y(x)-dY(x)X(x)) \cr 
&  =(0,0,[X,Y](x)) = \eps_1 \eps_2 [X,Y](x) \cr}
$$
by definition of the Lie bracket in Theorem 4.2.
\qed

\nin {\bf 14.5.} {\sl Remark on definition of the Lie bracket and 
the Jacobi identity.} 
The definition of the Lie bracket via Theorem 4.2 is not very conceptual.
Theorem 14.4 offers a more conceptual way of defining it: {\it
for two vector fields $X,Y \in \X(M)$, there exists a unique vector
field $[X,Y] \in \X(M)$ such that the group commutator 
$[\tilde{\eps_1 X}, \tilde{\eps_2 Y}]$ agrees with  
$\tilde {\eps_1 \eps_2 [X,Y]}$.} It is then easily seen that $[X,Y]$
depends $\K$-bilinearly on $X$ and $Y$ and that $[X,X]=0$.                                                   
One would like, then, to prove the Jacobi identity by intrinsic
arguments not involving chart computations. For this, it is necessary
to invoke the third order tangent bundle $T^3 M$ and the natural
group structure on the space $\X^3(M)$ of its sections (Theorem 28.2). 
Then the Jacobi identity can be proved in the same way as will be done
for the Lie algebra of a Lie group in Section 24.4.

\ssk
Essentially, this strategy of defining the Lie bracket and proving the
Jacobi identity is the one used in synthetic differential geometry,
 and it also corresponds
to the definition of the Lie algebra of a {\it formal group}
(see [Se65]). In the framework of formal groups, the Jacobi identity
is obtained by a third order computation from 
Hall's identity for iterated group commutators,
$$
[[u,v],w^u] [[w,u],v^w] [[v,w],u^v] = 1,
$$
(here, $x^y=yxy^{-1}$)
which is valid in any group (cf.\ [Se65,  p.\ LG 4.18]).
Compared to our framework, this proof rather corresponds to using 
the bundle $J^3 M$ 
 which is more complicated than $T^3 M$ since it does not have ``axes".
Similar arguments are used for the proof of the Jacobi identity
in synthetic differential geometry, cf.\
[Lav87, p.\ 74 ff] and [MR91, p.\ 187/88].

\vfill \eject 

%% file: Dgeo2.tex

\def\date{25.10.2006}
\def\title{IV. THIRD AND HIGHER ORDER GEOMETRY} 
\def\rightheadline{\hfil{\eightrm \title}\hfil\tenrm\folio}

\sectionheadline{
IV. Third and higher order differential geometry}

\bigskip
\nin
Having linearized the bundle $TF$ over $M$ by means of a linear
connection, our plan of attack is to linearize in a similar way all
higher tangent bundles $T^k F$ over $M$. 
It will turn out that there is an essentially canonical procedure to
do so, starting with a connection on $F$ and one on $TM$.
The procedure is canonical up to a permutation acting on $T^k F$;
in particular, for $k=2$, we get two versions of this procedure.
Their difference can be interpreted as the curvature of $L$.

\subheadline
{15. The structure of $T^k F$: Multilinear bundles}

\msk \nin
{\bf 15.1.} {\sl Higher order models of $F$:  projection and injection cube.}
For a vector bundle $p:F \to M$,
we consider the fiber bundle $T^k F \to M$ as a ``$k$-th order model of $F$.''
The case $F=TM$ deserves special attention: $T^k F = T^{k+1}M$ may be seen
as the model of order $k+1$ of $M$. 
The  bundle $T^k F$ comes with a natural
``$k+1$-dimensional cube of projections'': e.g., for
$TTF$, resp.\ for $T^3 M$,
this is the commutative diagram of projections
$$
\xymatrix@dl@!0@C+10pt{
&F \ar@{-}[rr]\ar@{-}'[d][dd] & & M \ar@{-}[dd]
\\
\epsilon_2 TF \ar@{-}[ur]\ar@{-}[rr]\ar@{-}[dd]
& & \epsilon_2 TM \ar@{-}[ur]\ar@{-}[dd]
\\
& \eps_1 TF \ar@{-}'[r][rr] & & \epsilon_1 TM
\\
TTF \ar@{-}[rr]\ar@{-}[ur] & & TTM \ar@{-}[ur]
}
\quad \quad \quad \quad
\xymatrix@dl@!0@C+10pt{
&\eps_0 TM \ar@{-}[rr]\ar@{-}'[d][dd] & & M \ar@{-}[dd]
\\
\eps_0 \epsilon_2 T^2M \ar@{-}[ur]\ar@{-}[rr]\ar@{-}[dd]
& & \epsilon_2 TM \ar@{-}[ur]\ar@{-}[dd]
\\
& \eps_0 \eps_1 T^2M \ar@{-}'[r][rr] & & \epsilon_1 TM
\\
\eps_0 \eps_1 \eps_2 T^3 M \ar@{-}[rr]\ar@{-}[ur] & & \eps_1 \eps_2
T^2M \ar@{-}[ur]
}
\eqno {(15.1)}
$$
where the notation 
 $\eps_i TM$ and $\eps_i TF$ has the same meaning as $T_{\eps_i} M$ and
$T_{\eps_i} F$ in Section 7.4.  In case
$F=TM$ we write also $F=\eps_0 TM$.
 Each vertice of the cube is the projection
of a vector bundle onto its base, but the composition of more than one
of such projections is just the projection of a fiber bundle onto its base.
In particular, $T^kF$ is a fiber bundle over $M$. The fiber over $x \in M$
will be denoted by $(T^kF)_x$.
For general $k$, the cube for $T^{k+1} M$ will have vertices
$\eps^\alpha T^{|\alpha|} M$ (notation for multi-indices
$\alpha \in \{0,1\}^{k+1}$  being as in Section 7.4, resp.\ MA.1),
 and the cube for
$T^k F$ will have vertices $\eps^\alpha T^{|\alpha|} M$ and
$\eps^\alpha T^{|\alpha|} F$ with $\alpha \in \{ 0,1 \}^k$; it contains
the $k$-dimensional cube for $T^k M$ as a sub-cube.
Diagram (15.1) can also be interpreted as a diagram of zero
 sections, just by reading the vertices as upward arrows.
In particular, inclusion of spaces of the second row from below  
in the top space defines  $k+1$ injections
$$
\eps_i TM \to T^k F \quad (i=1,\ldots,k), \quad \quad \quad F \to T^k F
\eqno (15.2) 
$$
which, for $x \in M$ fixed, define $k+1$ ``basic axes'' in $(T^k F)_x$ 
which will be denoted  by $\eps_i T_xM$, $F_x$. Altogether there are
$2^{k+1}-1$ ``axes in $(T^k F)_x$'' whose definition will be given below.

\msk \nin {\bf 15.2.} {\sl The axes-bundle.} For a multi-index 
$\alpha = (\alpha_0,\alpha_1,\ldots,\alpha_k) \in I_{k+1}=
\{ 0,1 \}^{k+1}$, we denote by
$$
F^\alpha := \Big\{ \matrix{ \eps^\alpha F  & {\rm if} & \alpha_0 = 1 \cr
\eps^\alpha TM & {\rm if} & \alpha_0 = 0\cr}
\eqno (15.3)
$$
a copy of $F$, resp.\ of $TM$ (for the moment, $\eps^\alpha$ is just
a label to distinguish the various copies of $F$ and of $TM$)
and we define the {\it axes-bundle} to be the direct sum (over $M$) of
all of these copies of $F$ and of $TM$:
$$
A^k F := \bigoplus_{\alpha \in I_{k+1} \atop \alpha >0} F^\alpha
\eqno (15.4)
$$
We let also
$$
A^{k+1} M := 
 \bigoplus_{\alpha \in I_{k+1} \atop \alpha >0} \eps^\alpha T M.
\eqno (15.5)
$$
These are  direct sums of vector bundles over $M$ and hence are itself
a vector bundles over $M$, in contrast to $T^k F$ and to $T^{k+1} M$
 which do not carry
a canonical vector bundle structure.

\msk\nin {\bf 15.3.} {\sl Isomorphism $T^k F \cong A^k F$ over a chart
domain.}
Assume $\phi_i:M \supset U_i \to V$ is a chart such that $\phi_i(x)=0$.
Then $T^k \phi_i:T^k M \supset T^k U_i \to T^k V$ is a bundle chart of
$T^k M$ defining a bijection of fibers
$$
(\Phi_i)_x:=(T^k \phi_i)_x:(T^k M)_x  \to (T^k V)_0.
$$
Taking another chart $\phi_j$ with $\phi_j(x)=0$, we have the
following commutative diagram of bijections, describing the transition
functions of the atlas of $T^k M$ induced from the one of $M$:
$$
\matrix{  &  & (T^k M)_x & & \cr
& {\buildrel (T^k \phi_j)_x \over \swarrow} & & 
{\buildrel (T^k \phi_i)_x \over \searrow} & \cr
(T^k V)_0 & & {\buildrel (T^k \phi_{ij})_0 \over \longrightarrow} & &
(T^k V)_0 \cr}
\eqno (15.6)
$$
The transition functions of $T^k M$,
$$
(\Phi_{ij})_0:=(T^k \phi_{ij})_0:  (T^k V)_0 = \oplus_{\alpha >0} \eps^\alpha
V_\alpha \to  \oplus_{\alpha >0} \eps^\alpha
V_\alpha
$$
 are not linear: they are given by the Tangent Map Formula (Theorem 7.5), with
$f=\phi_{ij}$, $f(0)=0$, 
$$
T^k f(0 + \sum_\alpha \eps^\alpha v_\alpha) 
= \sum_\alpha \eps^\alpha \sum_{\ell=1}^{|\alpha|}
\sum_{\lambda \in {\cal P}_\ell(\alpha)} b^\lambda(v_{\lambda^1},\ldots
v_{\lambda^\ell})
$$
where, for a partition $\lambda$ of $\alpha$ of length $\ell$, we have written
 $b^\lambda$ for the multilinear map $d^\ell f(0)$.
This shows that $(\Phi_{ij})_0$ is a {\it multilinear map} in the sense
of Section MA.5. Thus transition functions of $T^k M$ are multilinear
maps, and $T^k M$ is a {\it multilinear bundle} in the following sense:

\msk \nin {\bf 15.4.} {\sl Multilinear bundles.}
A {\it multilinear bundle (with base $M$, of degree $k$)}
 is a fiber bundle $E$ over
$M$ together with a family $(E_\alpha)_{\alpha \in I_k}$ of vector bundles
and
 a bundle atlas such that, in a bundle chart around $x \in M$,
the fiber $E_x$ carries a  structure of a multilinear space  
over $\K$ whose axes are given by natural inclusions $E_\alpha \to E$, 
and change of charts, restricted to fibers, is by isomorphisms
of multilinear spaces. For $k=1$, a multilinear bundle is simply a
vector bundle, and for $k=2$ it is bilinear bundle in the sense of
Section 9.6.

\Theorem 15.5. If $F$ is a vector bundle over $M$, then 
the iterated tangent
 bundle $T^k F $ carries a canonical structure of a multilinear
bundle over $M$. In particular, for all $p \in M$,
the general multilinear group $\Gm^{0,k+1}((T^k F)_p)$ of the fiber over $p$
is intrinsically defined, and there are canonical inclusions of axes
$F^\alpha \subset T^k F$.

\Proof.
The arguments are essentially the same as the ones used above in the
special case $F=TM$:
bundle charts of $F$ are of the form $\tilde g_i:F \supset \tilde U_i \to
V \times W$, with transition functions $\tilde g_{ij}(x,w)=
(\phi_{ij}(x), g_{ij}(x)w)$, where $g_{ij}(x):W \to W$ is a linear map.
We assume $\phi_i(p)=\phi_j(p)=0$, whence $\phi_{ij}(0)=0$.
Bundle charts of $T^k F$ are of the form $T^k \tilde g_{ij}$, and  
as an analog of (15.6) we have the
following commutative diagram of bijections:
$$
\matrix{  &  & (T^k F)_p & & \cr
& {\buildrel (T^k \tilde g_j)_p \over \swarrow} & & 
{\buildrel (T^k \tilde g_i)_p \over \searrow} & \cr
(T^k V)_0 \times T^k W & &
 {\buildrel (T^k \tilde g_{ij})_0 \over \longrightarrow} & &
(T^k V)_0 \times T^k W. \cr}
\eqno (15.7)
$$
For $k=1$, the chart expression has been calculated in Section 9.1.
In order to give the general chart expression,
let $\tilde U \subset F$ be a bundle chart domain with base
$U \subset M$ and let $W$ be the model space for the fibers
of $F$. Then $T^k(\tilde U) \subset T^k F$ is a bundle chart domain
for $T^k F$, and elements of $T^k F_p$ are represented in the form
$$
u =  \sum_{\alpha > 0}  \eps^\alpha v_\alpha
 \quad \quad  \eps^\alpha v_\alpha \in V_\alpha =
\Big\{ \matrix{ \eps^\alpha W & {\rm if} & \alpha_0 = 1 \cr
\eps^\alpha V & {\rm if} & \alpha_0 = 0.\cr}
\eqno (15.8)
$$
For $k=2$, we may also use the following two equivalent
notations (cf.\ Chapter 7):
$$
\eqalign{
u & = \eps^{001} w_{001} + \eps^{010} v_{010} + \eps^{011} w_{011} +
\eps^{100} v_{100} + \eps^{101} w_{101} + \eps^{110} v_{110} +
\eps^{111} w_{111} \cr
& =  
\eps_0 w_0 + \eps_1 v_1  +\eps_0 \eps_1 w_{01} + \eps_2 v_2 +
 \eps_0 \eps_2 w_{02} + \eps_1 \eps_2 v_{12} + 
\eps_0 \eps_1 \eps_2 w_{012} \cr}
\eqno (15.9)
$$
where we replaced the letter
$v$ by $w$ whenever the argument belongs to $W$. Now,
the general formula for the transition functions is gotten from
Theorem 7.5: we let
$$
\eqalign{
c^m_x & :=d^m \phi_{ij}(x),\cr
b^{m+1}_x(w,v_1,\ldots,v_m)& := \ldots :=
b^{m+1}_x(v_1,\ldots,v_m,w):=
\partial_{v_1}\cdots \partial_{v_m} g_{ij}(x,w). \cr}
$$
For $x=p$, we get                                 
$$
\eqalign{
T^k \tilde g_{ij}( \sum_{\alpha_0 =1}\eps^\alpha &  w_\alpha +
\sum_{\alpha_0 =0} \eps^\alpha  v_\alpha)  = 
\sum_{\alpha_0 =1} \eps^\alpha g_{ij}(0) w_\alpha +
\sum_{\alpha_0 =0}  \eps^\alpha  d\phi_{ij}(0) v_\alpha    \cr  
& + \sum_{m=2}^{k+1}
\Bigl(
\sum_{ \alpha_0 = 1 \atop |\alpha|=m}
\eps^\alpha    \sum_{\ell=2}^m \sum_{\beta^1+\ldots+\beta^\ell= \alpha \atop
\beta^1 < \ldots < \beta^\ell} b^\ell_0(w_{\beta^0},v_{\beta^1},\ldots,
v_{\beta^\ell}) 
\cr & \quad \quad \quad 
+ \sum_{\alpha_0 = 0 \atop |\alpha|=m} 
\eps^\alpha
\sum_{\ell=2}^m \sum_{\beta^1+\ldots+\beta^\ell= \alpha \atop
\beta^1 < \ldots < \beta^\ell}
c^\ell_0(v_{\beta^1},\ldots,v_{\beta^{\ell}})  \Bigr)
\cr}
\eqno (15.10)
$$
(For $m=k+1$, only terms of the first sum over $\alpha$ actually appear, and
if $V=TM$, then this is the expression for $T^{k+1}\phi_{ij}$,
with  $\epsilon_0$ being an additional infinitesimal unit.)
Since $g_{ij}(0)$ and $d\phi_{ij}(0)$ are invertible linear maps,
comparing with Section MA.5, we see that $T^k \tilde g_{ij}$ belongs
to the general multilinear group $\Gm^{0,k+1}(E)$ of the fiber $E =
\sum_{\alpha >0} V_\alpha$. 
From this our claim follows in the way already described in Section 9.3.
\qed 

\msk
\nin As for any multilinear space,  the {\it special multilinear
group} $\Gm^{1,k+1}(E)$ of the fiber $E$ is now well-defined.
In terms of the preceding proof this means that,  
for a fixed point $p \in M$, we may, ``without loss of generality'',
reduce to the case where $g_{ij}(0)=\id_W$ and $d\phi_{ij}(0)=\id_V$.
The special multilinear group will play the most important r\^ole in the
sequel since it is the group acting simply transitively on the space
of multilinearly related linear structures (cf.\ Section MA.7, MA.9).
More generally, all the groups $\Gm^{j,k+1}((T^k)_x)$ with $j=0,1,\ldots,k+1$
(see Section MA.5 and Theorem MA.6) are well-defined.

\vfill\eject 

\subheadline
{16. The structure of $T^k F$: Multilinear connections}

\msk \nin {\bf 16.1.} {\sl Definition of multilinear
connections and their homomorphisms.}
Recall that a linear structure $L$ on $T^k F$ over $M$ is simply
given by two smooth structure maps, $T^k F\times_M T^k F\to T^k F$ and
$\K \times T^k F \to T^k F$, defining a $\K$-module structure
on each fiber of $T^k F$ over $M$. 
A {\it multilinear connection on $T^k F$} is a linear structure $L$ that
is {\it multilinearly related} (i.e.\ conjugate
under $\Gm^{0,k}(T^k F)_x$) to all linear structures induced from charts.
(As in case $k=2$ it is clear from Theorem 15.4 that this notation is 
well-defined.) Assume $K$ is a multilinear connection  on $T^k F$ and
$L$ a multilinear connection on $T^k G$ and $f:F \to G$ a vector bundle
homomorphism.  
We say that $f$ is a {\it $K$-$L$
homomorphism} if $T^k f:T^k F \to T^k G$ is fiberwise linear (with
respect to $K_x$ and $L_{f(x)}$). Clearly, the usual functorial rules
hold, and thus vector bundles with multilinear connections (of a fixed order
$k$) form a category. For $k=0$ this is just the category of vector bundles,
and for $k=1$ it is the category of vector bundles with linear connection.

\Theorem 16.2. {\rm (The generalized Dombrowski Splitting Theorem.)}
For any vector bundle $F$ over $M$, the following objects are in bijection
to each other:
\item{(1)}  multilinear connections $L$ on $T^k F$
\item{(2)} isomorphisms $\Phi: A^k F \to T^k F$ of multilinear
bundles over $M$ such
that the restriction of $\Phi$ to each axis $F^\alpha$ is the identity
map on this axis (here and in the sequel,
 the corresponding axes $F^\alpha$ of $A^k F$ and
of $T^k F$ are identified with each other).
\ssk
\nin
The bijection is described as follows: given $\Phi$, the linear structure
$L=L^\Phi$ is just the push-forward of the canonical linear structure $L_0$
on $A^k F$ by $\Phi$ . Given $L$, we let
$$
\Phi:=\Phi^L:  A^k F =\bigoplus_{\alpha > 0} F^\alpha \to T^k F, \quad
(x;v_\alpha)_{\alpha > 0} \mapsto
\sum_{\alpha >0} v_\alpha.
$$
In this formula,
the sum on the right hand side is taken in the fiber $(T^k F)_x$ with
respect to the linear structure $L_x$.

\Proof. Assume first that 
$\Phi: A^k F \to T^k F$ is any isomorphism of multilinear
bundles over $M$. By definition of an isomorphism, the push-forward
$L^\Phi$ of $L_0$ is then a multilinear connection on $T^k F$. 

Conversely, assume a multilinear connection $L$ on $T^k F$ be given and define
$\Phi=\Phi^L$ as in the claim. In order to prove
that $\Phi$ is a diffeomorphism, we need the
chart representation of multilinearly related
 linear structures:
with respect to a fixed chart, both the fiber $E:=(T^k F)_x$ and the
fiber $(A^k F)_x$ are  represented in
the form $E = \oplus_{\alpha>0} V_\alpha$ with
$
V_\alpha
$
as in Equation (15.8),
$V$ being the model space for the manifold $M$ and $W$ the one for the
fiber $F_x$. Since $A^k F$ is a vector bundle, the linear structure
of the fiber agrees with the one induced from the chart, and the
linear structure induced by $L$ 
 is multilinearly related to this linear structure:
it is obtained by push-forward via a
map $f_{b}$ which is the chart representation of $\Phi$, where $f_b$ is
defined as in Section MA.5:
 $b = (b_\lambda)_{\lambda \in {\rm Part}(I)}$ 
is a family of multilinear maps
$$
b^\lambda:=b^\lambda_x:
 V_{\lambda^{(1)}} \times \ldots \times V_{\lambda^{\ell(\lambda)}}
\to V_\alpha,
\eqno (16.1)
$$
 $\lambda$ is a partition of $\alpha$, and 
$V_\beta = \eps^\beta V$ if $\beta_0=0$,
$V_\beta = \eps^\beta W$ if $\beta_0=1$, and
for $v=\sum_{\alpha >0} v_\alpha$,
$$
f_b(v)=v + \sum_{\vert \alpha \vert \geq 2}
\sum_{\lambda \in {\cal P}(\alpha) \atop \ell(\lambda) \geq 2} 
b^\lambda(v_{\lambda^{(1)}}, \ldots,v_{\lambda^{\ell(\lambda)}} ).
\eqno (16.2)
$$
(The map $f_b=f_{b_x}$, resp.\ their components $(b^\lambda)=(b^\lambda_x)$,
will be called the {\it Christoffel tensors of $L$ in the given chart};
for $k=1$, there is just one component, and it agrees with the classical
Christoffel symbols of a connection, cf.\ Section 10.4.)
We claim that the linear structure $L$ is smooth if and only if, for all
bundle charts, the Christoffel tensors
$b$ depend smoothly on $x$ in the sense that
$$
x \mapsto  b^\lambda_x(v_{\lambda^{(1)}}, \ldots,v_{\lambda^{\ell(\lambda)}} )
$$
is smooth for all choices of the $v_{\lambda^{(i)}}$. 
In fact, it is clear that the linear structure is smooth if it is related
to chart structures via smooth $f_b$'s, and conversely,
$f_b$ and hence $b$ can be recovered from the vector
addition: in a chart, addition of elements from the axes is given by
$\sum_{\alpha >0}^{(b)} v_\alpha = f_b(v)$,
where the supersript $(b)$ in the first sum means ``addition with
respect to the linear structure given by $b_x$". Hence the map $f_b$
defined by a smooth linear structure is also smooth.
We have proved that $\Phi$ is smooth. It is bijective with smooth
inverse: in fact, $f_{b_x}$ belongs to the group
$\Gm^{1,k}(E)$ and hence is invertible (Theorem MA.6), and its inverse
is of the form $f_{c_x}$ with another multilinear map $c_x$ which can
be calculated by the algorithm used in the proof of Theorem MA.6.
Since this algorithm consists of a sequence of compositions of multilinear
maps and multiplications by $-1$, $c_x(v)$ depends again smoothly on $x$,
and hence $\Phi^{-1}$ is smooth.

Finally, it is immediate from the definitions that the constructions
$L \mapsto \Phi^L$ and $\Phi \mapsto L^\Phi$ are inverse of each other.
\qed

\nin 
The isomorphism $\Phi=\Phi^L$ is called the {\it linearization map of $L$}.
In case $F=TM$ is the tangent bundle, we have $F^\alpha = \eps^\alpha TM$
for all $\alpha$, and the linearization map can be written
$$
\Phi^L_k: \bigoplus_{\alpha \in I_k,\alpha > 0} \eps^\alpha TM \to T^k M.
$$

\msk
\nin {\bf 16.3.} {\sl Multilinear differential calculus on manifolds
with multilinear connection.}
We can define a ``matrix calculus'' for multilinear maps between
higher order tangent bundles $T^k M$, $T^k N$ equipped with multilinear
connections, simply by using fiberwise the purely algebraic ``matrix
calculus'' from Section MA.12. Basically, the choice of a multilinear
connection on $T^k M$ is the analog of the choice of a basis in a
(say, $n$-dimensional) vector space $E$; the splitting map
$\Phi^L:A^k F \to T^k F$ corresponds to the induced isomorphism $\K^n \to E$,
the ``matrix'' is the induced homorphism $\K^n \to \K^m$, individual
``matrix coefficients'' are defined by composing with the various
base-depending
projections $\K^m \to \K$ and injections $\K \to \K^n$, and finally
there is a ``matrix multiplication rule'' describing the matrix coefficients
of the matrix of a composition. 

\ssk Let us describe these concepts in more detail.
Assume $K$ is a multilinear connection on $T^k M$ and $L$ a multilinear
connection on $T^k N$ (the case of general vector bundles is treated
similarly), and $f:M \to N$ is an arbitrary
 smooth map. We define its {\it
multilinear higher differentials (with respect to $K$ and $L$)} by
$d^L_K f:= (\Phi^{L})^{-1} \circ  T^k f \circ \Phi^K: A^k M \to A^k N$, i.e.,
by the following commuting diagram:
$$
\matrix{ T^k M & {\buildrel \Phi^K \over \leftarrow} &
A^k M = \bigoplus_{\alpha >0} \eps^\alpha TM \cr
T^k f \downarrow \phantom{T^k f}& & \phantom{x}\downarrow d_K^L f \cr
 T^k N & {\buildrel \Phi^{L} \over \leftarrow} & A^k N =
\bigoplus_{\alpha >0} \eps^\alpha TN \cr}.
\eqno (16.3)
$$
 As for matrices of linear maps, we have the functorial rules
$$
d_{L_2}^{L_3} g \circ d_{L_1}^{L_2} f = d_{L_1}^{L_3} (g \circ f), \quad
d_L^L(\id) = \id
\eqno (16.4)
$$
which are immediate from the definitions. In each fiber, the maps
$$
(d_K^{L})_x f: (A^k M)_x \to (A^k N)_{f(x)}
$$
are homomorphisms of multilinear spaces (in the sense of MA.5) --
in fact, with respect to chart structures, this is true
by Theorem 7.5, and hence it is true with respect to any other structure 
that is multilinearly related to chart structures. It follows that
 the ``matrix coefficients'', indexed by partitions $\Lambda$ and $\Omega$,
$$
\eqalign{
d^\Lambda_x f & :=
(d_K^L f)_x^\Lambda:=\pr_\alpha \circ T^k f \circ \iota_{\Lambda}:
\bigoplus_{i=1}^{\ell(\Lambda)} \eps^{\Lambda^i} TM \to \eps^\alpha TM
\cr
(d_K^L f)_x^{\Omega | \Lambda} & :=
\pr_\Omega \circ d^L_K  f \circ \iota_{\Lambda}:
\bigoplus_{i=1}^{\ell(\Lambda)} \eps^{\Lambda^i} TM \to
\bigoplus_{i=1}^{\ell(\Omega)} \eps^{\Omega^i} TN \cr}
\eqno (16.5)
$$
are fiberwise multilinear maps in the usual sense; 
here $\pr$ and $\iota$ denote the canonical
projections and injections in the direct sum bundle $A^k N$ defined by
$$
\pr_{\Omega}: A^k N \to \bigoplus_{i=1}^{\ell(\Omega)} \eps^{\Omega^i} TN,
\quad \quad 
\iota_\Lambda: \bigoplus_{i=1}^{\ell(\Lambda)} \eps^{\Lambda^i} TM \to T^k M .
$$
 The whole map is recovered from the matrix entries via
$$
T^kf(x+\sum_\alpha \eps^\alpha v_\alpha) =f(x)+
\sum_\alpha \eps^\alpha \sum_{\Lambda \in {\cal P}(\alpha)}
(d^\Lambda f)_x(v_{\Lambda^1},\ldots,v_{\Lambda^\ell}),
\eqno (16.6)
$$
where the sum on the left-hand side is taken with respect to the linear
structure $K$, and on the right-hand side with respect to $L$.
Homomorphisms of multilinear connections
in the sense of 16.1 are characterized by the property that
the linearization maps commute with the natural action 
$f_*: A^k M \to A^k N$,
 $\sum_\alpha v_\alpha \mapsto \sum_\alpha Tf(v_\alpha)$:
$$
\matrix{ T^k M & {\buildrel \Phi^K \over \leftarrow} &
A^k M = \bigoplus_{\alpha >0} \eps^\alpha TM \cr
T^k f \downarrow \phantom{T^k f}& & \phantom{x}\downarrow f_* \cr
 T^k N & {\buildrel \Phi^{L} \over \leftarrow} & A^k N =
\bigoplus_{\alpha >0} \eps^\alpha TN \cr}.
\eqno (16.7)
$$
This, in turn, means that the matrix $d_K^L f$
is a ``diagonal matrix'', i.e.\  all components 
$(d_K^L f)_x^\Lambda$ for $\ell(\Lambda)>1$ vanish, and for $\ell(\Lambda)=1$
they agree with $Tf$. 

\msk \nin {\bf 16.4.} {\sl Composition rules and
``Matrix multiplication".}  Let us return to Equation (16.5).
Assume $g:M_1 \to M_2$ and $f:M_2 \to M_3$ are smooth maps and
$L^i$, for $i=1,2,3$, are multilinear connections on $T^k M_i$.
The ``matrix coefficients" of $f \circ g$ with respect to
$L^3$ and $L^1$ can be calculated from those of $f$ with respect
to $L^3$, $L^2$ and those of $g$ with respect to $L^2,L^1$
by using the composition rule  for multilinear maps
(MA.19), resp.\ Prop.\ MA.11:
$$
 (d_{L_3}^{L_1} (f \circ g))^{\Lambda} =
\sum_{\Omega \preceq \Lambda} (d_{L_2}^{L_1} f)^\Omega \circ
(d_{L_3}^{L_2} g)^{\Omega | \Lambda}.
\eqno (16.8)
$$
The notation $d^\Lambda$
generalizes the usual differential $d$ since, taking
$M=V$ with the flat connection of a linear space, the matrix coefficient
$d^\Lambda f$ as defined here coincides with the usual higher differential
$d^\ell f$ with $\ell$ being the length of $\Lambda$, cf.\ Theorem 7.5.

\msk \nin {\bf 16.5.} 
{\sl On existence of multilinear connections.} For $k>1$,
the space of all multilinear connections on $T^k F$ is no longer an 
affine space over $\K$ in any natural way
 (since $\Gm^{1,k}(E)$ is no longer a vector group),
and hence it is more difficult to construct multilinear connections
than to construct linear ones. They may be defined in charts by brute
force by choosing arbitrary ``Christoffel
symbols'' $b^\Lambda$ which then have to behave under chart changes according 
to the composition rule in the general multilinear group (16.8). 
In the sequel we will see that there are canonical procedures to construct
multilinear connections out of linear ones.

\vfill \eject 

\subheadline{17. Construction of multilinear connections}

\msk \nin {\bf 17.1.} {\sl Deriving linear structures.}
Assume $L=(a,m)$ is a linear structure on a
fiber bundle $E$ over $M$, with structural maps $a:E \times_M E \to E$,
$m:\K \times E \to E$. Then $TL=(Ta,Tm)$ defines a linear
structure on $TE$ over $TM$, called the {\it derived linear 
structure} (this is even a linear structure
with respect to $T\K$, but we will only consider it as a linear structure 
with respect to $\K$). 
However, in general
there is no canonical way to see the derived linear structure 
$TL$ as a linear structure on the bundle $TE$ over the smaller base
space $M$.

\msk \nin {\bf 17.2.} {\sl Deriving multilinear connections.}
Assume we are given the following data: a multilinear connection
$L$ on $T^k F$ with corresponding isomorphism of linear bundles 
$\Phi: A^k F =\bigoplus_{\alpha > 0} F^\alpha \to T^k F$,
as well as linear connections $L'$ on $TM$ and $L''$ on $F$.
Then we can construct, in a canonical way, an isomorphism $T^{k+1} F
\cong A^{k+1}F$. This is obtained in two steps:

\ssk 
\item{(a)} $T^{k+1} F \cong T(A^k F)$. This obtained by deriving $\Phi$, 
which gives a  bundle isomorphism
$$                                          
T\Phi:T(\bigoplus_{\alpha >0} F^\alpha)=T(A^k F)  \to T(T^k F)=T^{k+1} F.
$$
\item{(a)} $T(A^k F) \cong A^{k+1} F$.
Since $A^k F$ is a direct sum of certain copies of $F$ and of $TM$,
we can equip $A^k F$ with the direct sum  of the connections $L'$ (on
the copies of $TM$) and $L''$ (on the copies of $F$) --
see Section 10.10.  This gives us
a linear structure on $T(A^k F)$ and an isomorphism of linear bundles
$T(A^k F) \cong A^{k+1} F$.

\ssk \nin Putting (a) and (b) together, we end up with a bundle isomorphism
$T^{k+1} F \cong A^{k+1} F$ which, by transport of structure, can be
used to define a linear structure on $T^{k+1} F$.
Let us make this isomorphism more explicit:
in the following equalities,
we use first the isomorphism $T\Phi$, then
the canonical isomorphism (10.8) $T(A \oplus_M B) \cong TA \oplus_{TM} TB$,
 the isomorphisms
$\Phi_1$ corresponding to $L'$ and to $L''$ and finally
the canonical isomorphism (10.9)
$(A \oplus_M TM)\oplus_{TM} (B \oplus_M TM) \cong A \oplus_M TM \oplus_M B$.
By  $\oplus^{TM}$ we denote the
 direct sum of vector bundles over $TM$ and by $\oplus = \oplus_M$
the direct sum of vector bundles over $M$,
and $e_1,\ldots,e_k$ is the canonical basis of $I_k=\{ 0,1 \}^k$:
$$
\eqalign{
T^{k+1} F & \quad {\buildrel T\Phi \over\cong} \quad 
 T (\bigoplus_{\alpha \in I_k \atop \alpha >0}  F^\alpha ) 
\quad {\buildrel {(10.8)} \over \cong} \quad 
 \bigoplus_{\alpha \in I_k \atop \alpha >0}^{TM} T F^\alpha 
\cr
& \quad {\buildrel \Phi_1 \over\cong} \quad
 \bigoplus_{\alpha \in I_k \atop \alpha >0}^{TM} 
(F^\alpha \oplus \eps_{k+1} TM \oplus
\eps_{k+1} F^\alpha ) \cr
& \quad {\buildrel {(10.9)} \over\cong} \quad 
 \eps_{k+1} TM \oplus \bigoplus_{\alpha \in I_k \atop \alpha > 0} 
F^\alpha \oplus \bigoplus_{\alpha \in I_k \atop \alpha > 0} 
\eps_{k+1} F^\alpha  
\quad {\buildrel {\rm def.} \over \cong} \quad 
\bigoplus_{\alpha \in I_{k+1} \atop \alpha > 0} F^\alpha . \cr}
\eqno (17.1)
$$
Let us write $D\Phi$ for the isomorphism (17.1) of bundles over $M$.
Suppressing the canonical isomorphisms (10.8) and (10.9), we may write
$$
D\Phi = T \Phi \circ (\times_{\alpha \in I_k} \Phi_1^{e_{k+1},\alpha}):
\bigoplus_{\alpha \in I_{k+1} \atop \alpha > 0} F^\alpha \to T^{k+1} F,
\eqno (17.2)
$$
where $\Phi_1^{e_{k+1},\alpha}$ is the splitting isomorphism $\Phi_1$
for $TF^\alpha$. By transport of structure, $D\Phi$ defines 
on $T^{k+1} F$ the structure of a linear bundle over $M$ which we denote
by $DL$. 

\Theorem 17.3.
The linear structure $DL$ defined on $T^{k+1} F$ by $D\Phi$ is 
a multilinear connection on $T^{k+1}F$. 

\Proof.  As seen in the proof
of the Generalized Dombrowski Splitting Theorem 16.2, in a chart, $\Phi$
is given by (with 
$v= \sum_{\alpha \in I_{k} \atop \alpha > 0} \eps^\alpha v_\alpha$)
$$
\eqalign{
\Phi(x+v) & =x+ v + 
\sum_{\alpha \in I_k \atop \alpha >0} \eps^\alpha 
\sum_{\lambda \in {\cal P}(\alpha)\atop \ell(\lambda) \geq 2}
 b_x^\lambda (v_{\lambda^{1}},
\ldots,v_{\lambda^{\ell(\lambda)}} ) \cr
& =x+ 
\sum_{\alpha \in I_k \atop \alpha >0} \eps^\alpha 
\sum_{\lambda \in {\cal P}(\alpha)}
 b_x^\lambda (v_{\lambda^{1}},
\ldots,v_{\lambda^{\ell(\lambda)}} ) \cr}
$$
with multilinear maps $b_x^\lambda$, depending smoothly on $x$, and 
$b_x^\lambda(u)=u$ if $\ell(\lambda)=1$.
The tangent map $T\Phi$ is calculated by taking $\eps_{k+1}$ as new
infinitesimal unit:
$$
\eqalign{
T \Phi  (x+\sum_{\alpha \in I_{k+1} \atop \alpha > 0} \eps^\alpha v_\alpha) &
=
\Phi (x+\sum_{\alpha \in I_{k} \atop \alpha > 0} \eps^\alpha v_\alpha) + \cr
& \quad \quad
\eps^{k+1} 
\Bigl(
d\Phi(x+\sum_{\alpha \in I_{k} \atop \alpha > 0} \eps^\alpha v_\alpha)
\cdot \sum_{\alpha \in I_{k}} \eps^\alpha v_{\alpha + e_{k+1}} \Bigr) \cr 
&
= x+
\sum_{\alpha \in I_k \atop \alpha >0} \eps^\alpha 
\sum_{\lambda \in {\cal P}(\alpha)} 
b_x^\lambda (v_{\lambda^{1}},\ldots,v_{\lambda^{\ell(\lambda)}} ) \, \, + \cr
& \quad \quad \sum_{\alpha \in I_k \atop \alpha >0} \eps^{\alpha + e_{k+1}}
\sum_{\lambda \in {\cal P}(\alpha)} 
\Bigl(   \sum_{i=1}^{\ell(\lambda)} 
b_x^\lambda (v_{\lambda^{1}}, \ldots, v_{\lambda^{i} + e_{k+1}}, \ldots,
v_{\lambda^{\ell(\lambda)}} ) 
\cr & \quad \quad\quad \quad\quad \quad 
 + \partial_{v_{e_{k+1}}}b_x^\lambda
(v_{\lambda^{1}},\ldots,v_{\lambda^{\ell(\lambda)}})\Bigr)
 \cr}
$$
(the partial derivative in the last term is with respect to $x$;
the dependence on the other variables is linear and hence deriving amounts
to replace the argument $v_{\lambda^i}$ by the direction vector
$v_{\lambda^i + e_{k+1}}$).
%
%
Next, we have to insert the chart representation of the isomorphisms $\Phi_1$:
we have to replace, for all $\beta \in I_k$,
$v_{\beta + e_{k+1}}$ by $v_{\beta + e_{k+1}} + d_x(v_\beta,v_{e_{k+1}})$
where, if $\beta_0 = 0$,  $d_x$ is the Christoffel tensor of the connection
$L'$ on $TM$, and, if $\beta_0 = 1$, $d_x$ is the Christoffel tensor of the
 connection $L''$ on $F$. 
The resulting formula for $D\Phi$ is:
$$
\eqalign{
D\Phi & (x+\sum_{\alpha \in I_{k+1} \atop \alpha > 0} \eps^\alpha v_\alpha)=
x+ 
\sum_{\alpha \in I_k \atop \alpha >0} \eps^\alpha 
\sum_{\lambda \in {\cal P}(\alpha)} 
b_x^\lambda (v_{\lambda^{1}},\ldots,v_{\lambda^{\ell(\lambda)}} ) + \cr
&   \sum_{\alpha \in I_k \atop \alpha >0} \eps^{\alpha + e_{k+1}}
\sum_{\lambda \in {\cal P}(\alpha)} 
\Bigl(   \sum_{i=1}^{\ell(\lambda)} 
b_x^\lambda (v_{\lambda^{1}}, \ldots, v_{\lambda^{i} + e_{k+1}}, \ldots,
v_{\lambda^{\ell(\lambda)}} ) +
 \cr
& \quad \quad 
b_x^\lambda (v_{\lambda^{1}}, \ldots, d_x(v_{\lambda_i},v_{e_{k+1}}),
 \ldots, v_{\lambda^{\ell(\lambda)}} )
+ \partial_{v_{e_{k+1}}}b_x^\lambda
(v_{\lambda^{1}},\ldots,v_{\lambda^{\ell(\lambda)}})   \Bigr)
 \cr}
\eqno (17.3)
$$
This  expression is again a sum of
multilinear terms and hence defines a new multilinear map
$f_B=f_{B_x}$ where $B_x$ depends smoothly on $x$. Thus $DL$
is a multilinear connection on $T^{k+1} F$.
\qed 

\nin For later use, let us spell out the formula in case $k=2$, writing
$b_x(w,v)$ for the Christoffel tensor of $L$ and $c_x(u,v)$ for the one
of $L'$ in the chart, und using the more traditional notation
$v =  \eps_0 w_0 + \eps_1 v_1 + \eps_2 v_2 +
\eps_0 \eps_1 w_{01} + \eps_0 \eps_2 w_{02} + \eps_1 \eps_2 v_{12} + 
\eps_0 \eps_1 \eps_2 w_{012}$ (see Eqn.\ (15.9))
 for elements in the fiber $(TTF)_x$:
$$
\eqalign{D\Phi_1(x,v) & = v +
\eps_0 \eps_1 b_x(w_0,v_1) + \eps_0 \eps_2b_x(w_0,v_2) +
\eps_1 \eps_2 c_x(v_1,v_2) + \cr
&  \eps_0 \eps_1 \eps_2 \bigl( b_x(w_0,v_{12}) + b_x(w_{01},v_2)
+ b_x(w_{02},v_1) + \cr
&  \quad \quad \quad  b_x(w_0,c_x(v_1,v_2))+b_x(b_x(w_0,v_2),v_1)+
\partial_{v_2} b_x(w_0,v_1) \bigr).\cr}
\eqno (17.4)
$$
If $F=TM$, then this formula describes a linear structure on $T^3 M$,
depending on two linear connections $L,L'$ on $TM$. In this case, unless
otherwise stated, we will choose $L=L'$, i.e.\ $b_x = c_x$.

\msk
\nin {\bf 17.4.} The preceding construction has the property that
the zero section $z:T^k F \to T^{k+1}F=T_{\eps_{k+1}} T^k F$ is a
$L$-$DL$-linear
homomorphism of linear bundles over $M$;
more precisely, by construction
$L$ and $DL$ agree on the image of the zero section:
$$
\matrix{
T^k F & {\buildrel z \over \to} & T^{k+1} F \cr
\Phi \uparrow {\phantom \Phi} & & \phantom{D \Phi} \uparrow D\Phi \cr
A^k F & \subset & A^{k+1} F \cr}
$$
(this may also be seen from the chart representation (17.4), where
the image of the zero section corresponds to points $v$ with
$v_{\beta+e_{k+1}}=0$ for all $\beta$.)

\msk \nin {\bf 17.5.} {\sl The sequence of derived linear structures.}
Assume $L$ is a linear connection on the vector bundle
$p:F \to M$ and $L'$ a linear connection on $TM$.
Then, using the preceding construction, we can derive $L$ and construct
a connection $L_2:=DL$ on $TTF$ over $M$, and so on:
we get a sequence of derived linear structures $L_k=D^{k-1}L$
 on $T^k F$ over $M$,
$k=2,3,\ldots$ with corresponding isomorphisms $\Phi_k=D^{k-1}\Phi$
of linear bundles over $M$.
In case of tangent bundles, $F=TM$,
 then (unless otherwise stated) we will take $L=L'$ in this
construction, and the construction is functorial in the following sense:
if $f:M \to N$ is affine for given connections $L$ on $TM$, $L'$ on $TN$,
then $T^k f:T^k M \to T^k N$ is an $L_k$-$L_k'$ homomorphism in the sense
of 16.1 (and conversely).
 This is seen by an easy induction using that the construction
involves only natural isomorphisms and the tangent functor. 
For instance, $\gamma:\K \to M$ is a geodesic if and only if, for some
$k>1$, $T^k \gamma :T^k \K \to T^k M$ is a $L_k$-$L_k'$ homomorphism,
where $L_k$ is the natural connection on $\K$ given by the canonical
identification  $T^k \K = A^k \K$.

\ssk
However, the ``operator $D$ of covariant derivative'' has not the 
property of the usual derivative $d$ that higher derivatives are
independent of the order of the $\eps_i$: 
in general, $\Phi_k$ and $L_k$ will depend on the order $\eps_1,\eps_2,
\ldots,\eps_k$ of the subsequent scalar extensions leading from $F$ to
$T^k F$, and we denote by $L_k^\sigma $, $\sigma \in \Sigma_k$,
the linear structure obtained with respect to the order $\eps_{\sigma(1)},
\ldots,\eps_{\sigma(k)}$. 
Then the push-forward of $L_k$ by the canonical action of $\sigma$ on
$T^k F$ is given by
$$
\sigma \cdot L_k = L^\sigma_k  .
\eqno (17.5)
$$
These remarks make it necessary to study in more detail the behavior of
the linear structures $L_k$ with respect to the action of the symmetric
group (see next chapter).

\msk \nin
{\bf 17.6.} {\sl Example: trivializable tangent bundles.}
We will see later that the tangent bundle of a Lie group is always
trivializable, and thus the following example can be applied to the
case of an arbitrary Lie group:
assume that the tangent bundle $TM$ is trivializable over $M$, i.e.\ that
there is a diffeomorphism $\Phi_0:M \times V \to TM$ which is
fiber-preserving and  linear in
fibers over $M$. We define a diffeomorphism
$\Psi_1: (M \times V) \times (V \times V)  \to TTM$
by the following diagram:
$$
\matrix{
M \times V \times V \times V & {\buildrel \Psi_1 \over \to} & TTM \cr
\downarrow & & \phantom{T\Phi} \uparrow T\Phi_0 \cr
M \times V \times TV & {\buildrel \Phi_0 \times \id \over \to} &
 TM \times TV \cr}
$$
Then $\Psi_1$
is a trivialization of the bundle $TTM$ over $M$. We can continue in this
way and trivialize all higher tangent bundles $T^k M$ over $M$ via
$$
\Psi_{k+1}:=T^k \Phi_0 \circ (\Psi_k \times \id)=T\Psi_k \circ (\Phi_0 \times
\id): \quad 
M \times \bigoplus_{\alpha \in I_{k+1} \atop \alpha >0} \eps^\alpha V
\to T^{k+1} M.
$$
Via push-forward,
this defines linear structures $L_k$ on all higher order tangent bundles
$T^k M$. By our assumption on $\Phi_0$, $L_1$ is just the vector bundle
structure on $TM$. 

\Lemma 17.7. Assume the tangent bundle $TM$ is trivializable, and retain
notation from {\rm Section 17.6}. Then
the linear structure
$L_2$ is a connection on the tangent bundle $TM$,
and  $L_k = D^{k-2} L_2$ for $k>2$, i.e.,
the linear structures $L_k$ agree with the derived connections of $L_2$
defined in the preceding section.

\Proof.
In order to prove that $L_2$ is a connection,
in a bundle  chart, $\Phi_0$ is represented by $\Phi_0(x,w)=(x,g_x(w))=
(x,g(x,w))$; then 
$$
T\Phi_0 ((x,w)+\eps(x',w')) = (x,g_x(w)) 
\, +  \, \eps \bigl(g_x(w'),d_1g(x,w)x'\bigr) 
$$
(cf.\ Section 9.1) where the last term depends bilinearly on $(w,x')$, proving
that $L_2$ is a connection.
Explicitly,
the Christoffel tensor of this connection is given by Equation (9.3):
$$
b_x(w,x')= g_x^{-1} d_1 g(x,w)x'
$$
(Note that the bilinear map $b_x$ is in general not symmetric.)
Next, it follows that $L_3 = D L_2$ since both
$L_3$ and $DL_2$ are constructed by deriving $L_2$ and composing with
 $\Psi_2$, resp.\ with $\Phi_2$, and $\Psi_2$ and $\Phi_2$ induce the
 same linear structure $L_2$ on $TTM$. By induction, this argument
 shows that $L_k=D^{k-1} L_2$.
\qed

\vfill \eject

\subheadline {18. Curvature} 

\msk \nin {\bf 18.1.} {\sl Totally symmetric multilinear connections.} 
The symmetric group $\Sigma_k$ acts on $T^k F$; in the notation (15.8) 
this action is by permuting the last $k$ coordinates
of $\{ 0,1 \}^{k+1}$. If $F=TM$, then $T^k F=T^{k+1}M$, and 
$\eps_0$ may be considered as another infinitesimal unit,
and the permutation group $\Sigma_{k+1}$ acts on all $k+1$ coordinates.
We say that a multilinear connection $L$ is {\it totally} or {\it very
symmetric}
if it is invariant under the action of $\Sigma_k$, resp.\ of
$\Sigma_{k+1}$. Equivalently, this means that $\Sigma_k$ acts linearly
with respect to $L$, or that $L$ is fiberwise
 conjugate  to all chart-structures under the action of
{\it very symmetric group} 
$\Vsm^{1,k}((T^kF)_x) = (\Gm^{1,k}(T^k F))^{\Sigma_k}$  defined
in Section SA.2. In terms of the Dombrowski Splitting Map $\Phi:A^k F \to
T^k F$, the condition is equivalent to saying that $\Phi$ is equivariant
with respect to the natural $\Sigma_k$-actions on both sides.

\msk \nin {\bf 18.2.} {\sl Curvature operators and curvature forms.}
Curvature is the obstruction for  a multilinear connection
$L$ to be very symmetric, or, in
other words, it measures the difference between $L$ and
 the linear structure $L^\sigma=\sigma. L$ obtained by push-forward
via $\sigma \in \Sigma_k$.
In order to formalize this idea, let us denote the canonical linear structure
on the axes-bundle $A^k F$ by $L^0$; it is $\Sigma_k$-invariant.
The linear structure $L$ is obtained as push-forward of $L^0$ via $\Phi$, i.e.\
$L=\Phi.L^0$. It follows that $\sigma. L = \sigma.\Phi.L^0=
\sigma \Phi \sigma^{-1}.L^0$, and hence
the linear structure $L^\sigma$ is given by push-forward via
$\Phi^\sigma = \sigma \circ \Phi \circ \sigma^{-1}$. 
It is proved in Prop.\ SA.3 that $\Phi^\sigma$ is again a multilinear
map, and thus we see that $\sigma.L$ is again a multilinear connection.
Now we define the {\it curvature operator}
$R:=R^\sigma:=R^{\sigma,L}$ by
$$
R^{\sigma,L}:= \Phi^\sigma  \circ \Phi^{-1}=
\sigma \circ \Phi \circ \sigma^{-1} \circ \Phi^{-1}: T^kF \to T^kF, 
\eqno (18.1)
$$
respectively its action on the axes-bundle $\tilde R:=\tilde R^\sigma:=
\tilde R^{\sigma,L}$ given by
$$
\tilde R^{\sigma,L}:= \Phi^{-1} \circ R \circ \Phi =
\Phi^{-1}  \circ \sigma \circ  \Phi \circ
 \sigma^{-1} = \Phi^{-1} \circ \Phi^{\sigma^{-1}}:
A^kF \to A^kF
\eqno (18.2) 
$$
Both operators are compositions of multilinear maps and hence are
multilinear automorphisms of $T^k F$, resp.\ of $A^k F$.
The curvature operator $R$ can be characterized as
the unique  element of $\Gm^{1,k}(T^kF)$ relating the linear structures
 $L$ and $L^\sigma$:
$$
R \cdot L  = L^\sigma.
\eqno (18.3)
$$
In fact,
$R \cdot L = \sigma \circ \Phi \circ \sigma^{-1} \circ \Phi^{-1} \cdot L = 
\sigma \circ \Phi \circ \sigma^{-1}  \circ L_0 =
\sigma \circ \Phi \cdot L_0 = \sigma \cdot L$.
Immediately from the definition we get the following ``cocycle relations'' :
for fixed $L$ and for all $\sigma,\tau \in \Sigma_k$,
$$
R^{\id} = \id, \quad 
R^{\sigma \tau} = \sigma R^\tau \sigma^{-1} \circ R^\sigma.
\eqno (18.4)
$$
In fact, it is also possible to interprete the curvature operators as
higher differentials
(``matrices'') in the sense of Section 16.3, and then the cocycle relations
are seen as a special case of the matrix multiplication rule (16.5):
from (18.2) it follows that 
$$
\matrix{ T^k F & {\buildrel \Phi^{\sigma^{-1}} \over \leftarrow} & A^k F \cr
T^k \id_F  \downarrow \phantom{T^k f}& & \phantom{x}\downarrow \tilde R  \cr
 T^k F & {\buildrel \Phi  \over \leftarrow} & A^k F \cr}
$$
commutes, which shows that
$$
\tilde R^\sigma = d_{\sigma^{-1}.L}^L (\id_F)
\eqno (18.5)
$$
is a higher differential in the sense of 16.3.
The matrix coefficients of the curvature operators are called the
{\it curvature forms},
$$
\tilde R^{\Omega | \Lambda} := \pr_\Omega \circ \tilde R \circ \iota_\Lambda:
  F_{\Lambda^1} \oplus_M \ldots \oplus_M F_{\Lambda^l} \to
F_{\Omega^1} \oplus_M \ldots \oplus_M F_{\Omega^r}.
\eqno (18.6)
$$
Of particular interest are the {\it elementary curvature forms}
$$
\tilde R^{\Lambda} :=\tilde R^{\uline \Lambda | \Lambda} : 
  F_{\Lambda^1} \oplus_M \ldots \oplus_M F_{\Lambda^l} \to
F_{\uline \Lambda}
$$
which are ordinary tensor fields, i.e.\ fiberwise multilinear maps
in the usual sense. For the longest partition 
$\Lambda = (e_1,\ldots,e_k)$ we get the {\it top curvature form}.
In particular, for $k=3$ one finds the classically known curvature tensor
as a top curvature form:

\Theorem 18.3. Assume $L$ is a linear connection on the vector bundle
$F$ over $M$ and $L'$ a torsionfree connection on $TM$. 
Let $k=2$ and $R = R^{\kappa,DL}$ where $\kappa: T^2 F \to T^2 F$ is
the canonical flip and $DL$ is the connection derived from $L$ as
described in {\rm Theorem 17.3}. Then $R$ does not
depend on the choice of $L'$. Moreover,  $R$ belongs to the subgroup
$ \Gm^{2,3}(TTF)$, i.e.\
 $\tilde R_x$ has only one non-trivial
component, namely the top curvature form
$$
\tilde R^\Lambda:
\eps_1 T M \times \eps_2 T M \times\eps_0 F \to \eps_0 \eps_1 \eps_2.F
$$
where $\Lambda = (e_1,e_2,e_3)$.
In other words, for $v = (v_\alpha) \in (A^2 F)_x$,
$$
\tilde R_x(v) = v + \tilde R_x^\Lambda(\eps_1 v_1,\eps_2 v_2)\eps_0 v_0.
$$
In a chart, with Christoffel tensor $b_x$ of $L$, the
top curvature form  $\tilde R_x^\Lambda$ is given by
$$
\tilde R_x^\Lambda(v_1,v_2)w_0 = 
 b_x(b_x(w_0,v_2),v_1) -  b_x(b_x(w_0,v_1),v_2) +
\partial_{v_2} b_x(w_0,v_1) - \partial_{v_1} b_x(w_0,v_2) .
$$

\Proof.  We use the chart
representation (17.4) of the map $\Phi_2=D\Phi$ for $L_2=DL$:
the chart formula for $L_2^\kappa$ is obtained from this formula
by exchanging $\eps_1 \leftrightarrow \eps_2$, $v_1 \leftrightarrow v_2$, 
$w_{01} \leftrightarrow w_{02}$.
Note that the resulting formula differs from the one of $\Phi_2$ only in
the $\eps_0 \eps_1 \eps_2$-component and in the $\eps_1 \eps_2$-component, 
and the difference in
the $\eps_1 \eps_2$-component vanishes 
if $c(v_1,v_2)=c(v_2,v_1)$, i.e.\ if $L'$ is torsionfree. Let us assume
that this is the case.
Then from Formula (17.4)  we get the following chart formula for the
difference between $L_2$ and $L_2^\kappa$: 
$$
\eqalign{
&\Phi_2(v)-\Phi_2^\kappa (v) =
\eps_0\eps_1 \eps_2 \cdot \cr
& \quad \Bigl( b_x(b_x(w_0,v_2),v_1) -  b_x(b_x(w_0,v_1),v_2) +
\partial_{v_2} b_x(w_0,v_1) - \partial_{v_1} b_x(w_0,v_2) \Bigr)
\cr}
\eqno (18.7)
$$
since $c_x$ is assumed to be symmetric. This formula proves all claims
of the theorem: if we
 denote by $A_x(v_1,v_2)w_0$ the right-hand side of (18.7) and by $f_A$ the
multilinear map given by $f_A(v)=v+A(v_1,v_2)w_0$, then
from the composition rule of multilinear maps it follows that
$f_A  \cdot L_2 = L_2^\kappa$; thus  (18.3) is satisfied,
and by uniqueness in (18.3) we get $R=f_A$, and $A$ is the only
non-trivial component of the curvature maps $R$, resp.\ $\tilde R$.
Formula (18.7) shows that it is indeed independent of $c$ and hence of $L'$.
\qed

\nin
The chart representation (18.7) of the top curvature form
 coincides with the usual expression
of the curvature tensor of a connection in a chart by its 
Christoffel tensor (cf.\ e.g.\ [La99, p.\ 232]), thus showing that our
interpretation of the curvature agrees with all other known concepts.

\msk \nin {\bf 18.4.} {\sl Case of the tangent bundle.} 
Assume we are in the case of the tangent bundle: $F=TM$; then we usually
choose $L=L'$. If $L$ has torsion, then also $L'$ has torsion, and
the assumptions of Theorem 18.3 are not satisfied. However, Formula (17.4)
still can be used to establish a relation between the classical curvature
tensor and the curvature forms of $L$; we will not go into the details
and leave apart
the topic of connections with torsion. Let us just add the remark that
in this case the $\Sigma_3$-orbit 
of the linear structure $L_3$ on $T^3 M$ has $6$ elements which we may
represent by a hexagon. None of the 6 structures is invariant under a
transposition, thus we may represent
the three transpositions  by symmetries
with respect to three axes that do not pass through the vertices.

\msk \nin {\bf 18.5.} {\sl Bianchi's identity.} Now assume that $F=TM$
is the tangent bundle and 
that $L$ is torsionfree. Then, by torsionfreeness, 
 $L$ and$ DL$ are invariant under
the transposition $(12)$, and the usual curvature tensor is the
top curvature form of $R=R^{(23)}$. Let us look at the $\Sigma_3$-orbit
of the linear structure $K:= DL$ on $T^3M$.
This orbit contains at most 3 elements since $K$ is stabilized by the
transposition $(12)$. We denote these elements by
$K^1 := K = K^{\id}=K^{(12)}$, $K^2 :=K^{(23)}=K^{(132)}$ and
 $K^3:=K^{(13)}=K^{(123)}$. From the ``matrix multiplication rule" (16.5),
  we get
$$
\id = d_{K^1}^{K^1} \id = 
d_{K^2}^{K^1} \id \circ d_{K^3}^{K^2} \id \circ d_{K^1}^{K^3} \id.
$$
Taking the top curvature forms, we get the {\it Bianchi identity} for the
curvature form $R=R^{(23)}$
$$
R + R^{(132)} + R^{(123)} = 0.
$$
In a similar way, we get the usual skew-symmetry of the curvature form
$R=R^{(23)}$. More generally, one can show that the top degree
curvature forms $R^{(23),D^k L}$ for the derived linear structures $L_{k+1}$
agree with the covariant derivatives of the curvature tensor in the
usual sense, and under suitable assumptions
one can prove algebraic relations for these. 

\msk
The following result is included for the sake of completeness. It relates
the curvature as defined here to the usual definition in terms of
covariant derivatives:

\Theorem 18.6.
Assume $\nabla$ is the covariant derivative associated 
to the linear connection $L$ on the vector bundle $F$.
Then, for all sections $X,Y$ of $TM$ and $Z$ of $F$ and all $x \in M$, 
the curvature is given by
$$
R_x(X(x),Y(x))Z(x) = ([\nabla_X,\nabla_Y] Z -
 \nabla_{[X,Y]} Z)(x) .
$$
In particular, the right hand side depends only on the values of $X,Y,
Z$ at $x$.

\Proof.  The right hand side has the following chart expression
(the calculation given in [La99, p.\ 232] can be repeated verbatim): 
$$
b_x(b_x(Z,Y),X) - b_x(b_x(Z,X),Y) +\partial_Y b_x(Z,X) - \partial_X
b_x(Z,Y)
\eqno (18.8)
$$
where $b_x$ is the Christoffel tensor of $L$. This
coincides with the chart expression for the curvature $R$
given by Eqn.\ (18.7). \qed

\msk
\nin
{\bf 18.7.} {\sl Structure equations.}
%
Assume $L,L'$ are two connections on $F$ with covariant
derivatives $\nabla,  \nabla'$ and curvature tensors
$R,R'$; let $A:=L'-L=\nabla'-\nabla$
(tensor field of type $(2,1)$).
For $X,Y \in \X(M)$, we write $A_X(Y):=A(X,Y)$.
Then the ``new" curvature is calculated from the ``old" one via
$$
\eqalign{
R'(X,Y)&=[(\nabla+A)_X,(\nabla+A)_Y] -(\nabla + A)_{[X,Y]} \cr
&
=R(X,Y) + [A_X,\nabla_Y] - [A_Y,\nabla_X] - A_{[X,Y]} + [A_X,A_Y].\cr}
$$
Thus $R'$ depends on $A$ in a quadratic way.
In particular, the set of flat connections can be seen as the zero
set of a quadratic map.  If we assume that $F=TM$ and that $\nabla$
is torsionfree, then
we can write $[\nabla_X,A_Y]-A_{ \nabla_X Y} = (\nabla A)(X,Y)$,
and the preceding equation reads 
$$
R'(X,Y)=R(X,Y) + (\nabla A)((X,Y)-(Y,X)) + [A_X,A_Y].
\eqno (18.9)
$$
If, moreover $L$ was the flat connection of a chart domain,
then $R=0$, and the tensor
 $A:TM \times W \to W$ can be interpreted as an $\End(W)$-valued
one-form (where $W$ is the model space for the fibers of $F$).
The exterior derivative of $A$ is $\dd A(X,Y)=(\nabla A)((X,Y)-(Y,X))$,
which gives the second term.
The third term is symbollically written ${1 \over 2}[A,A]$, and we get
$$
R' = \dd A + {1\over 2} [A,A]
\eqno (18.10)
$$
which is a version of the structure equations of E.\ Cartan.
If also $L'$ is flat, we get the condition
$\dd A + {1 \over 2} [A,A]=0$ which is another version of the
structure equations.

\vfill \eject 

\subheadline{19. Linear structures on jet bundles}

\msk \nin {\bf 19.1.} {\sl Jet-compatible
 multilinear connections.} 
We say that a multilinear connection $L$ on $T^k F$
 is {\it weakly symmetric} or {\it jet-compatible} if the
jet-bundle $J^k F=(T^k F)^{\Sigma^k}$ is a 
{\it linear} subbundle of the linear bundle $(T^k F,L)$.
Equivalently, $L$ is conjugate to all chart-structures under the
 {\it symmetric multilinear group} $\Sm^{1,k}((T^kF)_x)$ defined
in Section SA.2.
In terms of the Dombrowski Splitting Map $\Phi:A^k F \to
T^k F$, the condition is equivalent to saying that 
$$
\Phi((A^k F)^{\Sigma_k}) = (T^k F)^{\Sigma_k}.
$$
 A totally symmetric multilinear connection is weakly
symmetric (since $\Vsm(E)\subset \Sm(E)$) but the converse is not true:
for instance, every connection on $TM$ is weakly symmetric since the chart
representation of $\Phi$ is $\Phi(x;u,v,w)=(x,u,v,w+b_x(v,w))$, 
and this preserves
the diagonal $\{ (x,u,u,w)|x \in U, \, u,w \in V \}$.
 On the other hand, $\Phi$ is not always
totally symmetric (it is totally 
symmetric if and only if $b_x$ is a symmetric bilinear map).
More generally, a connection may be weakly symmetric even if it
has non-trivial curvature:

\Theorem 19.2.
 Assume $L$ is a linear connection on $F$ and $L'$ a linear
connection on $TM$.
Then the derived multilinear connections $L_k=
D^{k-1} L$ on $T^k F$ are weakly symmetric.

\Proof.
The claim is proved by induction on $k$. The case $k=1$ is trivial
since $\Sigma_1 = \{ \id \}$.

The case $k=2$ can be proved by using the explicit formula (17.4)
for $\Phi_2$ in a chart: the vector $v$ is fixed under $\Sigma_2$
iff $v_1=v_2$ and $w_{01}=w_{02}$; then (17.4) shows that also
$\Phi(v)$ is fixed under $\Sigma_2$.

Now assume that the property 
$\Phi_k((A^k F)^{\Sigma_k}) = (T^k F)^{\Sigma_k}$
 holds for the linear structure
$L_k$ on $T^k F$.
In order to prove that 
$\Phi_{k+1}((A^{k+1} F)^{\Sigma_{k+1}}) \subset (T^{k+1} F)^{\Sigma_{k+1}}$,
note first that the set $R = \Sigma_k \cup \{ (k,k+1) \}$
generates $\Sigma_{k+1}$,
where  $\Sigma_k$ is imbedded in $\Sigma_{k+1}$ as the subgroup permuting 
$\eps_1,\ldots,\eps_k$. Thus the fixed point spaces are
$$
\eqalign{
(A^{k+1} F)^{\Sigma_{k+1}} & = (A^{k+1} F)^{\Sigma_{k}} \cap
(A^{k+1} F)^{(k+1,1)}, \cr
(T^{k+1} F)^{\Sigma_{k+1}} & = (T^{k+1} F)^{\Sigma_{k}} \cap
(T^{k+1} F)^{(k+1,1)},  \cr}
$$
and it suffices to show that

\ssk
\item{(a)}
 $\Phi_{k+1}((A^{k+1} F)^{\Sigma_k}) =
(T^{k+1} F)^{\Sigma_k}$ and
\item{(b)}
$\Phi((A^{k+1} F)^{(k+1,k)}) = (T^{k+1} F)^{(k+1,k)}$.

\ssk \nin 
In order to prove (a),
recall from Formula (17.2) that $\Phi_{k+1}$ is defined as a composition
of two maps:
$\Phi_{k+1} = T\Phi_k \circ  \oplus_\alpha \Phi_1^{e_{k+1},\alpha}$.
We claim that both maps preserve the subspace of $\Sigma_k$-invariants, i.e.,
the following diagram commutes:
$$
\matrix{A^{k+1} F & {\buildrel \oplus_\alpha \Phi_1^{e_{k+1},\alpha} \over \to}
& T(A^k F) & {\buildrel T\Phi_k \over \to} & T(T^k F) \cr
\uparrow & & \uparrow & & \uparrow \cr
(A^{k+1} F)^{\Sigma_k}
 & {\buildrel \oplus_\alpha \Phi_1^{e_{k+1},\alpha} \over \to}
& (T(A^k F))^{\Sigma_k} & {\buildrel T\Phi_k \over \to} & T(T^k F)^{\Sigma_k}
 \cr}
$$
In fact, for the square on the right of the diagram this follows  by
applying the tangent functor to the
induction hypothesis since the action  of $\tau \in \Sigma_k$ on
$T^{k+1}F$ is simply given by the tangent map  $T\tau$ of the
action on $T^k F$.
The map on the left of the diagram commutes also with the $\Sigma_k$-action;
this follows from the fact that $\Sigma_k$ does not act on $\eps_{k+1}$, hence
simply interchanges the various copies of $\Phi_1$ used in the construction,
and it can also be directly seen from the chart representation
$$
\eqalign{
 \oplus_\alpha \Phi_1^{e_{k+1},\alpha}:
\sum_{\alpha \in I_{k+1}} \eps^\alpha v_\alpha & \mapsto
\sum_{\alpha_{k+1}=0} \eps^\alpha v_\alpha + \eps_{k+1} v_{k+1} + \cr
& \quad \quad
\sum_{\alpha_{k+1}=0} \eps^{\alpha + e_{k+1}} (v_{\alpha + e_{k+1}} +
d_x(v_\alpha,v_{e_{k+1}} ) ). \cr}
$$
Next, we will show that the space of invariants of the transposition $(k,k+1)$
is preserved. In order to prove this, we will use another representation of
the map $\Phi_{k+1}$, 
introducing first a more detailed notation for (17.2): we write
$
\Phi_j^{e_{i_j},\ldots,e_{i_1};\beta} 
$
for the map $\Phi_j:A^j F^\beta \to T^k F^\beta$ taken with respect to
the sequence $\eps_{i_1},\ldots,\eps_{i_j}$ of scalar extensions 
of the bundle $F^\beta$.
Then our recursion formula for $\Phi_{k+1} = \Phi_k^{e_{k+1},\ldots,e_1;e_0}$ 
reads
$$
\Phi_{k+1} = T \Phi_k^{e_k,\ldots,e_1;e_0} \circ
(\times_{\alpha \in I_k \atop \alpha > 0} \Phi_1^{e_{k+1},\alpha}).
\eqno (19.1)
$$
Applying twice this recursion formula, we get another
 recursion formula for $\Phi_{k+1}$, this time
in terms of $\Phi_{k-1}$ and $\Phi_2$:
$$
\eqalign{
\Phi_{k+1}& = T \Phi_k^{e_k,\ldots,e_1;e_0} \circ
(\times_{\alpha \in I_k \atop \alpha > 0} \Phi_1^{e_{k+1};\alpha}) \cr
& = T^2 \Phi_{k-1}^{e_{k-1},\ldots,e_1;e_0} \circ
(\times_{\alpha \in I_{k-1} \atop \alpha >0} T \Phi_1^{e_k;\alpha})
\circ 
(\times_{\alpha \in I_k \atop \alpha > 0} \Phi_1^{e_{k+1};\alpha}) \cr
& = T^2 \Phi_{k-1}^{e_{k-1},\ldots,e_1;e_0} \circ
(\times_{\alpha \in I_{k-1} \atop \alpha >0  } \Phi_2^{e_{k+1},e_k;\alpha}).
\cr}
\eqno (19.2)
$$
Now, $T^2\Phi_{k-1}$, being a second tangent map with respect to scalar
extensions $\eps_{k+1},\eps_k$, commutes with the
transposition $(k,k+1)$.
As seen for the case $k=2$ (beginning of our induction),
 $\Phi_2^{e_{k+1},e_k;\alpha}$ preserves
the subspace fixed under $(k,k+1)$. Altogether, it follows that also
$\Phi_{k+1}$ preserves the subspace fixed under $(k,k+1)$, as had to be
shown.
\qed

\msk \nin
%

\msk \nin {\bf 19.3.} {\sl  Multilinear differential calculus for jets. }
If $L$ is a jet-compatible connection on $T^k F$, then structures such
as the linearization map or the higher differentials from Chapter 16 can be
restricted to the subbundle $J^k F$. We give a short description of the
theory, where, for simplicity, in the following we only consider the
case of the tangent bundle, $F=TM$,
 leaving to the reader the generalization to arbitrary vector bundles.
The model for the jet bundle $J^k M$ is the subbundle of $A^k M$ fixed
under $\Sigma_k$, 
$$
(A^k M )^{\Sigma^l} = (\bigoplus_\alpha \eps^\alpha TM)^{\Sigma_k}=
 \bigoplus_{j=1}^k \delta^{(j)} TM
$$
with the $\delta^{(j)}$ as in (8.1). 
If $L$ is a jet-compatible connection
on $T^k M$, then it induces a bundle isomorphism
$$
J\Phi:  \bigoplus_{j=1}^k \delta^{(j)} TM \to J^k M,
$$
and $J\Phi$ is the same for all connections $L^\sigma$, $\sigma \in \Sigma_k$.
In particular, if $L=L_k=D^{k-1} K$ is derived from some linear connection
$K$ on $TM$, then $J\Phi$ depends only on $K$ and not on the order 
of the scalar extensions in the definition of the iterated tangent functor
$T^k$. 

Now let $f:M \to N$ be a smooth map between manifolds
with linear connections $K$ on $TM$ and $L$ on $TN$.
Restricting the covariant derivative $d_K^L f$ from Section 16.3
to a map $J_K^L f:J^k M \to J^k N$,
$$
\matrix{ J^k M & {\buildrel \Phi_k \over \leftarrow} &
\bigoplus_{j=1}^k \delta^{(j)} TM \cr
J^k f \downarrow \phantom{T^k f}& & \phantom{x}\downarrow J_K^L f \cr
 J^k N & {\buildrel \Phi_k \over \leftarrow} &
\bigoplus_{j=1}^k \delta^{(j)}  TN \cr},
$$
we get a jet version of the higher covariant derivatives
that no longer depends on the standard order of scalar extensions 
but only on the connections themselves. 
Restricting the expression (16.6) of $T^k f$ by its matrix coefficients
 to spaces of
$\Sigma_k$-invariants and recalling from Section MA.3 that
$\lambda \sim \mu$ iff $\ell(\lambda)=\ell(\mu)=:\ell$ and
$|\lambda^j|=|\mu^j|=:i_j$ for all $j$, we get
$$
J^k f(x+\sum_j \delta^{(j)} v_j)=
f(x)+\sum_{j=1}^k \delta^{(j)} \sum_{\ell=1}^j
\sum_{i_1 \leq \ldots \leq i_\ell \atop i_1 + \ldots + i_\ell = j}
(d^{(i_1,\ldots,i_\ell)} f)_x(v_{i_1},\ldots,v_{i_\ell})
$$
with the ``matrix coefficient"
$$
(d^{(i_1,\ldots,i_\ell)} f)_x(v_{i_1},\ldots,v_{i_\ell}) =
\sum_{\lambda: \forall r=1,\ldots,\ell:
\atop \vert \lambda^r \vert = i_r} (d^\lambda f)_x(v_{i_1},\ldots,v_{i_\ell}). 
$$
This formula should be compared with the one from Theorem 8.6 which
concerns the ``flat case'' (chart connection). One expects
that also in this general case there exist symmetry properties of
the matrix coefficients and combinatorial relations between different
coefficients, depending of course on the curvature operators of the
connections. This seems to be an interesting topic for future work.

\vfill \eject 
\subheadline {20. Shifts and symmetrization}

\msk \nin {\bf 20.1.} {\sl Totally symmetrizable multilinear connections.}
We say that a multilinear connection $L$ on $T^k F$
 is {\it (totally) symmetrizable} if all its curvature operators
$R^{\sigma,L}$ belong to the central vector group $\Gm^{k,k+1}(T^k F)$
of $\Gm^{1,k+1}(T^k F)$.
Equivalently, all curvature forms $\tilde R^\Lambda$ for $\ell(\Lambda) < k$
vanish, and only the top curvature form  possibly survives.

\Theorem 20.2. Assume that
 $L$ is a totally symmetrizable multilinear connection
on $T^k F$ and that $2,\ldots,k$ are invertible in $\K$. Then
all connections $\sigma \cdot L$, $\sigma \in \Sigma_k$, belong to
the orbit $\Gm^{k,k+1}(E) \cdot L$, which is an affine space over $\K$,
and their barycenter in this affine space,
$$
\Sym(L) := {1 \over k!} \sum_{\sigma \in \Sigma_k} \sigma \cdot L,
$$
is a very symmetric multilinear connection on $T^k F$. Moreover, if $L$
is invariant under $\Sigma_{k-1}$, then
$$
\Sym(L) = {1 \over k} \sum_{j=1}^k (12\ldots k)^j  \cdot L.
$$

\Proof. This follows by applying pointwise the corresponding purely
algebraic result (Cor.\ SA.9).
\qed

\Corollary 20.3.
Assume $L$ is a torsionfree connection on $TM$. Then the linear structure
$DL$ on $T^3 M$ is symmetrizable.

\Proof. According to Theorem 18.3, the only non-vanishing curvature form
of $DL$ is the top curvature form, which means that $DL$ is symmetrizable.
\qed

\nin
In the situation of the corollary,
another proof of Bianchi's identity (cf.\ Section 17.5) is obtained
by observing that the curvature of $\Sym(L)$ vanishes -- see Section
SA.10 for the purely algebraic argument.

\msk \nin {\bf 20.4.} {\sl Shift invariance.}
Besides the action of $\Sigma_k$, the iterated tangent bundle $T^k F$
carries another canonical family of endomorphisms: recall from Section 7.3 (D)
that, deriving the canonical action $\K \times F \to F$ of $\K$ on the
vector bundle $F$, we get an action $T^k \K \times T^k F \to T^k F$
which is trivial on the base $T^k M$ and
linear in the fibers of the vector bundle $T^k F \to T^k M$.
The chart representation is
$$
\sum_{\alpha \in I_k} \eps^\alpha r_\alpha \cdot \sum_{\alpha \in I_k} 
\eps^\alpha (v_\alpha + w_\alpha) = \sum_\alpha \eps^\alpha
(v_\alpha + \sum_{\beta + \gamma = \alpha} r_\beta w_\gamma)
\eqno (20.1)
$$
where $v_\alpha \in V$, $w_\alpha \in W$ and $v_0:=x \in U$.
(This generalizes (7.11).) In particular, for $j=1,\ldots,k$, the action of the
``infinitesimal unit'' $\eps_j$  on $T^k F$, denoted by
$$
S_{0j}:T^k M \to T^k M, \quad u \mapsto \eps_j.u, 
\eqno (20.2)
$$
will be called a {\it shift operator}. This terminology is explained by
the chart representation of $S_{0j}$:
letting $r_\alpha =1$ if $\alpha = e_j$ and $r_\alpha=0$ if $\alpha \not= e_j$,
 we see from (20.1) that $\eps_j$ acts in a fiber of the bundle $T^k F \to
T^k M$  by ``shifting'' the index,
$$
\eps_j \cdot  \sum_{\alpha \in I_k} 
\eps^\alpha (v_\alpha + w_\alpha)  = \sum_\alpha \eps^\alpha v_\alpha +
\sum_{\alpha_j =1} \eps^\alpha 
w_{\alpha - e_j}  
 = \sum_\alpha \eps^\alpha v_\alpha +
\sum_{\alpha_j=0} \eps^{\alpha + e_j} w_\alpha . 
\eqno (20.3)
$$
All terms $\eps^\alpha w_\alpha$ with $\alpha_j =1$ are annihilated since
$\eps_j^2 =0$. The action $T^k \K \times T^k F \to T^k F$
is $\Sigma_k$-invariant in the sense that, for all $\sigma
\in \Sigma_k$ and $r \in T^k \K$, $u \in T^k F$,
$\sigma(r.u)=\sigma(r).\sigma(u)$. 
In particular, the shift operators $S_{0j}$ are permuted among each other
by the $\Sigma_k$-action. 

But now let $F=TM=T_{\eps_0}M$.
Then the action $T^k \K \times T^{k+1}M \to T^{k+1} M$ is not invariant
under $\Sigma_{k+1}$. In particular, conjugating by the transposition
$(0i)$ for $i=1,\ldots,k$, we get further shift operators
$$
S_{ij} := (0i) \circ S_{0j} \circ (0i): T^{k+1}M \to T^{k+1} M.
$$
In a chart, they are represented by
$$
S_{ij}:T^{k+1} M \to T^{k+1} M, \quad (x;\sum_{\alpha>0}\eps^\alpha v_\alpha)
 \mapsto (x;
\sum_{\alpha_i=1 \atop \alpha_j =0}\eps^{\alpha + e_j} v_\alpha +
\sum_{\alpha_i=0} \eps^\alpha v_\alpha). 
$$
We say that $S_{ij}$ is a shift {\it in direction $j$ with basis $i$},
and if moreover $i>j$, we say that the shift is {\it positive}. 
For $k=1$, the two  shifts $S_{10}$ and $S_{01}$ are precisely the two almost
dual structures on $TTM$ (Section 4.7).

We say that a multilinear connection $L$ on $T^k F$
 is {\it (positively) shift-invariant}
if all (positive) shifts are linear operators with respect to $L$.
Equivalently,  the splitting map $\Phi:A^k F \to T^k F$ is equivariant under
the groups $\Shi(T^k M)$, resp.\ $\Shi_+(T^k M)$ defined by (SA.20). 

Note that shift 
invariance  is a natural condition: chart connections are shift invariant.
More generally, if $F$ is a $\K$-vector bundle over
$M$, then $T^k F$ is a $T^k \K$-vector bundle over $T^k M$, and if
$\tilde f:F \to F'$ is a $\K$-vector bundle morphism, then $T^k \tilde f:
T^k F \to T^k F'$ is a $T^k \K$-vector bundle morphism.
By $T^k \K$-linearity in fibers, we deduce that $T^k \tilde f$ commutes with
all shifts. Similarly, if $f:M \to N$ is smooth, then $T^{k+1}f$ commutes
with all shifts $S_{ij}$. Invariance of a multilinear connection under $S_{ij}$
means that the shifts act linearly with respect to the linear structure,
and the fibers are actually modules over the ring obtained by adjoining
to $\K$ an infinitsimal unit $\eps_j$ whose action is specified by $i$. 

\Theorem 20.5. Assume $L$ is a torsionfree connection on $TM$.
\item{(1)} The linear structure $DL$ on $T^3 M$ is invariant under all shifts.
\item{(2)} The linear structure $D^kL$ on $T^{k+2} M$ is invariant under
the shifts $S_{k+2,k+1}$ and $S_{k+1,k+2}$.

\Proof. (1)
We use the chart representation (17.4) for the map $\Phi_2 = D\Phi_1$:
$$
\eqalign{\Phi_2(x,v) & = v +
\eps_0 \eps_1 b_x(w_0,v_1) + \eps_0 \eps_2b_x(w_0,v_2) +
\eps_1 \eps_2 b_x(v_1,v_2) +
\cr
& \quad \quad \quad   \eps_0 \eps_1 \eps_2
 \bigl( b_x(w_0,v_{12}) + b_x(w_{01},v_2)+ b_x(w_{02},v_1) + t(w_0,
v_1,v_2) \bigr) \cr}
\eqno (20.4)
$$
with trilinear component 
$$
t(w_0,v_1,v_2)=
b_x(w_0,b_x(v_1,v_2))+b_x(b_x(w_0,v_2),v_1)+
\partial_{v_2} b_x(w_0,v_1) .
$$
Clearly, all bilinear components of $\Phi_2(x,\cdot)$ are conjugate
to each other, and hence the algebraic
Shift Invariance Condition (Prop.\ SA.13) is fulfilled.

\ssk
(2) The claim is proved by induction. For $k=1$ it follows from Part (1).
Now assume the claim is proved for $k>1$. We use the recursion formula (19.2)
 for $\Phi_{k+1}$:
$$
\Phi_{k+1}:
\matrix{A^{k+1} F & {\buildrel \oplus_\alpha \Phi_2^{e_{k+1},e_k;\alpha}
 \over \to}
& TT(A^{k-1} F) & {\buildrel TT\Phi_{k-1} \over \to} & TT(T^{k-1} F)=
T^{k+1}F \cr}
$$
The map  $TT\Phi_{k-1}$, like any second tangent map,
 commutes with both shifts $S_{k+1,k}$ and $S_{k,k+1}$.
By the case $k=1$ of the induction, 
 $\Phi_2^{e_{k+1},e_k;\alpha}$ also commutes with both
 shifts $S_{k+1,k}$ and $S_{k,k+1}$, whence the claim.
\qed

\nin As seen in the algebraic part (Theorem SA.14), the condition that $L$
be invariant under {\it all} shifts is a rather strong condition: it implies
in turn that $L$ is symmetrizable. -- Part (1) of
the preceding theorem does not carry over in this form to arbitrary
$k$: using the recursion formula (17.3), one can see that $D^2 L$
is in general no longer invariant under the shifts $S_{12}$ and $S_{21}$.
On the other hand, the statement (2) is not the strongest possible result:
combinatorial arguments, using again the recursion formula (17.3), show
that $D^k L$ is invariant under all shifts of the form $S_{i,k+2}$,
$i<k+2$. One would like to understand this situation in a conceptual way.

\Corollary 20.6. Under the assumptions of the preceding theorem, and
if $3$ is invertible in $\K$, then the linear structure
$\Sym(DL)$ on $T^3M$ is totally symmetric and invariant under all
shifts.

\Proof. By Cor.\ 20.3, $DL$ is indeed symmetrizable, and by Theorem 20.2,
$\Sym(DL)$ is then totally symmetric. We have to show that it
 is invariant under all shifts. 
In a chart, symmetrization amounts to replace in Equation (19.1)
the trilinear map $t$ by its symmetrized version ${1 \over 3}
(t + (123). t + (132). t)$, whereas the other components do not change
(recall that $b_x$ is already symmetric, by torsionfreeness).
Thus $\Sym(DL)$ still satisfies the shift invariance condition
SA.13. \qed

\Theorem 20.7.
Assume $L$ is a {\rm flat} linear structure on $TM$, i.e.,
 without torsion and such
that all curvature forms of $DL$ vanish. Then all linear structures
$D^k L $, $k\in \N$, are totally symmetric and invariant under all shifts.
       
\Proof.
The proof is by induction on $k$: for $k=1$, $DL$ is totally symmetric by
assumption and shift-invariant by Theorem 20.5.
For the induction step, we
use exactly the same arguments as in the proof of Theorem 19.2,
replacing the property of preserving the space of $\Sigma_k$-invariants
by full $\Sigma_k$-invariance.
In this way it is seen that $D^k L$ also is $\Sigma_{k+1}$-invariant,
i.e.\ all $D^k L$ are totally symmetric. By Theorem 20.5 they are invariant
under one shift (namely $S_{k+1,k}$), but since they are totally symmetric
and all other shifts are conjugate among each other under the permutation
group, $D^k L$ is then invariant under all shifts.
\qed

\msk \nin 
For a complete understanding of the structure of $T^k M$, in relation with
connections, it seems of basic importance to find a canonical construction of
 a sequence of totally symmetric
and shift-invariant multilinear structures coming from the linear structures
$D^k L$, $k=1,2,\ldots$. Unfortunately, already $D (\Sym(DL))$ is in general
 no longer totally symmetrizable in the sense of 20.1, and
thus an appropriate symmetrization procedure has to
 be more subtle than the one used in Theorem 20.2. 
We will solve this problem for Lie groups (Chapter 25) and symmetric spaces
(Chapter 30), where the solution is closely related to the {\it exponential
mapping of polynomial groups} (Chapter PG), and then
come back to the general problem in Chapter 31.
-- 
In the following chapter we recall a similar construction in a context
which is ``dual'' to the one treated here, namely in the context of
linear differential operators.
 
\vfill \eject 

\subheadline
{21. Remarks on differential operators and symbols}

\msk \nin
{\bf 21.1.} {\sl Overview over definitions of differential operators.}
In the literature several definitions of differential operators
can be found. They agree for finite-dimensional real manifolds
but may lead to quite different concepts in our much more general
setting. All have in common that a {\it linear differential
operator} between vector bundles $F \to M$ and $F' \to M$ should be
a $\K$-linear map
$$
D: \Gamma(F) \to \Gamma(F')
$$
from sections of $F$ to sections of $F'$ satisfying some additional
conditions such as:

\ssk
\item{(1)} 
$D$ is {\it support-decreasing}. Since there is no
hope to prove an analog of the classical Peetre theorem (cf., e.g.,
[Hel84, Th.\ II.1.4] or [KMS93, Ch.\ 19]) in our
setting, this cannot serve as a possible definition in the general case.
\item{(2)}
The value $Df(x)$ depends only on the $k$-jet of the section $f:M \to F$
at $x$. More precisely, let $J^k(M,F,{\rm sec})$ be the space whose fibers
over $M$ are given by
$k$-jets of sections of $F$ (if $F=M \times \K$ is the trivial bundle,
then $J^1(M,F,{\rm sec})$ essentially is the cotangent bundle $T^*M$).
Following [KMS93, p.\ 143], we say that $D$ is {\it local and of order $k$} 
if there is a map $D_k:J^k (M,F,{\rm sec}) \to F'$
such that $Df(x)=D_k(J^k_x f)$. In case of the trivial line
bundle, i.e.,  $F=F'=M \times
\K$, a differential operator of order one is thus seen as a section of 
the dual bundle of $T^*M$, and not as a section of $TM$. 
Already this example shows that, when trying to 
generalize this definition,
we inevitably run into difficulties related to double dualization.
\item{(3)} 
Following the elegant algebraic definition of differential operators given
e.g.\ in [ALV91, p.\ 85] or [Nes03, p.\ 125, p.\ 131], we may require that
there exists $k \in \N$ such that, for all smooth functions
$f_0,\ldots,f_k \in {\cal C}^\infty(M,\K)$, 
$$
[\ldots [[D,f_0],f_1]\ldots f_k] = 0,
$$
where, for an operator $A$ from sections to sections,
$[A,f]X=A(fX) - fA(X)$. Moreover, one would add some
smoothness condition. It is not hard to see that (2) implies
(3), whereas already in the infinite-dimensional Banach case
the converse fails in general 
(cf. remarks in Section 4.6). Hence it seems that this definition is
best suited for the finite-dimensional case.
\item{(4)} 
In the special case $F=F'$, we may assume that there exists
$k \in \N$ such that $D$ can (at least locally over a chart domain)
be written as a finite product
of operators $A \nabla_{X_1} \circ \ldots \circ \nabla_{X_k}$,
where $\nabla$ runs over the covariant derivatives associated to  connections 
on $F$ and $A$ is a field of endomorphisms (cf. [BGV92, p.64]).
In case $F=M \times \K$ is the trivial bundle over a finite-dimensional
manifold, this means that,
 in a chart, 
$D =g_0+\sum_{\alpha} g_\alpha \partial^\alpha$ is a sum a partial
derivatives multiplied by functions.
This is the most classical notion; it implies (2) and hence (3).

\ssk
\nin
All of these concepts are useful, and hence one should
distinguish between different classes of differential operators.
In any case, $k$ is the {\it degree} or {\it order of $D$},
and the degree defines a natural filtration on the space of
differential operators. 
The second approach seems to be the most conceptual one,
and probably a correct formalization of concepts of differential operators in
our context will need the use of (algebraic) double duals. However,
this is a topic for another work, and therefore the approach outlined
in the sequel will be closest to the most classical definition (4).
Things will be most clear if we start with differential operators on
the trivial line bundle, i.e.\ with scalar differential operators.

\msk \nin
{\bf 21.2.} {\sl Pure differential operators.}
A {\it pure scalar differential operator of order at most $k$} ($k\geq 1$)
 on a manifold $M$ is a smooth section
$D: M \to T^k M$
of the higher order tangent bundle $T^k M$ over $M$. The space of such sections
is denoted by $\X^k(M)$. (See Chapter 28 for a more systematic theory.)
 If $f:M \to \K$ is a smooth 
function, then $T^k f \circ D: M \to T^k \K$ is smooth, and we let
$$
Df:= \pr_k \circ T^k f \circ D: M \to \K,
\eqno (21.1)
$$
where $\pr_k:T^k \K \to \K$, $\sum_{\alpha \geq 0} a_\alpha
 \eps^\alpha  \mapsto a_{(1\ldots1)}$
is the projection onto the last factor. 
Thus $D$ induces a map
$$
D^{op}: {\cal C}^\infty(M) \to {\cal C}^\infty(M), \quad
f \mapsto Df
$$
which clearly is $\K$-linear. Let ${\cal D}_k(M) \subset
\End_\K({\cal C}^\infty(M))$ be the $\K$-span
of all such operators, called the
space of {\it linear scalar differential operators of degree at most $k$}.
In a chart, using Formula (7.19), the ``constant coefficient'' operator
$D=\sum_{\alpha >0} \eps^\alpha v_\alpha$ acts by
$$
(\sum_{\alpha >0} \eps^\alpha v_\alpha \cdot f)(x) =
\pr_k(T^k f(x+\sum_\alpha v_\alpha))=\sum_{\ell=1}^k
\sum_{\beta \in {\cal P}_\ell(1,\ldots,1)}
d^\ell f(x)(v_{\beta^1},\ldots,v_{\beta^\ell}) ,
\eqno (21.2)
$$
In other words,
$$
(\sum_{\alpha >0} \eps^\alpha v_\alpha)^{op} =
\sum_{\ell=1}^k
\sum_{\beta \in {\cal P}_\ell(1,\ldots,1)}
\partial_{v_{\beta^1}} \cdots \partial_{v_{\beta^\ell}}
$$
is a scalar differential operator in the most classical sense.
For instance, $\eps_1 v_1 + \eps_2 v_2$ acts as $\partial_{v_1}
\partial_{v_2}$ and $\delta^{(2)} v$ as $\partial_v^2$.
The chart formula also shows that, for all $\sigma \in \Sigma_k$ acting
on $T^k M$ in the canonical way, $\sigma \circ D$ and $D$ give
rise to the same operator; therefore,
if $2,\ldots,k$ are invertible in $\K$, then we could as well work
with sections of $J^k M$. Also, if $T^k z:T^k M \to T^{k+1} M$ is
the injection obtained by deriving some zero section $z=z_j:
M \to T_{\eps_j} M$, then $D$ and $z \circ D$ give rise to the same
operator on functions. Somewhat loosely, we thus may say that
a differential operator of order $k$ is also one of degree $k+1$.
(Taking in the definition of $Df$
instead of $\pr_k$ the projection onto the ``augmentation ideal"
of $T^k \K$, one gets another kind of operators satisfying an analog of
the Leibniz rule, called {\it expansions} in [KMS93, Section 37.6].)


The {\it composition} of two pure differential operators
$D_i:M \to T^{k_i} M$, $i=1,2$, is defined by
$$
D_2 \cdot D_1:= T^{k_1} D_2 \circ D_1:M \to T^{k_1} M \to T^{k_1+k_2} M.
\eqno (21.3)
$$ 
It is easily checked that this composition is associative, and that
$D \mapsto D^{op}$ is a ``contravariant" representation:
$$
\eqalign{
(D_2 \cdot D_1)^{op} f & = \pr_{k_1+k_2} \circ T^{k_1+k_2}f \circ
 (D_2 \cdot D_1) \cr
& = \pr_{k_1+k_2} \circ T^{k_1+k_2}f \circ T^{k_1} D_2 \circ D_1 \cr
& = \pr_{k_1+k_2} \circ T^{k_1} (T^{k_2}f \circ D_2) \circ D_1 \cr
& = \pr_{k_1} \circ T^{k_1} (\pr_{k_2} \circ T^{k_2}f \circ D_2) \circ D_1 \cr
& = (D_1^{op} \circ D_2^{op})f. \cr}
$$
It follows that ${\cal D}^k(M) \cdot {\cal D}^\ell(M) \subset
{\cal D}^{k+\ell}(M)$, and hence ${\cal D}^\infty(M):=\bigcup_{i=1}^\infty
{\cal D}^i(M)$ is an associative filtered algebra over $\K$. 

\msk \nin
{\bf 21.3.} {\sl Principal symbol of a scalar differential operator.}
Let $D:M \to T^k M$ be a section and
$$
\pr:T^k M \to \eps_1 TM \times \ldots \times \eps_k TM
$$
be the product of the canonical projections, and let 
$$
t:\eps_1 TM \times \ldots \times \eps_k TM \to 
\eps_1 TM \otimes \ldots \otimes \eps_k TM
$$
be the tensor product map (where $\otimes=\otimes_M$ is the algebraic
tensor product of vector bundles over $M$, i.e.\ given by tensor products
in fibers, without considering a topology). 
We say that
$$
\sigma^k(D):=
t \circ \pr \circ D: M \to \eps_1 TM \otimes \ldots \otimes \eps_k TM
$$
is the {\it principal symbol of $D$}. It is a section of the (algebraic)
bundle $\otimes^k TM$ over $M$.
Comparison with the chart representation (21.2) shows that the principal
symbol of $\sum_{\alpha >0} \eps^\alpha v_\alpha$ is
$v_1 \otimes \ldots\otimes v_k$, coming from the ``highest term''
$\partial_{v_1} \cdots \partial_{v_k}$.
If $D$ happens to be a section of $J^k M$, then
$\sigma^k(D)$ is a section of the algebraic bundle $S^k(TM)$ over $M$
($k$-th symmetric power of $TM$). 
If the ``operator representation'' $D \mapsto D^{op}$ is injective, then
one can show that the principal symbol map factors through a $\K$-linear
map
$$
\tilde \sigma^k: {\cal D}^k(M) \to \Gamma(S^k TM)
$$
(where $\Gamma(S^k TM)$ is the space of sections that locally are
finite sums of $t$ composed with smooth sections of 
$TM \times \ldots\times TM$), i.e.\
$\sigma^k(D)=\tilde \sigma^k(D^{op})$, and that
the sequence
$$ \matrix{
0&\to& {\cal D}^{k-1}(M)& \to & {\cal D}^k(M) & {\buildrel \tilde \sigma^k
\over \to}& \Gamma(S^k(TM))& \to 0 & \cr}
\eqno (21.4)
$$
is exact. Details are left to the reader.

\msk \nin
{\bf 21.4.} {\sl Correspondence between total 
symbols and differential operators (scalar case).}
For real finite dimensional manifolds 
it is known that, in presence of a connection on $TM$, one can associate
to a  differential operator a {\it total symbol} (which is a section
of $S^k TM$) in such a way that this correspondence becomes a bijection 
between operators and their total symbols.\footnote{$^1$}{
[ALV91, p.87]: ``Obviously, to one and only one symbol there 
correspond many differential operators. However, it is a remarkable
fact that in presence of a connection this correspondence can be
made one-to-one.'' On the same page, the authors add the remark:
``This procedure of restoring a differential operator from its symbol
is completely analogous to the basic procedure of quantum mechanics --
quantization.''}
For the case of second order operators, the construction is essentially
described in [Lo69] (see also [P62] and [Be00, I.B]);
 the general construction seems to be
folklore -- the few references we know about include 
[BBG98] and [Bor02] (where the construction is called {\it pr\'escription
d'ordre standard}). Let us briefly describe the
construction: clearly, the problem is equivalent to find $C^\infty(M)$-linear
sections of the exact sequences (21.4) for all $k$.
Using the covariant derivative $\nabla$ of $L$, for a function $f$ one defines
$\nabla^k f:\times^k TM \to \K$,
the ``$k$-th order Hessian of $f$'' (defined as usual in differential
geometry by $\nabla^1 f:=\dd f$, $\nabla^{j+1}f:=\nabla(\nabla^j f)$).
Then to a section $s$ of $S^k TM$ we associate a differential 
operator $D:=D(s)$:
if $s(x)=v_1 \otimes \ldots \otimes v_k$ with all $v_j \in T_x M$, let
$$
D f(x) = (\nabla^k f)(x;v_1,\ldots,v_k).
$$
One easily proves that $\tilde\sigma^k D(s)=s$, hence $s \mapsto D(s)$ gives
the desired section of the exact sequence (21.4).
This construction should be seen as an analog of the construction of linear
structures on higher order tangent bundles (Chapter 17) ; it
can be adapted to the case of differential operators
acting on sections of (finite dimensional) vector bundles over finite
dimensional real manifolds (cf.\ [Bor02]).

\msk \nin
{\bf 21.5.} {\sl Differential operators on general vector bundles.}
First of all, differential operators on a trivial bundle $M \times W$
may be defined as in the scalar case. 
For general vector bundles, the problem arises that we need a good class
of differential operators of degree zero (which can be avoided on
trivial vector bundles). The reason for this is that, although vector
bundles can be trivialized locally, the preceding definition of differential
operator over chart domains would no longer be chart-independent.
The smallest class of differential operators of degree zero which one
could take here are
fields of endomorphisms that arise precisely from chart changes.
Taking a somewhat bigger class, one arrives at
Definition (4) of a differential operator mentioned in Section 21.1.

Finally, in order to define differential operators acting from sections
of one vector bundles $F,$ to sections of {\it another} bundle $F'$, one needs 
to single out some subspace in 
$\Hom_\K(F_x,F_x')$  giving the differential operators
of degree zero; but there is no natural candidate for such a space
(in the general case).


\vfill \eject

\subheadline
{22.  The exterior derivative}

\nin {\bf 22.1.} {\sl Skewsymmetric derivatives.}
In the preceding chapters, we have focused on symmetric derivatives.
They depend on an additional structure (connection) and thus may be
considered as
``not natural''.  The skewsymmetric version of the derivative (the
exterior derivative), in contrast, is completely natural, and therefore
plays a dominant role in differential geometry\footnote{$^1$}{Cf., e.g.,
 [Sh97, p.\ 54]:  ``As usual in modern
differential geometry, we shall be concerned only with the {\it
skew-symmetric} part of the higher derivatives.''}.
Since the exterior derivative is so natural, we wish to give
also a definition which, in our context, is as natural as possible.
Recall (Section 4.5) that we consider a differential $k$-form $\omega$
as a smooth map $TM \times_M \ldots \times_M TM \to \K$ which is
multilinear and alternating in fibers. 
One could proceed  as in [La99, Ch.\ V], where the
exterior derivative of a differential form $\omega$ is defined in a chart by
$$
\eqalign{
\dd \omega(x;v_1,\ldots,v_{k+1}) & =
\sum_{\sigma \in \Sigma_{k+1}} \sgn(\sigma) \partial_{v_{\sigma(1)}}
\omega(\hat{v_{\sigma(1)}},v_{\sigma(2)},\ldots,v_{\sigma(k+1)})
\cr
& =\sum_{i=1}^{k+1} (-1)^i \partial_{v_i}
\omega(v_1,\ldots,\hat{v_i},\ldots,v_{k+1}).
\cr}
\eqno (22.1)
$$
Then one shows that 
evaluation of $\dd \omega$ on vector fields $X_1,\ldots,X_k$
is given by the formula
$$
\eqalign{
\dd \omega(X_1,\ldots,X_{k+1}) & =
\sum_{i=1}^{k+1} (-1)^i X_i \omega(X_1,\ldots,\hat{X_i},\ldots,X_{k+1}) +
\cr
&  \sum_{1 \leq i < j \leq k+1} (-1)^{i+j}
\omega([X_i,X_j],X_1,\ldots,\hat X_i,\ldots,\hat X_j,\ldots,X_{k+1}),
\cr}
\eqno (22.2)
$$
which implies that (22.1) is chart-independent. 
The usual properties of the exterior derivative then follow by
standard calculations.
A ``natural definition'' in our sense should not rely on chart
formulas, and it should be built only on the natural operators
$T$ (tangent functor) and $\alt$ (alternation operator). 
For the exterior derivative of one-forms this has been carried out
in Chapter 13. Here is the pattern for the general case:

\ssk
\item{(1)}
We define the notion of an {\it (intrinsic) $k$-multilinear
map} $f:T^k M \to \K$, and we define also {\it homogeneous} such
maps.
\item{(2)}
There is a canonical bijection between
tensor fields $\uline \omega:\times^k_M TM \to \K$ and 
homogeneous $k$-multilinear maps $\oline \omega:T^k M \to \K$.
\item{(3)}
If $f:T^k M \to \K$ is  $k$-multilinear,
then $\pr_2 \circ Tf:T^{k+1} M \to \epsilon \K$ is 
$k+1$-multilinear.
\item{(4)}
The alternation operator $\alt$ maps intrinsic $k$-multilinear
maps onto homogeneous ones.
\item{(5)}
Finally, we let $\dd \uline \omega := \uline{\alt (\pr_2 T(\oline \omega))}$.
In other words, the natural operator $\dd$ is defined by
the following commutative diagram:
$$
\matrix{
{\cal T}^{k,0}(M) & \cong & {\cal M}_h(T^kM) & \subset & {\cal M}(T^k M) \cr
\dd \downarrow \phantom{\dd} & & & & \phantom{T} \downarrow T \cr
{\cal T}^{k+1,0}(M) & \cong 
 & {\cal M}_h(T^{k+1}M) &  {\buildrel \alt \over \leftarrow} & 
{\cal M}(T^{k+1} M) \cr}
$$
where ${\cal T}$ denotes tensor fields and ${\cal M}$ (resp.
${\cal M}_h$) intrinsic multilinear (homogeneous) maps.
\item{(6)}
The  property $\dd \circ \dd = 0$ follows from the fact that
 $\alt \circ T \circ T = 0$ due to symmetry of second differentials.

\ssk \nin
In the following, we carry out this program.

\msk \nin
{\bf 22.2.} {\sl Multilinear maps on $T^k M$.}
We say that a smooth map $\omega:T^k M \to \K$ is {\it (intrinsically)
multilinear} if $\omega:T(T^{k-1} M) \to \K$ is linear in fibers,
for all the various ways $p_j:T^k M \to T^{k-1}M$, $j=1,\ldots,k-1$ of seeing
$T^k M$ as a vector bundle over $T^{k-1} M$.
Equivalently, for all $x \in M$, the restriction
$\omega_x:(T^k M)_x \to \K$ to the fiber of $T^k M$ over $x$ is
intrinsically multilinear in the sense of Section MA.15.
We denote, for fixed $x \in M$, by $H_j$ the  kernel of $p_j$
(tangent spaces $T_{0_x} (T^{k-1} M) \subset (T^k M)_x$ with first
$T=\eps_j T$).  
We say that $\omega$ is {\it homogeneous} if $\omega_x$ is homogeneous
in the sense of MA.15, i.e.\ the restriction of $\omega_x$ to intersections
$H_{ij}:=H_i \cap H_j$ ($i \not= j$) is constant.
In a chart, an intrinsically multilinear map is given by
$$
\omega(x+\sum_\alpha \eps^\alpha v_\alpha)=
\sum_{\ell=1}^k \sum_{\Lambda \in {\cal P}_\ell(1,\ldots,1)} b^\Lambda_x
(v_{\Lambda^1}, \ldots, v_{\Lambda^\ell}),
\eqno (22.3)
$$
where $b^\Lambda_x:V^\ell \to \K$ is a multilinear map in the usual sense,
and it is homogeneous if and only if the only non-zero term in (22.3) is
the one belonging to $\ell=k$, i.e., to the longest partition 
$\Lambda = (e_1,\ldots,e_k)$ of the full multi-index $(1,\ldots,1)$
(see Prop.\ MA.16).
This implies as in Prop.\ MA.16
 that we have a bijection between tensor fields and homogeneous
multilinear maps: 
$$
\eqalign{
 {\cal T}^{k,0}(M,\K) \to {\cal M}_h(T^k M,\K) , &  \quad
\omega \mapsto \oline \omega = \omega \circ \pr,
\cr
{\cal M}^h(T^k M,\K) \to  {\cal T}^{k,0}(M,\K) , & \quad
\oline \omega \mapsto \uline \omega :=
\oline \omega \circ \iota^L, \cr}
\eqno (22.4)
$$
where $\pr:T^k M \to TM \times_M \ldots \times_M \times TM$,
$v \mapsto (\pr_1(v),\ldots,\pr_k(v))$ is the vector bundle direct sum
of the various projections $\pr_i:T^k M \to TM$. 
The definition of the first map in (22.4)
 does not involve connections or charts, 
whereas the definition of the second map does -- the injection
 $\iota^L:\times^k TM \to T^k M$  obtained from the Splitting Map
$\Phi_k:\oplus_\alpha \eps^\alpha TM \to T^k M$ depends on a
connection $L$. 
The situation is summarized by the following diagram (cf.\ (MA.36)):
$$
\matrix{T^k M &  &  \cr \pr \downarrow \uparrow \iota^L 
& {\buildrel \oline \omega
\over \searrow} & \cr
\times^k_M TM &{\buildrel \uline \omega \over \to} & \K \cr}
$$
Later on, one may consider $\oline \omega$ and $\uline \omega$ as the
same object $\omega$.

\Lemma 22.3. If $\oline \omega:T^k M \to \K$ is $k$-multilinear, then
$\pr_2 \circ T\oline \omega:T^{k+1} M \to \epsilon \K$ is $k+1$-multilinear.

\Proof. This can be proved by a direct computation using (22.3). 
More intrinsically, the lemma is proved as follows.
Let the last tangent functor correspond to scalar extension with
$\epsilon_{k+1}$.
Then we have two kinds of tangent spaces in $T^{k+1} M$:
the first kind is given by spaces that are stable under $\epsilon_{k+1}$;
they correspond to projections $p_i$, $i=1,\ldots,k$, and
are tangent bundles of tangent spaces in $T^k M$.
The restriction of $Tf$ to such tangent spaces is the tangent map
of $f$ restricted to tangent spaces in $T^k M$ and hence is linear.
The second kind of tangent spaces corresponds to the
last projection $p_{k+1}$; but since $Tf$ is defined with respect
to $\epsilon_{k+1}$, the restriction of $Tf$ to such tangent spaces
is linear as it is a tangent map.
\qed 

\Lemma 22.4.
If $\oline \omega:T^k M \to \K$ is $k$-multilinear, then
$$
\alt  \oline \omega:T^k M \to \K,
\quad z \mapsto \sum_{\sigma \in \Sigma_k} \sgn(\sigma)\oline\omega(\sigma(z))
$$
is homogeneous $k$-multilinear.

\Proof. 
It is clear that alt$(\omega)$
 is $k$-multilinear since it is a sum of multilinear
functions. It remains to
prove homogeneity: notation being as in MA.15, 
we have to show that $\alt \oline \omega (v)=0$ for
$v \in H_{ij}$ with $i \not= j$. Let $\tau=(ij)$ be transposition between
$i$ and $j$; then $\tau$ acts trivially an $H_{ij}$ since $v=\sum_\alpha
\eps^\alpha v_\alpha$ where $\alpha_i=1=\alpha_j$ for all $\alpha$ contributing
effectively to the sum. It follows that $\tau v= v$ and hence
$$
\alt \oline \omega(v)=
\sum_{\sigma \in \Sigma_n} \sgn(\sigma) \oline \omega(\sigma v) =
\sum_{\sigma \in A_n} (\oline \omega(\sigma v) -
 \oline \omega(\sigma \tau v) ) = 0.
\qeddis

\msk \nin {\bf 22.5.} {\sl Exterior derivative.}
 The exterior derivative of
a $k$-form $\uline \omega:\times^k TM \to \K$ is now defined by
$$
\dd \uline \omega := \uline{ \alt_{k+1} \pr_2 T \oline \omega }
= \uline{\pr_2 \alt_{k+1}  T \oline \omega }.
$$
In fact, the same definition can be given
 for any tensor field $\uline \omega:T^k M \to \K$ (not
necessarily skew-symmetric), but then 
$\dd \uline \omega = \dd (\alt \uline\omega)$, 
so this gives nothing really new.
Similarly, we define the exterior derivative of $\oline \omega \in
{\cal M}_h(T^k M)$ by
$$
\dd \oline \omega:=\alt_{k+1} (\pr_2 T \oline \omega)
=\pr_2(\alt_{k+1} T \oline \omega)\in {\cal M}_h(T^{k+1}M).
$$

\Theorem 22.6. For all $k$-forms $\uline \omega$, we have
$\dd \dd \uline \omega = 0$.

\Proof. We show, equivalently, that $\dd \dd \oline \omega =0$.
The operator $\dd = \dd_k$ is defined as
$$
\dd_k = \alt_{k+1} \circ \pr_2 \circ T = \pr_2 \circ \alt_{k+1}\circ  T,
$$
and hence
$$
\dd_{k+1} \circ \dd_k = \alt_{k+2} \circ \pr_2 \circ T  \circ
\alt_{k+1} \circ \pr_2 \circ T =
\pr_4 \circ \alt_{k+2} \circ  T  \circ\alt_{k+1}\circ T.
$$
Now, we have $T \circ \alt_k = \alt_k \circ T$ for all $k$, because
 for all $\sigma \in \Sigma_k$, acting on $T^k M$, we have
$T(\oline \omega \circ \sigma)=T\oline \omega \circ T\sigma$, where
$T\sigma$ is the action of $\sigma$ on $T^{k+1} M$ if we inject
$\Sigma_k$ into $\Sigma_{k+1}$ in the natural way.
Hence $T \circ \alt_k = \alt_k \circ T$ for all $k$, and from this follows
$\pr_2 T \circ \alt_k = \alt_k \circ \pr_2 T$, and we get for the
operator $\dd_{k+1} \circ \dd_k$: 
$$
\dd_{k+1} \circ \dd_k =\pr_4 \circ \alt_{k+2} \circ T  \circ
\alt_{k+1}  \circ T \circ  \pr_2 \circ T 
=
\pr_4 \circ \alt_{k+2} \circ \alt_{k+1} \circ T \circ T .
$$
But, by symmetry of the second differential, $TT\oline\omega$
 is invariant under the permutation $(k+2,k+1) \in \Sigma_{k+2}$
 which corresponds to exchanging $\eps_{k+2}$ and $\eps_{k+1}$,
 and hence applying the alternation operator
$\alt_{k+2}$ annihilates the whole expression.
\qed 

\nin The {\it de Rham cohomology} is now defined as usual:
$$
H^k (M):=  \ker(\dd_{k+1})/\im(\dd_k).
$$

\ssk \nin {\bf 22.7.} {\sl Exterior derivative and covariant derivative.}
If $L$ is a linear connection of $TM$, we define the covariant derivative
of a tensor field $\uline \omega: \times^k TM \to \K$ by
$$
\nabla \uline \omega := \uline{\pr_2 T \oline \omega} \in {\cal T}^{k+1,0}
(M),
$$
where $\uline f :=f \circ \iota^L$ is defined as in (22.4) and
depends on the connection $L$ unless $f$ is homogeneous.

\Proposition 22.8.
For all $k$-forms $\uline \omega:\times^k M \to \K$, 
$$
\dd \uline \omega = \alt_{k+1}(\nabla \uline \omega).
$$

\Proof.
$$
\alt_{k+1}(\nabla \uline \omega)=\alt_{k+1}\uline{\pr_2 T \oline \omega}=
\uline{\alt_{k+1} \pr_2 T \oline \omega} =
\dd \uline \omega.
\qeddis

\Corollary 22.9. In a bundle chart, $\dd \omega$ is given by {\rm (22.1)},
and evaluation of $\dd \omega$ on vector fields $X_1,\ldots,X_k$
is given by {\rm (22.2)}.

\Proof. 
Let $\nabla$ be the covariant derivative induced by the chart, i.e.\
$\nabla = d$ is the ordinary derivative. Then
$$
(\nabla \omega)(x;v_1,\ldots,v_{k+1})  =
 \partial_{v_{k+1}}
\omega(x;v_{1},v_{2},\ldots,v_{k})
$$
Applying $\alt_{k+1}$ gives the first line of (22.1),
and it implies the second line since already $\omega$ was supposed to be
skew-symmetric.

In order to prove the formula for evaluation on vector fields,
we first have to check that the result of the right-hand side
at a point $p \in M$ depends only on the values of the $X_i$ at
$p$; then take $X_i = \tilde v_i$, the extension to a constant
vector field in a chart, and the formula reduces to the one
established before.
\qed

\msk
Of course, one may ask now: what about the Poincar\'e lemma and de Rham's
theorem? Clearly, these results
 do not carry over in their classical forms (not even
for the infinite-dimensional real case). Already
determining the solution space of the
``trivial'' differential equation $df=0$ for a scalar function $f$,
 as a first step towards the Poincar\'e lemma, is in general not possible
(in the $p$-adic case there are non-locally constant solutions,
see [Sch84] for an example),
and one may expect that there is some cohomological theory
giving a better insight into the structure of the
solution space of this ``trivial'' differential equation. 
For general $k$-forms, one should then combine such a theory with
the usual de Rham cohomology.

\vfill\eject

%% file: Dgeo3.tex

\def\date{25.10.2006}
\def\title{V. LIE THEORY} 
\def\rightheadline{\hfil{\eightrm \title}\hfil\tenrm\folio}

\sectionheadline{V. Lie theory}

\nin In this part we develop the basic ``higher order theory'' of
iterated tangent bundles of Lie groups and symmetric spaces.
It can be read without studying most of the preceding results --
see the ``Leitfaden'' (introduction).
 Strictly speaking, not even Theorem 4.4 on the Lie bracket
of vector fields is needed since we give, in Chapter 23, another,
independent definition of the Lie bracket associated to a Lie group, see below.

\subheadline
{23. The three canonical connections of a Lie group}

\msk \nin {\bf 23.1.} 
In the following theorem we will give another characterization of
the Lie bracket of a Lie group. One might as well use it as
{\it definition} of the Lie bracket; this has the advantage that
one can develop Lie theory without speaking about vector fields,
and in this way the theory becomes simpler and more transparent. 
In the following,
 $G$ is a Lie group with multiplication $m:G \times G \to G$
(cf.\ Chapter 5).
Then $TTG$ is a Lie group with multiplication $TTm$, and the fiber
$(TTG)_e$ is a normal subgroup. We let $\g:=T_e G$; then the three
``axes'' $\eps_1 \g$, $\eps_2\g$, $\eps_1 \eps_2 \g$ of $(TTG)_e$
are canonical, and the fiber of the axes-bundle over the origin is 
$(A^2 G)_e = \eps_1 \g \oplus
\eps_2 \g \oplus \eps_1 \eps_2 \g \cong \g \times \g \times \g$.
(Cf.\ Chapter 15 for this notation.)

\Theorem 23.2. {\rm (Canonical commutation rules in $(TTG)_e$.)}
The group commutator $[g,h]=ghg^{-1}h^{-1}$ of an element
 $g=\eps_1 v$ from the first and an element $h=\eps_2 w$ from the second
axis in $(TTG)_e$, where $v,w \in \g$,
 can be expressed by the Lie bracket $[v,w]$ in the following way:
$$
[\epsilon_1 v, \epsilon_2 w] =\epsilon_1 \epsilon_2 [v,w],
\eqno (23.1)
$$
and the third axis $\epsilon_1 \epsilon_2 \g$ is central in $(TTG)_e$, i.e.
 for $j=1,2$ and all $u,v \in \g$,
$$
[\epsilon_1 \epsilon_2 u, \epsilon_j v] = 0  .
\eqno (23.2)
$$

\Proof.
 First we prove (23.2). For any Lie group $H$,  the tangent space
$T_e H$ is an abelian subgroup of $TH$ whose group law is vector addition;
this follows from Equation (5.1). Now, taking $H:=TG$,
 $\epsilon_1 \epsilon_2 u$ and $\epsilon_j v$ both belong
to the fiber over the origin in $T_{\epsilon_j} TH$ and hence commute.
Moreover, we see that
$$
\eps_i v \cdot \eps_i w = \eps_i(v+w) \quad  (i=1,2), \quad \quad 
\eps_1 \eps_2 v \cdot \eps_1 \eps_2 w = \eps_1 \eps_2 (v+w).
\eqno (23.3)
$$
Next, we prove  the fundamental relation (23.1):
by definition, the Lie bracket of $v,w \in \g$ is given by
evaluating $[v^R,w^R]$ at $e$, where $v^R, w^R$
are the right-invariant vector fields extending $v,w$ (cf.\ our convention
concerning the sign of the Lie bracket from Section 5.3).
Recall that $v^R$ is given by left multiplication by $v$ in $TG$,
so the corresponding infinitesimal automorphism is left translation
$l_v:TG \to TG$, $u \mapsto v \cdot u$.
The vector field $v^R$ gives rise to two sections of $TTG$ over $G$,
with the notation from  Section 14.2 denoted by $\eps_i  v^R$ ($i=1,2$),
and each section gives rise to a diffeomorphism
$\tilde{\eps_i  v^R}$ of $TTG$ (Theorem 14.3).
These are nothing but the two versions of the tangent map $T  l_v:TTG \to
TTG$, and hence both can be written as left multiplications in $TTG$,
namely  
 $\tilde{\eps_1  v^R}=l_{\epsilon_1 v}$ and
$\tilde{\eps_2  v^R} =l_{\epsilon_2 v}$ with $\epsilon_j v \in (TTG)_e$,
$j=1,2$, defined as above. Now, applying first the definition of the
Lie bracket in $\g$ and then  Theorem 14.4, we get
$$
\eqalign{
\epsilon_1 \epsilon_2 [v,w] & 
=\epsilon_1 \epsilon_2 [v^R,w^R](e) \cr
&=[\tilde{\eps_1 v^R},\tilde{\eps_2 w^R}]_{\Diff(TTG)} \, \,  (e) 
=[l_{\epsilon_1 v},l_{\epsilon_2 w}]_{\Diff(TTG)} \, \, (e) \cr
&= (\epsilon_1 v)(\epsilon_2 w)(\epsilon_1 v)^{-1}(\epsilon_2 w)^{-1}
= [\epsilon_1 v,\epsilon_2 w]_{TTG}  \cr}
$$
where we have added as index the group in which the respective 
commutators are taken.
\qed

\ssk \nin We will often use (23.1) in the following form:
$$
\eps_1 v_1 \cdot \eps_2 w_2 = \eps_1 \eps_2 [v_1,w_2] \cdot
\eps_2 w_2 \cdot \eps_1 v_1 .
\eqno (23.4)
$$
Another way of stating (23.1) is by saying that the second order tangent
map of the commutator map
$$
c:=c_G:G\times G \to G,\quad (g,h)\mapsto [g,h]=ghg^{-1}h^{-1}
$$
is determined by
$$
T^2 c(\epsilon_1 v, \epsilon_2 w) =\epsilon_1 \epsilon_2 [v,w];
$$
this follows simply by noting that $c_{TTG}=TTc_G$.
One may also work with the conjugation map
$G\times G \to G$, $(g,h)\mapsto ghg^{-1}$; then 
the fundamental relation (23.1) can also be written in the 
following way
$$
\Ad(\eps_1 v) \eps_2 w = \eps_2 w \cdot \eps_1 \eps_2 [v,w] =
(\id + \eps_1 \ad(v))\eps_2 w .
$$
This relation  is equivalent to the
the fact  that the differential
of $\Ad:G \times \g \to \g$ is $\ad:\g\times \g \to \g$, $(X,Y) \mapsto
[X,Y]$ (cf.\ [Ne02b, Prop. III.1.6] for the real locally convex case).

\msk \nin
{\bf 23.3.} {\sl Left- and right trivialization of $TG$ and of $TTG$.}
So far we have avoided to use the well-known fact that the
tangent bundle $TG$ of a Lie group $G$ can be trivialized: the sequence
of group homomorphisms given by injection of $\g=T_eG$ in $TG$ and
projection $TG \to G$,
$$
\matrix{
0& \to & \g & \to & TG & \to & G & \to & 1 \cr},
$$
is exact and splits via the zero-section $z:G\to TG$.
Thus, group theoretically, $TG$ is a semidirect product of $G$ and of $\g$.
Explicitly, the map
$$
\Psi_1^R: \g \times G \to TG, \quad (v,g) \mapsto v \cdot g =
Tm(v,0_g)= T_e(r_g)v= l_v (0_g)
\eqno (23.5)
$$
is a smooth bijection with smooth inverse
$TG \to \g \times G$,
$u \mapsto (u \pi(u)^{-1},\pi(u))$.
The diffeomorphism $TG \cong \g \times G$ thus obtained
is called the {\it right trivialization of the tangent bundle}.
Similarly, the left trivialization $\Psi_1^L$ is defined.
The kernel of the canonical projection $TG \to G$ is $\g=T_e G$
on which $G$ acts via the {\it adjoint representation}
$$
\Ad: G \times \g \to \g , \quad (g,v) \mapsto \Ad(g)v=gvg^{-1}
$$
(product taken in $TG$; then as usual
 $\Ad(g):\g \to \g$ is the differential
of the conjugation by $g$ at the origin). From Formula (5.1) 
we see that the normal subgroup $\g$ of $TG$ is
abelian with product being vector addition. It follows that the group
structure of $TG$ in the right trivialization is given by
$$
(v,g) \cdot (w,h)= v\cdot g \cdot w \cdot h=v \cdot (gwg^{-1}) \cdot gh
=(v+\Ad(g)w,gh),
\eqno (23.6)
$$
which is the above mentioned realization of
$TG$ is a semidirect product.
Similarly, $TTG$ is described as an iterated semidirect product:
the right trivialization map of the group
 $TTG=T_{\eps_2} T_{\eps_1} G$ is obtained by
replacing in the right trivialization map for $TG$,
$\Psi^R_1:\g \times G \to TG$, $(\eps v,g) \mapsto  \eps v \cdot g$,
the group $G$ by $TG=T_{\eps_1} G$ and letting $\eps=\eps_2$.
We get a diffeomorphism
$$
\eqalign{
\Psi_2^R: \eps_1 \eps_2 \g \times \eps_2 \g \times \eps_1 \g \times  G   
& \to
TTG, \cr 
(\eps_2 \eps_1 v_{12},\eps_2 v_2,\eps_1 v_1,g) & \mapsto
 \eps_2 \eps_1 v_{12}\cdot \eps_2 v_2   \cdot \eps_1 v_1  \cdot g,
 \cr}
\eqno (23.7)
$$
where the products are taken in the group $TTG$. Similarly, we get the
left trivialization
$$
\eqalign{
\Psi_2^L: G \times \eps_1 \g   \times \eps_2 \g \times  \eps_1 \eps_2 \g 
& \to TTG, \cr
(g,\eps_1 v_1,\eps_2 v_2,\eps_2 \eps_1 v_{12}) & \mapsto
g \cdot \eps_1 v_1  \cdot \eps_2 v_2   \cdot  \eps_2 \eps_1 v_{12}. \cr}
\eqno (23.8)
$$
Clearly, both trivializations define  linear structures on
the fiber $(TTG)_e$  and hence on all other fibers. 
We will show below (Theorem 23.5) that both
linear structures indeed define  linear connections on $TG$.

\Theorem 23.4. With respect to the right trivialization $\Psi_2^R$,
 the group structure
of the subgroup $(TTG)_e \cong \epsilon_1 \epsilon_2 \g\times \epsilon_2 \g
\times \epsilon_1 \g$ is given  by the product
$$
\eqalign{
& \Psi_2^R(\eps_1 \eps_2 v_{12},\eps_2 v_2, \eps_1 v_1)  \cdot \Psi_2^R
(\eps_1 \eps_2 w_{12},\eps_2 w_2, \eps_1 w_1)  \cr
& = \Psi_2^R
(\eps_1 \eps_2 (v_{12}+w_{12}+[v_1,w_2]),\eps_2 (v_2+w_2), \eps_1 (v_1 +w_1))
\cr}
\eqno (23.9)
$$
and inversion
$$
(\Psi_2^R(\eps_1 \eps_2 v_{12},\eps_2 v_2, \eps_1 v_1))^{-1} =
\Psi_2^R(\eps_1 \eps_2 ([v_1,v_2] - v_{12}),-\eps_2 v_2, -\eps_1 v_1) .
\eqno (23.10)
$$

\Proof.  Using the Commutation Rules (23.1) -- (23.4), we get
$$
\eqalign{ \Psi_2^R
(\eps_1 \eps_2 v_{12}&,\eps_2 v_2, \eps_1 v_1)  \cdot 
\Psi_2^R(\eps_1 \eps_2 w_{12},\eps_2 w_2, \eps_1 w_1)  \cr
& =
\eps_1 \eps_2 v_{12} \cdot \eps_2 v_2 \cdot  \eps_1 v_1 \cdot 
\eps_1 \eps_2 w_{12} \cdot \eps_2 w_2 \cdot  \eps_1 w_1
\cr
&=
\eps_1 \eps_2 (v_{12}+w_{12}) \cdot \eps_2 v_2 \cdot 
 \eps_1 \eps_2 [v_1,w_2] \cdot
\eps_2 w_2 \cdot \eps_1 v_1 \cdot  \eps_1 w_1
\cr
 & =
\eps_1 \eps_2 (v_{12}+w_{12}+[v_1,w_2]) \cdot \eps_2 (v_2+w_2) \cdot
 \eps_1 (v_1 +w_1)) \cr
& = \Psi_2^R
(\eps_1 \eps_2 (v_{12}+w_{12}+[v_1,w_2]),\eps_2 (v_2+w_2), \eps_1 (v_1 +w_1))
\cr}
$$
and
$$
\eqalign{
(\Psi_2^R(\eps_1 \eps_2 v_{12},\eps_2 v_2, \eps_1 v_1))^{-1} & =
(\eps_1 \eps_2 v_{12} \cdot \eps_2 v_2 \cdot  \eps_1 v_1)^{-1} \cr
& = \eps_1 (-v_1) \cdot \eps_2(-v_2) \cdot \eps_1 \eps_2 (-v_{12}) \cr
& = \eps_2 (-v_2) \cdot \eps_1(-v_1) \cdot \eps_1 \eps_2 ([-v_1,-v_2]-v_{12})
 \cr
& = \eps_1 \eps_2 ([v_1,v_2]-v_{12}) \cdot \eps_2 (-v_2) \cdot \eps_1(-v_1) \cr
& = \Psi_2^R
(\eps_1 \eps_2 ([v_1,v_2]-v_{12}), \eps_2 (-v_2) , \eps_1(-v_1)) \cr}
$$
\qed

\nin The formula for the group structure of the whole group $TTG$ in the
right trivialization is easily deduced from the theorem: considering
$\Psi_2^R$ as an identification, it reads
$$
\eqalign{
& (\eps_1 \eps_2 v_{12},\eps_2 v_2, \eps_1 v_1,g)  \cdot 
(\eps_1 \eps_2 w_{12},\eps_2 w_2, \eps_1 w_1,h) =   \cr
&
(\eps_1 \eps_2 (v_{12}+\Ad(g)w_{12}+[v_1,\Ad(g)w_2]), \eps_2 (v_2+\Ad(g)w_2),
 \eps_1 (v_1 +\Ad(g)w_1),gh) \cr} 
\eqno (23.11)
$$
In particular, the group $G$ acts (by conjugation) 
{\it linearly} with respect to
the linear structure on $(T^2 G)_e$, and hence we can 
transport this linear structure
 by left- or right translations to all other fibers
$(T^2 g)_x$, $x \in G$. Clearly, this depends smoothly on $x$, and hence
we have defined two linear structures on $TTG$ which will be denoted by
$L^L$ and $L^R$, called the {\it left} and {\it right linear structures}.

\Theorem 23.5. The left and right linear structures $L^L$ and $L^R$ are
connections on $TG$, i.e.\ they are bilinearly related to all chart
structures. The difference $L^L - L^R$ is the tensor field of type
$(2,1)$ on $G$ given by the Lie bracket, and this is also the torsion
of $L^L$.

\Proof. Let us prove first that $L^L$ and $L^R$ are bilinearly related
to all chart structures. Recall from Section 17.6 that, for any manifold
$M$, a trivialization $TM \cong M \times V$ of the tangent bundle defines
a connection on $TM$. More specifically, in the present case, deriving
$\Psi_1^R:\g \times G \to TG$,
$$
T \Psi_1^R: T\g \times TG \to TTG, \quad (X,h) \mapsto X \cdot h,
$$
and composing again with $\Psi_1^R$, we get a trivialization of $TTG$:
$$
T(\Psi_1^R) \circ (\id \times \Psi_1^R ): 
\eps_1 \eps_2 \g \times \eps_2 \g \times \eps_1 \g \times G=
T\g \times (\g \times G) \to TTG,
\eqno (23.12)
$$ 
and the linear structure $L$ 
induced by this trivialization is in fact a connection (see Lemma 17.7).
We claim that $L=L^L$.
In fact, writing $\Psi_1 = Tm \circ (\iota \times z)$, where 
$\iota:\g \to TG$ and $z:G \to TG$ are the natural inclusions, 
we get $T\Psi_1 = TTm \circ(T \iota \times Tz)$, which means with
$T=T_{\eps_2}$,
$$
\eqalign{
T(\Psi_1^R)(  \eps_1 \eps_2 v_{12}, \eps_2 v_2,\Psi_1^R(\eps_1 v_1,g )) &=
(\eps_1 \eps_2 v_{12} \cdot \eps_1 v_1 ) \cdot (\eps_2 v_2 \cdot g) \cr
&=\eps_1 v_1 \cdot  \eps_2 v_2 \cdot \eps_1 \eps_2 v_{12} \cdot g. \cr}
$$
For $g=e$, this agrees with $\Psi_2^L$, and hence we have $L_e = (L^L)_e$. 
But both linear structures are invariant under left- and right translations
and hence agree on all fibers. We have proved that $L^L$ is a connection.
In the same way we see that $L^R$ is a connection.

We show that, on the fiber $(TTG)_e$, the difference between $L^L$ and
$L^R$ is the Lie bracket. In fact, using (23.1), we can calculate
the map relating both linear structures:
$$
\eqalign{(\Psi_2^L)^{-1} \circ \Psi_2^R 
(\eps_2 \eps_1 v_{12},\eps_1 v_2 ,\eps_1 v_1) & =
(\Psi_2^L)^{-1}(\eps_2 \eps_1 v_{12} \cdot \eps_2 v_2 \cdot \eps_1 v_1) \cr
&=
(\Psi_2^L)^{-1}( \eps_1 v_1  \cdot \eps_2 v_2 \cdot 
\eps_2 \eps_1 (v_{12} + [v_2,v_1]) ) \cr
&=(\eps_1 v_1,\eps_2 v_2,\eps_1 \eps_2 (v_{12} + [v_2,v_1])). \cr}
\eqno (23.13)
$$
Identifying $\eps_1 \g \times \eps_2 \g \times \eps_1 \eps_2 \g$ and
$\eps_1 \eps_2 \g \times  \eps_2 \g  \times \eps_1 \g$ in the canonical
way, this means that $L^L$ and $L^R$ are related via the map
$$
f_b(u,v,w) = (u,v,w+[u,v]),
\eqno (23.14)
$$
i.e.\ they are bilinearly related via the Lie bracket. 

In order to prove that $L^R - L^L$ is the torsion tensor of $L^R$, 
we have to show that $L^R = \kappa \cdot L^L$, where $\kappa:TTG\to TTG$
is the canonical flip. 
But, since $\kappa$ commutes with second tangent maps, $\kappa$ is
a group automorphism of $(TTG,TTm)$, whence
$$
\eqalign{
\kappa \circ \Psi_2^L(\eps_1 v_1,\eps_2 v_2, \eps_1\eps_2 v_{12}) & =
\kappa (\eps_1 v_1 \cdot \eps_2 v_2 \cdot \eps_1 \eps_2 v_{12}) \cr
& =
\kappa (\eps_1 v_1) \cdot \kappa(\eps_2 v_2) \cdot \kappa(\eps_1\eps_2 v_{12})
\cr
& = \eps_2 v_1 \cdot \eps_1 v_2 \cdot \eps_1 \eps_2 v_{12} \cr
&= \Psi_2^R(\eps_1 v_2,\eps_2 v_1,\eps_1 \eps_2 v_{12}) \cr
&= \Psi_2^R(\kappa(\eps_1 v_1,\eps_2 v_2,\eps_1 \eps_2 v_{12})) \cr}
\eqno (23.15)
$$
and hence $\Psi_2^R = \kappa \circ \Psi_2^L \circ \kappa$, which is
equivalent to $L^R = \kappa \cdot L^L$.
\qed                                                

\nin
{\bf 23.6.} {\sl Remarks on the symmetric connection.}
If ${1 \over 2} \in \K$, then the connection $L^S:={L^R + L^L \over 2}$
is called the {\it symmetric connection of $G$}; it is torsionfree and
invariant under left- and right-translations. Moreover, it is invariant
under the inversion map $j:G \to G$ because $j.L^L = L^R$. In fact,
this is the canonical connection of $G$ seen as a symmetric space,
cf.\ Chapter 26.

\Theorem 23.7. {\rm (Structure of $(J^2 G)_e$.)}
For $F=L,R,S$ (the latter in case ${1 \over 2} \in \K$), the bijections
$\Psi^F_2:\eps_1 \g \times \eps_2 \g \times \eps_1 \eps_2 \g \to (TTG)_e$
restrict to bijections
$$
J \Psi_2^F: \delta \g \times \delta^{(2)} \g \to (J^2 G)_e .
$$
These three maps coincide and are given by the explicit formula
$$
\Psi_2  (\delta^{(2)} v, \delta w)= 
\delta^{(2)} v  \cdot \delta w 
= \eps_1 \eps_2v \cdot \eps_1 w \cdot \eps_2 w,
$$
and the group structure of $(J^2 G)_e$ is given by
$$
\Psi_2(\delta^{(2)} v', \delta v ) \cdot 
\Psi_2(\delta^{(2)} w', \delta w)  = \Psi_2
(\delta^{(2)} (v'+w'+[v,w]), \delta (v+w) ).
$$
Inversion in $(J^2G)_e$ is multiplication by the scalar $-1$.

\Proof.
We have remarked above (Equation (23.15)) that the canonical flip
$\kappa$ acts on $(TTG)_e$ by
$\kappa (\eps_1 v_1 \cdot \eps_2 v_2 \cdot \eps_1 \eps_2 v_{12}) 
=  \eps_1 \eps_2 v_{12} \cdot  \eps_2 v_1 \cdot \eps_1 v_2$.
Now all claims follow from Formulae (23.7) -- (23.10) by observing
that $[v,v]=0$. \qed 

\msk \nin 
Finally, for the sake of completeness we add some results that are well-known
to hold in the real (say, Banach) case ; these results will not be needed
in the sequel.

\Lemma 23.8. 
The covariant derivatives associated to the three canonical connections
of the Lie group $G$ can be described as follows:
fix $x \in G$; for $v \in T_x G$ denote by
$v^L$, resp. $v^R$ the unique left- (resp. right-) invariant
vector field taking value $v$ at $x$.
Then, for all vector fields $X,Y \in \X(G)$,
$$
\eqalign{
(\nabla^L_X Y)(x) & =[Y,(X(x))^R](x), \quad \cr
(\nabla^R_X Y)(x) & =[Y,(X(x))^L](x), \cr
(\nabla^S_X Y)(x) & ={1 \over 2} [Y,(X(x))^L+(X(x))^R](x). \cr}
$$

\Proof. 
Note that $v^R$ is constant in the right trivialized picture, 
and hence $K^R \circ T v^R = 0$ for the connector of $L^R$.
It follows that
$$
\eqalign{
(\nabla^R_v Y)(x) & = K^R \circ TY \circ v^R (x) =
K^R \circ (TY \circ v^R - Tv^R \circ Y)(x) \cr
& =K^R \circ [Y,v^R](x)= [Y,v^R](x). \cr}
$$
Similarly for $\nabla^L$. The last equality follows by taking the
arithmetic mean of the preceding ones.
\qed

\Theorem 23.9. The curvature tensor of the connections $L^L$ and $L^R$
vanishes.

\Proof. 
Let $R^R$ be the curvature tensor of the connection $L^R$.
We use the ``classical" characterization of $R^R$ via 
$R^R(X,Y)Z = ([\nabla_X,\nabla_Y]-\nabla_{[X,Y]}) Z$ (Theorem 18.6).
Taking values at the origin $e$ and using Lemma 23.8, we see
that the right hand side equals
$(([\ad(X^R),\ad(Y^R)] - \ad[X^R, Y^R]) Z^R)(e)$ which is zero
since the space of right-invariant vector fields is a Lie algebra.
\qed

\nin Since $L^L$ and $L^R$ have torsion, the preceding
result should not be interpreted in the sense that $L^L$ and $L^R$ be flat. 

\msk \nin {\bf 23.10.} {\sl The adjoint Maurer-Cartan form.}
We denote by $\omega$ the tensor field of type $(2,1)$ on
$G$ given at each point by the Lie bracket (torsion of $L^R$). Using the
right trivialization of $TG$, it can also be seen as a
$\End(\g)$-valued one-form
$$
\omega:TG \cong G \times \g \to \End( \g), \quad (g,v) \mapsto \ad(v).
$$
We call this the {\it adjoint Maurer-Cartan form} since 
our $\End(\g)$-valued form is obtained from the usual $\g$-valued 
form by composing with the adjoint representation $\ad:\g \to \End(\g)$,
cf.\  e.g.\ [Sh97, p.\ 108].

\Theorem 23.11.
The adjoint Maurer-Cartan form satisfies the structure equation
$$
\dd \omega(X,Y) + [\omega_X,\omega_Y] = 0
$$
for all $X,Y \in \X(G)$ (where the bracket is a commutator
in $\End(\g)$). The curvature $R^S$ of the connection $L^S$
is given by taking pointwise triple Lie brackets:
$$
R^S(X,Y)Z ={1 \over 4} \dd \omega(X,Y)Z 
={1 \over 4} \omega(\omega(X,Y),Z) = {1 \over 4} [[X,Y],Z].
$$

\Proof. We have seen that the curvatures $R^L$ and $R^R$
of the connections $L^L$, resp. $L^R$ vanish. Thus
we can apply the structure equation (18.9), with the
difference term $A=\omega$, which leads to the first claim.

Next we apply the structure equation (18.10), but this time
with $R=R^L =0$ and $R'=R^S$, the symmetric curvature.
The difference term is now $A={1 \over 2} \omega$.
Thus 
$R^S={1\over 2}\dd \omega + {1 \over 8} [\omega,\omega]=
- {1 \over 8}[\omega,\omega]=
{1 \over 4} \dd \omega$.
\qed

\vfill \eject 

\subheadline{24. The structure of higher order tangent groups}

\nin We retain the assumptions from the preceding chapter:
 $G$ is a Lie group with multiplication $m:G \times G \to G$.
We are going to investigate the structure of the groups
$T^k G$ for $k \geq 3$.

\msk \nin {\bf 24.1.} {\sl Fundamental commutation relations for
elements of the axes.}
Let $\eps^\alpha v_\alpha$, $\eps^\beta w_\beta \in (T^k G)_e$ be
elements of axes in $(T_k G)_e$.
Then the commutator of these two elements in the group $(T^k G)_e$
is given by  the following fundamental
commutation rules: 
$$
[\eps^\alpha v_\alpha,\eps^\beta w_\beta]= \big\{ \matrix{
\eps^{\alpha +  \beta} [v_\alpha,w_\beta] & {\rm if} \, \alpha \perp \beta,
\cr 0 & {\rm else.} \cr}
\eqno (24.1)
$$
In fact, 
if $\supp(\alpha) \cap \supp(\beta) \not= \emptyset$,
then both elements commute:
$[\eps^\alpha v_\alpha,\eps^\beta v_\beta]=0$, 
because both belong to the tangent space $T_{\eps_i} (T^{k-1} G)$ at
the origin for some $\eps_i$ (namely for $i$ such that $e_i \in
\supp(\alpha) \cap \supp(\beta)$), and the group structure on the 
tangent space is just vector addition. This proves the second relation.
If $\alpha \perp \beta$, then we apply
 (23.1) with $\eps_1$ replaced by $\eps^\alpha$
and $\eps_2$ replaced by $\eps^\beta$, and we get the first relation.
In a unified way, (24.1) may be written
$$
[\eps^\alpha v_\alpha,\eps^\beta w_\beta]= 
\eps^{\alpha} \eps^{\beta} [v_\alpha,w_\beta]
\eqno (24.2)
$$
because $\eps^\alpha \eps^\beta=0$ if 
$\supp(\alpha) \cap \supp(\beta) \not= \emptyset$.
This can also be written
$$
\eps^\alpha v_\alpha \cdot \eps^\beta w_\beta=
\eps^{\alpha} \eps^{\beta} [v_\alpha,w_\beta] \cdot
\eps^\beta w_\beta \cdot \eps^\alpha v_\alpha =
\eps^\beta w_\beta \cdot \eps^\alpha v_\alpha \cdot
\eps^{\alpha} \eps^{\beta} [v_\alpha,w_\beta].
\eqno (24.3)
$$

\msk \nin {\bf 24.2.} {\sl Remark: $(T^k G)_e$ as a filtered group.}
Let $G_1 := (T^k G)_e$ and $G_j := \langle \eps^\alpha v_\alpha \vert \,
\vert \alpha \vert > j, \, v_\alpha \in \g \rangle$ 
be the ``vertical subgroup''
generated by elements from axes with total degree bigger than
$j$. Then it follows from (24.1) that
$$
\1 = G_{k} \subset G_{k-1} \subset \ldots \subset G_1
$$
is a {\it filtration} of $G_1$ in the sense of [Se65, LA 2.2] or
[Bou72, II.  4.4]. 
 Recall from [Se65, LA 2.3] or [Bou72] that the associated 
{\it graded group} carries a canonical structure of a Lie algebra.
This Lie algebra structure is compatible with ours, but in our
setting we can go further by constructing special bijections
(coming from the three canonical connections)
between the graded and the filtered objects.

\msk \nin {\bf 24.3.} {\sl Right and left trivialization of $T^k G$.}
Note that, according to our conventions, $T^k G$ is viewed as obtained
by successive scalar extensions via $\eps_1,\ldots,\eps_k$, i.e.\
$T^{k+1} G = T_{\eps_{k+1}} (T^k G)$.
The iterated right and left trivialization of $T^k G$ are defined inductively:
for $k=1,2$ they are given by (23.2) -- (23.4); then the right 
trivialization of $T^k G$ is given by appliying Formula (23.3) to
the right trivialization of $H:= T^{k-1} G$, with $TH=T_{\eps_k} H$.
For instance, the third order right trivialization is
$$
\eqalign{
\Psi_3^R : & (\eps_1 \eps_2 \eps_3 v_{123},\eps_3 \eps_2 v_{23},
\eps_3 \eps_1v_{13}, \eps_3 v_3, \eps_1 \eps_2 v_{12},
\eps_2 v_2,\eps_1 v_1,g) \mapsto \cr
& \quad 
\eps_1 \eps_2 \eps_3 v_{123} \cdot \eps_3 \eps_2 v_{23} \cdot 
\eps_3 \eps_1v_{13} \cdot  \eps_3 v_3 \cdot  \eps_1 \eps_2 v_{12} \cdot
\eps_2 v_2 \cdot \eps_1 v_1 \cdot g \cr}
$$
(products taken in $T^3 G$). The important feature is the order in which
the product in $T^3 G$ has to be taken: here it is the ``antilexographic''
order (opposite order of the lexicographic order).
An easy induction shows that the general formulae are
$$
\eqalign{
\Psi_k^R :  \bigoplus_{\alpha > 0} \eps^\alpha \g \times G \to
T^k G, & \quad 
 (\eps^\alpha v_\alpha)_{\alpha > 0}, g) \mapsto
\prod^{\downarrow}_{\alpha > 0} \eps^\alpha v_\alpha \cdot g, 
\cr
\Psi_k^L : \bigoplus_{\alpha > 0} \eps^\alpha \g \times G \to
T^k G, & \quad 
 (\eps^\alpha v_\alpha)_{\alpha > 0}, g) \mapsto
g \cdot \prod^{\uparrow}_{\alpha > 0} \eps^\alpha v_\alpha , \cr}
\eqno (24.4)
$$
where $\prod^\uparrow$ is the product of the elements labelled
by $\alpha >0$ in $(T^k G)_e$ taken with respect to the
lexicographic order of the index set
and $\prod^\downarrow$ is the product taken in
$(T^k G)_e$ in the opposite (``antilexicographic'') order.

\msk \nin {\bf 24.4.} {\sl The case $k=3$ and the Jacobi identity.}
As an application of our formalism, let us now give an
independent proof of the Jacobi identity (it is essentially equivalent
to the arguments used in [Se65]): we calculate both sides of
$\eps_3 v_3  \cdot  (\eps_2 v_2 \cdot \eps_1 v_1)  =
(\eps_3 v_3  \cdot \eps_2 v_2) \cdot \eps_1 v_1$ 
in $(T^3 G)_e$ and express the result in the left trivialization: on the one
hand,
$$
\eqalign{
&\eps_3  v_3  \cdot  (\eps_2 v_2 \cdot \eps_1 v_1)  =  \cr
&\eps_3 v_3 \cdot (\eps_1 v_1 \cdot \eps_2v_2 \cdot \eps_1 \eps_2 [v_2,v_1])=
\cr
&\eps_1 v_1 \cdot \eps_3 v_3 \cdot \eps_1 \eps_3 [v_3,v_1]   \cdot
\eps_2 v_2 \cdot \eps_1 \eps_2 [v_2,v_1]=   
\cr
&\eps_1 v_1 \cdot \eps_2v_2 \cdot \eps_3 v_3 \cdot \eps_1 \eps_3 [v_3,v_1]
\cdot \eps_2 \eps_3 [v_3,v_2] \cdot \eps_1 \eps_2 \eps_3 [[v_3,v_1],v_2] 
\cdot \eps_1 \eps_2 [v_2,v_1] =  
\cr
&\eps_1 v_1 \cdot \eps_2v_2 \cdot \eps_1 \eps_2 [v_2,v_1] 
\cdot \eps_3 v_3 \cdot \eps_1 \eps_3 [v_3,v_1]\cdot
\eps_2 \eps_3 [v_3,v_2] \cr
& \quad \quad \cdot  \eps_1 \eps_2 \eps_3 ([[v_3,v_1],v_2]+
[v_3,[v_2,v_1]])   =
\cr
&\Psi_3^L\Bigl( \eps_1 v_1 + \eps_2v_2 + \eps_1 \eps_2 [v_2,v_1] 
+ \eps_3 v_3 + \eps_1 \eps_3 [v_3,v_1] +
\eps_2 \eps_3 [v_3,v_2]    \cr
& \quad \quad   + \eps_1 \eps_2 \eps_3 ([[v_3,v_1],v_2]+ [v_3,[v_2,v_1]]
\Bigr). 
\cr}
$$
On the other hand, 
$$
\eqalign{
& (\eps_3 v_3  \cdot \eps_2 v_2) \cdot \eps_1 v_1  =
\eps_2 v_2 \cdot \eps_3 v_3  \cdot 
\eps_2 \eps_3 [v_3,v_2] \cdot \eps_1 v_1 = \cr
&  \eps_1 v_1 \cdot
\eps_2 v_2 \cdot \eps_1 \eps_2 [v_2,v_1] \cdot 
\eps_3 v_3 \cdot \eps_1 \eps_3 [v_3,v_1] \cdot
 \eps_2 \eps_3 [v_3,v_2] \cdot \eps_1 \eps_2 \eps_3 [[v_3,v_2],v_1] = \cr
& \Psi_3^L\Bigl( \eps_1 v_1 +
\eps_2 v_2 + \eps_1 \eps_2 [v_2,v_1] + 
\eps_3 v_3 + \eps_1 \eps_3 [v_3,v_1] +
 \eps_2 \eps_3 [v_3,v_2]  \cr
& \quad \quad + \eps_1 \eps_2 \eps_3 [[v_3,v_2],v_1]\Bigr). \cr}
 $$
Applying $(\Psi_3^L)^{-1}$ and comparing, we get
$$
[[v_3,v_1],v_2]+[v_3,[v_2,v_1]] = [[v_3,v_2],v_1]
$$
which is the Jacobi identity. Note that this proof does not make any use
of the interpretation of the Lie bracket via left or right
invariant vector fields, and note also that the expression of the
 Jacobi identity that we get here is the one defining the more general
class of {\it Leibniz algebras}, see Remark 24.12 below.  --
Another, rather tricky proof of the Jacobi identity
makes essential use of Theorem 6.2 by combining
algebraic and differential geometric arguments:
the tangent map of the adjoint action $\Ad_G:G \times \g \to \g$ is
$\Ad_{TG}:TG \times T\g \to T\g$. We know that $\Ad(g)$, for an element
 $g \in G$,
is a Lie algebra automorphism of $\g$, and hence $\Ad(h)$, for $h \in TG$,
is a Lie algebra automorphism of $T\g$. But Theorem 6.2 implies
that, moreover, $\Ad(h)$ is an automorphism of $T\g$ over the ring $T\K$.
Now, it is a simple exercise in linear algebra to show that
$\Aut_{T\K}(T\g)$ is a semidirect product of $\Aut_\K(\g)$ with the vector
group $\Der_\K(\g)$ of derivations. 
Putting things together, we see that, for all $X \in \g$,
$\Ad(\eps X)$ acts like a derivation from $\g$ to $\eps \g$,
 which gives again the Jacobi 
identity. This proof has the advantage that it is extremely close to
the well-known arguments from the theory of finite-dimensional real 
Lie groups.

\msk \nin {\bf 24.5.} {\sl Action of the symmetric group.}
The natural action of the symmetric group $\Sigma_k$ 
by automorphisms of         
$(T^k G)_e$ is described by the following formula: for $\sigma \in \Sigma_k$,
$$
\sigma \cdot \prod^\uparrow_\alpha  \eps^\alpha v_\alpha  
= \prod_\alpha^\uparrow \eps^{\sigma.\alpha} v_\alpha. 
\eqno (24.4)
$$
As a consequence, the push-forward $\Psi^{L,\sigma}=\sigma \circ \Psi^L
\circ \sigma^{-1}$ of $\Psi^L$ by $\sigma$ is given by
$$
\Psi^{L,\sigma}((\eps^\alpha v_\alpha)_{\alpha > 0}, g) =
g \cdot \prod^{\uparrow}_{\alpha > 0} \eps^{\sigma.\alpha} 
v_{\sigma.\alpha} .
\eqno (24.5)
$$
In fact, Formula (24.4) is proved by
exactly the same arguments as in case $k=2$ for
$\sigma = \kappa$ (cf.\ Equation (23.15)).
Re-ordering the terms and using the commutation rules (24.1), one
can calculate the curvature forms which are defined as in Chapter 18.
If $G$ is commutative, then by a change of variables the
left-hand side can be re-written $\prod^\uparrow \eps^\alpha
v_{\sigma^{-1} \alpha}$, and this can be interpreted by saying that
$\Psi^L$ is $\Sigma_k$-invariant. Clearly, then $\Psi^L$ and  
$\Psi^{L,\sigma}$ coincide.  However, if $G$ is non-commutative,
then we cannot re-write the formula in the way just mentioned
since the order would be disturbed.  In particular, for $k>2$,
$\Psi_k^L$ and $\Psi_k^R$ are no longer related via a permutation
since reversing the lexicographic order on $I_k$ is not induced by
an element of $\Sigma_k$ acting on $I_k$.

\msk \nin
{\bf 24.6.} {\sl Left and right group structures.}
We have defined two canonical charts $\Psi_k^F$, $F=R,L$,
 of the group $(T^k G)_e$, and
we want to describe the group structure in these charts. The model space
is $E=\bigoplus_{\alpha > 0} \eps^\alpha \g$, which is just
the tangent space $T_0((T^k G)_e)$ -- being a linear space, the tangent
space is canonically isomorphic to the direct sum of the axes.
Via $\Psi_k^F$ ($F=R,L$) we transfer the group structure to $E$. This
defines a product and an inverse,
$$
\eqalign{*^F: E \times E \to E, & \quad
 \Psi_k^F v  \cdot \Psi_k^F w  = \Psi_k^F(v *^F w), \cr
J^F: E \to E, & \quad
(\Psi_k^F (v))^{-1} = \Psi_k^F (J^F v), \cr}
\eqno (24.6)
$$
called the {\it right}, resp.\ {\it left group structure on $E$}. 
In Theorem 23.3 we have given explicit formulae for these maps in case
$k=2$. For the general case we use the same strategy: just apply the
Commutation Rules 24.1 along with $\eps^\alpha v \cdot \eps^\alpha w= 
\eps^\alpha(v+w)$ in order to re-order the product. Before coming to
the general result, we invite the reader to check, using the commutation
rules, that for $k=3$ we have the following:
with $v=\sum_{\alpha \in I_3 \atop \alpha > 0} \eps^\alpha v_\alpha$,
$w=\sum_{\alpha \in I_3 \atop \alpha > 0} \eps^\alpha w_\alpha$, we get
$$
\eqalign{
& \Psi^L_3(v)  \cdot \Psi^L_3(w)  \cr
& 
= \eps^{001}v_{001}\cdot \eps^{010}v_{010} \cdot \eps^{011} v_{011} \cdot
\eps^{100} v_{100} \cdot \eps^{101} v_{101} \cdot 
\eps^{110}v_{110} \cdot  \eps^{111} v_{111}
\cr
& 
 \cdot \eps^{001}w_{001}\cdot \eps^{010}w_{010} \cdot \eps^{011} w_{011} \cdot
\eps^{100} w_{100} \cdot \eps^{101} w_{101} \cdot 
\eps^{110}w_{110} \cdot  \eps^{111} w_{111} 
\cr
& 
= \eps^{001}(v_{001}+w_{001}) \cdot \eps^{010}(v_{010}+w_{010}) \cdot
 \eps^{011} (v_{011}+w_{011} + [v_{010},w_{001}]) \cdot
\cr
& 
 \eps^{100} (v_{100}+w_{100}) \cdot \eps^{101} (v_{101}+w_{101}+
[v_{100},w_{001}]) \cdot 
\eps^{110}(v_{110}+w_{110}+[v_{100},w_{010}])
\cr &  
 \cdot \eps^{111}  (v_{111}+w_{111}+[v_{110},w_{001}]+[v_{101},w_{010}] +
 [v_{100},w_{011}] + [[v_{100},w_{001}], w_{010}]) 
\cr
&  = \Psi_3^L \Bigl(
v+w+ \eps^{011}[v_{010},w_{001}]+\eps^{101}[v_{100},w_{001}]+
\eps^{110} [v_{100},w_{010}] + \cr
& \quad  \quad \quad
\eps^{111} ([v_{110},w_{001}]+[v_{101},w_{010}] +
 [v_{100},w_{011}] + [[v_{100},w_{001}], w_{010}]) \Bigr).
\cr}
$$
Similarly, we get for the inversion:  
$$
\eqalign{
& (\Psi_3^L  v)^{-1}  = 
(\eps^{001}v_{001}\cdot \eps^{010}v_{010} \cdot \eps^{011}v_{011}
 \cdot
\eps^{100} v_{100} \cdot \eps^{101} v_{101} \cdot 
\eps^{110}v_{110}\cdot  \eps^{111}v_{111})^{-1} 
\cr  
& = 
\eps^{111} (-v_{111}) \cdot \eps^{110} (-v_{110}) \cdot \eps^{101}(- v_{101}) 
\cdot \eps^{100} (-v_{100}) \cdot \eps^{011} (-v_{011}) \cdot
\cr
& 
\quad \quad   \eps^{010}(-v_{010}) \cdot \eps^{001} (-v_{001} ) 
\cr
& =  \eps^{011}(-v_{011}) \cdot
\eps^{100}(- v_{100}) \cdot \eps^{101} (-v_{101}) \cdot 
\eps^{110}(-v_{110}) \cdot  \eps^{111} (-v_{111} + [v_{100},v_{011}]))\cdot \cr
& \quad \quad
\eps^{010}(-v_{010})  \cdot \eps^{001} (-v_{001} ) 
\cr
& =
\eps^{010}(-v_{010}) \cdot
 \eps^{011}(-v_{011}) \cdot
\eps^{100}(- v_{100}) \cdot \eps^{101} (-v_{101}) \cdot 
\cr
& 
\quad \quad \eps^{110}(-v_{110}+[v_{100},v_{010}]) \cdot 
 \eps^{111} (-v_{111} + [v_{100},v_{011}]+[v_{101},v_{010}])) \cdot
 \eps^{001} (-v_{001} )  
\cr
& = 
\eps^{001} (-v_{001}) \cdot \eps^{010} (-v_{010}) \cdot
\eps^{011} (-v_{011}+[v_{010},v_{001}]) \cdot \eps^{100} (-v_{100}) \cdot 
\cr
& \quad \quad  \eps^{101} (-v_{101}+ [v_{100},v_{001}])  \cdot
\eps^{110} (-v_{110} + [v_{100},v_{010}]) \cdot \cr
& \quad \quad 
\eps^{111} (-v_{111}+[v_{110},v_{001}] + [v_{101},v_{010}]+[v_{100},v_{011}]-
[[v_{100},v_{010}],v_{001}]) 
\cr
& =  \Psi_3^L \Bigl(
-v + \eps^{011} [v_{010},v_{001}] + 
\eps^{101}[v_{100},v_{001}] + \eps^{110} [v_{100},v_{010}] +\cr
& \quad \quad 
\eps^{111} ([v_{110},v_{001}] + [v_{101},v_{010}]+[v_{100},v_{011}]-
[[v_{100},v_{010}],v_{001}]) \Bigr) .
 \cr}
$$

\Theorem 24.7. {\rm (Left and right product formulae for $(T^k G)_e$.)}
With respect to the left trivialization $\Psi_k^L$ of $T^k G$,  we have
$$
\sum_\alpha \eps^\alpha v_\alpha \, *^L \, \sum_\beta \eps^\beta w_\beta  =
\sum_\gamma \eps^\gamma z_\gamma
$$
with 
$$
z_\gamma  =  v_\gamma+w_\gamma + \sum_{m=2}^{|\gamma|}
\sum_{\lambda \in {\cal P}_m(\gamma)} [\ldots 
[[v_{\lambda^m},w_{\lambda^1}],w_{\lambda^2}],\ldots, w_{\lambda^{m-1}}].
$$
The left inversion formula is
$$
J^L(\sum_\alpha \eps^\alpha v_\alpha) =
\sum_\gamma \eps^\gamma (-v_\gamma+ \sum_{m=2}^{|\gamma|} (-1)^m
\sum_{\lambda \in {\cal P}_m(\gamma)}
[\ldots[[v_{\lambda^m},v_{\lambda^{m-1}}],v_{\lambda^{m-2}}], 
\ldots , v_{\lambda^1}] ),
$$
and the right  product is given by the formula
$$
(v *^R w)_\gamma =  v_\gamma+w_\gamma + \sum_{m=2}^{|\gamma|}
\sum_{\lambda \in {\cal P}_m(\gamma)}
[v_{\lambda^{m-1}}, \ldots ,[w_{\lambda^1},v_{\lambda^m}]].
$$
In these formulae, partitions $\lambda$ are always considered as
ordered partitions: $\lambda^1 <\ldots <\lambda^m$.

\Proof. According to (24.4), we have to show that
$$
\prod_{\alpha>0}^\uparrow  \eps^\alpha v_\alpha \cdot
\prod_{\beta>0}^\uparrow  \eps^\beta w_\beta =
\prod_{\gamma>0}^\uparrow  \eps^\gamma z_\gamma
$$
with $z_\gamma$ as in the claim.
Using the Commutation Rules 24.1, we will move all terms
of the form $\eps^\beta w_\beta$ to the left until we reach 
$\eps^\beta v_\beta$; then both terms give rise to a new term $\eps^\beta
(v_\beta + w_\beta)$. But every time we exchange $\eps^\beta w_\beta$ on its
way to the left with an 
element $\eps^\alpha v^\alpha$ with $\alpha > \beta$ and $\alpha \perp \beta$,
applying the commutation rule we get a term
$\eps^{\alpha + \beta} [v_\alpha,w_\beta]$.
Before doing the same thing with the next element from the second factor,
we may re-order the first factor: since $\alpha > \beta$ and 
$\alpha \perp \beta$, the term $\eps^{\alpha + \beta} [v_\alpha,w_\beta]$
commutes with all terms between position $\alpha$ and $\alpha + \beta$,
and hence it can be pulled to position $\alpha + \beta$ without giving
rise to new commutators. 

\msk
We start the procedure just described with the first  element
 $\eps^{0\ldots01} w_{0\ldots01}$ of the second factor, then take the
second element $\eps^{0\ldots 10} w_{0\ldots 10}$, and so on 
up to the last element. Each time we get
commutators of the type $\eps^{\alpha + \beta} [v_\alpha',w_\beta]$
with $v_\alpha'$ being the sum of $v_\alpha$ and all commutators produced
in the preceding steps.
This gives us precisely all iterated commutators of the form
$$
 [\ldots[v_{\mu^1},w_{\mu^2}]\ldots,
w_{\mu^{\ell}}]
$$
with $\mu^i \perp \mu^j$ ($\forall i \not= j$) and the
order conditions
$$
\mu^2 < \ldots < \mu^\ell, \quad \quad
\mu^1 > \mu^2, \, \quad
\mu^1+\mu^2 > \mu^3,\ldots,\mu^1+\ldots+\mu^{\ell-1} > \mu^\ell.
$$
These conditions imply $\mu^1 > \mu^j$ for all $j=1,\ldots,\ell-1$ because
the supports of the $\mu^i$ are disjoint and hence the biggest of the $\mu^j$
is bigger than the sum of all the others.
Hence $(\lambda^1,\ldots,\lambda^{\ell-1},\lambda^\ell):=
(\mu^2,\ldots,\mu^\ell,\mu^1)$ is a partition, and every partition is 
obtained in this way. 
Thus, by a change of variables $\lambda \leftrightarrow \mu$, we get
the formula for the left product.
Similarly for the right product.

\msk
The inversion formula is proved by the same kind of arguments: we write
$$
(\Psi_k^L(\sum_\alpha \eps^\alpha v_\alpha))^{-1} =
(\prod^\uparrow \eps^\alpha v_\alpha)^{-1} =
\prod^\downarrow \eps^\alpha (-v_\alpha)
$$
and start to re-order this product from the left, that is: we exchange the
first two terms on the left (they commute and hence this
 does not yet produce a commutator); then we pull the third element to
the first position (in case $k=2$ this produces already a commutator);
then we do the same with the fourth element, and so on. We collect 
the iteraded commutators and get precisely the formula for $J^L$ from
the claim.
\qed 

\bigskip
\nin Comparing our formulae
 with the ``usual" Campbell-Hausdorff formula
(say, for a nilpotent Lie group), one remarks that
in the usual formula,
inversion is just multiplication by $-1$, whereas here this is not the case.
 More generally, the powers $X^n$
with respect to the Campbell-Hausdorff  multiplication  correspond
to multiples $nX$, whereas here we have, e.g., the following
{\it (left) squaring  formula}: 
$$
(\Psi_k^L (v))^2 =
\Psi_k^L \Bigl(2v + \sum_{\gamma \atop |\gamma| \geq 2}
 \eps^\gamma (
\sum_{\lambda \in {\cal P}(\gamma) \atop \ell(\gamma) \geq 2}
 [\ldots[v_{\lambda^\ell(\lambda)},v_{\lambda^1}]\ldots, 
v_{\lambda^{\ell(\lambda)-1}}] )\Bigr).
$$
Note also that  our formulae for $J^L$ and $*^F$ ($F=L,R$)
are multilinear (resp.\
bi-multilinear) in the sense of Section MA.5. 
In the following we will give a conceptual explanation
of this fact by making the link with the theory of multilinear connections
from Part IV. (The reader who has not gone through that part may skip the
following theorem and its corollaries.)

\Theorem 24.8. Denote by
$\kappa \in \Sigma_k$ the permutation that reverses the order,
i.e.\ $\kappa(i)=k+1-i$, $i=1,\ldots,k$. (For $k=2$, this is 
the canonical flip, and for $k=3$, it is the transposition $(13)$.) Then 
left- and right trivialization $\Psi_k^L$ and $\Psi_k^R$ are multilinear
connections on $T^k G$ (in the sense of {\rm Chapter 16}),
 and, up to the permutation  $\kappa$, 
they coincide with the derived connections $D^{k-1} L^L$, resp.\
$D^{k-1} L^R$ defined in {\rm Chapter 17}.
 Put differently, the Dombrowski Splitting Maps 
$\Phi_k^F$ associated
to the multilinear connections $D^{k-1} L^F$,  $F=L,R$, are related
with the trivialization maps via
$$
\Phi_k^{F} = \Psi_k^{F,\kappa}
$$
where $\Psi_k^{F,\kappa}$ is as in {\rm Equation (24.5)}.

\Proof.
We proceed as in the proof of Theorem 23.5 where the case $k=2$ of
the claim has been proved.
Starting with the left-trivialization $\Psi_1^L:G \times \g \to TG$,
we have two possibilities to iterate the procedure:

\ssk
\item{(a)} we replace $G$ by $H:=TG$ and write the left-trivialization for 
$TG$: $TG \times T\g \to TTG$, compose again with $\Psi_1$ and obtain
$\Psi_2:G \times \g \times \g\times \g \to TTG$.
\item{(b)} we take the tangent map $T\Psi_1^L:TG \times T\g \to TTG$,
compose again with $\Psi_1$ and obtain another map $\tilde \Psi_2:
G \times \g \times \g\times \g \to TTG$.

\ssk \nin Both procedures yield essentially the same result. However,
the precise relation depends on the labelling of the copies of $TG$ we
use: as to convention (a) we agreed in Section 23.3 to write $H=T_{\eps_1} G$
so that $TTG=TH=T_{\eps_2}T_{\eps_1}G$.
Thus, in (b), we must take $T=T_{\eps_2}$, whence the copy $TG$ appearing
in the formula really is $T_{\eps_2} G$, but up to exchanging the order
of $\eps_1$ and $\eps_2$ both constructions are the same.
Now, when iterating the construction, in the $k$-th step,
we must take the {\it last} tangent
functor in (b) with respect to $\eps_1$, i.e.\ we must use the order of
scalar extensions $\eps_k,\ldots,\eps_1$ instead of the usual
order $\eps_1,\ldots,\eps_k$, and then constructions (a) and (b)
coincide.
According to Lemma 17.7, construction (b) corresponds to the construction
of the derived linear structure $D^{k-1} L^L$, where
due to the change of order the permutation
$\kappa$ appears in the final result. 
\qed


\Corollary 24.9.
Left- and right products $*^L$ and $*^R$
can be interpreted as higher covariant
differentials in the sense of {\rm Section 16.3}:
$*^F$ is the higher differential
of the group multiplication $m:G \times G \to G$
with respect to the linear structures derived from
 $L^F \times L^F$ and $L^F$, $F=L,R$, and inversion $J^F$
is the higher covariant
differential of group inversion.

\Proof. Since left, resp.\ right trivialization are the
linearizations associated to the connection $L$, 
the claim follows directly from the definition of the higher order
covariant differentials in Section 16.3.
\qed

\nin The corollary gives an {\it a priori} proof of the fact
that the left and right product formulae are bi-multilinear
(i.e.\ multilinear in the sense of MA.5 in both arguments) and that
$J^F$ is multilinear.
					     

\Theorem 24.10. {\rm (Restriction to jet groups.)}
Left and right trivialization of $T^k G$ induce, by restriction, 
two trivializations of the group $J^k G$; 
in other words, the restrictions 
$$
\Psi^F_k\vert_{E^{\Sigma_k}}:
 E^{\Sigma_k}=\oplus_{j=1}^k \delta^{(j)} \g \to (J^k G)_e, \quad
$$
$F=L,R$, are well-defined. 

\Proof.
For $k=2$, an elementary proof has been given in the preceding chapter
 (Theorem 23.7).
Let us give,
for $k=3$, an elementary proof in the same spirit:
we get from the explicit form of the left product
formula stated before Theorem 24.5, by letting
$v_1:=v_{001}=v_{010}=v_{100}$,
$v_2:=v_{011}=v_{101}=v_{110}$, $v_3:=v_{111}$, $v=\delta v_1 + 
\delta^{(2)} v_2 + \delta^{(3)} v_3$ (and similarly for $w$):
$$
\eqalign{
&\Psi_3^L(v)  \cdot \Psi_3^L(w) =
\eps^{001}(v_1+w_1) \cdot \eps^{010}(v_1+w_1) \cdot
 \eps^{011} (v_2+w_2 + [v_1,w_1]) \cdot
\cr
& \quad \quad \quad \quad \quad 
\eps^{100} (v_1+w_1) \cdot \eps^{101} (v_2+w_2+
[v_1,w_1]) \cdot 
\eps^{110}(v_2+w_2+[v_1,w_1]) \cdot 
\cr 
&   \quad \quad  \quad \quad \quad
\eps^{111}  (v_3+w_3+[v_2,w_1]+[v_2,w_1] +
 [v_1,w_2] + [[v_1,w_1], w_1])
\cr
& = \Psi_3^L\Bigl(v+w+ (\eps^{011}+\eps^{101}+
\eps^{110}) [v_1,w_1]  \cr
& \quad \quad  \quad \quad \quad
+\eps^{111} (2 [v_2,w_1] +
 [v_1,w_2] + [[v_1,w_1], w_1])\Bigr)
\cr
& = \Psi_3^L\Bigl(v+w+  \delta^{(2)} [v_1,w_1] + 
\delta^{(3)} (2 [v_2,w_1] +
 [v_1,w_2] + [[v_1,w_1], w_1])\Bigr)
\cr}
$$ 
which shows that the group structure of $(J^3 G)_e$ in the left trivialization
is given by
$$
\eqalign{
(\delta v_1 + & \delta^{(2)}   v_2+ \delta^{(3)} v_3) *^L 
 (\delta w_1 + \delta^{(2)} w_2+ \delta^{(3)} w_3) 
=\delta (v_1+w_1) + \cr
&  \delta^{(2)} (v_2+w_2 + [v_1,w_1]) 
+ \delta^{(3)} (w_3 + v_3 + [v_1,w_2] + 2 [v_2,w_1] +
[[v_1,w_1],w_1]) . \cr}
$$
Similarly, we have for the inversion:
$$
\eqalign{\Bigl(\Psi_3^L (\delta v_1 + 
\delta^{(2)} v_2 + \delta^{(3)} v_3)\Bigr)^{-1} & 
= \Psi_3^L(- \delta v_1 -
\delta^{(2)} v_2 + \delta^{(3)} (-v_3 +  [v_2,v_1]) ) \cr}.
$$
For general $k$, the combinatorial structure of the corresponding calculation
becomes rather involved. 
For this reason, we prefer to invoke Theorem 19.2 which in a sense
contains the relevant combinatorial arguments:
by Theorem 24.8, $\Psi^L_k$ and $\Psi_k^R$ are the Dombrowski splitting
maps associated to the connection $(D^{k-1} L^F)^\kappa$ ($F=L,R$).
According to Theorem 19.2, such connections are weakly
 symmetric (jet compatible), i.e.\
 they preserve the $\Sigma_k$-fixed subbundles.
\qed

\nin {\bf 24.11.} {\sl The left and right product formulae for
jets.} According to the preceding result,
the product formulae from Theorem 24.5 can be restricted to spaces
of $\Sigma_k$-invariants: there is a relation of the kind
$$
\sum_{i=1}^k \delta^{(i)} v_i *^L \sum_{j=1}^k \delta^{(j)} w_j =
\sum_{r=1}^k \delta^{(r)} z_r
$$
with $z_r \in \g$ depending on $v$ and $w$. For $k=3$ we have seen
in the preceding proof that
$$
\eqalign{
z_1 & = v_1 + w_1, \cr
z_2 & = v_2 + w_2 + [v_1,w_1], \cr
z_3 & = v_3 + w_3 + 2 [v_2,w_1] + [v_1,w_2] + [[v_1,w_1], w_1]. \cr}
$$
We will not dwell here on  the combinatorial details of the formula in
the general case, but content ourself with
the case of multiplying first order jets in $J^k G$: we claim that for
$v,w \in \g$,
$$
\delta v \, *^L \, \delta w  =
\delta (v+w) + \delta^{(2)} [v,w] + \ldots + \delta^{(k)}
[ [ v,w],\ldots w]. 
\eqno (24.8)
$$
In fact, $\delta v = \sum_{j=1}^k \eps_j v$, and
with the notation of Theorem 24.5, we have $v_\alpha = 0= w_\beta$ if
$\vert \alpha \vert , \vert \beta \vert > 1$.
Thus only partitions $\lambda$ with $\vert \lambda^i \vert \leq 1$
contribute to $z_\gamma$, i.e.\ we have
$$
z_\gamma = v_\gamma + w_\gamma + \sum_{\ell=2}^k \sum_{  
\lambda \in {\cal P}_\ell(\gamma)} [ [ v,w],\ldots w] .
$$
For $\vert \gamma \vert = 1$, there are no partitions of length
at least 2, whence $z_\gamma = v_\gamma + w_\gamma$.
This gives rise to the terms $\delta v + \delta w$.
For $\vert \gamma \vert = j > 1$, just one partition contributes to 
$z_\gamma$, giving rise to a term $\delta^{(j)} [[v,w],\ldots w]$. 
This proves (24.8). --
Note that,  since $[v,v]=0$,
$$
\Psi_k^L(\delta v) = \eps_1 v \cdot \ldots \cdot
\eps_k v =\eps_k v \cdot \ldots \cdot \eps_1 v = \Psi_k^R(\delta v),
$$
i.e.\ on $\delta \g$, left and right trivialization coincide;
for elements of the form $\delta^{(j)} v$ with $j >1$ this will no longer
be the case.

\ssk It would be interesting to investigate further the combinatorial
structures involved in the product and inversion formulae of $(J^k G)_e$.
In the case of arbitrary characteristic, these formulae are in a way the best
results one can get.
 In characteristic zero, one can use a ``better'' trivialization of
$(T^k G)_e$ given by the {\it exponential map}.

\msk \nin {\bf 24.12.} {\sl Remark on the synthetic approach.}
In the recent paper [Did06], Manon Didry has shown that, starting with
an arbitrary Lie algebra $\g$ over an arbitrary commutative ring $\K$,
the formula from Theorem 24.7 defines a group structure on the
$\K$-module $\oplus_{\alpha >0} \eps^\alpha \g$. 
Moreover, the symmetric group still acts by automorphisms, so that
the jet groups $\oplus_{j=1}^k  \delta^{(j)} \g$ are defined and form
a projective system. In some sense, these are the ``formal'' results
needed for a possible generalisation of Lie's Third Theorem. 
A surprising feature of
this ``synthetic'' approach is that it works not only for Lie algebras,
but more generally for {\it Leibniz algebras} (these algebras have
been introduced by J.-L.\ Loday; their product is not necessarily 
skew-symmetric, but the Jacobi-identity in its form given in Section
24.4 still holds). 
However, if the algebra is not a Lie algebra, then the symmetric
group in general no longer acts by automorphisms, so that the
groups we obtain are, in some sense, no longer ``integrable''.

\vfill \eject

\subheadline
{25. Exponential map and Campbell-Hausdorff formula for jet groups}

\msk
\nin {\bf 25.1.} {\sl $(T^k G)_e$ as a polynomial group.}
Recall from Chapter PG, Example PG.2 (2), that $(T^k G)_e$ satisfies
the properties of a {\it polynomial group}: there is chart (given by
left or right trivialization) such that, in this chart, group multiplication
is polynomial, and the iterated multiplication maps
$$
m^{(j)}: (T^k G)_e \times \ldots \times (T^k G)_e \to (T^k G)_e, \quad
(v_1,\ldots,v_j) \mapsto v_1 \cdots v_j
$$
are polynomial maps of degree bounded by $k$.

\Theorem 25.2.
Assume that $G$ is a Lie group over
 a non-discrete topological 
field $\K$ of characteristic zero. Then there exists a unique
map $\exp:(T^k \g)_0 \to (T^k G)_e$ such that:
\item{(1)} the representation of $\exp$ with respect to left trivialization,
$\exp_L:= (\Psi^L_k)^{-1} \circ \exp: (T^k \g)_0 \to (T^k \g)_0$, is a
polynomial map,
\item{(2)} $T_0 \exp = \id_{(T^k \g)_0}$,
\item{(3)} for all $n \in \Z$ and $v \in (T^k\g)_0$,
$\exp(n v)=(\exp(v))^n$.
\vskip 0mm
\nin The map $\exp$ is bijective and has a polynomial inverse $\log$.
Explicitly, the polynomials $\exp_L$ and $\log_L:=\log \circ \Psi^L_k$
are given by
$$
\eqalign{
(\exp_L(v))_\gamma & = v_\gamma + \sum_{j=2}^k {1 \over j!}
 \sum_{\lambda \in {\cal P}_j(\gamma)} [[v_{\lambda^j} ,v_{\lambda^{j-1}}],
\ldots , v_{\lambda^1} ]  \cr
(\log_L(v))_\gamma & = v_\gamma + 
\sum_{j=2}^k {(-1)^{j-1} \over j}
 \sum_{\lambda \in {\cal P}_j(\gamma)} [[v_{\lambda^j} ,v_{\lambda^{j-1}}],
\ldots , v_{\lambda^1} ] .\cr}
$$

\Proof.
Existence and uniqueness of $\log$ and $\exp$ follow from Theorem PG.6.
The explicit formulae are obtained by identifying the terms $\psi_p(x)$
and $\psi_{p,p}(x)$ appearing in the explicit formulae from Theorem PG.6 :
by definition, $\psi_1(v)=v$, and comparing the Left Inversion Formula
from Theorem 24.7 with Formula (PG.8) for the inversion,
$v^{-1} = \sum_j (-1)^j \psi_j(v)$, we get by induction, using that
$\psi_{p,m}=0$ for $m<p$ (i.e., $\psi_p$ vanishes at $0$ of order at least
$p$), 
$$
(\psi_j(v))_\gamma = (\psi_{j,j}(v))_\gamma =
 \sum_{\lambda \in {\cal P}_j(\gamma)} [[v_{\lambda^j} ,v_{\lambda^{j-1}}],
\ldots , v_{\lambda^1} ].
$$
Inserting this in the formulae from Theorem PG.6, we get the claim.
\qed

\msk \nin {\bf 25.3.}
It is clear that $\exp$ can also be expressed by a polynomial with respect 
to the right trivialization. It would be interesting to have also an
``absolute formula'' for $\exp$, i.e., a formula that does not refer to
any trivialization. For $k=2$, such a formula is given by
$$
\eqalign{
\exp (\eps_1 v_1 + \eps_2 v_2 + \eps_1 \eps_2 v_{12}) & 
=  \eps_1 v_1 \cdot \eps_2 v_2 
 \cdot \eps_1 \eps_2 (v_{12} + {1 \over 2} [v_2,v_1])  \cr
& = 
 \eps_2 {v_2 \over 2} \cdot \eps_1 v_1 
 \cdot \eps_1 \eps_2 v_{12}  \cdot 
 \eps_2 {v_2 \over 2}  \cr
& =
 \eps_2 {v_2 \over 2} \cdot \eps_1 {v_1 \over 2}
 \cdot \eps_1 \eps_2 v_{12} \cdot \eps_1 {v_1 \over 2} \cdot 
 \eps_2 {v_2 \over 2}, \cr}
$$
where the last expression appears to be fairly symmetric (it is a special
case of the formula for the canonical connection of a symmetric space,
see Chapter 26).  It is not at all obvious how to
 generalize this formula to the case of general $k$. 

\Theorem 25.4.
The exponential map commutes with Lie group automorphisms $\phi:(T^k G)_e
 \to (T^k G)_e$
in the sense that $ \phi \circ \exp  = \exp \circ T_0 \phi$,
 and it can be extended in a $G$-invariant way to a trivialization
of $T^k G$, again denoted by $\exp: G \times (\oplus_\alpha \eps^\alpha
\g) \to T^k G$. This trivialization is a totally symmetric multilinear
connection on $T^k G$. The restriction of $\exp$ to a map
$(J^k \g)_0 \to (J^kG)_e$ defines an exponential map of $(J^k G)_e$.

\Proof.
One easily checks that both $\exp$ and 
$ \phi \circ \exp  \circ (T_0 \phi)^{-1}$ satisfy Properties (1), (2), (3)
of an exponential map and hence agree by the uniqueness statement of
Theorem 25.2. Applying this to the inner automorphisms
$(T^k c_g)_e$ (where $c_g:G \to G$ is conjugation by $g \in G$), we
see that $\exp$ can be transported in a well-defined way to any fiber
of $T^k G$ over $G$, thus defining the trivialization map of $T^k G$.
On the fiber over $e$, it is multilinearly related to all chart structures
(by the explicit formula from Theorem 25.2, which clearly is multilinear),
and hence the same is true in any other fiber. Thus $\exp$ is a multilinear
connection. This connection is totally symmetric because the symmetric
group $\Sigma_k$ acts by automorphisms of $T^k G$, and we have just seen
that the exponential map is invariant under automorphisms.
Therefore the restriction of $\exp$ to $(J^k G)_e$ is well-defined, and
the uniqueness statement of Theorem PG.6 implies that it is then
``the'' exponential map of the polynomial group $(J^k G)_e$.
\qed

\Theorem 25.5.
With respect to the exponential connection from the preceding theorem,
the group structure of $(T^k G)_e$ and the one of $(J^k G)_e$ is
given by the Campbell-Hausdorff multiplication with respect to the
nilpotent Lie algebras $(T^k \g)_0$, resp.\ $(J^k \g)_0$.

\Proof. This follows from 
 Theorem PG.8 which states that
the group structure with respect to the exponential map is given by the
Campbell-Hausdorff formula.
\qed

\msk \nin Recall that  the classical Campbell-Hausdorff formula for Lie groups
is given by
$$
\eqalign{
& X*Y  :=\log(\exp(X) \cdot \exp(Y)) =
X + \cr 
& 
\sum_{k,m \geq 0 \atop p_i + q_i > 0}
{(-1)^k \over (k+1) (q_1+\ldots q_k +1)} \cdot 
 {(\ad X)^{p_1} (\ad Y)^{q_1} \cdots (\ad X)^{p_k} (\ad Y)^{q_k}
(\ad X)^m \over p_1 ! q_1 ! \cdots p_k ! q_k ! m !} Y,\cr}
$$
and the first four terms are:
$$
X*Y =
X +Y + {1 \over 2} [X,Y] + {1 \over 12}([X,[X,Y]] + [Y,[Y,X]]) + {1 \over 24}
[X,[Y,[X,Y]]] + \ldots
$$
(cf.,  e.g., [T04] for the fourth term).
In particular, we may apply this to the nilpotent Lie algebra
$$
(J^k \g)_0 =
\delta \g \oplus \delta^{(2)} \g \oplus \ldots \oplus \delta^{(k)} \g
\subset \g \otimes_\K J^k \K
$$
where $\g$ is an arbitrary $\K$-Lie algebra.
Elements are multiplied by the rule
$$
\delta^{(i)} \delta^{(j)} = \pmatrix{ i+j \cr i \cr}\delta^{(i+j)}, \quad \quad
[\delta^{(i)} X,\delta^{(j)} Y ] = \pmatrix{ i+j \cr i \cr}
\delta^{(i+j)} [X,Y].
$$
We re-obtain for $k=2$ the result from Theorem 23.7:
$$
(\delta v_1 + \delta^{(2)} v_2) * (\delta w_1 + \delta^{(2)} w_2) 
=\delta (v_1+w_1) + \delta^{(2)} (v_2+w_2 + [v_1,w_1]) ,
$$
and for $k=3$ we have:
$$
\eqalign{
(\delta  v_1 + \delta^{(2)}  v_2 & + \delta^{(3)}   v_3) * 
 (\delta w_1 + \delta^{(2)} w_2+ \delta^{(3)} w_3) 
= \cr
& \delta (v_1+w_1) + \cr
&   \delta^{(2)} (v_2+w_2 + [v_1,w_1])+  \cr
& 
 \delta^{(3)} 
(w_3 + v_3 + {3 \over 2}[v_1,w_2] + {3 \over 2}[v_2,w_1] +
{1 \over 2}([v_1,[v_1,w_1]]+[w_1,[w_1,v_1]])) . \cr}
$$
Note that this formula
makes sense whenever $2$ is invertible in $\K$.
Compare also with the formula from the proof of Theorem 24.10 which is
valid in arbitrary characteristic.  

\msk \nin {\bf 25.6.} {\sl Example:
Associative continuous inverse algebras.}
An associative topological
algebra $(V,\cdot)$ over $\K$ is called a {\it continuous
inverse algebra (CIA)} if  the
set $V^\times $ of invertible elements is open in $V$ and inversion
$i:V^\times \to V$ is continuous. We assume that $V$ has a unit element $\1$. 
It is easily proved that multiplication
and inversion are actually smooth, and hence $G=V^\times $ is a Lie group.
We consider $V$ as a global chart of $G$ and identify $(T^k G)_e$ with
$\1+(T^k V)_0$ where $(T^k V)_0 = \oplus_{\alpha >0} \eps^\alpha V$.
The group $(T^k G)$ is just the group of invertible elements in the
associative algebra $T^k V = V \otimes_\K T^k \K$.
The inverse of $\1+\eps v$ is $\1-\eps v$, and hence we get the group
commutator
$$
[\1+\eps_1 v,\1+\eps_2 w] =(\1+\eps_1 v)(\1+\eps_2 w)(\1-\eps_1 v) (\1
-\eps_2 v) 
=\1 + \eps_1 \eps_2 (vw-wv),
$$
which implies by the commutation relations (23.1) that the Lie algebra
of $G$ is $V$ with the usual commutator bracket. 

\Proposition 25.7. Let $G=V^\times$ be the group of invertible elements
in a continuous inverse algebra. Then
the group $(T^k G)_e = (T^k V)_0^\times$ is a polynomial group with respect
to the natural chart given by the ambient $\K$-module $(T^k V)_0$, and
the exponential map with respect to this polynomial group structure is given
by the usual exponential formula
$$
\exp(\sum_{\alpha > 0} \eps^\alpha v_\alpha) = \sum_{j=0}^k {1 \over j!}
(\sum_{\alpha > 0} \eps^\alpha v_\alpha)^j.
$$

\Proof. The group $(T^k G)_e$ is polynomial in the natural chart since
the algebra $(T^k V)_0$ is nilpotent of order $k$, and hence the degree
of the iterated product maps $m^{(j)}$ is bounded by $k$. 
Clearly, there is just one homogeneous component of $m^{(j)}$ and  we have
$\psi_j(x)=\psi_{j,j}(x)=x^j$. Inserting this in the formulae from Theorem
PG6, we get the usual formulae for the exponential map and the logarithm.
\qed

\msk
\nin {\bf 25.8.}
{\sl Comparing the exponential map with right trivialization.}
Expanding all terms in the formula for $\exp$ from the preceding
proposition, we get
$$
\eqalign{
& \exp (\sum_{\alpha > 0} \eps^\alpha v_\alpha) 
 = \sum_\alpha \eps^\alpha 
\sum_\ell {1 \over \ell!} \sum_{\Lambda\in {\cal P}_\ell(\alpha)}
\sum_{\sigma \in \Sigma_\ell} v_{\sigma.\Lambda} \cr
& =  \sum_\ell {1 \over \ell!} 
 \sum_{\Lambda\in {\rm Part}_\ell(I_k)}  \eps^{\uline \Lambda}
\sum_{\sigma \in \Sigma_\ell} v_{\sigma.\Lambda} 
\cr
& = \1 + \sum_\alpha  \eps^\alpha v_\alpha + 
\sum_{\Lambda : \atop
\ell(\Lambda)=2} \eps^{\uline \Lambda} \, \, {v_{\Lambda^1} v_{\Lambda^2} +
v_{\Lambda^2} v_{\Lambda^1} \over 2}  +
\sum_{\Lambda: \atop \ell(\Lambda)=3} 
\eps^{\uline \Lambda} \, \, {v_{\Lambda^1} v_{\Lambda^2} 
v_{\Lambda^3} + \ldots \over 6}  + \ldots
\cr}
\eqno (25.3)
$$
On the other hand, 
the trivialization map $\Psi_k^R$ is in the chart $T^k V$ described by
$$
\eqalign{
\Psi_k^{R} & (\sum_{\alpha >0} 
\eps^\alpha v_\alpha)  = \prod_\alpha^\downarrow
 (\1 + \eps^\alpha v_\alpha) =
\1 + \sum_{\alpha} \eps^\alpha \sum_{\lambda \in P(\alpha)} v_{\lambda^1}
\cdot \ldots \cdot v_{\lambda^\ell} \cr
& = e + \sum \eps^\alpha v_\alpha + 
\sum_{\ell(\Lambda)=2} \eps^{\uline \Lambda} v_{\Lambda^1} v_{\Lambda^2} + 
\sum_{\ell(\Lambda)=3} \eps^{\uline \Lambda} v_{\Lambda^1} v_{\Lambda^2} 
v_{\Lambda^3} + \ldots,
\cr}
\eqno (25.4)
$$
where partitions are considered as ordered partitions.
Comparing (25.3) and (25.4),  we see that
$\Psi_k^R$ is defined by making one particular
choice among the $\ell !$ terms in the sum in $\exp$
and forgetting division by $\ell !$, and in turn $\exp$ is obtained
by symmetrizing each term of $\Psi_k^R$. Here, ``symmetrization''
is to be understood with respect to the natural linear structure of
$T^k V$. In the case of a general Lie group, it seems much more difficult to
understand the passage from the right (or left) trivialization to the
exponential connection in terms of a suitable symmetrization procedure.

\msk \nin {\bf 25.9.} {\sl Relation with an exponential map of $G$.}
We say that a Lie group $G$ {\it has an exponential map} if the following
holds: for every $v \in \g$ there exists a unique Lie group homomorphism
$\gamma_v:\K \to G$ such that $\gamma_v'(0)=v$ and such that the map $\exp_G$
defined by
$$
\exp_G:\g \to G, \quad v \mapsto \gamma_v(1)
$$
is smooth. In general, a Lie group will not have an exponential map;
but if it does, then the higher differentials of $\exp$ are given
by the exponential map $\exp=:\exp^{(k)}$ from Theorem 25.2, i.e.
$$
(T^k \exp_G)_o = \exp^{(k)}: (T^k \g)_0 \to (T^k G)_e.
\eqno (25.5)
$$
In fact, it follows immediately from the definitions that $(T^k \exp_G)_0$
has the properties (1), (2), (3) from Theorem 25.2 and hence agrees with
$\exp^{(k)}$.
It is instructive to check (25.5) directly in the case of a CIA, using
 the Derivation Formula (7.19): let $p^j(x)=x^j$ be the $j$-th
power in $V$. Then $d^i p_j(0)=0$ if $i\not=j$ and 
$$
d^j p_j(0)(v_1,\ldots,v_j)=\sum_{\sigma \in \Sigma_j} 
v_{\sigma(1)} \cdots v_{\sigma(j)}
$$
is the total polarization of $p_j$. Thus 
$$
T^k \exp(0) (v)  = \sum_\alpha \eps^\alpha \sum_{\ell=1}^{\vert \alpha 
\vert} \sum_{\Lambda\in {\cal P}_\ell(\alpha)}
 {1 \over \ell !} d^\ell p_\ell(0)
(v_\Lambda)
$$
 agrees with (25.3).

\vfill \eject 

\subheadline
{26. The canonical connection of a symmetric space}
                                                              
\nin
{\bf 26.1.} {\sl Abelian symmetric spaces.}
We return to the study of symmetric spaces and start by some general
(in fact, purely algebraic)
remarks on {\it abelian} symmetric (or reflection) spaces (Chapter 5).
In the real, finite-dimensional case, connected abelian symmetric spaces
are always products of tori and vector spaces (cf.\ [Lo69]),
 hence are of group type (cf.\ Section 5.11). In our general context,
it is certainly not true that all abelian symmetric spaces are of group
type, but we will single out a property (the ``unique square root property'')
which permits to define a compatible group structure.
--  
Recall first that a symmetric (or reflection)
 space is called {\it abelian} if the group $G(M)$
generated by all $\sigma_x \circ \sigma_y$, $x,y \in M$, is abelian.
In particular, if $o \in M$ is a base point, the operators $Q(x)=\sigma_x \circ
\sigma_o$, $x \in M$, commute among each other. 
For instance, if $(V,\cdot)$ is a commutative group, then
$$
\mu(u,v)=uv^{-1}u=u^2 v^{-1}
\eqno (26.1)
$$
 defines on $V$ the structure of an abelian symmetric space,
and $Q(v)x=v^2 x$ is translation by $v^2$.
We say that a group $M$ (resp.\ a symmetric space $M$ with base point $o$)
{\it  admits unique square roots} if,  for all
$v \in M$ there is a unique element $\sqrt v \in M$ such that
$(\sqrt v)^2 = v$ (resp.\  $\mu(\sqrt v,o)=v$, which is the same
as $Q(\sqrt v)o=v$). 

\Proposition 26.2.
There is a bijection between
\item{(1)} abelian groups $(V,\cdot)$ admitting unique square roots, and
\item{(2)} abelian reflection spaces with base point $(V,\mu,o)$ admitting
unique square roots.
\vskip 0mm
\nin The bijection is given by associating to $(V,\cdot)$ the symmetric space
structure {\rm (26.1)}, and by associating to $(V,\mu,o)$ the group structure
given by
$$
u \cdot v =  \mu( \sqrt u, \mu(e,v)) = Q(\sqrt u) v = Q(\sqrt u) 
Q(\sqrt v) o .
\eqno (26.2)
$$

\Proof. In this proof and on the following pages,
 the ``fundamental formula'' (5.6) will be frequently used without further
comments. -- 
We have already seen that (26.1) is an (abelian) reflection space
 structure; the square roots are unique since the squaring operation is
the same as the one for the group. Let us prove that (26.2) defines an
abelian group law on $V$. We have $u \cdot v=v \cdot u$ since the operators
$Q(\sqrt u)$ and $Q(\sqrt v)$ commute by our assumption. Clearly, $v\cdot o
 =o \cdot v=v$. In order to establish associativity, note first that
$$
\sqrt{Q(x)y} = Q(\sqrt x)\sqrt y.
\eqno (26.3)
$$
In fact, we have $Q(\sqrt x)^2 = Q(Q(\sqrt x)o)=Q(x)$ and hence
$$
Q(Q(\sqrt x)\sqrt y)o=Q(\sqrt x)^2 Q(\sqrt y)o=Q(x)y,
$$
 whence (26.3). Using this, associativity follows:
$$
\eqalign{
u \cdot (v \cdot w) &= Q(\sqrt u) Q(\sqrt v) w, \cr
(u \cdot v) \cdot w & =Q(\sqrt{Q(\sqrt u) v}) w =
Q(Q(\sqrt{\sqrt u}) \sqrt v) w \cr
& = Q(\sqrt{\sqrt u})^2 Q(\sqrt v) w = Q(\sqrt u) Q(\sqrt v) w. \cr}
$$
Similarly, we see that the inverse of $x$ is $x^{-1}=Q(x)^{-1}x$,
and thus $(V,\cdot,o)$ is an abelian group; it admits unique square roots
since the squaring operation is the same as the one of $(V,\mu,o)$.
Finally, it is clear that both constructions are inverse to each other.
\qed

\nin
In the sequel, abelian groups will mainly be written additively, so that
$v^2$ corresponds to $2v$ and $\sqrt v$ corresponds to $
{1 \over 2} v$. Then (26.1) and (26.2) read
$$
\mu(u,v)=2u-v, \quad \quad \quad  u + v = \mu({1 \over 2} u, \mu(0,v)) =
Q({u \over 2})v = Q({u \over 2})Q({v \over 2})o.
\eqno (26.4)
$$
Now let $(M,\mu)$ be an arbitrary symmetric space over $\K$.
Recall that the higher tangent bundles $(T^k M,T^k \mu)$ are again
symmetric spaces. The purpose of this chapter is to prove the following
theorem:

\Theorem 26.3.
Assume $(M,\mu)$ is a symmetric space over $\K$.
There exists a linear connection $L$ on $TM$ which is uniquely determined
by one of the following two equivalent properties:
\item{(a)}
$L$ is invariant under all symmetries $\sigma_x$, $x \in M$.
\item{(b)} For all $x \in M$,
the symmetric space structure on $(TTM)_x$ is the canonical flat 
structure associated to the $\K$-module $((TTM)_x,L_x)$:
for all $u,v \in (TTM)_x$, $TT\mu(u,v)=2u-v$.
\ssk
\nin In particular, the fibers of $TTM$ over $M$ are abelian.
The torsion of $L$ vanishes.

\Proof.
Let us first prove uniqueness. Assume  two linear connections
$L,L'$ with Property (a) given. 
Then $A:=L-L'$ is a tensor field of type $(2,1)$ that is again
invariant under all symmetries. From $A_x \sigma_x=-\id_{T_xM}$ it
follows then that $-A_x(u,v)=A_x(-u,-v)=A_x(u,v)$ and
$A_x(u,v)=0$ for all $x \in M$ and $u,v \in T_x M$.
Hence $A=0$. The same argument also shows that the torsion 
of $L$ must vanish. 

\ssk
Now assume $L$ is a linear structure on $TM$ satisfying (b). 
Then multiplication by $2$ in the fiber $(TTM)_x$ is given by
$2u=TT\mu(u,0_x) = Q(u).0_x$; 
it is a bijection whose inverse is denoted by $2^{-1}$.
Thus $(TTM)_x$ has unique square roots, and therefore
addition in the fiber $(TTM)_x$ is necessarily given by
(26.4): $u + v = TT\mu(2^{-1} u, TT\mu(0,v)) =
Q({u \over 2}) Q({v \over 2}).0_x$ (here $Q$ is taken with respect to
the base point $0_x$).
Thus the addition is uniquely determined by the symmetric space
structure of $(TTM)_x$. But then also the bijection
$$
\eqalign{
\eps_1 T_xM \times \eps_2 T_xM \times \eps_1\eps_2 T_xM &  \to (TTM)_x, \cr
 (\eps_1 v_1 ,\eps_2 v_2,\eps_1\eps_2 v_{12}) & \mapsto
\eps_1 v_1 +\eps_2 v_2+\eps_1\eps_2 v_{12} \cr}
$$
depends only on the symmetric space structure.
 Since the axes are $\K$-modules
in a canonical way, it follows that also the action of scalars on
$(TTM)_x$ is uniquely determined. Summing up, there is at most one
linear structure on $TTM$ satisfying (b), and moreover it follows
 that such a linear structure will
be invariant under all automorphisms of $M$ since it is defined by
the symmetric space structure alone (more generally,
it will depend functorially on $(M,\mu)$).

Existence of a linear structure $L$ satisfying (b)
can be proved in several ways. In the following we give a proof by using
charts, which shows also that the linear structure is a connection.
The proof relies on the following

\Lemma 26.4.
If $U$ is a chart domain of $M$ around $x$ with model space $V$, then for all
$z,w \in (TTM)_x =V \times V \times V$,
$$
TT\mu(z,w)=2_b z -_b w,
$$
(notion for bilinearly related structures being as in {\rm Section BA.1})
 where 
$$
b:=b_x:= {1 \over 2} d^2 \sigma_x(x): V \times V \to V
$$
is given by the ordinary second differential of the symmetry
$\sigma_x$ at $x$ in the chart $U \subset V$.

\Proof. (Cf. [Pos01, p.\ 64--66] for a coordinate version of the
following proof in the  finite-dimensional real case.)
First of all, if $V_1,V_2,V_3$ are $\K$-modules and
$b\in \Bil(V_1\times V_2,V_3)$, then the symmetric structure
with respect to the linear structure $L^b$ is given by the
``Barycenter Formula'' (BA.4), with $r=-1$:
$$
2_b(u_1,u_2,u_3) -_b (v_1,v_2,v_3)=
(2u_1 -v_1,2u_2-v_2,2u_3-v_3+ 2 
 b(u_1-v_1,u_2-v_2)).
\eqno (26.5)
$$
Now choose a chart $U$ around $o$. From the general
expression of the second tangent map of $f$,
$$
TTf(x,u_1,u_2,u_3)=(f(x),df(x)u_1,df(x)u_2,df(x)u_3+d^2 f(x)(u_1,u_2)),
$$
for $f=\mu$, we get the chart formula for $TT\mu:TT(U \times U) \to TTM$
 at the point $(x,x)$:
$$
\eqalign{
 TT\mu& ((x,x),(u_1,v_1),(u_2,v_2),(u_3,v_3))   \cr
& =\Bigl(\mu(x,x),d\mu(x,x)(u_1,v_1),d\mu(x,x)(u_2,v_2), \cr
& \quad \quad \quad \quad \quad \quad \quad d\mu(x,x)(u_3,v_3) +
d^2m(x,x)((u_1,v_1),(u_2,v_2)) \Bigr)  \cr
&= \Bigl(x,2v_1-u_1,2v_2-u_2,  2v_3-u_3 + 
d^2\mu(x,x)((u_1,v_1),(u_2,v_2))\Bigr) \cr}
$$
Comparing with (26.5), we see that
the claim will be proved if we can show that
$$
d^2\mu(x,x)((u_1,v_1),(u_2,v_2)) = d^2 \sigma_x (x) (u_1-v_1,u_2 - v_2).
\eqno (26.6)
$$
In order to prove (26.6), we
assume that, in our chart, $x$ corresponds to the origin of $V$
and  write  Taylor expansions of order 2 at the origin for $\sigma_x=
\sigma_0$ and for $\mu$,
$$
\eqalign{
\sigma_0(tv) & =-tv + {t^2 \over 2} d^2 \sigma_0(0) (v,v) + t^2 O(t), \cr
\mu(tu,tv) & = t(2u-v)+{t^2 \over 2} d^2 \mu(0,0)((u,v),(u,v)) +t^2 O(t) \cr
& = t(2u-v)+ t^2 c(u,v) +t^2 O(t) \cr}
$$
with $c(u,v):={1 \over 2} d^2\mu(0,0)((u,v),(u,v))$
(this is not bilinear as a function of $(u,v)$ but is 
a quadratic form in $(u,v)$).
The defining identity (M2) of a symmetric space, $p=\mu(x,\mu(x,p))$, implies
$$
\eqalign{
tv & = \mu(0,\sigma_0(tv))=
\mu(0,t(-v+{t \over 2} d^2 \sigma_0(0) (v,v) + t  O(t))) \cr
& =
t(v-{t \over 2} d^2 \sigma_0(0) (v,v) - t O(t))
+t^2  c((0,-v+O(t))) + t^2 O(t) \cr
&=tv + t^2(c(0,-v) - {1 \over 2} d^2\sigma_0(0)(v,v)) + t^2 O(t).\cr}
$$
Comparing quadratic terms and letting $t=0$ we get
$$
d^2 \sigma_0(0)(v,v) = 2
c((0,-v)) = d^2 \mu(0,0) \, ((0,-v),(0,-v)).
$$
This is a special case of (26.6). The general formula now follows:
from the identity $\mu(p,p)=p$ we get
$$
\partial_{v,v} \mu(x,x)={\mu(x+tv,x+tv)-\mu(x,x) \over t}\vert_{t=0} =
{tv \over t}|_{t=0} = v
$$
and hence
$$
d^2 \mu(0,0)((u,w),(v,v))=
\partial_{(u,w)} \partial_{(v,v)} \mu(0,0)=
\partial_{(u,w)} v = 0.
$$
We write $(u,v)=(v,v)+(u-v,0)$, and so on, and get
$$
\eqalign{
d^2\mu(0,0)((&u_1,v_1),(u_2,v_2))  \cr
& =
d^2\mu(0,0)((v_1,v_1)+(u_1-v_1,0),(v_2,v_2)+(u_2-v_2,0)) \cr
& =
d^2\mu(0,0)((u_1-v_1,0),(u_2-v_2,0)) =
d^2 \sigma_0(0)(u_1-v_1,u_2-v_2) \cr}
$$
whence (26.6), and the claim is proved. \qed

\nin
Returning to the proof of the existence part of the theorem, using a chart
$U$ as in the lemma, we define a
 linear structure on $(TTF)_x = T_xM \times T_xM\times T_xM$
by $L_x:=L^{b_x}$ with $b_x={1 \over 2} d^2 \sigma_x(x)$. By the lemma,
this linear structure satisfies Property (b) from the theorem.
We have already seen that there is at most one linear structure on $(TTM)_x$
satisfying (b), and therefore $L_x$ does not depend on the choice of the
chart. Thus $(L_x)_{x \in M}$ is a well-defined linear structure on $TTM$,
and it is a connection since, again by the lemma, it is bilinearly related
to all chart structures via bilinear maps $b_x$ depending smoothly on $x$.
We have already remarked that  this connection
 is invariant under all automorphisms and hence satisfies also
Property (a).
\qed

\Corollary 26.5.
If $(\phi_i,U_i)$ is a chart of $M$ around $x$ such that
$\phi_i(x)=0$ and $\sigma_x(y)=-y$ for all $y \in U$, then
in this chart the symmetric space structure on $(TTM)_x$
is simply the usual flat structure of the vector space
$V \times V \times V$.

\Proof. We have
$d^2 (-\id)(x)=0$ and hence $b_x=0$.
\qed

For example, this may applied to the exponential map of
a symmetric space, if it exists (cf.\ Chapter 31, or see [BeNe04]), or
to Jordan coordinates of a symmetric space with twist ([Be00]).
(In fact, it it is well-known that
in the real finite dimensional case and for any
torsionfree connection one can find a chart such that
the Christoffel symbols at a given point $x$ vanish and
hence the linear structure of $(TTM)_x$ is the canonical flat
one from the chart.
 This generalizes to all contexts in which  one has an inverse function
theorem at disposition.)

\msk \nin {\bf 26.6.} {\sl Dombrowski splitting and the associated
spray.}
The Dombrowski splitting (Theorem 10.5) for the canonical connection of
a symmetric space is given by
$$
\eqalign{
 \Phi_1:   \eps_1  TM \times_M \eps_2  TM \times_M \eps_1\eps_2 TM
 & \to TTM, \cr
  (x; \eps_1 v_1, \eps_2 v_2, \eps_1 \eps_2 v_{12})  \mapsto
 \eps_1 v_1  & + \eps_2 v_2  + \eps_1 \eps_2 v_{12} \cr
& = 
 \eps_2 v_2  + (\eps_1 v_1  + \eps_1 \eps_2 v_{12}) \cr
& 
= Q_x(\eps_2 {v_2 \over 2}) Q_x(\eps_1 {v_1 \over 2})\eps_1 \eps_2 v_{12} 
\cr
& = Q_x(\eps_2 {v_2 \over 2}) Q_x(\eps_1 {v_1 \over 2}) Q_x(\eps_1 \eps_2
{v_{12} \over 2}).0_x,
\cr}
\eqno (26.7)
$$ 
where $Q_x(y)=s_y s_x$ is the quadratic map in $TTM$ with respect to the
base point $x=0_x$ (origin in $(TTM)_x$). Restricting to invariants under the
canonical flip $\kappa$, we get a bijection of vector bundles over $M$
$$
\eqalign{
\delta TM \times_M \delta^{(2)} TM  & \to J^2 M, \cr
(x;\delta u,
\delta^{(2)} w) &\mapsto u+\kappa(u)+w =
Q_x(\eps_2 {u \over 2}) Q_x(\eps_1 {u \over 2}) Q_x(\eps_1 \eps_2
{w \over 2}).0_x \cr}
$$
which for $w=0$ gives the spray
$$
TM \to J^2 M, \quad (x;u) \mapsto u+\kappa(u) =
Q_x(\eps_2 {u \over 2}) \eps_1 u =
Q_x(\eps_1 {u \over 2}) \eps_2 u.
$$
The last expression can be rewritten, with $z:M \to TM$ being
the zero section, 
$$
u+\kappa(u)= TT\mu({u \over 2}, - \kappa(u))
= T \sigma_{u\over 2} (-\kappa(u)) =
 T \sigma_{u\over 2} Tz(-u)
$$
since $TT\mu(v,\cdot)=T\sigma_v$ and $\kappa(v)=Tz(v)$.
Compare with [Ne02a, Th.\ 3.4].

\Corollary 26.7.
The covariant derivative associated to the canonical connection
of a symmetric space $M$ can be expressed in the following way:
for $x \in M$ and $v \in T_x M$ let $l_x(v):=\tilde v =
{1 \over 2}\circ Q(v) \circ z\in \X(M)$ be the canonical
vector field extending $v$ (see proof of
{\rm Proposition 5.9.}).
Then we have, for all $X,Y \in \X(M)$ and $x \in M$,
$$
(\nabla_X Y)(x) = [l_x(X(x)),Y](x).
$$

\Proof. Chose a chart such that $x$ corresponds to the origin $0$
and let $v:=X(0)$. Then the value of
$$
T \tilde v = {1 \over 2} TQ(v) \circ Tz = {1 \over 2} T\sigma_v
\circ T\sigma_0 \circ Tz
$$
at the origin is given, in the chart, by
$$
d \tilde v(0) w = -{1 \over 2} d^2 \mu(0,0)((v,0),(w,0)) =
{1 \over 2} d^2 \sigma_0(0)(v,w)
$$
where for the last equality we used (26.6). Hence
$$
\eqalign{
[l_0(X(0)),Y](0) & =
[\tilde v,Y](0)= dY(0) v - d\tilde v(0) Y(0) \cr
& =
 dY(0) v + {1\over 2} d^2 \sigma_0(0) (v,Y(0)) \cr
 & =  dY(0) v + b_x(v,Y(0)) \cr
& = (\nabla_X Y)(0) \cr}
$$
according to  Lemma 12.2.
\qed

\nin In [Be00], the canonical covariant derivative of a 
(real finite-dimensional)
symmetric space $G/H$ has been defined by the formula from the preceding 
lemma.

\vfill \eject
\subheadline{27. The higher order tangent structure of symmetric spaces}

\nin {\bf 27.1.}
In this chapter we study the structure of tangent bundles $T^k M$
of a symmetric space for $k \geq 3$.
 In order to simplify notation, if $v \in (T^k M)_o$, we will write $Q(v)$
for the quadratic map  $Q_o(v):T^k M \to T^k M$.
If the base point $o \in M$ is fixed and if we use the notation
$$
Q:=Q^{(M)}: M\times M \to M, \quad (x,y)\mapsto Q(x)y=\sigma_x\sigma_o(y)
$$
for the quadratic map of $M$ w.r.t.\ $o$, then the quadratic map of
$T^k M$ w.r.t.\ $0_o$ is given by $Q^{(T^k M)} = T^k (Q^{(M)})$.

\Theorem 27.2.
Assume $(M,\mu)$ is a symmetric space over $\K$ with base point $o \in M$.
For elements from the axes of $(T^3 M)_o$, the following ``fundamental
commutation relations'' hold: for all $u,v,w \in T_oM$ and $i,j\in
\{ 1,2,3\}$,
$$
\eqalign{
[Q(\eps_i u),Q(\eps_i v)]\eps_j w &= \eps_j w, \cr
[Q(\eps_i u),Q(\eps_j v)]\eps_j w &= \eps_j w, \cr
[Q(\eps_1 u),Q(\eps_2 v)]\eps_3 w &= \eps_3 w +\eps_1 \eps_2 \eps_3
[u,v,w],             \cr}
$$
where the commutator is taken in the group $G(T^3M)$
and $[u,v,w]$ is the Lie triple product of $(M,o)$. 

\Proof.
The first two relations simply take account of the fact that the
fibers of $TTM=T_{\eps_i} T_{\eps_j}M$
 are abelian symmetric spaces (Theorem 26.3), i.e.\ the quadratic operators
 commute on the fiber.

\ssk
In order to prove the third relation, recall from Theorem 6.8 that
the group $\Diff_{T\K}(TM)$ is a semidirect product of $\Diff_\K(M)$ with
 the vector group
$\X(M)$.
It follows that the symmetric space automorphism group $\Aut_{T\K}(TM)$
can  be written as a semidirect product of  $\Aut_\K(M)$ and a vector
group which is precisely the subspace $\g=\Der(M) \subset \X(M)$ of
{\it derivations of $M$}, introduced in Section 5.7. 
(Indeed, the $T\K$-smooth extension $\tilde X$
of a derivation $X:M \to TM$ 
is again a homomorphism of symmetric spaces, and since it it is invertible,
it is an automorphism. We call automorphisms of the form $\tilde X:TM \to 
TM$ {\it infinitesimal automorphisms}. Compare with the closely related
definition of infinitesimal automorphisms in [Lo69].)

\Lemma 27.3.
For all $v \in T_o M$, the map $Q(v):TM \to TM$ is an infinitesimal
automorphism. More precisely, $Q({v \over 2})$ is the infinitesimal
automorphism induced by  the derivation $\tilde v$ introduced in 
{\rm Equation (5.9)}. In other words, we have
 $Q({v \over 2})=\tilde {\eps \tilde v}$.

\Proof. The automorphism $Q(v)$ of $TM$ is smooth over $T\K$ since
it is defined in terms of $T\mu$, and it preserves fibers since for all
$u \in TM$,
$$
\pi (Q(v)u) = \pi(T\mu(v,T\mu(0_o,u))) =
\mu(\pi(v),\mu(o,\pi(u))) = \mu(o,\mu(o,\pi(u)))=\pi(u).
$$ 
Thus Theorem 6.8 implies that $Q(v)$ is an infinitesimal automorphism.
It is induced by the vector field
$$
p \mapsto Q(v)0_p=2 \tilde v(p).
\qeddis 

\nin
Using the lemma, we can give another proof (not depending on the preceding 
chapter) of the first two relations from the theorem:
on $TTM$, the quadratic operator is given by taking the tangent map of
the infinitesimal automorphism from the preceding lemma,
$$
Q(\eps_1 {v \over 2}) = \tilde {\eps_1 \tilde v},
$$
and similarly for $\eps_1$ and $\eps_2$ exchanged.  Then
Theorem 14.4 on the Lie bracket implies
$$
[Q(\eps_1 {u \over 2}), Q(\eps_2 {v \over 2})] (w)=
[\tilde {\eps_1 \tilde u},\tilde {\eps_2 \tilde v}](w)
=\tilde{\eps_1 \eps_2 [\tilde u,\tilde v]}(o) + w = 0_o + w = w
$$
since $[\tilde u,\tilde v]$ vanishes at $o$ (Prop.\ 5.9).
In a similar way, again with Theorem 14.4, we get the third relation:
the triple bracket $[[\tilde u, \tilde v],\tilde w]$ is given by
$$
\eps_1 \eps_2\eps_3 [[\tilde u, \tilde v],\tilde w]=
[[Q({\epsilon_1 u \over 2}), Q({\epsilon_2 v \over 2})],
Q({\eps_3 w\over 2})] |_{\eps_1 \eps_2 \eps_3 TM}.
$$
Evaluating at the base point, using that the vector field
$[\tilde u,\tilde v]$ vanishes at $o$ (Prop.\ 5.9) and that
$Q({\eps_3 w \over 2}).0 = \eps_3 w$, we get
$$
\eqalign{
\eps_1 \eps_2 \eps_3 [u,v,w] & =
\eps_1 \eps_2 \eps_3 [[\tilde u, \tilde v],\tilde w].0 \cr
& =
\eps_1 \eps_2 [\tilde u,\tilde v] \cdot \eps_3 w - \eps_3 \tilde w(0) =
[Q({\epsilon_1 u \over 2}), Q({\epsilon_2 v \over 2})].\eps_3 w -\eps_3 w, \cr}
$$
which had to be shown.
\qed

\nin
Another way of stating the preceding relations goes as follows: let us
 introduce (with respect to a fixed base  point $o \in M$) the map
$$
\eqalign{
\rho:=\rho_M&: M\times M\times M  \to M, \cr
 (x,y,z) & \mapsto 
[Q(x),Q(y)]z = Q(x)Q(y)Q(x)^{-1}Q(y)^{-1}z =
\sigma_x \sigma_o \sigma_y \sigma_x \sigma_o \sigma_y. z \cr}
$$
%
%
Then, writing $w=\sum_\alpha \eps^\alpha w_\alpha$, etc., we have 
$$
(T^3 \rho)_{o,o,o}(u,v,w)=w + \eps^{111} [u_{100},v_{010},w_{001}].
\eqno (27.1)
$$
This follows by noting that $T^3 \rho_M = \rho_{T^3M}$.

\msk \nin {\bf 27.4.} {\sl Fundamental commutation relations for elements
in the axes of a symmetric space.}
For elements
 $\eps^\alpha u_\alpha$,
 $\eps^\beta v_\beta$, $\eps^\gamma w_\gamma$ of the axes
of $(T^k M)_o$ with $k \geq 3$
we have the following ``commutation relations":
$$
[Q(\eps^\alpha {u_\alpha \over 2}),Q(\eps^\beta {v_\beta \over 2})] \eps^\gamma
w_\gamma = \Big\{ \matrix{
\eps^\gamma w_\gamma + \eps^{\alpha +  \beta + \gamma}
 [u_\alpha,v_\beta,w_\gamma] & 
 {\rm if} \, \alpha \perp \beta, \beta \perp
\gamma , \alpha \perp \gamma,
\cr \eps^\gamma w_\gamma & {\rm else,} \cr}
\eqno (27.2)
$$
where the bracket on the left hand side denotes the
commutator in the group $\Diff(T^kM)$.
In fact,
the general case is reduced to the case $k=3$; if $\alpha,\beta,
\gamma$ are not disjoint, then we may further reduce to the case $k=2$
in which the fibers are abelian symmetric spaces (Theorem 26.3), whence
the claim in this case; in the other case we apply Theorem 27.2.

\msk \nin {\bf 27.5.} {\sl The derived multilinear connections on  $T^k M$.}
Starting with the canonical connection $L$ on $TM$, we may define
a sequence of derived multilinear connections on $T^k M$ (Chapter 17).
We will give an explicit formula for the derived linear structure on $T^3M$:
let $L_1:=L$ be the canonical connection of the symmetric space $M$
and $L_k:=D^{k-1} L$ be the sequence of derived linear structures,
and $\Phi_k:A^{k+1}M \to T^{k+1}M$ be the corresponding linearization map.
The map $\Phi_1:A^2 M \to T^2 M$ is given by (26.7).
We are going to calculate $\Phi_2:A^3 M \to T^3 M$.
First of all, we have to calculate the tangent map
$T\Phi_1$. It is a bundle map over the base space
$TM=T_{\eps_3}M$; a typical point in the base is $\eps_3 v_3 \in T_x M$.
In the formula for $T\Phi_1$, the $Q$-operators in $T^2 M$ are replaced by
the corresponding ones in $T^3M$, but taken with respect to
the base point $\eps_3 v_3 \in T_x M$. 
The factor ${1 \over 2}$ in (26.7) really stands for the map
$$
({1 \over 2})_{TTM}: TTM \to TTM, \quad v \mapsto {v \over 2}.
\eqno (27.3)
$$
When taking its derivative with respect to $\eps_3$, the base point $\eps_3
v_3$ is fixed. In the following, the sign $+$ without index 
stands for $+_x$ (addition in fiber of $TTM$ over $x$), 
and $Q$ without index stands for $Q_{0_x}$;
the index $\eps_3 v_3$ labels the corresponding operations
in the fiber over $\eps_3 v_3$.
In the following calculation, we pass from
one base point to another via the ``isotopy formula''
$$
Q_x(y)=s_y s_o s_o s_x = Q_o(y) Q_o(x)^{-1}.
\eqno (27.4)
$$
 Then we have
$$
\eqalign{
& \Phi_2 (x;\eps_1 v_1,\eps_2 v_2,\eps_3 v_3,\eps_1 \eps_2 v_{12},\eps_2
\eps_3 v_{13},\eps_2\eps_3 v_{23})  \cr
&= (\eps_2 v_2 + \eps_3 v_3 + \eps_2 \eps_3 v_{23}) +_{\eps_3 v_3}
 (\eps_1 v_1 + \eps_3 v_3 + \eps_1 \eps_3 v_{13}) \cr
& \quad \quad +_{\eps_3 v_3}
 (\eps_1 \eps_2 v_{12} + \eps_3 v_3 + \eps_1 \eps_2 \eps_3  v_{123}) 
\cr
&= Q_{\eps_3 v_3}(\eps_2 {v_2 \over 2} + \eps_3 v_3 + \eps_2 \eps_3 {v_{23}
\over 2}) \circ 
 Q_{\eps_3 v_3}(\eps_1 {v_1 \over 2} + \eps_3 v_3 + \eps_1 \eps_3 {v_{13}
\over 2}) \circ
\cr
& \quad \quad 
 Q_{\eps_3 v_3}(\eps_1\eps_2 {v_{12} \over 2} + \eps_3 v_3 + \eps_1 \eps_2
\eps_3 {v_{123} \over 2}).0_{\eps_3 v_3} \cr
& =
Q_{\eps_3 v_3}\Bigl(Q(\eps_2 {v_2 \over 4})  Q(\eps_3 {v_3 \over 2})
  Q(\eps_2 \eps_3 {v_{23} \over 4}).0_x\Bigr) \circ 
 \cr
& \quad \quad  Q_{\eps_3 v_3}\Bigl(Q(\eps_1 {v_1 \over 4})  Q(\eps_3 {v_3 \over 2})
  Q(\eps_1 \eps_3 {v_{13} \over 4}).0_x\Bigr) \circ \cr
& \quad \quad 
Q_{\eps_3 v_3}\Bigr(Q(\eps_1\eps_2 {v_{12} \over 4})  Q(\eps_3 {v_3 \over 2})
  Q(\eps_1 \eps_2 \eps_3 {v_{123} \over 4}).0_x\Bigr). 0_{\eps_3 v_3} \cr
&=
Q(\eps_3{v_3 \over 2}) Q\Bigl(Q(\eps_2 {v_2 \over 4})  
  Q(\eps_2 \eps_3 {v_{23} \over 4}).0_x\Bigl) \, \,
Q\Bigl(Q(\eps_1 {v_1 \over 4})  
  Q(\eps_1 \eps_3 {v_{13} \over 4}).0_x\Bigr) \cr
& \quad \quad
Q\Bigl(Q(\eps_1\eps_2 {v_{12} \over 4}) 
  Q(\eps_1 \eps_2 \eps_3 {v_{123} \over 4}).0_x\Bigr). 0_x \cr
&=
Q(\eps_3{v_3 \over 2}) Q(\eps_2 {v_2 \over 2})  
  Q(\eps_2 \eps_3 {v_{23} \over 2}) 
Q(\eps_1 {v_1 \over 2})  
  Q(\eps_1 \eps_3 {v_{13} \over 2})
Q(\eps_1\eps_2 {v_{12} \over 2}) 
  Q(\eps_1 \eps_2 \eps_3 {v_{123} \over 2}). 0_x \cr}
$$
where the simplifications in the last lines are due to the
``fundamental formula''  and the fact
that $Q$-operators having an $\eps_i$  for argument in common commute.
The formula for $\Phi_k$ with general $k$ is obtained in a similar
way; we state the result without proof: if,
 as in Theorem 24.8, $\kappa$ denotes the order-reversing permutation
$\kappa (j)=k+1-j$, then
$$
\Phi_k(x;(v_\alpha)_{\alpha >0}) =
(\prod_{\alpha >0}^\uparrow  Q(\eps^{\kappa \cdot \alpha} 
{v_{\kappa \cdot \alpha} \over 2})) .0_x
\eqno (27.5)
$$
where $\prod^\uparrow_\alpha A_\alpha = A_{0\ldots01} \circ \ldots \circ
A_{1\ldots1}$ is the composition of operators $A_\alpha$
in the order given by the lexicographic order of the index set.

\Theorem 27.6.
The curvature of the canonical connection of a symmetric space agrees
with its Lie triple product (up to a sign): for $u,v,w \in T_oM$,
$$
R_o(u,v)w = - [u,v,w].
$$

\Proof.
We are going to use the characterization of the curvature tensor
by Theorem 18.3. (See [Lo69] and [Be00, I.2.5] for different
arguments in the real finite-dimensional case,
 rather using covariant derivatives.)
Let $\tau$ be the transposition $(23)$.
 The formula for the connection
$\tau \cdot L_2$ and for its Dombrowski splitting $\Phi_s^\tau$ is
obtained from the formula for $\Phi_2$ by exchanging 
$\eps_2 \leftrightarrow \eps_3$, $v_1 \leftrightarrow v_3$,
$v_{12} \leftrightarrow v_{13}$. Using the commutation rules, one sees
that as a result one gets for  $\Phi_s^\tau$
  a formula of the same type as for $\Phi_s$, but beginning
with $A\circ  B := Q(\eps_2{v_2 \over 2}) Q(\eps_2 {v_3 \over 3})$ instead of
$B \circ A =Q(\eps_3{v_3 \over 2}) Q(\eps_2 {v_2 \over 2}) $.
Using the operator formula $AB = [A,B] BA$, where $[A,B]$ is described
by Theorem 27.2 (the transposition denoted there by $\kappa$ is now $\tau$),
 we see that
the difference $\Phi_s^\tau - \Phi_s$
 (with respect to an arbitrary chart representation) belongs to the
axis $\eps_1 \eps_2 \eps_3 V$ and is given by
$$
-  [Q({\epsilon_1 v_1 \over 2}), Q({\epsilon_2 v_2 \over 2})].\eps_3 v_3.
$$
On the one hand, by  Formula (18.7),
 this term is equal to
$\eps_1 \eps_2 \eps_3 R(v_1,v_2)v_3$; on the other hand, by Theorem 27.2, it
is equal to $\eps_1 \eps_2 \eps_3 [v_1,v_2,v_3]$, whence the claim.
\qed

\msk \nin
{\bf 27.7.} {\sl The structure of $J^3 M$.} 
The preceding calculations permit to give an explicit description of the
symmetric space $(J^3 M)_o$ with respect to 
the symmetrized linear structure $\Sym(DL)$  from Corollary 20.3.
We will see that this ``third order model'' is the first non-trivial
one, and it contains information that is equivalent to the Lie triple 
system at $o$.
In fact, one finds that the restrictions of $\Sym(DL)$ and of the derived
linear structure $DL$ (corresponding to the map $\Phi_2$ calculated above)
coincide on $J^3 M$, and
the result of the calculation is given by the following formula,
where $\Phi_2$ is considered as an identification,
$$
\eqalign{
\mu(\delta v_1 + \delta^{(2)} v_2+ \delta^{(3)} v_3, 
 \delta w_1 + \delta^{(2)} w_2+ & \delta^{(3)} w_3)) 
 = \delta (2v_1-w_1) + \cr
&  \delta^{(2)} (2v_2-w_2) + \cr
& \delta^{(3)} (2v_3 - w_3 - ([v_1,w_1,v_1] + [w_1,v_1,w_1])) . \cr}
$$
In terms of the quadratic map, this reads
$$
\eqalign{
Q_o(\delta v_1 +  \delta^{(2)}  v_2+ \delta^{(3)} v_3, 
 \delta w_1 +  & \delta^{(2)} w_2+  \delta^{(3)} w_3)) 
 =  \delta (2v_1+w_1) +\cr
&  \delta^{(2)} (2v_2+w_2) + \cr
& \delta^{(3)} (2v_3 +w_3 + ([v_1,w_1,v_1] - [w_1,v_1,w_1])) . \cr}
$$
In particular, 
$$
Q(\delta v) = \tau_v + \delta^{(3)} ([v,\cdot,v] - [\cdot,v,\cdot])
$$
is ordinary translation by $v$ with respect to the linear structure plus
an extra term which is the sum of a linear and a homogeneous quadratic map.
Note that the Lie triple system of $(J^3 M, 0_o)$ is
$(J^3 \m)_0 = \delta \m \oplus \delta^{(2)} \m \oplus \delta^{(3)} \m$
 with product
$$
\eqalign{
& [\delta v_1 + \delta^{(2)} v_2 + \delta^{(3)} v_3,
\delta u_1 + \delta^{(2)} u_2 + \delta^{(3)} u_3,
\delta w_1 + \delta^{(2)} w_2 + \delta^{(3)} w_3]  = \cr
& (0,0,\delta^3 [v_1,u_1,w_1])
 = (0,0,6 \, \delta^{(3)} [v_1,u_1,w_1]) \cr}
$$
where $\m$ is the Lie triple system of $(M,o)$.
It is nilpotent, but in general non-abelian, and it contains (if $6$
is invertible in $\K$) all information
on the original Lts $\m$.
 Thus the space $(J^3 M)_o$ may be considered
as a ``faithful model" of the original symmetric space $M$ in the sense
that it contains all infinitesimal information on $M$.

\msk \nin {\bf 27.8.} {\sl On  exponential mappings.}
We say that a symmetric space $M$ {\it has an exponential map} if, for every
$x \in M$ and $v \in T_x M$, there exists a unique homomorphism of symmetric
spaces $\gamma_v:\K \to M$ such that $\gamma(0)=x$, $\gamma_v'(0)=v$, and such
that the map $\Exp_x:T_x M \to M$, $v \mapsto \gamma_v(1)$ is smooth
(cf.\ [BeNe05, Ch.\ 2]). 
It is known that every real finite-dimensional symmetric space has an
exponential map, and moreover this exponential map agrees with
 the exponential map
of the canonical connection (see [Lo69]). In our general setting, not every 
symmetric space will  have an exponential map, but we will show in Chapter
30 that, if $\K$ is a field of characteristic zero,
 there is a canonical ``exponential jet''
$$
\Exp^{(k)}_x: T_0((T^k M)_x) \to (T^k M)_x
$$
having the property that, if $M$ has an exponential map, then
 $\Exp^{(k)}_x = T^k(\Exp_x)_0$.

\vfill\eject

%% file: Dgeo4.tex

\def\date{25.10.2006}

\def\title{VI. DIFFEOMORPHISM GROUPS} 
\def\rightheadline{\hfil{\eightrm \title}\hfil\tenrm\folio}

\sectionheadline{VI. Diffeomorphism groups and the exponential jet }

\nin
In this final part of the main text we bring together ideas from
Part V (Lie theory) and Part IV (Higher order geometry)
by considering the group $G=\Diff_\K(M)$ as a sort of Lie group
and relating its higher order ``tangent groups'' $T^k G =
\Diff_{T^k \K}(T^k M)$ to classical notions such as flows and
geodesics.

\subheadline
{28. Group structure on the space of sections of $T^k M$}

\nin {\bf 28.1.} {\sl Groups of diffeomorphisms.}
We denote by $\Diff_{T^k \K}(T^k M)$ the group of
diffeomorphisms $f:T^k M \to T^k M$ which are smooth over the
ring $T^k \K$. 
By Theorem 7.2, there is a natural injection
$$
\Diff_\K(M) \to \Diff_{T^k \K}(T^k M), \quad f \mapsto T^k f. 
$$

\msk
\nin {\bf 28.2.} {\sl The space of sections of $T^k M$.}
We denote by  $\X^k(M):=\Gamma(M,T^k M)$
the space of smooth sections of the bundle $T^k M \to M$.
A chart $\phi:U \to V$ induces a bundle chart $T^k \phi:T^k U \to T^k V$,
and with respect to such a chart, a section $X:M \to T^k M$ over $U$
is represented by
$$
X(x) = x + \sum_{\alpha >0} \eps^\alpha X_\alpha(x)
\eqno (28.1)
$$
with (chart-dependent) vector fields $X_\alpha:U \to V$.
The inclusions of axes $\iota_\alpha:TM \to T^k M$ induce inclusions
$$
\X(M) \to \X^k(M), \quad X \mapsto \iota_\alpha \circ X = \eps^\alpha X.
$$
We call such sections {\it (purely) vectorial}.

\Theorem 28.3.
\item{(1)}
Every section $X:M \to T^k M$ admits a unique extension to a 
diffeomorphism $\tilde X:T^k M \to T^k M$ 
which is smooth over the ring $T^k \K$. Here, the term ``extension'' 
means that $\tilde X \circ z = X$, where $z:M \to T^k M$ is the zero section. 
In a chart representation as above, this extension is given by
$$
\eqalign{
\tilde X(x+\sum_\alpha \eps^\alpha v_\alpha) &
 =x+ \sum_\alpha \eps^\alpha \Bigl(v_\alpha + X_\alpha(x) + \cr
& \quad \quad 
\sum_{\ell=2}^{\vert \alpha \vert} \sum_{\Lambda \in {\cal P}_\ell(\alpha)} 
\sum_{j=1,\ldots,\ell} (d^{\ell-1} X_{\Lambda^j})(x)
(v_{\Lambda^1},\ldots,\hat v_{\Lambda^j},
\ldots, v_{\Lambda^\ell})\Bigr)  
\cr}
$$
(where ``$\hat v$" means that the corresponding term is to be omitted).
If $X$ is a section of $J^k M$ over $M$, then the diffeomorphism $\tilde X$
preserves $J^k M$.
\item{(2)}
There is a natural group structure on $\X^k(M)$, defined by the formula
$$
X \cdot Y := \tilde X \circ Y .
$$
In a chart representation, the group structure of $\X^k(M)$ is given by
$$
\eqalign{
(X \cdot Y)(p) & =p+ \sum_\alpha \eps^\alpha
\Bigl(X_\alpha(p) + Y_\alpha(p) + \cr
& \quad  \sum_{\ell=2}^{|\alpha|} \sum_{\Lambda \in
{\cal P}_\ell(\alpha)} \sum_{j=1}^\ell (d^{\ell-1} X_{\Lambda^j})(p)
 \bigl( Y_{\Lambda^1}(p),
\ldots,\hat{Y_{\Lambda^j}(p)}, \ldots,Y_{\Lambda^\ell}(p) \bigr)\Bigr).
\cr}
$$
The sections of $J^k M$ form a subgroup of $\X^k(M)$, and the imbeddings
$\X(M) \to \X^k(M)$, $X \mapsto \eps^\alpha X$ are group homomorphisms.
\item{(3)}
The map $X \mapsto \tilde X$ is an injective group homomorphism from
$\X^k(M)$ into $ \Diff_{T^k\K}(T^kM)$. Elements of $\Diff_{T^k \K}(T^k M)$
permute fibers over $M$, and hence there is a well-defined projection
$$
\Diff_{T^k\K}(T^k M) \to \Diff_\K(M), \quad  F \mapsto f:=\pi \circ F \circ z .
$$
The injection and projection thus defined give rise to
a splitting exact sequence of groups
$$
\matrix{
0 &\to& \X^k(M)&\to&\Diff_{T^k\K}(T^k M)& \to& \Diff_\K(M)& \to & 1 \cr},
$$
and hence $\Diff_{T^k\K}(T^k M)$
 is a semidirect product of $\Diff_\K(M)$ and $\X^k(M)$.

\Proof. (1)
We prove uniqueness and existence of the extension of $X:M \to T^k M$
by using a chart representation $X:U \to T^k V$ of $X$.
The statement to be proved is: 

\ssk
\item{(A)} {\it
Let $V$ and $W$ be $\K$-modules and $f: U \to T^k W$ a smooth map (over
$\K$) defined on an open set $U \subset V$. Then there exists a 
unique extension of $f$ to a map $\tilde f:T^k U \to T^k W$ which is smooth
over the ring $T^k \K$.}

\ssk \nin 
We first prove uniqueness.
In fact, here the arguments from the proof of Theorem 7.5 apply word by word --
there the uniqueness statement was proved under the assumption that
 $f(U) \subset W$ (leading to the uniqueness statement in Theorem 7.6).
 But all arguments (i.e., essentially, the second order
Taylor expansion) work in the same way for a $T^k \K$-smooth
map  taking values in $T^k W$, not only in $W$.  

\ssk
Next we prove the existence part of (A). 
By Theorem 7.2, $T^k f:T^k V \to T^k(T^k W)$
is smooth over the ring $T^k \K$. We will construct a canonical smooth
map $\mu_k:T^k(T^k W) \to T^k W$ such that $\tilde f:=\mu_k \circ
T^k f$ has the desired properties. 
Recall that $T^k \K$ is free over $\K$ with basis $\eps^\alpha$,
$\alpha \in I_k$. The algebraic tensor product $T^k \K \otimes_\K T^\ell \K$ 
is again a commutative algebra over $\K$, isomorphic to $T^{k+\ell}\K$
after choice of a partition $\N_{k+\ell} = \N_k \dot \cup \N_\ell$. 
For any commutative algebra $A$, the product map $A \otimes A \to A$ is
an algebra homomorphism, so we have,  since $T^k \K$ is commutative,
a  homomorphism of $\K$-algebras
$$
T^k m: T^{2k} \K \cong T^k \K \otimes_\K T^k \K \to T^k \K .
$$
 Explicitly, with respect
to the $\K$-basis $\eps^\alpha \otimes \nu^\beta$, $\alpha,\beta \in I_k$
(where the $\nu_i$ obey the same relations as the $\eps_i$),
 it is given by
$$
T^k m( \sum_{\alpha, \beta \in I_k } t_{\alpha,\beta} \,
 \eps^\alpha \otimes \nu^\beta) =
\sum_{\gamma \in I_k} \eps^\gamma
\sum_{\alpha , \beta \in I_k\atop
\alpha + \beta = \gamma} t_{\alpha,\beta}.
$$
For any $\K$-module $V$, the ring homomorphism $T^km:T^{2k}\K \to T^k \K$
induces a homomorphism of modules
$$
\eqalign{
\mu_k:
T^{2k} V = V \otimes_\K T^{2k}\K & \to T^k V = V \otimes_\K T^k \K, \cr
x + \sum_{\alpha, \beta \in I_k\atop
(\alpha,\beta) \not= 0 } \eps^\alpha \nu^\beta v_{\alpha,\beta} & \mapsto 
x + \sum_{\gamma \in I_k} \eps^\gamma
\sum_{\alpha , \beta \in I_k\atop
\alpha + \beta = \gamma} v_{\alpha,\beta}. \cr}
$$
This map is $T^{2k}\K$-linear (recall that
any ring homomorphism $R \to S$ induces an $R$-linear map 
$V_R \to V_S$ of scalar extensions; here $R=T^{2k}\K$, $S=T^k \K$).
It is of class $C^0$ since it is a finite composition of certain sums
and projections which are all $C^0$, and hence $\mu_k$ is of class $C^\infty$
over $T^{2k} \K$.
There are two imbeddings of $T^k \K$ as a subring of $T^{2k}\K \cong
T^k \K \otimes T^k \K$ having
image $1 \otimes T^k \K$, resp.\  $T^k \K \otimes 1$; thus $\mu_k$
is also smooth over these rings.
Therefore, if we let
$$
\tilde f:= \mu_k \circ T^k f: T^k U \to T^{2k} W \to T^k W,
$$
$\tilde f$ is smooth over $T^k \K$ (which we may identify with 
$T^k \K \otimes 1 \subset T^{2k} \K$). Let us prove that it is
an extension: for all $x \in U$, 
$$
\tilde f(x) = \mu_k (T^k f(x)) = \mu_k (f(x)) = f(x)
$$
where the second equality holds since $T^k f(x)=T^k f(x+0)=f(x)$ and the
third equality holds since $\mu_k(x+\sum_\alpha \eps^\alpha v_{\alpha,0})=
x+\sum_\alpha \eps^\alpha v_{\alpha,0}$, i.e.\ 
$\mu_k \circ \iota = \id_{T^k M}$, where $\iota:T^k V \to T^{2k} V$ is the
imbedding induced by $T^k \K \to T^{2k}\K$, $r \mapsto r \otimes 1$.
Thus claim (A) is completely proved.

\msk Now we prove part (1) of the theorem.
The extension $\tilde X$ of $X$ is defined, locally, by using (A) in a chart
representation; by the uniqueness statement of (A), this extension
does not depend on the chart and hence is uniquely and globally
defined. 
(See also Remark 28.4 for variants of this proof.)
The chart formula given in the claim is simply Formula
(28.2), where $f=X=\pr_1+\sum_\alpha \eps^\alpha X_\alpha$ is decomposed into 
its (chart-dependent) components according to Formula (28.1).
All that remains to be proved is that $\tilde X$ is bijective with smooth
 inverse.
This will be done in the context of the proof of part (2).

\msk (2)
Clearly, $X \cdot Y$ is again a smooth section of $T^k M$, and hence the
product is well-defined. The uniqueness part of (1)
shows that $\tilde{X \cdot Y}=\tilde X \circ \tilde Y$ since both sides
are $T^k \K$-smooth extensions of $\tilde X \circ Y$.  This property
implies associativity:
$$
(X \cdot Y)\cdot Z = \tilde{X \cdot Y} \circ Z = \tilde X \circ \tilde Y \circ
Z = \tilde X \circ (Y \cdot Z) = X \cdot (Y \cdot Z).
$$
Moreover, for the zero section $z:M \to T^k M$ we have $\tilde z = \id_{T^k M}$
and hence $z$ is a neutral element.

All that remains to be proved is that $X$ has an inverse in $\X^k(M)$.
Equivalently, we have to show that $\tilde X$ is bijective and has
an inverse of the form $\tilde Z$ for some $Z \in \X^k(M)$.
This is proved by induction on $k$. In case $k=1$, $\tilde X$
is fiberwise translation in tangent spaces and hence is bijective
with $(\tilde X)^{-1} = \tilde{-X}$. For the case $k=2$, see the
explicit calculation of $(\tilde X)^{-1}$ in the proof of Theorem 14.3.
For the proof of the induction step for general $k$, we give two
slightly different versions:

\ssk
(a) First version. Note that the chart formula for $X \cdot Y$ from the 
claim is a simple consequence of the chart formula for $\tilde X(v)$ from
Part (1) (take $v=Y(x)$). This formula shows that left multiplication
by $X$, $l_X:\X^k(M) \to \X^k(M)$, $Y \mapsto X \cdot Y$,
is an affine-multilinear self-map of the multilinear space $\X^k(M)$
in the sense  explained in Section MA.19.
 Moreover, $l_X$ is unipotent, hence regular, and by
the affine-multilinear version of Theorem MA.6 (cf.\ MA.19)
it is thus bijective.
Now let $Z:=(l_X)^{-1}(0)$ where $0:=z$ is the zero section. 
Then $X \cdot Z=l_X \circ (l_X)^{-1}(0)=0$.
Similarly, $Y \cdot X = 0$, where $Y:=(r_X)^{-1}(0)$ is gotten from
the bijective right multiplication by $X$. It follows that $X$ is invertible
and $X^{-1}=Z=Y$.

\ssk
(b) Second version. 
Note first that, if $\eps^\alpha X$, $\eps^\alpha Y$ with $X,Y \in \X(M)$
are purely vectorial, then $\eps^\alpha X \cdot \eps^\alpha Y = 
\eps^\alpha (X+Y)$, as follows from the explicit formula in Part (2).
Thus $\eps^\alpha X$ is invertible with inverse $\eps^\alpha (-X)$.
One proves by induction that every element $X \in \X^k(M)$ can
be written as a product of purely vectorial sections (Theorem 29.2):
$X = \prod_\alpha^\uparrow \eps^\alpha X_\alpha$, and hence $X$ is  
invertible.

\msk (3) We have already remarked that $X \mapsto \tilde X$ is
a group homomorphism, and injectivity follows from the extension
property. Next we show that the projection is well-defined. More generally,
we show that a $T^k \K$-smooth map $F:T^k M \to T^k N$ maps fibers to
fibers. In fact, letting $f:=F \circ z:M \to T^k N$, we are in the
situation of Statement (A) proved above. It follows that $F=\tilde f$
is, over a chart domain, given by the formula from Theorem 7.5, which
clearly implies that $F(x+\sum_{\alpha >0} \eps^\alpha v_\alpha) \in
F(x) + \sum_{\alpha >0} \eps^\alpha W$, and hence $F$ maps fibers to fibers.
In particular, the projection $\Diff_{T^k \K}(T^k M) \to \Diff_\K(M)$ is
well-defined.
Let us prove that the imbedded image $\X^k(M)$ is the
kernel of the natural projection $\Diff_{T^k \K}(T^k M) \to \Diff_\K(M)$,
$F \mapsto f$. In fact, $X:=F \circ z:M \to T^k M$ is a section of
$T^k M$ iff $f=\id_M$. But then, by the uniqueness statement in
Part (1), $F=\tilde X$; conversely, for any $X \in \X^k(M)$, 
$\tilde X$ clearly belongs to the kernel of the projection.
Finally, it is clear that $f \mapsto
T^k f$ is a splitting of the projection 
$\Diff_{T^k \K}(T^k M) \to \Diff_\K(M)$.
\qed

\msk \nin {\bf 28.4.} {\sl Comments on Theorem 28.3 and its proof.}
For the real finite-dimensional  case  of Theorem 28.3 (2), see
[KM87, Theorem 4.6] and [KMS93, Theorem 37.7], where the proof is
carried out in the ``dual" picture, corresponding to function-algebras.
The proof given there has the advantage that
the problem of bijectivity (i.e., the existence of inverses)
 is easily settled by using the Neumann 
series in a nilpotent associative algebra, but it does not carry over 
to our general situation. Moreover, in a purely real theory it is less
visible why the extension of a section $X$ to a diffeomorphism $\tilde X$
is so natural. 

\ssk
Comparing with the proof of [KMS93, Theorem 37.7], the proof of the
existence part in (1) may be reformulated as follows.
We define  $\tilde X = \mu_{2k,k} \circ T^k X$ where
$\mu_{2k,k}:T^{2k} M \to T^k M$ is the natural 
map which, in a chart representation, is given by $\mu_k:T^k(T^k V) \to T^k
V$. Note that, for $k=1$, $\mu_{2,1}:TTM \to TM$ can directly be
defined by
$$
\mu_{2,1}= a \circ (p_1 \times_M p_2): TTM \to TM \times_M TM \to TM,
$$
where $p_i:TTM \to TM$, $i=1,2$, are the canonical projections and 
$a:TM \times_M TM \to TM$ the addition map of the vector bundle $TM \to M$. 
In a chart,
$$
\mu_{2,1}(x+\eps_1 v_1 + \eps_2 v_2 + \eps_1\eps_2 v_{12})=
x + \eps(v_1 + v_2)
$$
which corresponds to the formula for $\mu_1:TTV \to TV$. 



\vfill \eject 

\subheadline{29. The exponential jet for vector fields}

\nin {\bf 29.1.} {\sl The ``Lie algebra of $\X^k(M)$".}
As explained in the introduction, heuristically we may think of 
$G:=\Diff_\K(M)$ as a Lie group with Lie algebra $\g=\X(M)$.
Then the group $\Diff_{T^k\K}(T^k M)$ takes the r\^ole of $T^k G$ and the group
$\X^k(M)$ the one of $(T^k G)_e$. The Lie algebra of $T^k G$ is then
$T^k \g = \g \otimes_\K T^k \K$ and therefore the one of $\X^k(M)$
has to be $(T^k \g)_0 = \oplus_{\alpha >0} \eps^\alpha \g$. It is the
space of sections of the {\it axes-bundle}
$A^k M=\bigoplus_{\alpha \in I_k,\alpha >0} \eps^\alpha TM$. 
The space of  sections of $A^k M$  will be denoted by
$$
{\cal A}^k(M):= \Gamma(M,A^k M) = 
\bigoplus_{\alpha \in I_k,\alpha >0} \eps^\alpha \X(M).
\eqno (29.1)
$$
We will construct  canonical bijections ${\cal A}^k(M) \to
\X^k(M)$, taking the r\^ole of left-trivialization,
 resp.\ of an exponential map.

\Theorem 29.2. 
There is a bijective ``left trivialization'' of the group $\X^k(M)$,
given by
$$
{\cal A}^k(M) \to \X^k(M), \quad
(X_\alpha)_\alpha \mapsto \prod^\uparrow_{\alpha \in I_k \atop \alpha >0} 
\eps^\alpha X_\alpha,
$$
where the product is taken in the group $\X^k(M)$, with order corresponding
to the lexicographic order of the index set $I_k$.
Then the product in $\X^k(M)$ is given by 
$$
\prod^\uparrow_\alpha \eps^\alpha X_\alpha \cdot 
\prod^\uparrow_\beta \eps^\beta Y_\beta = 
\prod^\uparrow_\gamma \eps^\gamma Z_\gamma
$$
 with $Z_\gamma$ given by
the ``left product formula" from {\rm Theorem 24.7.} 

\Proof. 
First of all, note that Theorem 14.4 implies the following commutation
relations for purely vectorial elements in the group $\X^k(M)$:
for all vector fields $X_\alpha$, $Y_\beta \in \X(M)$,
$$
[\eps^\alpha X_\alpha, \eps^\beta Y_\beta] = \eps^{\alpha + \beta}
[X_\alpha,Y_\beta],
\eqno (29.2)
$$
where the last bracket is taken in the Lie algebra $\X(M)$.
Now we define the map
$\Psi:{\cal A}^k(M) \to \X^k(M)$, $(X_\alpha)_\alpha \mapsto 
\prod^\uparrow_{\alpha \in I_k \atop \alpha >0} \eps^\alpha X_\alpha$
as in the claim and show that it is a bijection. 
We claim that the image of $\Psi$ 
is a subgroup of $\X^k(M)$. In fact, this follows from the relations (29.2):
re-ordering the product of two ordered products is carried out
by the procedure described in the proof of Theorem 24.7 and leads to the
left-product formula mentioned in the claim. 
Thus the image of $\Psi$ is the subgroup of $\X^k(M)$ generated by the purely
vectorial sections, with group structure as in the claim.
But this group is in fact the whole group $\X^k(M)$:
using the explicit chart formula for the product in $\X^k(M)$ 
(Theorem 28.3 (2)), it is seen that $\Psi$ is a unipotent multilinear
map between multilinear spaces over $\K$, whence is bijective (Theorem MA.6).
\qed

\msk 
\nin In the same way as in the theorem, 
a ``right-trivialization map" ${\cal A}^k(M) \to \X^k(M)$ is defined.
Next we prove the analog of Theorem 25.2 on the exponential map:

\Theorem 29.3. Assume  that $\K$ is a non-discrete
topological field of characteristic zero and $M$ is a manifold over $\K$.
Then there is a unique map
$\exp_k: {\cal A}^k(M) \to \X^k(M)$ such that:
\item{(1)}
In every chart representation, $\exp_k$ is a polynomial map (of degree
at most $k$); equivalently, $\exp_k^L:=\Psi^{-1} \circ \exp_k:{\cal A}^k
(M) \to {\cal A}^k(M)$, with $\Psi$ as in the preceding theorem, is
a polynomial map,
\item{(2)}
for all $n \in \Z$ and $X \in \Gamma(M,A^k M)$,
$\exp(nX)=(\exp X)^n$,
\item{(3)}
on the axes, $\exp_k$ agrees with the inclusion maps, i.e., if $X:M \to TM$
is a vector field, then
$\exp(\eps^\alpha X)= \eps^\alpha X$ is the purely vectorial section
corresponding to $\eps^\alpha X$.

\Proof.
As for Theorem 25.2, this in an immediate
consequence of Theorem PG.6, applied to the polynomial group
$\X^k(M)$.
\qed

\Theorem 29.4.
The map $\exp_k$ from the preceding theorem 
commutes with the canonical action of the symmetric
group $\Sigma_k$ and hence restricts to a bijection 
$$
\exp_k: ({\cal A}^k(M))^{\Sigma_k} =
\oplus_{j=1}^k \delta^{(j)} \X(M) =
\Gamma (\oplus_{j=1}^k \delta^{(j)} TM)
\to (\X^k(M))^{\Sigma_k} = \Gamma(J^k M).
$$
In particular, we get an injection
$$
\X(M) \to G^k(M), \quad X \mapsto \tilde{\exp_k(\delta X)}
$$
such that $\tilde{\exp_k(\delta nX)} = (\tilde{\exp_k(\delta X)})^n$
for all $n \in \Z$. 

\Proof.
This is the analog of Theorem 25.4, and it is proved in the same way.
\qed

\msk
\nin {\bf 29.5.} {\sl Smooth actions of Lie groups.}
Assume $G$ is a Lie group over $\K$, acting smoothly on $M$ via 
 $a:G \times M \to M$.
Deriving, we get an action
$T^k a:T^kG \times T^k M \to T^k M$ of $T^k G$ on $T^k M$.
Then $T^k a(g,\cdot)$ is smooth over $T^k \K$ and hence
 belongs to the group $\Diff_{T^k\K}(T^k M)$.
Therefore we may consider $T^k a$ and $(T^k a)_e$ as group homomorphisms
$$
T^k a: T^k G \to \Diff_{T^k\K}(T^k M),
 \quad {\rm resp.} \quad (T^k a)_e:(T^k G)_e \to
\X^k(M).
$$
In particular, for $k=1$ we get
a homomorphism of abelian groups $\g = T_e G \to \X(M)$, associating to
$X \in \g$ the ``induced vector field on $M$", classically often denoted
by $X_* \in \X(M)$. 
This map is a Lie algebra homomorphism since $T^2 a$ is a group homomorphism
and the Lie bracket on $\g$, resp.\ on $\X(M)$, can be defined in similar
ways via group commutators (Theorem 14.5 and Theorem 23.2).
Taking a direct sum of such maps, we get a Lie algebra homomorphism
$(T^k \g)_0 \to {\cal A}^k(M)$, and from the uniqueness property of the
exponential map we see that the following diagram commutes:
$$
\matrix{
(T^k G)_e & {\buildrel T^k a \over \to} & \X^k(M) \cr
\exp \uparrow \phantom{\exp} & &
\phantom{\exp_k} \uparrow \exp_k \cr
(T^k \g)_0 & \to & {\cal A}^k (M) \cr}
$$
Summing up, everything behaves as if $a:G \to \Diff(M)$ were
a honest Lie group homomorphism.

\msk \nin {\bf 29.6.} {\sl Relation with the flow of a vector field.}
The diffeomorphism $\tilde{\exp_k(\delta X)}$ of $T^k M$ (Theorem 29.4) can be
related to the $k$-jet of the flow of $X$, if flows exist.
Let us explain this.
We assume that $\K=\R$ and $M$ is finite-dimensional (or a Banach manifold).
Then any vector field $X \in \X(M)$ admits a flow map
$(t,x) \mapsto \Fl^X_t(x)$, defined on some neighborhood of
$\{ 0 \} \times M$ in $\R \times M$ (cf.\ [La99, Ch.\ IV]). 
In general, $\Fl^X_1$ is not defined as a diffeomorphism on the 
whole of $M$, and this obstacle prevents us from defining  
the exponential map of the group $G=\Diff(M)$ by $\exp(X):=\Fl^X_1$. 
Nevertheless, this is what one would like to do, and one may give some
sense to this definition by fixing some point $p \in M$ 
and by considering the subalgebra
$$
\X(M)_p := \{ X \in \X(M) \vert \, X(p)=0, \, \forall Y \in \X(M):
[X,Y](p) = 0 \}
$$
of vector fields vanishing at $p$ of order at least $2$. Then it is  
easily seen from the proof of the existence and uniqueness theorem 
for local flows that, for all $X \in \X(M)_p$, $\exp(X):=\Fl^X_1$
is defined on some open neighborhood of $p$,
and moreover that $T_p(\exp (X)) = \id_{T_p M}$ (cf., e.g., [Be96,
Appendix (A2)]). 
Thus the group $\Diff(M)_p$ of germs of local diffeomorphisms that are 
defined at
$p$ and have trivial first order jet there, admits an exponential map
given by $\exp = \Fl_1:\X(M)_p \to \Diff(M)_p$. 
We claim that, under these assumptions, the $k$-jet of $\Fl_1^X$ at $p$
is given by the map from Theorem 29.4:
$$
(\tilde{\exp_k(\delta X)})_p = (T^k (\Fl_1^X))_p,
$$
or, in other words, that the following diagram commutes:
$$
\matrix{
\Diff(M)_p & {\buildrel T^k \over \to} & G^k(M) & {\buildrel \ev_p \over \to}
& {\rm Am}^{k,1}(T^k M)_p \cr
\Fl_1 \uparrow \phantom{\exp} & & & & \phantom{\ev_p} \uparrow \ev_p \cr
\X(M)_p & {\buildrel X \mapsto \delta X \over \to} & {\cal A}^k(M) &
{\buildrel \exp_k \over \to} & G^k(M) \cr}
$$
where ${\rm Am}^{k,1}(T^k M)_p$ denotes the affine-multilinear group of the
fiber $(T^k M)_p$.
The proof is by using the uniqueness of the map $\exp$ from Theorem
PG.6: we contend that
the map $X \mapsto (T^k (\Fl^{X}_1))_p$ also has the properties
(1) and (2) from Theorem PG.6 and hence coincides with the exponential
map constructed above by using Theorem PG.6. Now, Property (1) follows by
observing that
$$
T^k ({\Fl}^{nX}_1) = T^k (({\Fl}^X_1)^n) =(T^k ({\Fl}^X_1))^n,
$$
and Property (2) follows from the fact that $T_p({\Fl}_1^X) = \id_{T_pM}$
(here we need the fact that $X$ vanishes at $p$ of order at least two,
cf.\ [Be96, Appendix (A.2)]) and hence our candidate for the exponential
mapping acts trivially on axes, as it should. 
Summing up, for any $X \in \X(M)_p$ the $k$-jet of $\exp(X)$ at $p$
may be constructed in a purely ``synthetic way", i.e., without using
the integration of differential equations.
If $X$ does not belong to $\X(M)_p$, then there still is a relation,
though less canonical, with the $k$-jet of the flow of $X$, see the
next subsection.

\msk \nin {\bf 29.7.} {\sl Relation with ``Chapoton's formula''.}
If $X$ is a vector field on $V=\R^n$, identified with the corresponding map
$V \to V$, then the local flow of $X$ can be written
$$
{\Fl}_t^X(v) = v + \sum_{j \geq 1} {t^j \over j!} X^{\star j}(v) 
$$
where 
$$
(X \star Y)(x) = (\nabla_X Y)(x) = dY(x) \cdot X(x)
$$
denotes the product on $\X(V)$ given by the canonical flat connection of 
$V$, and $X^{\star j+1}=X \star X^{\star j}$ are the powers with respect
to this product -- see [Ch01, Prop.\ 4]. Thus, letting $t=1$, 
 the exponential map associated to the Lie algebra of
vector fields on $V$ can be written as a formal series
$$
\exp(X) = \id_V + \sum_{j \geq 1}  {t^j \over j!} X^{\star j} .
$$
This formal series may be seen as the projective limit of our
$\exp_k(\delta X):V \to J^k V$ for $k \to \infty$ (cf.\ Chapter 32).
In fact, 
according to Theorem PG.6, we have the explicit expression of 
the exponential map of the polynomial group $\X^k(V)$ by
$\exp(\delta X)=\sum_{j=1}^k 
{1 \over j!} \psi_{j,j}(\delta X)$ where 
$\psi_{j,j}(X)$ is gotten from the homogeneous term of degree $j$
in the iterated multiplication map $m^{(j)}$. The chart formula
in the global chart $V$ for the product (Theorem 28.3 (2)) permits
to calculate these terms explicitly, leading to the formula
$$
\exp_k(\delta X) = \id_V + \sum_{j \geq 1}  {\delta^{(j)} 
\over j!} X^{\star j} .
$$

A variant of these arguments is implicitly contained  in
[KMS93, Prop.\ 13.2] where the exponential mapping of the ``jet groups"
is determined: for a vector field $X$ vanishing
at the origin,
$\exp(j_0^k X)= j_0^k{\Fl}_1^X$.

\vfill \eject 

\subheadline{30. The exponential jet of a symmetric space}

\nin  
{\bf 30.1.}
As an application of the preceding results, we will construct a canonical
isomorphism $\Exp_k:A^k M \to T^k M$ (called the {\it exponential jet},
cf.\ Introduction) for  symmetric spaces over
$\K$. In order to motivate our construction,
we start by considering the group case, where
such a construction has already been given in Chapter
25 (cf.\ also Section 29.5). We relate those results to the
context explained in the preceding chapters.

\msk
\nin {\bf 30.2.} {\sl Case of a Lie group.}
We continue to assume that  $\K$ is a non-discrete topological field
of characteristic zero.
We construct the exponential jet
$\Exp_k:A^k G \to T^k G$ as follows:
for $X \in \g$, let $X^R \in \X(G)$ be the right-invariant vector field
such that $X^R(e)=X$. The ``extension map"
$$
l_R:\g \to \X(G), \quad X \mapsto l_R(X):=X^R
$$
is injective. Denote by   
$$
(T^k l_R)_0:(T^k \g)_0 \to (T^k \X(G))_0 = {\cal A}^k(G)
$$
its scalar extension (simply a direct sum of copies of $l_R$).
Then the exponential
map $\Exp=\Exp_k:(T^k \g)_0 \to (T^k G)_e$ of the polynomial group
$(T^k G)_e$ (Theorem 25.2) is recovered from these data via 
$$
\Exp_k(X)= (\exp ((T^k l_R)  X))(e): \quad \quad
\matrix{(T^k G)_e &  {\buildrel \ev_e \over \leftarrow}  & \X^k(G) \cr
\Exp \uparrow \phantom{\Exp} &  & \phantom{\exp_k} \uparrow \exp_k \cr
(T^k \g)_0 &{\buildrel (T^k l_R)_0 \over \to} & 
{\cal A}^k(G) \cr}
$$
In fact, this is seen by the arguments given in Section 29.5, applied to the
action of $G$ on itself by left translations.

In the case of finite dimension over $\K=\R$ (or for general real Banach
Lie-groups), the group $G$ does admit an exponential map $\Exp:\g \to G$ in
the sense explained in Section 25.9, and it is well-known that then 
the exponential map of $G$ is given by
$\exp(X)=\Fl^{X^R}_1(e)$, i.e., by the commutative diagram
$$
\matrix{
G & {\buildrel \ev_e \over \leftarrow} & \Diff(G) \cr
\Exp \uparrow \phantom{\Exp} & & \phantom{\Fl_1} \uparrow \Fl_1 \cr
\g & {\buildrel l_R \over \to} & \X(G)_c \cr}
$$
(where $\X(G)_c$ is the space of complete vector fields on $G$).
The preceding definition is obtained by
applying, formally (i.e., as if $\Diff(G)$ were a Lie group),
 the functor $(T^k)_0$ to this diagram.

\msk
\nin {\bf 30.3.} 
{\sl Exponential jet of a symmetric space.}
Assume $M$ is a symmetric space over $\K$ with base point $o \in M$.
Recall from the proof of Prop.\ 5.9 that every tangent vector $v \in T_o M$
admits a unique extension to a vector field $l(v)$ (denoted by $\tilde v$ 
in Chapter 5) such that $l(v)$
is a derivation of $M$ and $\sigma_o \cdot l(v) = - l(v)$.
The ``extension map'' $l:=l_o:T_o M \to \X(M)$ is linear, and its higher order
tangent map
$(T^k l)_0:(A^k M)_o \to {\cal A}^k (M)$ is simply a direct sum 
$\oplus_{\alpha>0} \eps^\alpha l$ of copies of this map.

\Theorem 30.4. Let $(M,o)$ be a symmetric space over $\K$ and assume that
$\K$ is non-discrete topological field of characteristic zero.
Then the map $\Exp:=\Exp_k:(A^k M)_o \to (T^k M)_o$ defined by
$$
\Exp_k(v)= (\exp_k (T^k l(v)))(o): \quad \quad
\matrix{(T^k M)_o  & {\buildrel \ev_o \over \leftarrow}  & \X^k(M) \cr
\Exp \uparrow \phantom{\Exp} &  & \phantom{\exp_k} \uparrow \exp_k \cr
(A^k M)_o &{\buildrel (T^k l)_0 \over \to} & 
{\cal A}^k(M) \cr}
$$
is a unipotent
multilinear isomorphism of fibers over $o \in M$ commuting with
the action of the symmetric group $\Sigma_k$.
Moreover, $\Exp_k$ satisfies the ``one-parameter subspace property'' from
{\rm Section 27.8}: for all $u \in (A^k M)_o$, the map
$$
\gamma_u: \K \to (T^k M)_o, \quad t \mapsto \Exp_k(t u)
$$
is a homomorphism of symmetric spaces such that $\gamma_u'(0)=u$.
 In particular, the property
$\Exp(n u)= (\Exp(u))^n$ holds for all $n \in \Z$ (where the powers
in the symmetric space $(T^k M)_o$ are defined as in {\rm Equation (5.5)}).

\Proof. By its definition, $\Exp_k$ is a composition of multilinear maps
and hence is itself multilinear. In order to prove that $\Exp_k$ is
bijective,  it suffices, according to Theorem MA.6, to show that the restriction
of $\Exp_k$ to axes is the identity map. But this is an immediate
consequence of the corresponding property of the map $\exp_k$
(Theorem 29.3, (3)): let $v = \eps^\alpha v_\alpha$ be an element of the axis
$\eps^\alpha V$; then $X:=T^k l( \eps^\alpha v_\alpha)=
\eps^\alpha l(v_\alpha)$ is a purely vectorial section, hence $\exp_k$
applied to this section is simply given by inclusion of axes. Evaluating
at the base point, we get
$$
(\exp_k(\eps^\alpha l(v_\alpha)))(o)= \eps^\alpha l(v_\alpha)(o) = \eps^\alpha
v_\alpha
$$
since $\ev_o \circ l = \id$. 
The map $\Exp$ commutes with the $\Sigma_k$-action since so do all maps
of which it is composed. 

Finally, we prove the one-parameter subspace property.
In a first step, we show that, for all $n \in \Z$,
$ \Exp(n u) = \Exp(u)^n$.
We write $u=x+\sum_\alpha \eps^\alpha v_\alpha$; then, using that
$\exp_k(nX)=(\exp_k(X))^n$,
$$
\eqalign{
\Exp(nu) & = \exp_k(n  \sum_\alpha \eps^\alpha l(v_\alpha)).o =
(\exp(\sum_\alpha \eps^\alpha l(v_\alpha)))^n.o \cr & =
(\exp(\sum_\alpha \eps^\alpha l(v_\alpha)).o)^n = (\Exp u)^n, \cr}
$$
where the third equality is due to a general property of symmetric spaces:
if $(N,p)$ is a symmetric space with base point and $g \in \Aut(N)$ such
that $\sigma_p \circ g \circ \sigma_p = g^{-1}$, then
$(g.p)^n = g^n.p$. (This is proved by an easy induction, using the definition
of the powers in a symmetric space given in Equation (5.5).)
This remark is applied to $(N,p)=((T^kM)_o, 0_o)$ and $g = 
\exp_k(\sum_\alpha \eps^\alpha l(v_\alpha))$. 

Now let $\gamma(t):=\Exp(t u)$. Then, in any chart representation,
$\gamma(t)=tu$ modulo higher order terms in $t$ (since $\Exp$ is unipotent,
as already proved), and hence $\gamma'(0)=u$.
We have to show that $\gamma:\K \to T^k M$ is a homomorphism of symmetric
spaces. As we have just seen,  $\gamma(n)=(\gamma(1))^n$.
By power associativity in symmetric spaces (cf.\ Equation (5.7)), and denoting
by $\mu$ the binary product in the symmetric spaces $\K$, resp.\ $T^k M$,
this implies $\gamma(\mu(n,m)) =\gamma(2n-m)=(\gamma(1))^{2n-m}=\mu(\gamma(m),
\gamma(n))$, and thus the restriction of $\gamma$ to $\Z$ is an ``algebraic
symmetric space homomorphism''.
 But $\gamma:\K \to (T^k M)_o$ is a polynomial map between
two symmetric spaces with polynomial multiplication maps, and hence
by the ``polynomial density argument'' 
used several times in Section PG, it follows
that also $\gamma$ must be a homomorphism. 
\qed

As in Step 1 of the proof of Theorem PG.6, it is seen that $\gamma_u$
is the only {\it polynomial} one-parameter subspace with $\gamma_u'(0)=u$.
We do not know wether this already implies that $\gamma_u$ is the only
{\it smooth} one-parameter subspace with this property.

\msk
In the finite-dimensional case over $\K=\R$, it is known
that the exponential map of a symmetric space can be defined by
$$
\Exp_o(v)=\Fl_1^{l(v)}(o) = \ev_o \circ \Fl_1 \circ l: \quad \quad
\matrix{
M & {\buildrel \ev_o \over \leftarrow} & \Diff(M) \cr
\Exp \uparrow \phantom{\Exp} & & \phantom{\Fl_1} \uparrow \Fl_1 \cr
T_o M & {\buildrel l \over \to} & \X(M)_c \cr}
$$
(see [Lo69], [Be00, I.5.7]).
As in the case of Lie groups, the definition of $\Exp_k$ from the preceding
theorem is formally obtained by 
applying the functor $(T^k)_0$ to this diagram.

\vfill \eject

\subheadline{31. Remarks on the exponential jet of a general connection}

\nin {\bf 31.1.} {\sl The exponential jet.}
As explained in the Introduction (Section 0.16),
 for real finite dimensional or Banach
manifolds $M$ equipped with a (torsionfree) linear connection on $TM$,
we can define a canonically associated bundle isomorphism $T^k \Exp:
A^k M \to T^k M$, called the {\it exponential jet (of order $k$)}.
For a general base field or -ring $\K$, there is no exponential map
$\Exp_x$, but still it is possible to
construct a bundle isomorphism $A^k M \to
T^k M$ having all the good properties of the exponential jet and
coinciding with $T^k \Exp$ in the real finite-dimensional or Banach case.
In the following, we explain the basic ideas of two different approaches
to such a construction; details will be given elsewhere.

\msk \nin {\bf 31.2.} {\sl Vector field extensions.}
The definition of the exponential jet for Lie groups and symmetric spaces
in the preceding chapter suggests the following approach:
in the real finite-dimensional  case, every connection $L$ gives rise
to a, at least locally, on a star-shaped neighborhood $U$
of $x$ defined ``vector field extension map'' $l:=l_x:T_x M \to 
\X(U)$, assigning to $v \in T_x M$ the {\it adapted vector field} $l(v)$
which is defined by taking as value $(l(v))(p) \in T_p M$ the 
parallel transport of $v$ along the unique geodesic segment joining
$x$ and $p$ in $U$. Then clearly $l$ is linear,
and $l(v)(x)=v$. Moreover, one may check that the
covariant derivative of $L$ is recovered by $(\nabla_X Y)(x) = [l_x(X(x)),
Y](x)$. The exponential map $\Exp_x$ is given, as in the preceding chapter,
by $\Exp_x(v)=(\Fl_1^{l(v)})(x)$. 

For general base fields and rings $\K$, one may formalize the properties
of such ``vector field extension maps''. One does not really need 
the whole vector field $l(v) \in \X(U)$, but only its $k$-jet at $x$.
Then one may try to define the exponential jet in a ``synthetic way''
by using the exponential maps of the polynomial groups $\Gm^{1,k}(T^k M)_x$
and  ${\rm Am}^{1,k}(T^k M)_x$ acting on the fiber over $x$.

\msk
\nin {\bf 31.3.} {\sl  The geodesic flow.}
Another ``synthetic'' approach to the exponential jet is motivated as follows:
 in the classical case, if the connection is  
complete, the geodesic flow $\Phi_t: TM \to TM$ is the flow of a vector
field $S:TM \to T(TM)$ which is called the {\it spray of the connection}
and which is equivalent to the (torsionfree part of) the connection itself
(cf.\ Chapter 11). 
Then the exponential map $\Exp_x:T_x M \to M$ at a point $x \in M$
is recovered from the flow $\Fl_1^S$ via 
 $\Exp_x(v)=\pi (\Fl_1^S(v))$, where $\pi:TM \to M$ is the
canonical projection:
$$
\matrix{TM & {\buildrel \Fl_1^S \over \to} & TM \cr
\uparrow & & \downarrow \cr
T_x M & {\buildrel \Exp_x \over \to } & M \cr}
$$
Applying the $k$-th order tangent functor $T^k $ to this diagram and
 restricting to the fiber over $0_x$, we express the exponential
jet via the $k$-jet of $\Fl_1^S$ and canonical maps:
$$
\matrix{T^{k+1}M & {\buildrel T^k \Fl_1^S \over \to} & T^{k+1}M \cr
\uparrow & & \downarrow \cr
T^k(T_x M) & {\buildrel T^k (\Exp_x)
 \over \to } & T^k M \cr}
$$
As explained in Chapter 29, under some conditions,
the map $T^k \Fl_1^S$ can be expressed in our general framework via
$\exp(\delta S)$, and thus we should get an expression
of the exponential jet in terms of maps that all can be defined in
a ``synthetic'' way.
However, the problem in the preceding diagram is that we cannot simply
restrict everything to the fiber over $x \in M$,  because
$\Exp_x$ does not take values in some ``fiber over $x$''.
One may try to overcome this problem by deriving once more, or, in other
words, by introducing one more scalar extension, thus imbedding 
$A^k M$ into $T^{k+2} M$ and not into $T^{k+1}M$ and defining
$$
\Exp_k:= \pr \circ \, \tilde{\exp(\delta S)} \circ \iota: \quad \quad 
\matrix{
T^{k+1} (TM) & {\buildrel \tilde{\exp(\delta S)} \over \to} 
& T^{k+1} (TM) \cr
\iota\uparrow\phantom{x} & &\phantom{xx} \downarrow \pr \cr
A^k M  & {\buildrel \Exp_k \over \to} & T^k M \cr}
$$
where the imbedding $\iota$ and the projection $\pr$ are defined by
suitable choices of $k$ infinitesimal units among the
$k+2$ infinitesimal units $\eps_0,\ldots,\eps_{k+1}$, and  the diffeomorphism
 $\tilde{\exp(\delta S)}:T^{k+2}M \to T^{k+2}M$ is associated to
the spray $S:T_{\eps_0}M \to T_{\eps_1}T_{\eps_0}M$ 
of the connection $L$ in the way described by Theorem 29.4.

\msk \nin
{\bf 31.4.} {\sl Jet of the holonomy group.}
A major interest of the construction of the exponential jet, either in
the way indicated in 31.2 or in the one from 31.3, should be the possibility
to describe the $k$-jet of the {\it holonomy group} of $(M,L)$ at $x$.
It should be given by parallel transport along ``infinitesimal rectangles'',
which can be developed at arbitrary order into its ``Taylor series''.
All these are interesting topics for further work.

\vfill \eject

\subheadline{32. From germs to jets and from jets to germs}

\msk \nin {\bf 32.1.} {\sl Infinity jets (the functor $J^\infty$).}
The jet projections $\pi_{\ell,k}:J^k M \to J^\ell M$ ($k \geq \ell$) form a
projective system, i.e., if $k_1 \geq k_2 \geq k_3$, then
$\pi_{k_3,k_1}=\pi_{k_3,k_2} \circ \pi_{k_2,k_1}$ (cf.\ Equation (8.6)).
Therefore we can form the projective limit
$J^\infty M:= \lim_{\leftarrow} J^k M$.
We do not define a topology on $J^\infty M$. As usual for projective limits,
we have projections
$$
\pi_k: J^\infty M \to J^k M, \quad u=(u_0,u_1,u_2,\ldots) \mapsto u_k
$$
(where $u_j \in J^j M$ with $u_j=\pi_{j,i}(u_i)$). In particular, 
$$
\pi_0: J^\infty M \to M, \quad u=(u_0,u_1,u_2,\ldots) \mapsto u_0
$$
is the projection of the ``algebraic bundle $J^\infty M$'' onto its base
$M$. 
Clearly,  smooth maps $f:M \to N$ induce maps $J^\infty f: J^\infty M \to 
J^\infty N$, and the functorial rules hold. 
For instance, for every Lie group $G$ we get an abstract group $J^\infty G$.
Note that it is not possible to define a functor $T^\infty$ in the same
way since the projections associated to $T^k$ form a cube and not
an ordered sequence of projections. 

\msk \nin {\bf 32.2.} {\sl Formal jets.}
There are several possible ways to define homomorphisms between
fibers $(J^k M)_x \to (J^k N)_y$, $x \in M$, $y \in N$,
$k \in \N \cup \{ \infty \}$. Our definition will be such that, 
for $k = \infty$, homomorphisms in a fixed fiber essentially
 correspond to ``formal power series''.
First, for $k \in \N$ and $x \in M$ and $y \in N$, we let
$$
J_x^k(M,N)_y := \Hom ((J^k M)_x,(J^k N)_y)),
$$
where ``$\Hom$'' means the space of all algebraic homomorphisms 
$f:(J^k M)_x \to (J^k N)_y$ in the
strong sense that $f$ comes from a totally symmetric and shift-invariant
homomorphism of multilinear spaces $\tilde f:(T^k M)_x \to (T^k N)_x$.
(One would also get a category if one only requires that $\tilde f$ be
weakly symmetric in the sense of Section SA.2. But the condition of
total symmetry and shift-invariance is necessary if one wants to integrate
jets to germs, see below.)
A {\it formal jet (from the fiber at $x \in M$ to the fiber at $y \in N$)}
is the projective limit 
$f_\infty :(J^\infty M)_x \to (J^\infty N)_y$ 
given by a family $f_k \in J_x^k(M,N)_y$, $k \in \N$, such that
$\pi_{\ell,k} \circ f_k = f_\ell \circ \pi_{\ell,k}$. 
We denote by
 $J^\infty_x(M,N)_y$ the set of formal jets from $x$ to $y$.
For instance, if $h:M \to N$ is a smooth map such that $h(x)=y$,
then $(J^\infty h)_x:(J^\infty M)_x \to (J^\infty N)_y$ is a formal jet.
It is not required that all formal jets are obtained in this way --
see Section 32.5 below.

\msk \nin {\bf 32.3.} {\sl Germs.}
For manifolds $M,N$ over $\K$ and $x \in M$, $y \in N$,
 we denote by ${\cal G}_x(M,N)_y$ the space of germs of locally defined
 smooth maps $f$ from $M$ to $N$ such that $f(x)=y$, defined as usual
by the equivalence relation $f \sim g$ iff the restrictions of $f$ and $g$ 
agree on some
neighborhood of $x$. Then the $k$-jet $(J^k f)_x$ of a representative $f$
of $[f]=f/\sim$ does not depend on the choice of $f$ because, in a chart,
$f|_U=g|_U$ implies, by the ``Determination Principle'' (Section 1.2),
$f^{[k]}|_{U^{[k]}} = g^{[k]}|_{U^{[k]}}$ (notation from 
Chapter 1) and hence $d^k f(x)=d^k g(x)$ for all $x \in U$.
Therefore, for all $k \in \N \cup \{ \infty \}$, we have a well-defined 
{\it $k$-th order evaluation map}
$$
\ev_x^k: {\cal G}_x(M,N)_y  \to J^k_x(M,N)_y, \quad [f] \mapsto (J^k f)_x
$$
assigning to a germ at $x$ its $k$-th order jet there.
The kernel of $\ev_x^k$ is the space of germs vanishing of order $k+1$ at $x$;
if $k=\infty$ and $\ev_x^\infty([f])=0$, we say that {\it $f$ vanishes of
infinite order at $x$}.

\Proposition 32.4. Assume the integers are invertible in $\K$. Then,
for all $k \in \N$ and $x \in M$, $y \in N$, the $k$-th order evaluation
map $\ev_x^k$ is surjective. 

\Proof. 
Let $f \in J^k_x(M,N)_y$, i.e., $f$ is given by a totally symmetric and 
shift-invariant homomorphism $\tilde f:(T^k M)_x \to (T^k N)_y$.
We  choose bundle charts at $x$ and at $y$  and represent $\tilde f$ in the
usual form of a homomorphism of multilinear maps
$$
\tilde f(\sum_\alpha \eps^\alpha v_\alpha) = \sum_\alpha \eps^\alpha 
\sum_{\ell=1}^
{|\alpha|} \sum_{\Lambda \in
{\cal P}_\ell(\alpha)} b^\Lambda(v_{\Lambda^1},\ldots,v_{\Lambda^\ell}).
$$
As shown in Theorem SA.16, total symmetry and shift-invariance imply that
$b^\Lambda$ depends essentially only on $\ell=\ell(\Lambda)$: there exist
multilinear maps $d^\ell:V^\ell \to V$ (in the usual sense) such that
$$
\tilde f(\sum_\alpha \eps^\alpha v_\alpha) = \sum_\alpha \eps^\alpha 
\sum_{\ell=1}^
{|\alpha|} \sum_{\Lambda \in
{\cal P}_\ell(\alpha)} d^\ell(v_{\Lambda^1},\ldots,v_{\Lambda^\ell}).
$$
Then we define, still with respect to the fixed charts, on a neighborhood 
of $x$, a map
$$
F(x+h):=y + \sum_{\ell=1}^k {1 \over \ell !} d^\ell(h,\ldots,h).
$$
Then $d^\ell F(x) = d^\ell$ and hence, by Theorem 7.5., $(T^k F)_x = \tilde f$,
whence $(J^k F)_x = f$. Therefore, the germ $[F] \in {\cal G}_x(M,N)_y$
satisfies $\ev_x^k([F])=f$. 
\qed

\nin In particular, for $k \in \N$, the evaluation maps from germs of curves to $J^k M$,
$$
\eqalign{
\ev_0^k:{\cal G}_0(\K,M)_x & \to J_0^k(\K,M)_x \cong (J^k M)_x, \cr
[\alpha] & \mapsto (J^k \alpha)_0 \cong (J^k \alpha)_0 (\delta 1) \cr}
$$
are surjective.

\msk \nin {\bf 32.5.}
{\sl From jets to germs.} 
The problem of ``integrating jets to germs'' amounts to construct some sort of
inverse of the ``infinite order
evaluation map'' $\ev_x^\infty: {\cal G}_x(M,N)_y  \to J^\infty_x(M,N)_y$.
However, in general this map is not invertible -- there are obstructions both
to surjectivity and to injectivity.

\ssk
(a) {\sl Failure of surjectivitiy. Convergence.}
Let us call the image of $\ev_x^\infty$ in $J^\infty_x(M,N)_y$
the {\it space of integrable jets}. 
For a general topological field or ring $\K$ it seems impossible to give 
some explicit description of this space. But if, e.g., 
$\K$ is a complete valued field, then we have the notion of {\it convergent
power series} (see, e.g., [Se65]),
 and we may define the notion of {\it convergent jets}.
Then the convergent jets are integrable.
For instance, if $G$ is a Lie group, then the sequence of exponential 
maps $\exp_k:(J^k \g)_0 \to (J^k G)_e$ defines a formal jet
$\exp_\infty \in J_0^\infty(\g,G)_e$ with inverse $\log_\infty$.
 It is of basic interest to know
wether these jets are integrable or not. Our construction of $\exp_k$ and
$\log_k$
is explicit enough (Theorem 25.2) to provide estimates assuring
convergence in suitable situations. 

\ssk
(b) {\sl Failure of injectivity. Analyticity.}
The kernel of the evaluation map $\ev_x^\infty$ in $J^\infty_x(M,N)_y$
is the space of germs of 
 smooth maps vanishing of infinite order at $x$ and hence
measures, in some sense, the ``flabbiness'' of our manifolds.
In order to control this kernel one needs some sort of rigidity,
e.g., the assumption that our manifolds are analytic (where again we
need additional structures on $\K$ such as above).

\ssk
Summing up, under the above mentioned assumptions, the 
construction of an inverse of 
$\ev_x^\infty$ would look as follows: given an infinite jet $f \in
J^\infty_x(M,N)_y$,
we pick a chart around $x$ and write $f = \sum_{j=0}^\infty \delta^{(j)} d^j$
(i.e., $f=\lim_{\leftarrow} f_k$, $f_k =\sum_{j=0}^k \delta^{(j)} d^j$),
then, if this jet is convergent, we let in a chart
$$
F(x+h)= \sum_{j=0}^\infty {1 \over j!} d^j (h,\ldots,h),
$$
which is ``well-defined up to terms vanishing of infinite order at $x$''.

\msk
\nin {\bf 32.6.} {\sl Analytic structures on Lie groups and symmetric spaces.}
It is a classical result that finite-dimensional real Lie groups admit a 
compatible analytic Lie group structure. Recently, this result has been extended
to the case of Lie groups over ultrametric fields by H.\ Gl\"ockner [Gl03c].
However, the proof uses quite different methods than in the real case.
It seems possible that, with the methods presented in this work, a common
proof for both cases could be given, including also the case of symmetric spaces.
The first step should be to prove that, in finite dimension over a complete valued field,
 the exponential jet $\exp_\infty$ of a Lie
group or a symmetric space (Chapter 30)
is in fact convergent and defines the germ
of a diffeomorphism (with inverse coming from $\log_\infty$);
 then to define, with respect
to some fixed chart, an exponential map at the origin by this convergent jet;
to use left translations to define an atlas of $G$ and finally to show that this
atlas is analytic since transition functions are given by the Campbell-Hausdorff
series (limit of Theorem 25.5).

\vfill \eject  

\def\title{LIMITATIONS} 
\def\rightheadline{\hfil{\eightrm \title}\hfil\tenrm\folio}

\sectionheadline
{Appendix L. Limitations}

\nin
The prize for generalizing a theory like 
differential geometry is of course a limitation of the number
of properties and results that will survive in the more
general framework. Thus,
in our general category of manifolds over topological fields or
rings $\K$, some
constructions which the reader will be used to from
the theory of real manifolds
(finite-dimensional, Banach, locally convex or other)
are not possible. In particular, this concerns {\it dual}
and {\it tensor product constructions} which would permit to construct
new vector bundles from old ones.
In this appendix we give a short overview over the kinds of limitations
that we have to accept when working in the most general category. 

\msk \nin {\bf L.1.} {\sl Cotangent bundles and dual bundles.}
Formally, we can define a cotangent bundle $T^*M$ 
following the pattern described in  Chapter 3,
just by taking the topological dual of the tangent space
$T_x M$ as fiber. However, in general there is no reasonable
topology on topological duals which could serve to turn
$T^* M$ into a smooth manifold, and similarly for 
general dual bundles of vector bundles over $M$.
Moreover, in order to define a manifold structure on $T^*M$ it
is not enough just to have the structure of a topological
vector space on the dual of the model space --
see [Gl03b] for a counterexample. It seems that
the most natural framework where dual bundles always exist is
the category of (real, complex or ultrametric) {\it Banach
manifolds}, i.e.\ all model spaces are Banach spaces
(cf.\ [La99], [Bou67]).

\msk \nin {\bf L.2.} {\sl Tensor products.}
In our approach, the use of multilinear bundles allows us to
develop the theory of connections and related notions without
having to introduce tensor products.\footnote{$^1$}{We agree with the remark
 in the introduction
of [La99] (loc.\ cit.\ p.\ vi): ``An abuse of multilinear algebra in standard
treatises arises from an unnecessary double dualization and
an abusive use of the tensor product."}
However, in specific situations one may wish to have a tensor
product of vector bundles at disposition. In fact, it is possible to
define {\it topological tensor products} in the category
of topological $\K$-modules be inquiring the usual universal
property for {\it continuous} bilinear maps
(generalizing the so-called {\it projective topology} from the
real locally convex case, cf. [Tr67, Ch.\ 43]), and doing this
construction fiberwise, we get an abstract bundle over the base $M$
whose fibers are topological vector spaces.
But in general
it is not possible to turn this bundle into a manifold
(see [Gl03b] for a discussion of the real locally convex case).
Moreover, the topological tensor product of topological modules
 does not have the good properties known from the algebraic
tensor product: tensor products in this sense are in general {\it not}
associative ([Gl04b]); this fails already for (non-locally convex)
real topological vector spaces.
For all of these reasons, we will never use topological tensor
products in the present work.
As for dual bundles, the most natural category where tensor products
of bundles make no problem seems to be the category of Banach manifolds
(cf.\ [La99]).

\msk \nin
{\bf L.3.} {\sl The main categories.}
Here is a short (and by no means complete)
list of the main categories of  manifolds, resp.\ of
 their model spaces. We indicate what
constructions and results are still available,
starting from general categories and ending with most specific ones.
Details can be found in the quoted references.
For more information on related topics we refer, in particular, to work of H.
Gl\"ockner (cf.\ References).

\bigskip
\item{(1)}
{\sl The case where  $\K$ is a field.}
Locality -- see Section 2.4.

\msk
\item{(2)}
{\sl
Metrizable topological vector spaces over $\K=\R$ or an ultrametric field.}
Smoothness can be tested by composition with smooth maps
$\K^k \to V$ from finite dimensional spaces into the model space
[BGN04, Theorem 12.4].

\msk
\item{(3)}
{\sl Locally convex vector spaces over $\K=\R,\C$.}
Manifolds in our sense are the same as those in the sense of
Michal-Bastiani [BGN04, Prop. 7.4].

\itemitem{--} Abundance of smooth functions (at least locally): 
the representation of the space
of vector fields by differential operators, 
$$
\X(U) \to
\End_\K({\cal C}^\infty(U)), \quad  X \mapsto L_X,
$$
 is injective for
chart domains $U \subset M$; thus the Lie
bracket can be defined via $L_{[X,Y]}=-[L_X,L_Y]$ (or taking the
opposite sign).
\itemitem{--} Solutions of the differential equation $df=0$
are locally constant:
 uniqueness results for  flows  (if they exist), etc.
\itemitem{--} Fundamental theorem of differential and integral
calculus for curves: Poinca\-r\'e Lemma (cf.\ [Ne02b, Lemma II.3.5]).

\msk
\item{(4)}
{\sl Complete locally convex vector
spaces over $\K=\R,\C$.} One can speak about analytic maps;
analytic manifolds can be defined, and they are smooth in our sense
[BGN04, Prop. 7.20].
For general complete ultrametric fields this holds e.g.\ in the category of
Banach manifolds [Bou67].

\msk
\item{(5)}
{\sl Banach manifolds over $\K=\R,\C$ or
a complete ultrametric field.}

\itemitem{--} Cotangent, dual 
 and hom-bundles exist in the given category. 
\itemitem{--} Tensor products of vector bundles exist and have the
usual properties.
\itemitem{--} The group $\Gl(V)$ of continuous linear bijections of
 the model space is a Lie
group; hence the frame bundle of a manifold exists, and
natural bundles are associated bundles of the frame bundle
(in the sense of [KMS93]). 
\itemitem{--} Inverse function and implicit function theorems
(cf.\ [La99] in the real case and [Gl03a] for the ultrametric case,
including much more general results); 
their consequences for the theory of submanifolds

\msk
\item{(6)}
{\sl Banach manifolds over $\K=\R$.}
Existence and uniquenes theorem for ordinary differential
equations and all of its consequences:
 existence of flows, geodesics, parallel transports;
Frobenius theorem;
exponential map: existence and local diffeomorphism.
(See [La99] for all this.)
Connected symmetric spaces are homogeneous (cf. [BeNe05]).

\msk
\item{(7)}
{\sl Finite dimensional  manifolds over locally compact fields.}

\itemitem{--} Existence of compact objects.
\itemitem{--} Existence of volume forms (cf. [Bou67]),
Haar measure on Lie groups.

\msk
\item{(8)}
{\sl Finite dimensional real manifolds.}

\itemitem{--} Existence of partitions of unity: global existence
of metrics and connections and all consequences.
\itemitem{--} Equivalence of all definitions of linear differential
operators mentioned in Chapter 21:
``algebraic definition" of tangent spaces and of vector fields; tensor fields
are the same as differential operators of degree zero.

\vfill \eject

\def\title{GENERALIZATIONS} 
\def\rightheadline{\hfil{\eightrm \title}\hfil\tenrm\folio}

\sectionheadline
{Appendix G. Generalizations}

\nin 
Generalizations of the present work are possible in two directions:
one could treat more general classes of manifolds, or one
could leave the framework of manifolds and use more formal concepts.
In either case, it seems that the main features of our approach
-- use of scalar extensions and use of ``geometric multilinear algebra'' --
could carry over and play a useful r\^ole.

\msk \nin
{\bf G.1.} {\sl General $C^0$-concepts.}
As to generalizations in the framework of manifolds, we just discuss
here the manifolds associated to a general {\it $C^0$-concept}
in the sense of [BGN04].
Differential calculus, as explained in Chapter 1, can be generalized
as follows:
$\K$ is a field or unital ring with topology and $\cal M$ a class
of $\K$-modules with topology; for each open set $U$ and each
$W \in {\cal M}$, a set $C^0(U,W)$ of continuous mappings is given,
such that several natural conditions are satisfied, among them the
{\it determination principle} (DP) from Section 1.1 which makes
a differential calculus possible.
We require also that there
is a rule defining a topology on the direct product of spaces
from the class $\cal M$ (and making them a member of $M$) such that
the ring operations (addition and multiplication)
are of class $C^0$. A typical example is the class of finite dimensional
vector spaces over an infinite field $\K$, equipped with their
Zariski-topology, and $C^0$ meaning ``rational":
in this case, $\K$ is not a topological field
in the usual sense, and the members of $\cal M$ are not topological
vector spaces in the usual sense. Nevertheless, differential calculus,
the theory of manifolds and differential geometry
can be developed almost word by word as in the present work.
This is true, in particular, for Theorem 6.2 and its proof.
The only difference is that direct products of manifolds now
carry a topology which, in general, is finer than the product
topology; this is true also for local product situations:
bundles are locally direct products and hence their topology, even
locally, is in general finer then the local product topologies.
We leave the details to the interested reader.

\msk \nin
{\bf G.2.} {\sl Synthetic differential geometry.}
Among the more formal approaches to differential geometry we would
like to mention the so-called ``synthetic differential geometry''
(cf. [Ko81] or [Lav87]) and  the theory of its models such as
the so-called ``smooth toposes'' (over $\R$) (cf.\ [MR91]). 
 In [MR91, p.\ 1--3], the authors give three main reasons 
for generalizing the theory of manifolds by the theory of smooth toposes:

\ssk
\item{(1)} the category of smooth manifolds is not cartesian closed
(spaces of mappings between manifolds are not always manifolds),
\item{(2)} the lack of finite inverse limits in the category of manifolds
(in particular, manifolds can not have ``singularities''),
\item{(3)} the absence of a convenient language to deal explicitly and
directly with structures in the ``infinitely small".

\ssk \nin
The ``naive'' approach by R. Lavendhomme [Lav87] focuses mainly on Item (3),
whereas A. Kock in [Ko81] emphasizes  Item (1), quoting F.\ W.\ Lawvere
(loc.\ cit.\ p.\ 1): ``...it is necessary that the mathematical world picture
involve a cartesian closed category $\cal E$ of smooth morphisms between
smooth spaces.'' Of course, this  statement is more philosophical than
mathematical in nature; but we feel that our approach proves at least
that Items (1) and (3) are not as strongly coupled among each other as
one used to believe up to now -- the infinitely small has a very natural
and useful place in our approach, and this works in a category
based on conventional logic and usual set-theory
 and not needing any sophisticated model theory.\footnote{${}^1$}{cf.\
 [Lav87, p. 2]: ``...nous nous sommes strictement
limit\'es au point de vue ``na\"\i f'' ... sans tenter d'en d\'ecrire
des ``mod\`eles'' qui eux vivent dans des topos hors d'atteinte \`a
ce niveau \'el\'ementaire.''}
Nevertheless, it is certainly still interesting to generalize results
to categories that behave well with respect to (1) and (2) --
here it seems possible, in principle, to develop a theory of ``smooth
toposes over $\K$'' (based on differential calculus over $\K$ as used
here) and to transpose much of the material presented here
to this context which then would unify differential geometry
and algebraic geometry over topological base fields and -rings.

\msk \nin
{\bf G.3.} {\sl Other categorial approaches.}
Just as the concept of a polynomial is a formal generalization
of the more primitive concept of a polynomial mapping (see
[Ro63] for the general case; cf. also [BGN04, Appendix A]),
it seems possible to generalize the differential calculus
from [BGN04] in a formal way such that it becomes meaningful
for any commutative ring or field, including finite ones.
A ``categorial'' generalisation of ordinary differential
calculus  has already been given
by L.D. Nel
[N88]; however, instead of taking Condition (0.2) as starting
point, he works with the condition
$f(y)-f(x)=\Phi(x,y)\cdot (y-x)$ which 
is well-adapted for differential calculus in {\it one} variable,
but causes much trouble
in the further development of differential calculus in several
variables. 
Since the arguments proving our Theorem 6.2 are of purely
categorial nature, we believe that the scalar extension functor
should exist and play a central r\^ole in a categorial theory
based on (0.2).

\vfill \eject 

%% file: Dgeo5.tex

\def\date{25.10.2006}
\def\title{APPENDIX: MULTILINEAR GEOMETRY} 
\def\rightheadline{\hfil{\eightrm \title}\hfil\tenrm\folio}

\sectionheadline{Appendix: Multilinear Geometry}

This appendix is completely independent of the main text and 
contains all algebraic definitions and results that are used there.
Chapter BA on bilinear algebra is the special case $k=2$ of Chapter MA
on multilinear algebra,
and the reader who is interested in general results may directly
start by looking at Chapter MA. However, the case $k=2$ is of independent
interest since it corresponds to {\it linear} connections (Chapter 10).
-- 
In the following, $\K$ is any commutative ring with unit $1$.

\subheadline {BA. Bilinear algebra}

\nin {\bf BA.1.} {\sl Bilinearly related structures.}
Assume 
 $V_1,V_2,V_3$  are three $\K$-modules together with a bilinear map
 $b: V_1 \times V_2 \to V_3$. 
We let $E:=V_1 \times V_2 \times V_3$. Then
there is a structure of a $\K$-module on $E$,
depending on $b$, and given by
$$
\eqalign{
(u,v,w)+_b (u',v',w') & := (u+u',v+v',w+w'+b(u,v')+b(u',v)),
\cr
\lambda_b (u,v,w) & := (\lambda u,\lambda v, \lambda w + \lambda (\lambda-1)
 b(u,v)). \cr}
\eqno {\rm (BA.1)}
$$
One may either check directly the axioms of a $\K$-module, or
use the following observation: the map
$$
f_b:E \to E, \quad
(u,v,w) \mapsto (u,v,w+b(u,v))
\eqno {\rm (BA.2)}
$$
is a bijection  whose inverse is given by
 $f_b^{-1}(u,v,z)=(u,v,z-b(u,v))$.
In fact, it is immediately seen that
$$
f_b \circ f_c = f_{b+c}, 
\quad f_0 = \id_E,
$$
i.e.\ the additive group ${\rm Bil}(V_1 \times V_2,V_3)$ of bilinear maps
from $V_1\times V_2$ to $V_3$ acts on $E$.
This action is not linear, and a straighforward calculation
shows that the linear structure
$$
L^b := (E, +_b, \cdot_b)
\eqno {\rm (BA.3)}
$$
is simply the push-forward of the original structure (corresponding to $b=0$)
via $f_b$, i.e.
$$
\lambda_b x = f_b(\lambda \, f_b^{-1}(x)), \quad
x +_b y = f_b(f_b^{-1}(x) + f_b^{-1} (y)),
$$
where $+$ means $+_0$.
We call two linear structures $L^b,L^c$ with $b,c \in {\rm Bil}(V_1 \times V_2,
V_3)$ {\it bilinearly related}, and
we denote by $\brs (V_1,V_2;V_3)$ the space of all $\K$-module structures
on $E$ that are bilinearly related to the original structure. 
From Formula (BA.1) we get the following expressions for the
negative, the difference and the barycenters with respect to $L^b$:
$$
\eqalign{
&-_b(u,v,w)  =(-u,-v,-w+2 b(u,v)), \cr
& (u,v,w)-_b (u',v',w')  = 
(u-u',v-v',w-w'+ 2b(u',v')-b(u,v')-b(u',v)), \cr
&(1-r)_b (x,y,z) +_b r_b(x',y',z')  = \cr
& \quad \quad ((1-r)x+rx',(1-r)y+ry' , 
(1-r)z+rz'+r(r-1) b(x-x',y-y')) .\cr}
\eqno {\rm (BA.4)}
$$
Translations and Dilations now depend on $b$. We will
take account of this by adding $b$ as a subscript in the
notation: for $x=(x_1,x_2,x_3) \in E$ and $b \in \Bil(V_1 \times V_2,V_3)$,
the {\it translation by $x$ with respect to $b$} is
$$
\tau_{b,x}:E   \to E, \quad 
 (u,v,w)  \mapsto (x_1,x_2,x_3)+_b(u,v,w) , 
\eqno {\rm (BA.5)}
$$
and the {\it dilation with ratio $r \in \K$ with respect to $b$} is
$$
r_b:E \to E, \quad
x \mapsto r_b x = (rx_1,rx_2,rx_3 + (r^2-r) b(x_1,x_2)).
\eqno {\rm (BA.6)}
$$

\msk \nin {\bf BA.2.} {\sl Intrinsic objects of bilinear geometry.} We have 
remarked above that the set of transformations of $E$ given by
$$
\Gm^{1,2}(E) := \{ f_b | \, \,  b \in \Bil(V_1 \times V_2,V_3) \}
\eqno {\rm (BA.7)}
$$
is a vector group, isomorphic to the $\K$-module $\Bil(V_1 \times V_2,V_3)$.
It will be called the {\it special bilinear group}. As already remarked,
 it acts
simply transitively on the set $\brs(V_1,V_2;V_3)$  of bilinearly related
structures which in this way
is turned into an affine space over $\K$.
The ``original'' linear structure (which corresponds to $b=0$)
shall not be preferred to the others.
By {\it bilinear geometry} we mean the study of all geometric
and algebraic properties of $E$ that
are common for all bilinearly related linear structures, or,
equivalently, that are invariant under the special bilinear group.
Slightly more generally, we may consider also all properties that are invariant
under the {\it general bilinear group}
$$
\Gm^{0,2}(E) := \{ f = (h_1,h_2,h_3) \circ f_b | \, \, h_i \in \Gl(V_i),
\, i=1,2,3, \, \, f_b \in \Gm^{1,2}(E) \},
\eqno {\rm (BA.8)}
$$
i.e.\ under all bijections of $E$ having the form 
$$
f(u,v,w)=(h_1(u),h_2(v),h_3(w+ b(u,v)))
$$
with $b$ as above and invertible linear maps $h_i$. In fact, the group
$\Gm^{0,2}(E)$ also acts by push-forward
on the space of bilinearly related structures.
The stabilizer of the base structure $L^0$ is the group
$H:=\Gl(V_1)\times \Gl(V_2) \times \Gl(V_3)$. The group $H$ acts affinely
on the affine space $\brs(V_1,V_2;V_3)$ and hence acts linearly with respect
to zero vector $L^0$ in this space. The situation is the same as for
usual affine groups, and hence we have an exact splitting sequence of groups
$$
\matrix{
0 & \to & \Gm^{1,2}(E) & \to & \Gm^{0,2}(E) & \to & \Gl(V_1)\times \Gl(V_2)
\times \Gl(V_3) & \to & 1 \cr}.
\eqno{\rm (BA.9)}
$$ 
Now we say that an object 
on $E$ is {\it intrinsic} if its definition
is independent of $b$. Here are some examples:

\ssk
{\sl (A) Injections, axes.} The origin $(0,0,0) \in E$
is invariant under $\Gm^{1,2}(E)$, and so are the three ``axes''
$V_1 \times 0 \times 0$, $0 \times V_2 \times 0$, $0 \times 0 \times V_3$
(referred to as ``first, second, third axis"; they are submodules
isomorphic to, and often identified with, $V_1$, $V_2$, resp. $V_3$).
Also, the submodules $V_1 \times 0 \times V_3$ and $0 \times V_2 \times V_3$
together with their
inherited module structure are intrinsic, but the ``submodule
$V_1 \times V_2$" (which should correspond to
$V_1\oplus_b V_2 \oplus_b 0$) is not intrinsic: in fact,
$$
f_b(V_1 \times V_2 \times 0)=
V_1\oplus_b V_2 \oplus_b 0 = \{ (u,v,b(u,v)) | \, u \in V_1, v \in V_2 \}
\eqno {\rm (BA.10)}
$$
clearly depends on $b$ (it is the graph of $b$).
Thus we have the following intrinsic diagram of inclusion of ``axes"
and ``planes'':
$$
\matrix{
V_1 & \to & E^{13}:=V_1 \times V_3 & &  \cr
  & \nearrow & & \searrow & \cr
V_3 & & & & E=V_1 \times V_2 \times V_3 \cr
 & \searrow & & \nearrow &  \cr
V_2 & \to & E^{23}:=V_2 \times V_3 & & \cr}
$$

\ssk
{\sl (B) Projections, fibrations.}
It follows immediately from (BA.1) that the two projections
$$
\pr_j: V_1 \times V_2 \times V_3 \to V_j
\eqno {\rm (BA.11)}
$$
for $j=1,2$ are invariant under $f_b$ in the sense that
$\pr_j \circ f_b = \pr_j$, and hence they are
linear for every choice of linear structure $L^b$.
 (This is not true for the third
projection!)  Thus also
$$
\pr_{12}:=\pr_1 \oplus \pr_2: V_1 \times V_2 \times V_3 \to
V_1 \times V_2, \quad (u,v,w) \mapsto (u,v)
\eqno {\rm (BA.12)}
$$
is an intrinsic linear map, and the following three exact
sequences are intrinsically defined:
$$
\matrix{
0 & \to & V_3 & \to & V_1 \times V_2 \times V_3 & {\buildrel \pr_{12}
\over \longrightarrow} & V_1 \times V_2 & \to 0 
\cr
0 & \to &  V_2 \times V_3 & \to & V_1 \times V_2 \times V_3 
& {\buildrel \pr_{1} \over \longrightarrow} & V_1 & \to 0 
\cr
0 & \to & V_1 \times V_3 & \to & V_1 \times V_2 \times V_3 & {\buildrel \pr_{2}
\over \longrightarrow} & V_2 & \to 0 
\cr}
\eqno {\rm (BA.13)}
$$
The fibers of $\pr_j$ are invariant under $\Gm^{1,2}(E)$, and moreover
they are affine subspaces of $E$, for any choice of linear structure $L^b$.
But more is true: the induced affine structure on the fibers
 does not depend on $L^b$ at all.

\Proposition BA.3. For all $z \in E$ , the spaces
$$
E_z^{13}: = z +_b E^{13}, \quad
E_z^{23}: = z +_b E^{23}
$$
are intrinsic affine subspaces in the sense that both the sets and
the structures of  an affine space over $\K$ induced from the ambient
space $(E,L^b)$ are independent of $b$.

\Proof. Let $b$ be fixed. The map $f_b$ preserves the subset
$$
z +_0 E^{13} = \{ (z_1+v,z_2,z_3+w) | \, v \in V_1, w \in V_3 \}
$$
on which it acts {\it linearly}: with respect to the origin
$(0,z_2,0)$ it is described, in the obvious matrix notation, by
$$
\pmatrix{\1 & b(\cdot,z_2) \cr 0 & \1 \cr}.
$$
Therefore the set and the linear structure on this set are independent of $b$.
The claim for $E_z^{13}$ now
follows since $E_z^{13}$ and $z +_0 E^{13}$ are the same sets,
and similarly for $E_z^{23}$.
\qed

\nin The proposition shows that $E$ has an
{\it intrinsic double fibration} by affine spaces.
Every fiber contains one distinguished element (intersection
with an axis), and taking this element as origin, the fiber
has a canonical linear structure (e.g., the origin
in $V_{(u,v,w)}^{13}$ is $(0,v,0)$, etc.) This implies, for instance, that
the  maps $E \to E$, $(u,v,w) \mapsto (\lambda u,v,\lambda w)$ and
$(u,v,w) \mapsto (u,\lambda v,\lambda w)$ have an intrinsic interpretation,
and so has their composition $(u,v,w) \mapsto (\lambda u,\lambda v,
\lambda^2 w)$.

\ssk {\sl (C) Parallelism.}
The notion of being ``parallel to an axis'' is not always intrinsic:
it is so for the parallels to the third axis, that is, fibers of $\pr_{12}$,
which we call {\it vertical spaces}, written, for $z=(z_1,z_2,z_3) \in E$,
$$
E_z^3:= E_z^{13} \cap E_z^{23} = z +_b  V_3 =
\{ (z_1,z_2,z_3+w) | \, w \in V_3 \}.
\eqno {\rm (BA.14)}
$$
But the parallel to, say, $0 \times V_2 \times 0$ through
$(u,0,0)$ is the ``horizontal space''
$$
H_{u,b}:=(u,0,0)+_b V_2 = 
\{ (u,0,0) +_b (0,v,0)  = (u,v,b(u,v)) | \, v \in V_2 \}
\eqno {\rm (BA.15)}
$$
which clearly depends on $b$.

\bigskip \centerline{\bf Special cases}

\msk \nin {\bf BA.4.} {\sl
The special case $V_2=V_3$ and shifts.} In the setup of Section BA.1,
consider the special case $V_2=V_3$. Then a direct computation shows that
the {\it shift} 
$$
S = \pmatrix{1 & & \cr & 0 & 0 \cr  & 1 & 0 \cr}:
V_1 \times V_2 \times V_3 \to V_1 \times V_2 \times V_3, \quad
(u,v,w) \mapsto (u,0,v)
\eqno {\rm (BA.16)}
$$
commutes with all $f_b \in \Gm^{1,2}(E)$, and hence $S$
is an ``intrinsic linear map'', i.e.\ it is linear with respect
to all $L^b$.
It stabilizes the fibers $E_z^{23}$, and the restriction
to $E_z^{23}$, with respect to the origin $(\pr_1(z),0,0)$, is
a two-step  nilpotent linear map given in matrix form by
the lower right corner of the matrix $S$. It defines
on $E_z^{23}$  an
intrinsic structure of a module over the dual numbers $\K[\epsilon]$.
If $V_1 = V_3$, we get in a similar way a shift $S'$.


\msk \nin {\bf BA.5.} {\sl
The special case $V_1=V_2$ and torsion.}
In case $V_1=V_2$, 
the symmetric group $\Sigma_2 = \{ \id, (12) \}$ acts on
$E=V_1 \times V_1 \times V_3$, on $\Bil (V_1 \times V_1,V_3)$ and on the
space of bilinearly related linear structures, where 
the transposition $(12)$ acts
by the {\it exchange map} or {\it flip}
$$
\kappa:V_1 \times V_1 \times V_3 \to V_1 \times V_1 \times V_3, \quad
(u,v,w) \mapsto (v,u,w).
$$ 
We say that $L^b$ is {\it torsionfree} if $b$ is symmetric;
this means that the exchange map is a linear automorphism
with respect to $L^b$. 
In general, the push-forward of $L^b$ by $\kappa$ is given by
$\kappa \cdot L^b = L^{\kappa \cdot b}$, where
$$
\kappa \cdot b(u,v)=b(v,u) .
$$
The difference $L^b-L^c$ in the affine space $\brs(V_1,V_2;V_3)$
corresponds to the difference $b-c$; thus the
 difference $L^b - \kappa \cdot L^b$  corresponds
to the skew-symmetric bilinear map $t(u,v):=b(u,v)-b(v,u)$,
called the {\it torsion of the linear structure $L^b$}.

\msk \nin {\bf BA.6.} {\sl Quadratic spaces.}
Assume $V_1=V_2$ and $b$ is symmetric. Then
 the $\kappa$-fixed space
$\{ (v_1,v_1,v_3) \vert \, v_1 \in V_1,v_3\in V_3\}$ is a submodule of
$(E,L^b)$, isomorphic to $V_1 \times V_3$ with $\K$-module structure
$$
\eqalign{
(v_1,v_3)+_b(v_1',v_3') & =(v_1+v_1',v_3+v_3'+b(v_1,v_1') +b(v_1',v_1)) \cr
\lambda(v_1,v_3) & =(\lambda v_1,\lambda v_3 + \lambda(\lambda-1) b(v_1,v_1)).
\cr}
$$
Note that this structure is already determined by the quadratic map
$q(v)=b(v,v)$, even in case $2$ is not invertible in $\K$, and thus
for any quadratic map $q:V_1 \to V_3$ we may define a linear structure
on $V_1 \times V_3$ given by
$$
\eqalign{
(v_1,v_3)+_q(v_1',v_3') & =(v_1+v_1',v_2+v_2'+q(v_1+v_1') - q(v_1') -
q(v_1)) \cr
\lambda(v_1,v_3) & =(\lambda v_1,\lambda v_3 + \lambda(\lambda-1) q(v_1)),
\cr}
$$
leading to the notion of a {\it quadratic space} (we will not 
give here an axiomatic definition).
The third axis still corresponds to an intrinsic subspace in $V_1 \times
V_3$, but the first and second axis do not give rise to intrinsic subspaces:
they are replaced by the {\it diagonal imbedding}
$$
V_1 \to V_1 \times V_1 \times V_3, \quad v \mapsto 
(v,0,0)+_b(0,v,0) =(v,v,q(v))
$$
which clearly depends on $q$ and hence on $b$. 


\subheadline{Intrinsically bilinear maps}

\Proposition BA.7. Let $E,V_i$, $i=1,2,3$ be as above and 
$W$ be a linear space, i.e.\ a $\K$-module. For a map $\omega:E \to W$ 
the following properties are equivalent:
\item{(1)} There exists a linear map $\lambda:V_3 \to W$ and a 
$\K$-bilinear map $\nu:V_1 \times V_2 \to W$ such that
$$
\omega(u,v,w)=\nu(u,v) + \lambda(w).
$$
\item{(2)} For all $z  \in E$, the restrictions
$$
\omega_z^1:=\omega\vert_{E_z^{13}}: E_z^{13} \to W, \quad 
\omega_z^2:=\omega\vert_{E_z^{23}}: E_z^{23} \to W
$$
to the fibers of the canonical double fibration
are linear. 

\Proof.
Assume (1) holds. Let $z=(0,v,0)$; then
$\omega_z^1(u,w)=\nu(u,v)+\lambda(w)$ is linear in $(u,w)$.
Conversely, assume $\omega$ satisfies (2).
We let $\lambda:=\omega\vert_{V_3}$ (restriction to third axis)
and $\nu:= \omega \circ \iota_{12}$, where
$\iota_{12}:V_1 \times V_2 \to E$ depends on the linear structure
$L^0$. Then, by (2), $\nu(u,v)=\omega(u,v,0)$ is linear in $u$  and
in $v$, i.e.\ it is $\K$-bilinear. 
We decompose 
$$
(u,v,w)=((0,0,w)+_{(0,0,0)}(0,v,0))+_{(0,v,0)} (u,v,0),
$$
where the first sum is taken in the fiber $E_0^{13}$ and the second
sum in the fiber $E_z^{23}$ with $z=(0,v,0)$. 
Using Property (2), we get
$$
\eqalign{
\omega(u,v,w) & =\omega((0,0,w)+_{(0,0,0)}(0,v,0)) + \omega(u,v,0) \cr
& =
\omega(0,0,w) + \omega(0,v,0) + \omega(u,v,0) \cr
& =\lambda(w) + \nu(v,0)+\nu(u,v)=\lambda(w) + \nu(u,v) \cr}
$$
as had to be shown.
\qed

\ssk
\nin If $\omega$ satisfies the equivalent conditions of the preceding
proposition, we say that $\omega$ is {\it (intrinsically) bilinear.}
For instance, the third projection with respect to a fixed linear structure,
$\pr_3^b:E \to V_3$, is intrinsically bilinear (it corresponds to $W=V_3$,
$\nu =0$, $\lambda = \id_{V_3}$).

\Corollary BA.8. In the situation of the preceding proposition are
equivalent:
\item{(1)} $\lambda = 0$.
\item{(2)} For all $z\in E$,
the restriction of $\omega$ to the third axis is constant,
$$
\forall w \in V_3 : \quad \omega(z_1,z_2,z_3+w)=\omega(z_1,z_2,z_3).
\qeddis

\nin If $\omega$ satisfies the properties of the corollary,
then we say that $\omega$ is {\it homogeneous bilinear}, and we let
${\cal M}_h(E,W)$ the space of homogeneous bilinear maps $E \to W$.
The preceding statements show that
$$
\Bil(V_1 \times V_2;W) \to {\cal M}_h(E,W), \quad
\nu \mapsto \omega:=\nu \circ \pr_{12}
$$
is a bijection. Its  inverse $\omega \mapsto \omega \circ \iota_{12}$ 
is not canonical since the injection $\iota_{12}$ depends on the choice
of a  linear structure on $E$.

\msk \nin
{\bf BA.9.} {\sl The skew-symmetrization operator.}
Assume we are in the special case
 $V_1=V_2$. Then to every intrinsically bilinear map
$\omega:E \to W$ we can associate
a {\it homogeneous} intrinsically bilinear map
$$
\alt \omega := \omega - \omega \circ \kappa
$$
In fact, $\omega$ and $\omega\circ \kappa$ are intrinsically bilinear, hence so
is $\omega-\omega \circ \kappa$, and $\alt \omega$  is homogeneous: writing
$\omega$ in the form  $\omega(u,v,w)=\nu(u,v)+\lambda(w)$ as in the preceding
paragraph, we get
$$
\alt \omega (u,v,w)=\omega(u,v,w) - \omega(v,u,w)=\nu(u,v) - \nu(v,u).
$$
which clearly is homogeneous bilinear.
The bilinear map $V_1 \times V_1 \to V_3$ corresponding to $\alt \omega$ is
$$
\nu_{\alt \omega} (u,v) = \omega(u,v,w) - \omega(v,u,w).
$$
For instance, if $\omega=\pr_3^b$, then $\alt \omega$ is the
torsion of $L^b$.

\bigskip
\centerline{\bf Relation with the tensor product}

\msk \nin {\bf BA.10.} {\sl The tensor space.}
It is possible to translate everything we have said so far
into usual {\it linear} algebra by using {\it tensor
products.}
For arbitrary $\K$-modules $V_1,V_2,V_3$, we define the {\it tensor space}
$Z:= (V_1 \otimes V_2) \oplus V_3$ and the {\it big tensor space}
$\hat Z:=V_1 \oplus V_2 \oplus V_3 \oplus (V_1 \otimes V_2)$
together with the canonical maps
$$
\eqalign{
V_1 \times V_2 \times V_3 \to Z, & \quad (v_1,v_2,v_3) \mapsto v_1 \otimes v_2
+ v_3, \cr
V_1 \times V_2 \times V_3 \to \hat Z, &  \quad
(v_1,v_2,v_3) \mapsto v_1+v_2+v_3+ v_1 \otimes v_2.  \cr}
\eqno {\rm (BA.17)}
$$
We claim that the linear structure
on $Z$, resp.\ on $\hat Z$
 is intrinsic, i.e. independent of the linear structure
$L^b$ on $V_1 \times V_2 \times V_3$: in fact,
let $B: V_1 \otimes V_2 \to V_3$ be the linear map corresponding to
$b:V_1 \times V_2 \to V_3$ and let
$$
F_b:=\pmatrix{ \1 & 0 \cr B & \1 \cr}: Z \to Z.
$$
Then the diagram
$$
\matrix{(V_1 \times V_2) \times V_3 & \to & Z \cr
f_b \downarrow & & \phantom{F_b} \downarrow F_b \cr
(V_1 \times V_2) \times V_3 & \to & Z \cr}
\eqno {\rm (BA.18)}
$$
commutes.
Since $F_b$ is a linear automorphism of $Z$, it follows that
the linear structure on $Z$ does not depend on $b$, and similarly for
$\hat Z$. Thus $Z$ and $\hat Z$
are a sort of universal linear space for all
bilinearly equivalent structures.
(The space $\hat Z$ is faithful in the sense that different 
elements of $\Gm^{0,2}(E)$ induce different automorphisms
of $\hat Z$, whereas for $Z$ this property holds only with respect to the
smaller group $\Gm^{1,2}(E)$.)
Note that the splitting of $Z$ into a direct sum is not
intrinsic, but the subspace $0 \oplus V_3$ is intrinsic:
the sequence
$$
0 \to V_3 \to Z \to V_1 \otimes V_2 \to 0
\eqno {\rm (BA.19)}
$$
is intrinsic. (If $V_1=V_2$, then, with respect to
 symmetric $b$'s, then we may define $Z':=S^2 V_1 \oplus V_3$ and
 we get an intrinsic sequence $V_3 \to Z' \to S^2 V_1$.)
 As we have seen, in presence of a fixed linear structure
$L^b$, this sequence splits. Conversely, a splitting
$\beta:V_1 \otimes V_2 \to Z$ of this sequence induces a linear map
$B:=\pr_{V_3} \circ \beta: V_1 \otimes V_2 \to V_3$
(where the projection is taken w.r.t. the splitting given
by $b=0$); then, if $b$ is the bilinear map corresponding to $B$,
it is seen that the splitting $\beta$ is the one induced by
$L^b$. Summing up,
splittings of (BA.19) are in one-to-one
correspondence with bilinearly related structures.
Once again this shows that the space of bilinearly related structures
carries the structure of an affine space over $\K$.
Moreover, intrinsically bilinear maps
$f:V_1 \times V_2 \times V_3 \to W$
correspond  to  linear maps $F:Z \to W$; homogeneous ones
correspond to linear maps that vanish on $W$.
In particular, the ``third projection'' corresponds to the
linear map $Z \to V_3$, and the projection $Z \to V_1 \otimes V_2$
corresponds to $f_\omega$, where $\omega:V_1 \times V_2 \to V_1\otimes V_2$
is the tensor product map.

\vfill \eject

\subheadline
{MA. Multilinear algebra}

For a correct generalization of the bilinear theory from the
preceding chapter to the multilinear case we have to change notation:
the $\K$-modules $V_1,V_2$ from BA.1 play a symmetric r\^ole
and should be denoted by
$V_{01},V_{10}$, whereas $V_3$ belongs to a ``higher'' level and shall be
denoted by $V_{11}$. The multi-indices $01,10,11$ in turn correspond to
non-empty subsets of the $2$-element  set $\{ 1,2 \}$. 
In the general $k$-multilinear
 case, our $\K$-modules will be indexed by the lattice of non-empty subsets
of the $k$-element set $\{ 1,\ldots ,k \}$.

\msk \nin {\bf MA.1.} {\sl The index set.}
For a finite set $M$, we denote by $2^M$ its power set,
and taking for $M$ the $k$-element set
$$
\N_k :=\{ 1,\ldots,k \},
\eqno{\rm (MA.1)}
$$
we define our {\it index set} 
$I:= I_k := 2^{\N_k}$.
Often we will identify $I_k$ with $\{ 0,1 \}^k$ via the bijection
$$
2^{\N_k} \to \{ 0 ,1 \}^k, \quad A \mapsto \alpha:=\1_A ,
\eqno {\rm (MA.2)}
$$
where $\1_A$ is the characteristic function of a set $A \subset \N_k$.
The inverse mapping assigns to $\alpha \in \{ 0,1 \}^k$ its
{\it support}
$$
A:=\supp(\alpha):=\{ i \in \N_k | \, \alpha_i =1 \} .
$$
The natural total order on $\N_k$ induces a natural total
order ``$\leq$" on $I_k$ which
is the lexicographic order. For instance, $I_3$ is ordered as
$$
\eqalign{
I_3 & = ( \emptyset, \{ 1 \}, \{ 2 \}, \{ 1,2 \}, \{ 3 \},
\{ 3,1 \}, \{ 3,2 \}, \{ 3,2,1 \} ) \cr
& \cong ((000),(001),(010),(011),(100),(101),(110),(111)),\cr}
\eqno {\rm (MA.3)}
$$
where, for convenience, we
prefer to write elements $\alpha \in
\{ 0,1 \}^3 $ in the form $(\alpha_3,\alpha_2,\alpha_1)$ and not
$(\alpha_1,\alpha_2,\alpha_3)$. 
There are other useful total orderings of $I_k$, but
the natural order has the advantage that it is compatible with the
natural inclusion $\N_k \subset \N_{k+1}$ and hence is best adapted
to induction procedures. 
In terms of multi-indices $\alpha \in \{ 0,1 \}^k$, the basic
set-theoretic operations are interpreted as follows:

\ssk
\item{(a)}
the cardinality of $A$ corresponds to 
$\vert \alpha \vert :=  \sum_i \alpha_i$;
\item{(b)}
the singletons $A = \{ i \}$ correspond to 
the ``canonical basis vectors'' $e_i$, $i=1,\ldots,k$;
\item{(c)}
inclusion $\alpha \subseteq \beta$ is defined by
$\supp(\alpha) \subseteq \supp(\beta)$ (note that this partial order
is compatible with the total order:
$\alpha \subseteq \beta$ implies $\alpha \leq \beta$);
\item{(d)}
disjointness will be denoted by
$$
\alpha \perp \beta  \quad :\Leftrightarrow \quad  
\supp(\alpha) \cap \supp(\beta)=\emptyset
 \quad \Leftrightarrow \quad  \sum_i \alpha_i \beta_i =0, 
\eqno {\rm (MA.4)}
$$
which is equivalent to saying that the componentwise sum
$\alpha + \beta$ is again an element of $\{0,1 \}^k$.

\msk \nin {\bf MA.2.} {\sl Cubes of $\K$-modules.}
A {\it $k$-dimensional cube of $\K$-modules} is given by a family
$V_\alpha$, $\alpha \in I_k \setminus 0$, of $\K$-modules. We define
the {\it total space}
$$
E := \bigoplus_{\alpha \in I_k \setminus 0} V_\alpha.
$$
Elements of $E$ will be denoted by $v=\sum_{\alpha>0} v_\alpha$ or
by $v=(v_\alpha)_{\alpha >0}$ with $v_\alpha \in V_\alpha$.
The subspaces $V_\alpha \subset E$ are called the {\it axes of $E$}.
For every $\alpha \in I$, $\alpha > 0$, we let
$$
\tilde V_\alpha := \bigoplus_{\beta \subseteq \alpha} V_\beta, \quad \quad
H_\alpha := \bigoplus_{\alpha \subseteq \beta} V_\beta.
\eqno {\rm (MA.5)}
$$
If $\alpha \subseteq \alpha'$, then 
$\tilde V_\alpha \subseteq \tilde V_{\alpha'}$, and 
the lattice of subsets of $I_k$ corresponds thus to a
 commutative cube of inclusions. For instance, for $k=3$,
$$
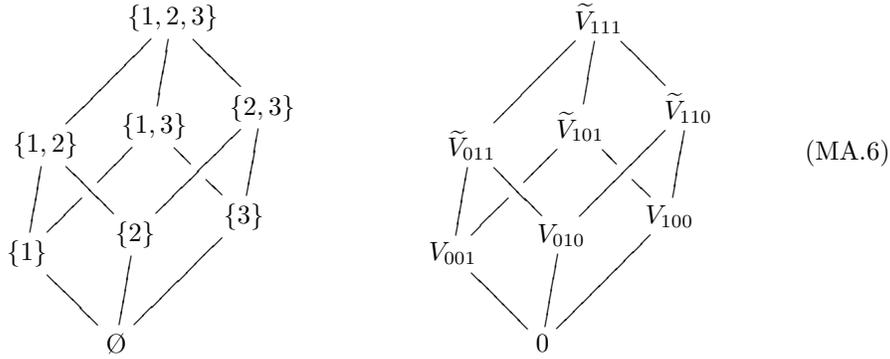
\xymatrix@dl@!0@C+10pt{
&\{ 3 \} \ar@{-}[rr]\ar@{-}'[d][dd]
 & & \emptyset \ar@{-}[dd]
\\
 \{ 2,3 \} \ar@{-}[ur]\ar@{-}[rr]\ar@{-}[dd]
& & \{ 2 \} \ar@{-}[ur]\ar@{-}[dd]
\\
&  \{ 1, 3 \} \ar@{-}'[r][rr]
& & \{ 1 \}
\\
\{ 1,2,3 \} \ar@{-}[rr]\ar@{-}[ur]
& & \{ 1,2 \} \ar@{-}[ur]
}
\quad \quad \quad \quad
\xymatrix@dl@!0@C+10pt{
& V_{100}  \ar@{-}[rr]\ar@{-}'[d][dd]
 & & 0 \ar@{-}[dd]
\\
\tilde V_{110} \ar@{-}[ur]\ar@{-}[rr]\ar@{-}[dd]
& & V_{010} \ar@{-}[ur]\ar@{-}[dd]
\\
& \tilde V_{101} \ar@{-}'[r][rr]
& & V_{001}
\\
\tilde V_{111} \ar@{-}[rr]\ar@{-}[ur]
& & \tilde V_{011} \ar@{-}[ur]
}
\eqno {\rm (MA.6)}
$$

\msk \nin {\bf MA.3.} {\sl Partitions.} 
A {\it partition} of a finite set $A$ is given by a subset
$\Lambda = \{ \Lambda^1,\ldots,\Lambda^\ell \}$ of the power set of $A$
such that all $\Lambda^i$ are non-empty and
 $A$ is the disjoint union of the $\Lambda^i$:
$A = \dot \cup_{i=1,\ldots,\ell} \Lambda^i$. We then say that 
$\ell:=\ell(\Lambda)$
 is the {\it length} of the partition.  We denote by ${\cal P}(A)$ the
set of all partitions of $A$ and by ${\cal P}_\ell(A)$ 
(resp.\ by ${\cal P}_{\geq \ell}(\alpha)$)
the set of all partitions of $\alpha$ of length $\ell$ (resp.\ 
at least $\ell$):
$$
{\cal P}_\ell (A) = 
\{ \Lambda = \{ \Lambda^1,\ldots,\Lambda^\ell \} \in 2^{(2^A)} \vert \,
\quad \dot \cup_{i=1,\ldots,\ell} \, \Lambda^i =A; \, \emptyset 
\notin \Lambda \}.
\eqno{\rm (MA.7)}
$$
Finally, for a finite set $M$, we  denote by   
$$
{\rm Part}(M):= \bigcup_{A \in 2^M \atop A \not= \emptyset}
{\cal P}(A) \subset 2^{2^M}
\eqno{\rm (MA.8)}
$$
the set of all partitions of non-empty subsets of $M$, and for a partition
$\Lambda$ of some set $A \subset M$, the set $A$ 
is called the {\it total set of the partition $\Lambda$}, denoted by 
$$
\uline \Lambda := A =  \bigcup_{\nu \in \Lambda} \nu.
\eqno{\rm (MA.9)}
$$
If we use multi-index notation, this corresponds to $\uline \lambda := \alpha =
\sum_i \lambda^i$.

Now let $M=\N_k$, $A \subset M$.
 Using the total order on $I_k = 2^{\N_k}$, every partition
of length $\ell$ can also be written as an ordered $\ell$-tupel:
$$
\Lambda = (\Lambda^1,\ldots,\Lambda^\ell), \quad 
\emptyset < \Lambda^1 < \ldots < \Lambda^\ell \leq A.
$$
In  terms of multi-indices,
 a {\it partition of a multi-index $\alpha \in \{ 0,1 \}^k$}
is described as
$\lambda = \{ \lambda^1,\ldots,\lambda^\ell \} \subset I_k$
such that $\sum_{i=1}^\ell \lambda^i = \alpha$ and all $\lambda^i > 0$, or,
using ordered sequences:
$$
\lambda = (\lambda^1,\ldots,\lambda^\ell) \in (I_k)^\ell, \quad 
0< \lambda^1 < \ldots < \lambda^\ell, \quad
\sum_{i=1}^\ell \lambda^i = \alpha.
\eqno {\rm (MA.10)}
$$
The {\it trivial partition} $\{ A \}$
 of a set $A$, i.e.\ the one of length one, 
is identified with $A$ itself. 
Sometimes it is useful to represent a partition by an $\ell \times k$-matrix 
with rows $\lambda^1,\ldots,\lambda^\ell$. So, for instance,
$$
{\cal P}(\{ 1,2,3 \}) = \Bigl\{ \{ 1,2,3 \}, \{ \{ 1 \}, \{ 2,3 \} \},
\{ \{ 2 \}, \{ 1,3 \} \}, \{ \{ 1,2 \}, \{ 3 \} \},
\{ \{ 1 \}, \{ 2 \}, \{ 3 \} \Bigr\}
$$
corresponds to 
$$
\eqalign{
{\cal P}((111))& =\Bigl\{ \pmatrix{1 & 1 & 1 \cr}, 
\pmatrix{0 & 0 & 1 \cr 1 & 1 & 0 \cr},
\pmatrix{0 & 1 & 0 \cr 1 & 0 & 1 \cr}, \pmatrix{0 & 1 & 1 \cr 1 & 0 & 0 \cr},
\pmatrix{0&0&1 \cr 0&1&0 \cr 1&0&0 \cr} \Bigr\}.
\cr}
\eqno {\rm (MA.11)}
$$
The set ${\cal P}(A)$ can again be totally ordered, but we will not use this
order. 

\ssk 
The symmetric group $\Sigma_k$ acts on $ \N_k$, on $I_k$ and on
${\rm Part}(\N_k)$ in the usual way.
Then $\alpha,\beta \in I_k$ are conjugate under this action iff
$|\alpha|=|\beta|$, and two partitions $\lambda,\mu \in {\cal P}(I_k)$  
are conjugate under $\Sigma_k$ iff
$$
\ell(\mu)=\ell(\lambda)=:\ell \,\quad {\rm and} \, \quad
\forall i=1,\ldots,\ell: \, \vert \lambda^i \vert = \vert \mu^i \vert.
\eqno {\rm (MA.12)}
$$ 
For instance, in (MA.11) the second, third and fourth partition are
equivalent under $\Sigma_3$. Note that,
when writing partitions as ordered sequences, the action of $\Sigma_k$ does
not preserve the order. In (MA.11),
the lines of the fourth matrix had to be exchanged
after letting act the transposition (23) on the third matrix.

\msk \nin {\bf MA.4.} {\sl Refinements of partitions.} 
We say that a partition $\Lambda$ is a {\it refinement} of another
partition $\Omega$, or that $\Omega$ is {\it coarser than $\Lambda$},
 and we write $\Omega \preceq \Lambda$, if
$\Lambda$ and $\Omega$ are partitions of the same set $A$ and every
set $L \in \Lambda$ is contained in some set $O \in \Omega$:
$$
\Omega \preceq \Lambda \quad : \Leftrightarrow \quad 
\uline \Omega = \uline \Lambda, \, \, \forall L \in \Lambda:
\exists O \in \Omega: \, L \subseteq O.
\eqno{\rm (MA.13)} 
$$
This means that the equivalence relation on $\uline \Lambda$ induced
by $\Lambda$ is finer than the one induced by $\Omega$.
For instance, the last partition in Equation (MA.11) is a refinement of all the
preceding ones, and  the first partition is coarser than all the others.
The relation $\preceq$ defines a partial order on ${\cal P}(\alpha)$;
for instance, the three middle partitions in (MA.11) cannot be compared
among each other. In any case, if $\Lambda$ is a partition of a set $A$
with
$|A|=j$, we have a finest and a coarsest partition of this set:
$$
\uline \Lambda = A = \{ a_1,\ldots, a_j \} \, \preceq \Lambda \preceq \,
\oline \Lambda := \Bigl\{  \{a_1 \},\ldots,\{ a_j\} \Bigr\}=
\Bigl\{ \{ a \} | \, a\in A \Bigr\}.
\eqno{\rm (MA.14)}
$$ 
If $\Lambda$ is finer than $\Omega$, we
define, for $O \in \Omega$, the {\it $\Lambda$-induced partition
of $O$}, by
$$
O | \Lambda := \{ L \in \Lambda | \, L \subset O \} \in {\cal P}(O).
\eqno{\rm (MA.15)} 
$$
Another way of viewing refinements is via {\it partitions of partitions}:
a partition $\Lambda$ is a set and hence can
again be partitioned. In fact, every partition $\Omega$ that is coarser
than $\Lambda$ defines a partition of the partition
$\Lambda$ via
$$
\Lambda = \dot \cup_{O \in \Omega} \, \,  O| \Lambda ,
$$
and in this way we get a canonical bijection between the set of all partitions
$\Omega$ that are coarser than $\Lambda$ and the set ${\cal P}({\cal P}
(\Lambda))$ of all partitions of the partition $\Lambda$.

\msk \nin {\bf MA.5.} {\sl Multilinear maps between cubes of
$\K$-modules.}
Given two cubes of $\K$-modules $(V_\alpha),(V_\alpha')$ with total
spaces $E,E'$, for every partition $\Lambda$ of an element $\alpha \in I_k$,
we assume a $\K$-multilinear map
$$
b^\Lambda: V_{\Lambda^{1}} \times \ldots \times V_{\Lambda^{\ell(\Lambda)}}
\to V_\alpha '
\eqno {\rm (MA.16)}
$$
be given. Then a {\it multilinear map} or a {\it homomorphism of cubes of
$\K$-modules} is a map of the form
$$
f:=f_b:E \to E', \quad
v=\sum_{\alpha > 0}  v_\alpha  \mapsto \sum_{\alpha > 0}
\sum_{\Lambda \in {\cal P}(\alpha)} b^\Lambda(v_{\Lambda^{1}},\ldots,
v_{\Lambda^{\ell(\Lambda)}} ).
\eqno {\rm (MA.17)}
$$
Written out more explicitly,
$$
\eqalign{
f_b(v) & = \sum_{\alpha \in I_k \atop  \alpha >0} 
\sum_{\ell=1}^{\vert \alpha \vert}
\sum_{\Lambda \in {\cal P}_\ell(\alpha)} b^\Lambda(v_{\Lambda^{1}},\ldots,
v_{\Lambda^{\ell}} ).
\cr}
$$
To get a more concise notation, we let
$$
V_\Lambda := V_{\Lambda^1} \times \ldots \times  V_{\Lambda^\ell}
\eqno {\rm (MA.18)}
$$
and write $b^\Lambda : V_\Lambda \to V_{\uline \Lambda}'$.
Note that $\Lambda = \{ \Lambda^1,\ldots,\Lambda^\ell \}$ is a set,
hence we will not distinguish between the various versions of 
$V_\Lambda$ or $b^\Lambda$ obtained by permuting the factors.
 If one whishes, one may
single out the ``ordered version'' $\Lambda^1 < \ldots < \Lambda^\ell$,
but as long as we distinguish the various spaces $V_\alpha$ by their
index (even if they happen to be isomorphic), this will not be
necessary. We write $\Hom(E,E')$ for the set of all multilinear maps
$f=f_b:E \to E'$ and $\End(E):=\Hom(E,E)$.
The family 
$b = (b^\lambda)_{\lambda \in {\rm Part}(I)}$ of multilinear maps
is called a {\it multilinear family},
and we denote by ${\rm Mult}(E,E')$ the set of all multilinear families on $E$
with values in $E'$.
We introduce the following conditions on multilinear maps and on the
corresponding multilinear families:

\ssk
\item{(reg)} 
We say that $f_b:E \to E'$ is {\it regular} if, whenever
$\ell(\Lambda)=1$ (i.e.\ $\Lambda = \{ \alpha \} = \{ \Lambda^1 \}$
 is a trivial partition), the
(linear) map $b^\Lambda: V_\alpha \to V_\alpha'$ is bijective.
We denote by
$$
\Gm^{0,k}(E):=\{ f_b \in \End(E) | \, b \in {\rm Mult}(E), \,
\ell(\Lambda)=1 \Rightarrow b^\Lambda \in \Gl_\K(V_{\lambda^1}) \}
$$
the set of all regular endomorphisms of $E$.
\item{(up)} 
We say that an element $f_b \in \Gm^{0,k}(E)$ is {\it unipotent}
if, whenever
$\Lambda = \{ \alpha \}$ is a trivial partition, we have
 $b^\Lambda=\id_{V_\alpha}$. We denote by
$$
\Gm^{1,k}(E) :=\{ f_b \in \Gm^{0,k}(E) | \, b \in \Mult(E), \, 
\ell(\Lambda)=1 \Rightarrow b^\Lambda = \id_{V_{\lambda^1}} \}
$$
the set of all unipotent endomorphisms of $E$.
\item{(si)} 
Let $1 <  j \leq k$. 
We say that $f_b$ {\it is singular of order $j$ (at the origin)} or
just: $f_b$ {\it is of order $j$} if $b^\Lambda =0$ whenever $\ell(\Lambda)
\leq j$,  and we denote by $\Mult_{>j}(E)$ the set of all
multilinear families on $E$ that are singular of order $j$.
\item{(ups)}
 For $1<j \leq k$, we denote by
$$
\Gm^{j,k}(E):=\{ f_b \in \Gm^{1,k}(E) | \, 
\ell(\Lambda)=2,\ldots,j \Rightarrow b^\Lambda = 0 \}
$$
the set of all unipotent endomorphisms that are of the form
$f = \id_E  + h_b$ with $h_b$ singular of order $j$. 

\Theorem MA.6. 
\item{(1)}
Cubes of $\K$-modules together with their homomorphisms 
form a category (which we call the {\rm category of $k$-multilinear spaces}).
\item{(2)}
A homomorphism $f_b \in \Hom(E,E')$ is invertible if and only if
it is regular, and then $(f_b)^{-1}$ is again a homomorphism.
 In particular, $\Gm^{0,k}(E)$ is a group, namely it is
the automorphism group of $E$ in the category of $k$-multilinear spaces.
\item{(3)} We have a chain  of  subgroups
$$
\Gm^{0,k}(E)\supset \Gm^{1,k}(E) \supset
\ldots \supset \Gm^{k-1,k}(E)\supset \Gm^{k,k}(E) =\{\id_E \}.
$$

\Proof. (1) It is fairly 
obvious that the composition of homomorphisms is again a homomorphism.
However, let us, for later use, spell out the composition rule in more
detail: 
assume $g:E'\to E''$ and $f:E \to E'$ are multilinear maps;
for simplicity we will denote the multilinear families corresponding
to $g$ and to $f$ by $(g^\Lambda)_\Lambda$ and $(f^\Lambda)_\Lambda$.
Then the composition rule 
is obtained simply by using multilinearity and ordering terms
according to the space in which we end up:
let $w_\alpha=
\sum_{\nu \in {\cal P}(\alpha)} f^\nu(v_{\nu_1},\ldots,v_{\nu_\ell})$, then
$$
\eqalign{g(f(v)) &= g(\sum_\alpha w_\alpha)
 \cr
&=\sum_\alpha \sum_{\Omega \in {\cal P}(\alpha)} g^\Omega(w_{\Omega^1},\ldots,
w_{\Omega^r})  
\cr
&=\sum_\alpha \sum_{\Omega \in {\cal P}(\alpha)} 
\sum_{\Lambda_1 \in {\cal P}(\Omega^1) \ldots \atop
\Lambda_r \in {\cal P}(\Omega^r) }
g^\Omega\Bigl(f^{\Lambda_1}(v_{\Omega_1^1},\ldots,v_{\Omega^1_{\ell_1}}),\ldots,
f^{\Lambda_r}(v_{\Omega^r_1},\ldots,v_{\Omega^r_{\ell_r}})\Bigr)  
\cr
&= \sum_\alpha \sum_{\Lambda \in {\cal P}(\alpha)} 
\sum_{\Omega \preceq \Lambda}
g^\Omega \Bigr(f^{ \Omega^1|\Lambda }(v_{\Omega_1^1},\ldots,
v_{\Omega^1_{\ell_1}} ),
\ldots, f^{ \Omega^r|\Lambda}(v_{\Omega^r_1},\ldots,v_{\Omega^r_{\ell_r}} )\Bigl)
\cr}
$$
where $r=\ell(\Omega)$ and $\ell_j = \ell(\Lambda_j)$ and 
$\Lambda_j = \Omega^j\vert\Lambda$ (cf.\ Eqn.\ (MA.15); this implies
$\Lambda =\cup_{j=1}^r \Lambda_j$). It follows that $g \circ f$ corresponds 
to the multilinear family
$$
(g \circ f)^\Lambda(v_\Lambda) = 
\sum_{\Omega \preceq \Lambda}
g^\Omega \Bigl(f^{ \Omega^1|\Lambda }(v_{\Omega_1^1},\ldots,
v_{\Omega^1_{\ell_1}} ),
\ldots, f^{ \Omega^r|\Lambda}(v_{\Omega^r_1},\ldots,v_{\Omega^r_{\ell_r}})\Bigr).
\eqno ({\rm MA.19})
$$
Finally, the identity map clearly is a homomorphism, and hence we have
defined a category.  

(2) 
If $f_b:E \to E'$ is a homomorphism, then the definition of $f_b$
shows that $f_b$ maps the $\alpha$-axis into the $\alpha$-axis:
$f_b(V_\alpha) \subset V_\alpha'$, and the induced map
$V_\alpha \to V_\alpha'$ agrees with $b^\alpha:=b^\Lambda$
 for $\Lambda = \{\alpha \}$.
Therefore, if $f_b$ is bijective, then so is $b^\Lambda$ whenever
$l(\Lambda)=1$, and thus $f_b$ is regular.
 In order to prove the converse, we define, for every
homomorphism $f=f^b:E \to  E'$, the map
$$
T_0 f:= \times_\alpha b^\alpha:
 E \to E', \quad \sum_\alpha v_\alpha \mapsto \sum_\alpha
b^\alpha (v_\alpha),
\eqno ({\rm MA.20})
$$
where as above $b^\alpha := b^{\{ \alpha \} }$. This is again a homomorphism,
acting on each axis in the same way as $f$. From the definition of 
multilinear maps
one immediately gets the functorial rules $T_0(g \circ f)=T_0g \circ T_0 f$
and $T_0(T_0 f)=T_0 f$. 
Now assume that $f:=f_b$ is regular. This is equivalent to saying that $T_0 f$
is bijective. Then clearly $(T_0 f)^{-1}$ is again a homomorphism
and $T_0 (T_0 f)^{-1} = (T_0 f)^{-1}$. Let
 $g:= (T_0f)^{-1} \circ f: E \to E$; then 
$g$ is an endomorphism with
$T_0 g = \id_E$, i.e.\ $g$ is regular and moreover unipotent, i.e.\
$g \in \Gm^{1,k}(E)$.
Thus we are done if we can show that $\Gm^{1,k}(E)$ is a group.

Let us, more generally, 
 prove by a descending induction on $j=k,k-1,\ldots, 1$ that all $\Gm^{j,k}(E)$
are groups:
it is clear that $\Gm^{k,k}(E)=\{\id_E \}$ is a group.
For the induction step with $j\geq 1$, note first that $\Gm^{j,k}(E)$
is stable under composition. We need the lowest term of the composition
of elements
  $f_a,f_b \in \Gm^{j,k}(E)$: define $c$ by $f_c=f_b \circ f_a$; then we have
$c^\lambda = a^\lambda + b^\lambda$ whenever $l(\lambda)=j+1$. It follows that
$g:=f_{-b} \circ f_b \in \Gm^{j+1,k}(E)$ for all $f_b \in \Gm^{j,k}(E)$.
But by induction $\Gm^{j+1,k}(E)$ is a group, and hence $f_b$ is invertible
with inverse $g^{-1} \circ f_{-b} \in \Gm^{j,k}(E)$.
 (In particular, this argument shows 
that $\Gm^{k-1,k}(E)$ is a vector group.) 
Finally,
(3) has just been proved.
\qed

\nin
As seen in the proof, (2) implies that $f:E \to E'$ is invertible if and
only if $T_0 f$ is invertible. If we interpret $T_0 f$ as the ``total
differential of $f$ at the origin'', then this may be seen as a kind of
{\it inverse function theorem for homomorphisms of multilinear spaces}. 

\msk \nin
{\bf MA.7.} {\sl The general and the special multilinear group.} 
The automorphism group $\Gm^{0,k}(E)$ is called the {\it general multilinear
group of $E$}, and the group $\Gm^{1,k}(E)$ is called 
 the {\it special multilinear group of $E$}. 
The latter acts ``simply non-linearly'' on $E$ in the sense that
the only group element that acts linearly on $E$ is the identity. Therefore
it plays a more important r\^ole in our theory than the full automorphism
group. The preceding proof shows that
we have a (splitting) exact sequence of groups
$$
\matrix{
1& \to & \Gm^{1,k}(E)& \to & \Gm^{0,k}(E) & {\buildrel T_0 \over \to} &
\times_{\alpha>0} \Gl(V_\alpha) & \to & 1 \cr}
\eqno ({\rm MA.21})
$$
showing that $\Gm^{1,k}(E)$ is normal in the automorphism group, and the
quotient is a product of ordinary linear groups.
For more information on the structure of $\Gm^{1,k}(E)$, see Sections
MA.12 -- MA.15.\ and Chapter PG.

\msk \nin
{\bf MA.8.} {\sl Example: the case $k=3$.}
For $k=3$, elements of $\Gm^{1,3}(E)$ are
described by six bilinear maps and one trilinear map:
$$
b=(b^{001,010}, \, b^{001,100}, \, b^{010,100}, \,
b^{001,110}, \, b^{010,101}, \, b^{011,100}, \, \,
b^{001,010,100}),
$$
and thus
$$
\eqalign{ & f_b(v) =v+
b^{001,010}(v_{001},v_{010}) + b^{001,100}(v_{001},v_{100}) + b^{010,100}
(v_{010},v_{100}) +\cr
& \quad   \quad
 b^{001,110}(v_{001},v_{110}) + b^{010,101}(v_{010},v_{010})
+ b^{011,100}(v_{011},v_{100}) + \cr
& \quad \quad b^{001,010,100}(v_{001},v_{010},v_{100}),
\cr}
$$
where the term in the second and third 
line describes the component in $V_{111}$.
Composition $f^b \circ f^a = f^c$ of such maps is described by the
following formulae: if $\ell(\lambda)=2$, then we have simply
$c^\lambda = b^\lambda + a^\lambda$. If $\ell(\lambda)=3$, i.e.,
$\lambda = ((001),(010),(100))$, then there are three
non-trivial partitions $\Omega$ that are strictly coarser than
$\Lambda$, and we have:
$$
\eqalign{
& c^\lambda(u,v,w)  =a^\lambda(u,v,w) + b^\lambda(u,v,w) + \cr
& \quad b^{001,110}(u,a^{010,100}(v,w)) + 
b^{010,101}(v,a^{001,100}(u,w)) +
b^{011,100}(a^{001,010}(u,v),w). \cr}
$$
Inversion of a unipotent $f=f_b$ is described as follows:
if $\ell(\lambda)=2$, then $(f^{-1})^\lambda = - f^\lambda$,
and the coefficient belonging to $\lambda=((100)(010)(001))$ is
$$
\eqalign{
& (f^{-1})_\lambda (u,v,w)  = -b^\lambda(u,v,w) +   \cr
& \quad b^{001,110}(u,b^{010,100}(v,w)) + b^{010,101}(v,b^{001,100}(u,w)) +
b^{011,100}(b^{001,010}(u,v),w). \cr}
$$

\msk \nin {\bf MA.9.} {\sl Multilinear connections
 and the category of multilinear spaces.}
Since $\Gm^{1,k}(E)$ acts non-linearly on $E$,
a multilinear space $E$ does not have any distinguished linear structure.
We say that two linear structures on $E$ are {\it $k$-multilinearly related}
if they are conjugate to each other under the action of
$\Gm^{1,k}(E)$, and we say that a linear structure on $E$ is a {\it
(multilinear) connection on $E$} if it is 
$k$-multilinearly related to the initial
$\K$-module structure $L^0$ defined by the bijection
$E \cong \oplus_{\alpha>0} V_\alpha$. As in Section BA.1, we denote by
$L^b=(E, +_b, \cdot_b)$
the push-forward of the original linear structure on $E$ via $f_b$;
then $\Gm^{1,k}(E)$ acts simply transitively on the space of
all multilinear connections via $(f^c,L^b) \mapsto (f^c \circ f^b) \cdot L^0$.
Recall that, for $k=2$, $\Gm^{1,k}(E)$ is a vector group, and hence
in this case the space of connections is an affine space over $\K$. 
For general $k$,
the explicit formulae for addition $+_b$ and multiplication by scalars 
$\cdot_b$ in the $\K$-module $(E,L^b)$
 are much more complicated than the formulae given in Section BA.2 for
the case $k=2$,  and we will not try to write them  out in full detail.
Let us just make the following remarks.
It is easy to write the sum of elements taken from axes:
since $f_b$ acts trivially on each axis, this is the ``$b$-sum''
$$
(b\sum)_\alpha  v_\alpha := f_b(\sum_\alpha  v_\alpha).  
\eqno {\rm (MA.22)}
$$
Moreover, in the sum $v+_b v'$ for arbitrary elements $v,v' \in E$,
 components with $|\alpha|=1$ are simply of the form
$v_\alpha + v_\alpha'$; components with $|\alpha|=2$ are of the form
$v_\alpha + v_\alpha' + b^{\beta,\gamma}(v_\beta,v_\gamma')+
 b^{\beta,\gamma}(v_\beta',v_\gamma)$ where $\alpha = \beta + \gamma$
is the unique decomposition with $|\beta|=|\gamma|=1$ and $\beta < \gamma$.

\msk
\centerline{\bf Intrinsic geometric structures}

\msk \nin
We say that a structure on $E$ is {\it intrinsic} if it is invariant
under the action of the general multilinear  group $\Gm^{0,k}(E)$.

\msk \nin {\bf MA.10.} {\sl Intrinsic structures: subgeometries.}
The group $\Gm^{0,k}(E)$ preserves axes and acts linearly on them;
 hence the axes are
intrinsic subspaces, not only as sets, but also as sets with their induced
$\K$-module structure.
For any $\alpha \in I_k$, the set
$\tilde V_\alpha := \bigoplus_{\beta \subseteq \alpha} V_\alpha$ (cf.\
Eqn.\ (MA.5)) is invariant under the action of
$\Gm^{1,k}(E)$ and hence is an intrinsic subspace of $E$; but the action of 
$\Gm^{1,k}(E)$ on these subspaces is in general not trivial.
For $k=3$, diagram (MA.6) describes the cube of these intrinsic subspaces.
Note that, if $|\alpha|=j$, then $\tilde V_\alpha$ is a $j$-linear
space in its own right which is represented by the
$j$-dimensional subcube having bottom point $0$ and
top point $\tilde V_\alpha$.
For instance, if $k=3$, the three bottom faces of the cube represent
three bilinear subspaces of $E$. 
In a similar way, if $\alpha,\beta,\alpha + \beta \in I$, then the direct sum
$V_\alpha \oplus V_\beta \oplus V_{\alpha + \beta}$ is an invariant
bilinear subspace of $E$.
More generally,
 for any partition $\Lambda$, we can define an invariant subspace
$E_\Lambda$ which is a  cube of $\K$-modules with bottom
$\Lambda^1,\ldots,\Lambda^\ell$.

\msk \nin {\bf MA.11.} {\sl Intrinsic structures: projections, hyperplanes.} 
Let 
$$
p_\beta:E \to \tilde V_\beta, \quad \sum_\alpha v_\alpha \mapsto
\sum_{\alpha \subset \beta} v_\alpha
\eqno {\rm (MA.23)}
$$
be the projection onto the invariant subspace 
$\tilde V_\beta$. Then a straightforward check shows 
that $p_\beta$ is intrinsic in the sense that $p_\beta \circ f_b = f_b \circ
p_\beta$ for all $f_b \in \Gm^{1,k}(E)$:
$$
\matrix{
E & {\buildrel f_b \over \to} & E \cr
p_\beta \downarrow \phantom{p_\beta}& &p_\beta \downarrow \phantom{p_\beta}\cr
\tilde V_\beta & {\buildrel f_b \over \to} & \tilde V_\beta \cr}
\eqno {\rm (MA.24)}
$$
In particular, the factors of the decomposition
$$
E = \tilde V_\beta \oplus K_\beta, \quad \quad 
K_\beta:=\ker(p_\beta)= \bigoplus_{\alpha \not\subseteq \beta} V_\alpha
=\sum_{\alpha \perp \beta} H_\alpha
$$
are invariant under $\Gm^{1,k}(E)$; but of course $\Gm^{1,k}(E)$
 does in general not
act trivially on the factors, nor are the ``parallels" of $\tilde V_\beta$
and $K_\beta$ invariant under $\Gm^{1,k}(E)$.
We say that the subgeometry $\tilde V_\beta$ is {\it maximal} if
$\vert \beta \vert = k-1$ and {\it minimal} if $\vert \beta \vert = 1$.
(For $k=2$ every proper subgeometry is both maximal and minimal.)

\ssk
\item{(1)} Maximal subgeometries: fix $j \in \{ 1,\ldots,k\}$ and let
 $\beta:=e_j^*:
 = \sum_{i\not= j} e_i$ correspond to the complement $\N_k \setminus
 \{ j \}$. The corresponding projection is
$$
p_{e_j^*}:E \to \tilde V_{(1\dots101\ldots1)} = \bigoplus_{\alpha \in I \atop
\alpha_j =0 } V_\alpha, \quad \sum_{\alpha >0} v_\alpha \mapsto
\sum_{\alpha_j=0} v_\alpha.
\eqno {\rm (MA.25)}
$$
The kernel is called the {\it $j$-th hyperplane}:
$$
H_j:=K_{e_j^*} = \ker(p_j)=
\bigoplus_{\alpha \in I \atop \alpha_j = 1} V_\alpha =
 H_{e_j}
\eqno {\rm (MA.26)}
$$
We claim that $\Gm^{1,k}(E)$ acts trivially on $H_{j}$. Indeed,
let $v=\sum_{\alpha_j = 1} v_\alpha \in H_{e_j}$.
Then $f_b(v)=v + \sum_\lambda b^\lambda(v_{\lambda^1},\ldots,v_{\lambda^\ell})=
v$  since $\lambda^i_j=1$ for $i=1,\ldots,\ell$ and
hence the $\lambda^i$ cannot form
a non-trivial partition and the sum over $\lambda$ has to be empty.
Hence $H_{e_j}$ together with its linear structure is invariant under $f_b$
and thus it is an intrinsic subspace of $E$.
In the same way we see that all fibers
 $y + H_{e_j}$, for any $y \in E$, are intrinsic
affine subspaces of $E$.
There is a distinguished origin in $y+H_{e_j}$, namely $\sum_{\alpha_j=0}
 y_\alpha$, and hence $y+H_{e_j}$ carries an intrinsic linear structure.
 (In the context of higher order tangent bundles $T^k M$,
  these hyperplanes correspond to
tangent spaces $T_u (T^{k-1} M)$ contained in $E=(T^k M)_x$.)

\item{(2)} Minimal subgeometries: $\beta = e_j$. Then the projection onto
$\tilde V_{e_j}=V_{e_j}=:V_j$ is
$$
\pr_j:=  p_{e_j}: E\to V_j, \quad v = \sum v_\alpha \mapsto v_{e_j}.
\eqno {\rm (MA.27)}
$$
The group $\Gm^{1,k}(E)$ acts trivially on the axes $V_j$ and hence
$p_{e_j} \circ f_b =p_{e_j}$; thus not only the kernel, but
all fibers of $p_{e_j}$ are $\Gm^{1,k}(E)$-invariant
subspaces.

\item{(3)} Intersections. Clearly, intersections of intrinsic subspaces
are again intrinsic subspaces. For instance, if $i\not= j$,
$$
H_{ij}:=H_i \cap H_j = \ker( p_{e_i^*} \times p_{e_j^*}) =
 \bigoplus_{\alpha \in I \atop \alpha_i = 1 = \alpha_j}
V_\alpha = H_{e_i + e_j}
\eqno {\rm (MA.28)}
$$
is an intrinsic subspace, and so is the ``vertical subspace"
$\bigoplus_{\alpha: \vert \alpha \vert \geq 2} V_\alpha$, kernel
of the projection
$$
\pr:=\pr_1 \times \ldots \times \pr_k: E \to V_1 \times \ldots \times V_k,
\quad v \mapsto (\pr_1(v),\ldots,\pr_k(v)).
\eqno {\rm (MA.29)}
$$

\ssk \nin
Summing up, multilinear spaces have a non-trivial ``intrinsic
incidence geometry", and it seems that the complementation map
of the lattice $I_k$ corresponds to some duality of this incidence
geometry.  It should be interesting to develop this topic
more systematically, in particular in relation with Lie groups and
symmetric spaces. 


\bigskip
\centerline{\bf Complements on the structure of the general multilinear group}

\msk \nin
{\bf MA.12.} {\sl ``Matrix coefficients'' and ``matrix multiplication''.}
The choice of a multilinear
connection in a multilinear space should be seen as an
analog of choosing a base in a free (say, $n$-dimensional)
 $\K$-module $V$, and the analog of 
the map $\K^n \to V$ induced by a basis  is the map
$$
\Phi^b: \bigoplus_{\alpha >0} V_\alpha 
\to E, \quad v=(v_\alpha)_\alpha \mapsto f_b(v)
\eqno ({\rm MA.30})
$$
(in differential geometry, $\Phi^b$ is the {\it Dombrowski splitting} or
{\it linearization map}, cf.\ Theorem 16.3).
Then the ``matrix" $D^b_c(f)$ of a homomorphism $f:E \to E'$
with respect to linear structures $L^b$ and $L^c$ is defined by
$$
\matrix{
\bigoplus_{\alpha >0} V_\alpha & {\buildrel \Phi^b \over \to} & E \cr
D^b_c(f) \downarrow  \phantom{xxxx} & & \phantom{x} \downarrow  f \cr
\bigoplus_{\alpha >0} V_\alpha' & {\buildrel \Phi^c \over \to} & E' \cr}
\eqno ({\rm MA.31})
$$
We have the usual ``matrix multiplication rules" $D^a_c(g \circ f)=
D^b_c(g) \circ D^a_b(f)$, $D^a_b(\id)=\id$.
Just as the individual matrix coefficients
of a linear map $f$ with respect to a basis are defined by
$a_{ij}:=\pr_i \circ f \circ \iota_j$, where $\pr_i$ and $\iota_j$
are the projections, resp.\ injections associated with the basis vectors,
we define for a homomorphism $f:E \to E'$, with respect to fixed connections
in $E$ and in $E'$,
$$
f^{\Omega | \Lambda}:= \pr_\Omega \circ f \circ \iota_\Lambda: V_\Lambda \to
V_\Omega', \quad \quad 
\matrix{
E & {\buildrel f \over \to} & E' \cr
\iota_\Lambda \uparrow \phantom{xx} & & \phantom {xx} \downarrow \pr_\Omega \cr
V_\Lambda & {\buildrel f^{\Omega|\Lambda} \over \to} & V_\Omega' \cr}
\eqno {\rm (MA.32)}
$$
where $\iota_\Lambda:V_\Lambda \to E$ is inclusion and
$\pr_\Omega:E' \to V_\Omega'$ is projection with respect to the fixed
connections. Not all matrix coefficients are interesting -- we are mainly
interested in the following:

\ssk
\item{(1)}
We say that a matrix coefficient $f^{\Omega|\Lambda}$ is {\it effective}  
if $\Omega \preceq \Lambda$. 
In this case let us write
$\Omega = \{ \Omega^1,\ldots,\Omega^r \}$ if $\ell(\Omega)=r$, $r \leq 
\ell(\Lambda)$, 
$\Omega^i = \{ \Omega^i_1,\ldots,\Omega^i_{\ell_i} \}$. Then 
the effective matrix coefficient is explicitly given by
$$
\eqalign{
f^{\Omega | \Lambda } : V_\Lambda 
 & \to V_\Omega =V_{\Omega^1}\times \ldots \times V_{\Omega^r},
\cr
(v_{\Lambda^1},\ldots,v_{\Lambda^\ell}) & \mapsto
(b^{ \Omega^1| \Lambda }(v_{\Omega_1^1},
\ldots,v_{\Omega^1_{\ell_1}}),
\ldots, b^{ \Omega^r|\Lambda }(v_{\Omega^r_1},\ldots,v_{\Omega^r_{\ell_r}} ))
\cr}
\eqno{\rm (MA.33)} 
$$
\item{(2)} 
We say that a matrix coefficient $f^{\Omega|\Lambda}$ is {\it elementary}  
if $\ell(\Omega)=1$ and $\Omega = \uline \Lambda$. (Since $\uline \Lambda
\preceq \Lambda$, elementary coefficients are effective.)
With respect to the linear structures $L^0, (L^0)'$,
the elementary matrix coefficients are just
 the multilinear maps $b^\Lambda:V_\Lambda \to V_{\uline \Lambda}'$
defining $f=f^b$:
$$
b^\Lambda = f^{\uline \Lambda | \Lambda} : \,
 V_\Lambda \to V_{\uline \Lambda}'
\quad \quad 
\matrix{
E & {\buildrel f_b \over \to} & E' \cr
\iota_\Lambda \uparrow \phantom{xx} & & \phantom {xx} \downarrow \pr_\alpha \cr
V_\lambda & {\buildrel b^\Lambda \over \to} & V_\alpha' \cr}
\eqno {\rm (MA.34)}
$$
where $\alpha = \uline \Lambda$.
In this situation, we may also use the notation $f^{\Omega|\Lambda} =:
b^{\Omega | \Lambda}$, $f^\Lambda:=b^\Lambda$.

\Proposition MA.13. {\rm (``Matrix multiplication.")}
Assume $E,E',E''$ are multilinear spaces with fixed multilinear connections
and $g:E'\to E''$ and $f:E \to E'$ are homomorphisms.
Let $\Omega \preceq \Lambda$.
Then the effective matrix coefficient  $(g \circ f)^{\Omega|\Lambda}$ is
given by
$$
(g \circ f)^{\Omega|\Lambda} = \sum_{\Omega \preceq \nu \preceq \Lambda}
g^{\Omega \vert \nu} \circ  f^{\nu | \Lambda}
$$
In particular, for elements of the group $\Gm^{1,k}(E)$ we have the composition
rule
$$
(g \circ f)^\Lambda = g^\Lambda + f^\Lambda + 
 \sum_{r=2}^{\ell(\Lambda)-1}
\sum_{\Omega \prec \Lambda \atop \ell(\Omega)=r} 
g^\Omega \circ  f^{\Omega | \Lambda}.
$$

\Proof.
Using the terminology of matrix coefficients, Formula (MA.19) from
the proof of Theorem MA.6 can be re-written
$$
(g \circ f)^\Lambda =
\sum_{\Omega \preceq \Lambda} g^\Omega \circ  f^{\Omega | \Lambda}  
=  \sum_{r=1}^{\ell(\Lambda)}
\sum_{\Omega \preceq \Lambda \atop \ell(\Omega)=r} 
g^\Omega \circ  f^{\Omega | \Lambda}.
\eqno{\rm (MA.35)} 
$$
The proposition follows from this formula.
\qed

\Corollary MA.14. The subgroup  $\Gm^{k-1,k}(E) \subset \Gm^{1,k}(E)$
is a central vector group. More precisely,
if $f \in \Gm^{k-1,k}(E)$, $g \in \Gm^{1,k}(E)$,
then 
$$
(f \circ g)^\Lambda = f^\Lambda + g^\Lambda = (g \circ f)^\Lambda .
$$

\Proof. 
Note that $f$ belongs to $\Gm^{k-1,k}(E)$ if and only if:
 $f^\Lambda = 0$ whenever
$\ell(\Lambda) \in \{ 2,\ldots,k-1 \}$. Thus the sum in the composition rule
reduces to two terms, and we easily get the claim.
\qed

\msk \nin
Under suitable assumptions, one can show that
 $\Gm^{k-1,k}(E)$ is precisely the center of $\Gm^{1,k}(E)$.
The cosets of $\Gm^{k-1,k}(E)$ in $\Gm^{1,k}(E)$ are easily described:
$f_c = f_a \circ f_b$ with $f_b \in \Gm^{k-1,k}(E)$ iff
$$
\ell(\Lambda)<k \Rightarrow c^\Lambda=a^\Lambda,\quad \quad
\ell(\Lambda)=k \Rightarrow c^\Lambda=a^\Lambda + h
$$
with some $k$-multilinear map $h:V^k \to V$.

\subheadline
{Intrinsic multilinear maps}

\msk \nin {\bf MA.15.} {\sl Intrinsic multilinear maps.}
Let $V'$ be an arbitrary $\K$-module and $E$ a $k$-linear space.
 A map $f:E \to V'$ is called
{\it (intrinsically) multilinear} if, for all $y\in E$ and $j=1,\ldots,k$,
the restrictions to the intrinsic affine spaces (hyperplanes)
$f|_{y+H_j}:y+H_j \to V'$ are  $\K$-linear.
We say that $f$ is {\it homogeneous multilinear} if, for all $i \not= j$,
 the restriction of $f$ to the affine spaces
$y + H_{ij}$ is constant (for notation, cf.\ Eqns.\ (MA.26), (MA.28)).
We denote by ${\cal M}(E,V')$, resp.\ by ${\cal M}_h(E,V')$
 the space of (homogeneous) intrinsically multilinear maps
from $E$ to $\K$ and by $\Hom(V_1,\ldots,V_k;V')$ the space 
of multilinear maps $V_1 \times \ldots \times V_k \to V'$
in the usual sense (over $\K$).

\Proposition MA.16. 
\item{(1)} A map $f:E \to V'$ is intrinsically multilinear if and only if
it can be written, with respect to a fixed linear structure $L=L^0$, in the
form
$$
f(\sum_\alpha v_\alpha)= \sum_{\ell=1}^k \sum_{\Lambda \in {\cal P}_\ell(\N_k)}
b^\Lambda(v_{\Lambda^1},\ldots,v_{\Lambda^\ell}).
$$
with multilinear maps (in the ordinary sense) $b^\Lambda:V_\Lambda \to V'$.
Then $f$ is homogeneous if and only if the only non-vanishing term belongs
to $\ell=k$, i.e., 
$$
f(\sum_\alpha v_\alpha) = b(v_{e_1},\ldots,v_{e_k})
$$
with multinear $b:V_1 \times \ldots \times V_k \to V'$.
\item{(2)}
The map
$$
\Hom(V_1,\ldots,V_k;V') \to {\cal M}_h(E,V'), \quad
f \mapsto \oline f:= f \circ \pr
$$
is a well-defined bijection.

\Proof. 
(1) Assume first that $f$ is of the form given in the claim. We decompose
$v=\sum_{\alpha_j = 1} v_\alpha + \sum_{\alpha_j=0} v_\alpha= u+w$ so that the
first term $u$
 belongs to $H_j$. Keeping the second term constant, the dependence
on the first term is linear: in fact, since $\Lambda$ is a partition of the
full set $\N_k$, each term
$b^\Lambda(v_{\Lambda^1},\ldots,v_{\Lambda^\ell})$
involves exactly one argument $v_{\Lambda^i} \in H_j$, and hence depends
linearly on this argument. 

\ssk
Conversely, assume $f:E \to V'$ is intrinsically multilinear. Then we prove
by induction on $k$ that $f$ is of the form given in the claim.
For $k=1$ the claim is trivial, and for $k=2$ see Prop.\ BA.7.
By induction, the restriction of $f$ to all maximal subgeometries
$\tilde V_{(1\ldots 101 \ldots 1)}$ (which are $k-1$-cubes of $\K$-modules,
see Section MA.11) is given by the formula from the claim.
Thus, decomposing $v=u+w \in E$ as above, we write
$f(v)=f_u(w)$, apply the induction hypothesis to $f_u$ with coefficients
depending on $u$. But the dependence on $u$ must be linear, and this yields
the claim for $f$.

\ssk
(2) This is a direct consequence of (1) since $b = f \circ \pr$.
\qed

\ssk \nin
In order to define an inverse of the map $f \mapsto \uline f$, we fix a linear
structure $L^b$ and let
$$
\iota^b:V_1 \times \ldots \times V_k \to E, \quad
(v_{e_1},\ldots,v_{e_k}) \mapsto (b\sum)  v_{e_i}
$$
be the inclusion map for the longest partition
 (sum with respect to the linear structure $L^b$). Then the inverse map is 
given by
$$
{\cal M}_h(E,V') \to \Hom(V_1,\ldots,V_k;V'), \quad
f \mapsto \uline f:= f \circ \iota^b.
$$
Summing up, the following
diagram encodes  two different ways of seeing the same object $f$:
$$
\matrix{ E& &  \cr \pr \downarrow \uparrow \iota^b & {\buildrel
\oline f \over \searrow} & \cr
V_1\times \ldots \times V_k & {\buildrel \uline f \over \to} & V' \cr}
\eqno{\rm (MA.36)}
$$
Note that the map $f \mapsto \uline f$ can be extended to a map
$$
{\cal M}(E,V') \to \Hom(V_1,\ldots,V_k;V'), \quad
f \mapsto \uline f:= f \circ \iota^b.
$$
which does depend on $L^b$. 
Finally, as in Section BA.9, in case that all $V_\alpha$ are canonically
isomorphic (see next chapter), one can define a skew-symetrization operator
from multilinear maps to homogeneous ones.

\subheadline
{\bf Tensor space and tensor functor}

\nin {\bf MA.17.} {\sl The tensor space.}
For a partition $\Lambda$ we let
$$
\hat V_\Lambda := V_{\Lambda^1}\otimes \ldots \otimes V_{\Lambda^\ell}
$$
and we define the {\it (big) tensor space of $E$}  
$$
\hat E:= \bigoplus_{\Lambda \in \Part(I)} \hat V_\Lambda =
\bigoplus_{\alpha >0} \bigoplus_{\ell=1}^{\vert \alpha \vert}
 \bigoplus_{\Lambda \in {\cal P}_\ell(\alpha)}  \hat V_\Lambda.
$$
There is a canonical map
$$
E \to \hat E, \quad (v_\alpha)_\alpha \mapsto
\sum_\Lambda v_{\Lambda^1}\otimes \ldots \otimes v_{\Lambda^\ell}.
$$
Next, to a homomorphism $f:E \to E'$ we associate a linear map
$\hat f:\hat E \to \hat E'$ in the following way:
to every multilinear map $b^\Lambda:V_\Lambda \to V_\alpha$, 
using the universal property of tensor products, we associate
a linear map  $\hat b^\Lambda:\hat V_\Lambda \to V_\alpha$.
More generally, whenever $\Omega \preceq \Lambda$, 
to the matrix coefficient 
$f^{\Omega\vert \Lambda}:V_\Lambda \to V_\Omega$ we associate a linear map
$\hat f^{\Omega\vert
\Lambda}:\hat V_\Lambda \to \hat V_\Omega$:
we can write 
$$
f^{\Omega\vert \Lambda} = b^{\Omega^1\vert \Lambda} \times \ldots \times
b^{\Omega^r\vert \Lambda}: 
\times V_{\Omega^i \vert \Lambda} 
\to \times V_{\uline {\Omega^i \vert \Lambda}}
$$
and define
$$
\hat f^{\Omega\vert \Lambda}:= 
\hat b^{\Omega^1\vert \Lambda} \otimes \ldots \otimes
\hat b^{\Omega^r\vert \Lambda}:
\otimes V_{\Omega^i \vert \Lambda} 
\to \otimes V_{\uline {\Omega^i \vert \Lambda}}.
$$
We define all other components of $\hat f:\hat E \to \hat E$ to be
zero. In other words, we simply forget all non-effective matrix
coefficients; since $f$ is determined by the elementary coefficients, this
is no loss of information. 

\Proposition MA.18. The preceding construction is functorial:
$$
\hat \id_E= \id_{\hat E},\quad \quad \quad 
\hat{g \circ f} =\hat g \circ \hat f.
$$

\Proof. The first equality is trivial. The second one follows from
the matrix multiplication rule MA.13:
the component $(\hat g \circ \hat f)^{\Omega | \Lambda}$
 of $\hat g \circ \hat f$ from $\hat V_\Lambda$ to $\hat V_\Omega$
is calculated by usual matrix multiplication:
$$
(\hat g \circ \hat f)^{\Omega | \Lambda} = \sum_\nu
\hat g^{\Omega | \nu} \circ \hat f^{\nu | \Lambda},
$$
where the sum is over all partitions $\nu$; but, by the preceding definitions,
only terms with $\Omega \preceq \nu \preceq \Lambda$ really contribute
to the sum, and hence, by MA.13, we get
the component $(\hat{g \circ f})^{\Omega\vert \Lambda}$.
\qed

\nin
We call the functor $E \to \hat E$ the {\it tensor functor}.
This functor is an equivalence of categories: 
$f$ can be recovered from $\hat f$ since all $\hat b^\Lambda$ are
components of $\hat f$.
Finally, let us add the remark that in the special case of $E$ of
tangent or bundle type (see next chapter), we may define a smaller
tensor space, essentially by taking fixed subspaces of the permutation
group,   which is functorial with respect to totally symmetric
homomorphisms; it is this smaller tensor space that plays the role of
the ``higher osculating bundle of a vector bundle" in [Po62].

\subheadline
{The affine-multilinear group}

\nin {\bf MA.19.} Just as the usual general (or special)
 {\it linear} group can be extended
to the usual general (or special) {\it affine} group, the groups $\Gm^{j,k}(E)$
may be extended to the corresponding {\it affine multilinear groups}
${\rm Am}^{j,k}(E)$. We give the basic definitions and leave details to
the reader:
assume that
for every $\alpha \in I_k$ and $\beta \subseteq \alpha$ and 
$\Lambda \in {\cal P}(\beta)$,  a multilinear map
$$
C^{\Lambda;\alpha}: V_{\Lambda^1} \times \ldots \times 
V_{\Lambda^\ell} \to V_\alpha
$$
is given. (Note that $\beta$ is recovered from $\Lambda$ as $\beta =
\uline \Lambda$ and hence can be suppressed in the notation.)
We say that $f:E \to E$ is {\it affine-multilinear} if it is of the form
$$
f(\sum_\alpha v_\alpha) =
\sum_\alpha \sum_{\beta \subset \alpha} \sum_{\Lambda \in {\cal P}(\beta)}
C^{\Lambda;\alpha}(v_{\Lambda^1},\ldots,v_{\Lambda^\ell})
$$
for some family $(C^{\Lambda;\alpha})$.
Note that for $\beta = \emptyset$, $C^{\Lambda;\alpha}$ is a constant belonging
to $V_\alpha$; in particular, for $k=1$, we get the usual affine group of
$E=V_1$. For $\beta = \alpha$, $C^{\Lambda;\alpha}=:b^\Lambda$
contributes to the ``multilinear part'' of $f$ which is defined by
$$
F: E \to E, \quad 
\sum_\alpha v_\alpha \mapsto 
\sum_\alpha \sum_{\Lambda \in {\cal P}(\alpha)}
b^{\Lambda}(v_{\Lambda^1},\ldots,v_{\Lambda^\ell}).
$$
Clearly, $F$ is a multilinear map in the sense of Section MA.5.
We say that $f$ is {\it regular} (resp.\ {\it unipotent}) if
so is $F$. 
Then the following analog of Theorem MA.6 (2) holds:
{\it 
An affine-multilinear map $f$ is invertible iff it is regular,
and then the inverse of $f$ is again affine-multilinear.}
(The proof is left to the reader.)
Thus we have defined groups ${\rm Am}^{0,k}(E)$, resp.\
${\rm Am}^{1,k}(E)$, of regular (resp.\ unipotent)
affine-multilinear maps. The structure of these groups is rather
interesting and will be investigated elsewhere.

\vfill \eject
\subheadline{SA. Symmetric and shift invariant multilinear algebra}

\msk \nin
Fibers of higher order tangent bundles are cubes of $\K$-modules
having the additional property that all or at least some of the
$V_\alpha$ are canonically isomorphic as $\K$-modules; for this reason
we call such tubes {\it of tangent type}. This situation gives
rise to interesting new symmetries: the main actors are the
{\it permutation operators} which act horizontally, on the rows of
the cube, like the castle in a game of chess, and the {\it shift operators}
which act diagonally, like the bishops. The commutants in the general
multilinear group of these actions are the {\it very symmetric}, resp.\
the {\it shift-invariant} subgroups. The intersection of these two 
groups is (in characteristic zero)
the unit group of a semigroup of truncated polynomials on $V$.

\msk \nin {\bf SA.1.} {\sl Cubes of tangent type and of bundle type.}
We say that a cube of $\K$-modules {\it is of tangent type} if all $V_\alpha$
are ``canonically" isomorphic to each other and hence to a given
model space denoted by $V$.
Formally, this is expressed by assuming that
a family $\id_{\alpha,\beta}:V_\alpha \to V_\beta$ of
$\K$-linear ``identification isomorphisms" be given, satisfying the
compatibility conditions $\id_{\alpha,\beta} \circ \id_{\beta,\gamma} =
\id_{\alpha,\gamma}$ and $\id_{\alpha,\alpha} = \id_{V_\alpha}$.
 We say that the cube {\it is of (general) bundle
type} if there are two $\K$-modules $V$ and $W$ such that each
$V_\alpha$ is canonically isomorphic to $V$ or to $W$; more
precisely, we require that all $V_\alpha$ with $\alpha_1 = 1$ are
canonically isomorphic to $W$ and all other $V_\alpha$ are canonically
isomorphic to $V$. (In a differential geometric context, we then often use
also the the notation $\alpha=(\alpha_k,\ldots,\alpha_1,
\alpha_0)$.) Thus for a cube of bundle type the
axes of the cube (MA.6) are represented by
$$
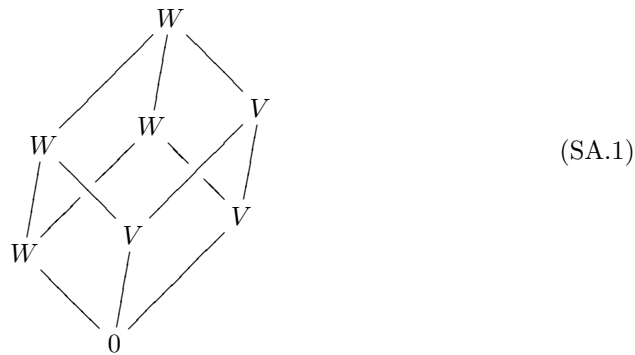
\xymatrix@dl@!0@C+10pt{
&V \ar@{-}[rr]\ar@{-}'[d][dd] & & 0 \ar@{-}[dd]
\\
V \ar@{-}[ur]\ar@{-}[rr]\ar@{-}[dd] & & V \ar@{-}[ur]\ar@{-}[dd]
\\
& W \ar@{-}'[r][rr] & & W
\\
W \ar@{-}[rr]\ar@{-}[ur] & & W \ar@{-}[ur]
}
\eqno {\rm (SA.1)}
$$
For simplicity, in the following text we treat only the case of
a cube of tangent type. In the differential geometric part, we apply
these results as well to the general bundle type, leaving the 
slight modifications in the statements to the reader.

\msk \nin {\bf SA.2.} {\sl Action of the permutation group $\Sigma_k$}.
If $E$ is of tangent type,
 fixing for the moment the linear structure $L^0$, there
is a linear action of the permutation group $\Sigma_k$ on $E$ such that,
for all $\sigma \in \Sigma_k$, $\sigma :V_\alpha \to V_{\sigma \cdot \alpha}$
agrees with the canonical identification
isomorphism:
$$
\sigma:E \to E, \quad 
\sigma(\sum_\alpha v_\alpha) = \sum_\alpha \id_{\alpha,\sigma(\alpha)}v_\alpha
\eqno({\rm SA.2})
$$
(Strictly speaking, in order to define such an action of $\Sigma_k$ 
we do not need the assumption that $E$ is of tangent type, but only
the weaker assumption that $V_\alpha \cong V_\beta$ whenever
$|\alpha| = |\beta|$.)
By our definition, the $\Sigma_k$-action is linear with respect to 
the linear structure $L^0$ on $E$.
It will not be linear with respect to all $L^b$, or, put differently,
it will not commute with the action of $\Gm^{1,k}(E)$.
We let $E^\sigma = \{ v \in E \vert \, \sigma v = v \}$ and
denote by
$E^{\Sigma_k} = \{ v \in E \vert\, \forall \sigma \in \Sigma_k: \,
\sigma v = v \}$
the fixed point space of the action of $\Sigma_k$. Then we define subgroups
of $\Gm(E)^{0,k}$ by
$$
\eqalign{
\Vsm^{j,k}(E) &:= \{ f^b \in \Gm^{j,k}(E) | \, \forall \sigma \in \Sigma_k: \,
f_b \circ \sigma = \sigma \circ f_b \}, \cr
\Sm^{j,k}(E) &:= \{ f^b \in \Gm^{j,k}(E) | \, \,
f^b(E^{\Sigma_k}) = E^{\Sigma_k} \}, \cr}
\eqno({\rm SA.3})
$$
called the {\it group of very symmetric elements},
resp.\ the {\it group of symmetric elements of $\Gm^{j,k}(E)$}. It is clear
that $\Vsm^{j,k}(E) \subset \Sm^{j,k}(E)$; we will see that 
in general this inclusion is strict.
We say that a linear structure $L^b$ is {\it very} or {\it totally
 symmetric} if
the action of $\Sigma_k$ is linear with respect to $L^b$, and that it is
{\it (weakly) symmetric} if all $E^\sigma$, $\sigma \in \Sigma_k$, are
linear subspaces of $E$ with respect to $L^b$.
Clearly, $L^b$ is (very) symmetric if so is $f^b$.
	       
\Proposition SA.3. The group $\Gm^{j,k}(E)$ is normalized by the action
of $\Sigma_k$. More precisely, 
for every multilinear family $b=(b^\Lambda)_\Lambda$ and every 
$\sigma \in \Sigma_k$,
$$
\sigma \circ f_b \circ \sigma^{-1} = f_{\sigma \cdot b},
$$
where $\sigma \cdot b$ is the multilinear family
$$
(\sigma \cdot b)^\Lambda := 
\sigma \circ b^{\sigma^{-1} \cdot \Lambda} \circ 
(\sigma^{-1} \times \ldots \times \sigma^{-1}).
$$
It follows that $\sigma \cdot L^b$ is again a multilinear connection on $E$;
more precisely, $\sigma\cdot L^b=L^{\sigma \cdot b}$.

\Proof. By a
direct calculation, making the change of variables $\beta=\sigma .\alpha$
and $\Omega:=\sigma.\Lambda$, we obtain
$$
\eqalign{\sigma f_b \sigma^{-1} (\sum_\alpha v_\alpha)&=
\sigma f_b(\sum_\alpha \id_{\alpha,\sigma^{-1}(\alpha)}(v_\alpha)) \cr
&= 
\sigma\Bigl(\sum_\alpha \sum_{\Lambda \in {\cal P}(\alpha)} 
b^{\Lambda} \bigl(\id_{\sigma.\Lambda,\Lambda}(v_{\sigma(\Lambda)} )
\bigr) \Bigr) \cr
&= 
\sum_\alpha \sum_{\Lambda \in {\cal P}(\alpha)} 
\id_{\alpha,\sigma.\alpha}
b^{\Lambda} (\id_{\sigma.\Lambda,\Lambda}(v_{\sigma(\Lambda)} )) \cr
&= 
\sum_\beta \sum_{\Omega \in {\cal P}(\beta)} 
\id_{\sigma^{-1}.\beta,\beta}
b^{\sigma^{-1} \Omega} 
(\id_{\Omega,\sigma^{-1}\Omega}(v_{\Omega} )) \cr
&= f_{\sigma.b}(\sum_\beta v_\beta). \cr}
$$
Finally, since $\sigma \cdot L^0 =L^0$, it follows that
$$
\sigma \cdot L^b = \sigma \cdot f^b \cdot L^0 = 
(\sigma \circ f^b  \circ \sigma^{-1}) \cdot L^0=
f^{\sigma \cdot b} L^0 = L^{\sigma \cdot b}
\qeddis

\Corollary SA.4. The following are equivalent:
\item{(1)} $f_b \in \Vsm^{1,k}(E)$ 
\item{(2)} $\forall \sigma \in\Sigma_k$, $L^b=\sigma \cdot L^b$
\item{(3)}  the multilinear family $b$
is $\Sigma_k$-equivariant in the sense that, for all $\alpha \in I_k$, all
$\Lambda \in {\cal P}(\alpha)$ and all $\sigma \in \Sigma_k$, 
we have $(\sigma \cdot b)^\Lambda = b^\Lambda$, i.e.,
$b^{\sigma \cdot \Lambda}=\sigma \circ b^\Lambda \circ 
(\sigma^{-1} \times \ldots \times \sigma^{-1})$:
$$
\matrix{
V_{\Lambda^1} \times \ldots \times V_{\Lambda^\ell} &
{\buildrel b^\Lambda \over \to} & V_\alpha \cr
\sigma \times \ldots \sigma \downarrow \phantom{
\sigma \times \ldots \sigma} & & \phantom{\sigma} \downarrow \sigma \cr
V_{\sigma(\Lambda^1)} \times \ldots \times V_{\sigma(\Lambda^\ell)} &
{\buildrel b^{\sigma \cdot \Lambda} \over \to} & V_{\sigma \cdot \alpha} \cr}
\eqno{\rm (SA.4)}
$$
\qed

\nin The condition from the corollary implies that the $b^\Lambda$ are
{\it symmetric multilinear maps}. Let us explain this.
We assume that $f_b$ is very symmetric, i.e.\ (SA.4) commutes for
all $\Lambda$ and $\sigma$, and 
let us fix $\Lambda \in {\cal P}(\alpha)$ such
 that, for some $i\not= j$, we have $|\Lambda^i| = |\Lambda^j|$.
Then there exists $\sigma$ such that $\sigma(\Lambda^i)=\Lambda^j$ and
$\sigma(\Lambda^m)=\Lambda^m$ for all other $m$, whence
$\sigma(\Lambda)=\Lambda$ and hence $\sigma(\alpha)=\alpha$.
The action of $\sigma$ on $V_\Lambda$ 
is determined by exchanging $V_{\Lambda^i}$
and $V_{\Lambda^j}$ and fixing all other $V_{\Lambda^m}$.
Therefore (SA.6) can be interpreted in this situation 
by saying that $b^\Lambda$ is symmetric with respect to exchange of arguments
from $V_{\Lambda^i}$ and $V_{\Lambda^j}$.
In other words,
 the multilinear map $b^\Lambda$ has all possible
permutation symmetries in the sense that it
is invariant under the subgroup
$\{ \sigma \in \Sigma_k \vert \, \sigma \cdot \Lambda = \Lambda \}$ of $\Sigma_k$.
For instance, if $k=2$, $f_b$ is very symmetric iff $b^{01,10}$ ``is a
symmetric bilinear map $V\times V \to V$'',
and if $k=3$ (cf.\ Section MA.8), then $f_b$ is very symmetric iff

\ssk
\item{(a)} the three bilinear maps $b^{001,110},b^{101,010}, b^{011,100}$
are conjugate to each other:
$(23) \cdot b^{101,010}=b^{011,100}= (13) \cdot b^{110,001}$,
\item{(b)} the three bilinear maps $b^{001,010},b^{001,100}, b^{010,100}$
are conjugate to each other:
$(23) \cdot b^{001,010}=b^{001,100}= (12) \cdot b^{010,100}$
 {\it and} are symmetric as bilinear maps,
\item{(c)} the trilinear map $b^{001,010,100}$ is a symmetric trilinear map.

\ssk
\nin Next we will give a sufficient condition for $f_b$
to be weakly symmetric. For $k=2$, all
$f_b \in \Gm^{1,2}(E)$ are weakly symmetric (in the notation of BA.1,
$f_b(v,v,w)=(v,v,w+b(v,v))$), and for $k=3$, the reader may check by 
hand that $f_b$ is weakly symmetric if
the preceding condition (a) holds in combination with

\ssk
\item{(b')} the three bilinear maps $b^{001,010},b^{001,100}, b^{010,100}$
are conjugate to each other:
$(23) \cdot b^{001,010}=b^{001,100}= (12) \cdot b^{010,100}$.

\ssk
\nin Our sufficient criterion will be based on the simple observation
that $E^{\Sigma_k} = \cap_{\tau \in T} E^\sigma$
where $T$ is any set of generators of $\Sigma_k$; for instance, we may
choose $T$ to be the set of transpositions.

\Proposition SA.5.
Let $f_b \in \Gm^{1,k}(E)$,  $\tau \in \Sigma_k$ and consider the following
conditions $(1)$ and $(2)$.
\item{(1)}  $f_b(E^\tau)=E^\tau$
\item{(2)} For all $\alpha \in I_k$ 
such that $\tau \cdot \alpha \not=\alpha$ and for
all partitions $\Lambda \in {\cal P}(\alpha)$, the condition        
$(\tau \cdot b)^\Lambda =b^\Lambda$ holds.
\ssk                                        
\nin Then $(2)$ implies $(1)$. It follows that,
if $f_b$ satisfies {\rm Condition (2)} for all elements $\tau$
of some set $T$ of generators of $\Sigma_k$,
then $f_b \in \Sm^{1,k}(E)$.

\Proof. 
By definition of $E^\tau$,
$$
v= \sum_\alpha v_\alpha \in E^\tau \, \Leftrightarrow \,
\forall \alpha \in I_k: \,\id_{\alpha,\tau \cdot \alpha}(v_\alpha) = 
v_{\tau \cdot \alpha}.
$$
It follows that $f^b(v)$ belongs to $E^\tau$ if and only if
$\tau((f^b v)_\alpha) =\id_{\alpha,\tau \cdot \alpha}
(f^b v)_{\tau \cdot \alpha}$ for all $\alpha$, i.e.
$$
\eqalign{
\id_{\alpha,\tau \cdot \alpha}(v_\alpha) + 
\sum_{\Lambda \in {\cal P}(\alpha)} b^\Lambda
(v_{\Lambda^1},\ldots,v_{\Lambda^\ell})) & =
v_{\tau \cdot \alpha} + \sum_{\Omega \in {\cal P}(\tau \cdot \alpha)}
b^\Omega(v_{\Omega^1},\ldots,v_{\Omega^\ell}) \cr
& = v_{\tau \cdot \alpha} +
\sum_{\Lambda' \in {\cal P}(\alpha)}
b^{\tau \cdot \Lambda'}(v_{(\tau \cdot \Lambda')^1},\ldots,v_{
(\tau \cdot \Lambda')^\ell})
\cr}
$$
where we made the change of variables $\Lambda':=\tau^{-1} \cdot \Omega$.
If we assume that $v\in E^\tau$, then this condition is
equivalent to: for all $\alpha \in I_k$,
$$
\sum_{\Lambda \in {\cal P}(\alpha)} \id_{\alpha,\tau \cdot \alpha}(b^\Lambda
(v_{\Lambda^1},\ldots,v_{\Lambda^\ell})) =
\sum_{\Lambda' \in {\cal P}(\alpha)}
b^{\tau \cdot \Lambda'}(v_{(\tau \cdot \Lambda')^1},\ldots,v_{
(\tau \cdot \Lambda')^\ell})
\eqno{\rm (SA.5)}
$$
If $\tau \cdot \alpha = \alpha$, then, since $\id_{\alpha,\alpha}= 
\id_{V_\alpha}$, this condition is automatically
fulfilled, and if  $\tau \cdot \alpha \not= \alpha$, then
Condition (2) guarantees that it holds. Thus (SA.5) holds for all
$\alpha$, and (1) follows. The final conclusion is now immediate.
\qed

\nin One can show that,
if in the preceding claim that $\tau$ is a transposition,
then  (1) and (2) are equivalent.

\msk \nin {\bf SA.6.} {\sl Curvature forms.}  
Let $L=L^b$ be a multilinear connection on $E$. We have seen that
$L$ and $\sigma \cdot L$ are multilinearly related (Prop.\ SA.3), and hence
there exists a unique element $R:=R^{\sigma,L} \in \Gm^{1,k}(E)$ such that
$$
R^{\sigma,b} \cdot L =\sigma \cdot L.
\eqno{\rm (SA.6)}
$$
This element is called the {\it curvature operator} and is given by
$$
R^{\sigma,L}= \sigma \circ f_b \circ \sigma^{-1}\circ (f_b)^{-1} 
=  f_{\sigma \cdot b} \circ (f_b)^{-1} .
\eqno{\rm (SA.7)}
$$
In fact, $\sigma \circ f_b \circ \sigma^{-1}
\circ (f_b)^{-1} \cdot L^b = \sigma \cdot f_b \cdot  \sigma^{-1} L^0 = 
\sigma \cdot L^b$ since
$L^0$ is very symmetric, i.e.\ $\sigma \cdot L^0 = L^0$.
Directly from the definition of the curvature operators, we get the
{\it cocycle relations}
$$
R^{\sigma \circ \tau}=\sigma \circ R^\tau \circ \sigma^{-1}
 \circ R^\sigma, \quad \quad
R^{\id}=\id.
\eqno{\rm (SA.8)}
$$
The elementary matrix coefficients of $R$ with respect to $L$,  
$$
R^\Lambda:=(R^{\sigma,L})^\Lambda:V_\Lambda \to V_{\uline\Lambda}
\eqno{\rm (SA.9)}
$$
are called {\it curvature forms (of $\sigma,L$)}. 

\msk \nin {\bf SA.7.} {\sl Symmetrizability.}
We say that a multilinear connection $L=L^b$ (resp.\ an endomorphism 
$f_b \in \Gm^{1,k}(E)$) is {\it symmetrizable} if all
its curvature operators $R^{\sigma,L}$, $\sigma \in \Sigma_k$,
 belong to the central vector group  $\Gm^{k-1,k}(E)$, i.e.\
$(\sigma \circ f_b \circ \sigma^{-1}\circ (f_b)^{-1})^\Lambda =0$
whenver $\ell(\Lambda)<k$, or, equivalently,
$$ 
\forall \Lambda : \, \ell(\Lambda)<k \, \Rightarrow
\, (\sigma \cdot b)^\Lambda = b^\Lambda. 
$$

\Lemma SA.8.
Assume $L$ is symmetrizable. Then the orbit
${\cal O}_L:=\Gm^{k-1,k}(E).L$ of $L$ under the action of the vector
group $\Gm^{k-1,k}(E)$ is an affine space over $\K$, isomorphic to 
$\Gm^{k-1,k}(E)$.  This affine space  is stable under the
action of the group $\Sigma_k$ which acts by affine maps on ${\cal O}_L$.

\Proof.
All orbits of the action of the vector group $\Gm^{k-1,k}(E)$ on the 
space of connections have
trivial stabilizer and hence are affine spaces. Let us prove that ${\cal O}_L$
is stable under the action of $\sigma \in \Sigma_k$:
for $f_c \in \Gm^{k-1,k}(E)$,
$$
\sigma . f_c.L = (\sigma \circ f_c \circ \sigma^{-1}) . \sigma.L=
f_{\sigma\cdot c} . R^{\sigma,L}.L
$$
belongs to the orbit $\Gm^{k-1,k}(E).L ={\cal O}_L$ (since $R^{\sigma,L}
\in \Gm^{k-1,k}(E)$ by assumption of symmetrizability, and 
$f_{\sigma\cdot c} \in \Gm^{k-1,k}(E)$ by Prop.\ SA.3).
Next, we show that the last expression depends affinely on $c$:
indeed, conjugation by $\sigma$ is a {\it linear} automorphism of
the vector group $\Gm^{k-1,k}(E)$
(the action $c^\Lambda \mapsto (\sigma.c)^{\Lambda}$, for the only non-trivial
component $\Lambda =\{ \{1\},\ldots,\{ k\} \}$, is the usual action
of the permutation group on $k$-multilinear maps $V^k \to V$, which
is linear), and hence $f_c.L \mapsto f_{\sigma.c}.L$ is a linear
automorphism of the linear space $({\cal O}_L,L)$. 
Composing with the translation by $R^{\sigma,L}$, we see
that $f_c.L \mapsto R^{\sigma,L}.f_{\sigma.c}.L=\sigma.(f_c.L)$ is affine.
\qed

\Corollary SA.9.
Assume that $L=L^b$ is symmetrizable and that the integers $2,\ldots,k$ are
invertible in $\K$.
Then the barycenter of the $\Sigma_k$-orbit of $L$,
$$
L_{\Sym} := {1 \over k!} \sum_{\nu \in \Sigma_k} \nu.L \in {\cal O}_L
\eqno{\rm (SA.10)}
$$
is well-defined and is a very symmetric linear structure on $E$.
Assume that, in addition, $L$ is invariant under $\Sigma_{k-1}$. Then
$$
L_{\Sym} = {1 \over k} \sum_{j=1}^k (12\cdots k)^j.L.
\eqno{\rm (SA.11)}
$$

\Proof.
For any finite group $\Gamma$, acting affinely on an affine space over
$\K$, under the suitable assumption on $\K$, the barycenter of $\Gamma$-orbits
is well-defined and is fixed under $\Gamma$ since affine maps are
compatible with barycenters. This can be applied to the situation of
the preceding lemma in order to prove the first claim. 
Moreover, if $\Sigma_{k-1}.L=L$, then
$\Sigma_k.L=\Z/(k).L$, where $\Z/(k)$ is the cyclic group generated by
the permutation $(12\cdots k)$, and we get (SA.11).
\qed

\nin
Summing up,  $P:={1 \over k!} \sum_{\nu \in \Sigma_k} \nu$
may be seen as a well-defined projection operator on the space of all
symmetrizable linear structures, acting affinely on orbits of the type
of ${\cal O}_L$ and having as image the space of all very symmetric linear
structures. The fibers of this operator are linear spaces, isomorphic to
the kernel of the symmetrization operator on multilinear maps 
$V^k \to V$.

\msk \nin {\bf SA.10.} {\sl Bianchi's identities.}
Assume $L$ satisfies the assumptions leading to Equation (SA.11).
Then $R^{L_{\Sym}}=0$ since $L_{\Sym}$ is very symmetric;
on the other hand,
$$
0=R^{L_{\Sym}} = {1 \over k} \sum_{j=1}^k R^{(12\cdots k)^j.L} =
{1 \over k} \sum_{j=1}^k (12\cdots k)^j.R^{L}.
\eqno{\rm (SA.12)}
$$
Thus the sum of the cyclically permuted $k$-multilinear maps 
corresponding to $R^L$ vanishes; for $k=3$ this is Bianchi's identity.
Also, if $\sigma$ is a transposition, then $R^\sigma$ has the
familier skew-symmetry of a curvature tensor: this follows from the
cocylce relations
$$
\id = R^{\id} = R^{\sigma^2} = \sigma \circ R^\sigma \circ \sigma
\circ R^\sigma;
$$
but under the assumption of SA.9, the composition of the curvature
operators is just vector addition in the vector group $\Gm^{k-1,k}(E)$,
whence
$$
0 = \sigma \circ R^\sigma \circ \sigma + R^\sigma,
\eqno{\rm (SA.13)}
$$
expressing that the multilinear map 
$R^\sigma$ is skew-symmetric in the $i$-th and $j$-th component 
if $\sigma=(ij)$.

\subheadline{Shift operators and shift invariance}

\msk \nin
{\bf SA.11.} {\sl Shift of the index set.}
For $i \not= j$, the 
{\it (elementary) shift of the index set $I_k \cong 2^{\N_k}$
(in direction $j$ with basis $i$)} is defined by
$$
s:=s_{ij}:2^{\N_k} \to 2^{\N_k}, \quad
s(A) = \Big\{ \matrix{A & {\rm if} \, i \notin A,
 \cr
A \cup \{ j \} & {\rm if} \, i \in A, j \notin A , 
 \cr \emptyset & {\rm if} \, i,j \in A .\cr}
\eqno{\rm (SA.14)}
$$
In terms of multi-indices,
$$
s:=s_{ij}:I_k \to I_k, \quad
s(\alpha) = \Big\{ \matrix{\alpha & {\rm if} \, \alpha_i =0,
 \cr
\alpha +e_j & {\rm if} \, \alpha_i=1, \alpha_j=0, 
 \cr 0 & {\rm if} \, \alpha_i=1, \alpha_j=1 .\cr}
\eqno{\rm (SA.15)}
$$
Then $s\circ s(A)=A$ if $i\in A$ and $s \circ s(A)=\emptyset$ else;
hence $s \circ s$ is idempotent and may be considered as a ``projection
with base $i$''.

\msk \nin {\bf SA.12.} {\sl Shift operators.}
We continue to assume that $E$ is of tangent type. Then, with respect
to the linear structure $L^0$, there is a unique linear map
$S:=S_{ij}:E \to E$ corresponding to the elementay shift $s$ of the
index set, called an {\it (elementary) shift operator}
$$
S_{ij}:E \to E, \quad
S_{ij}( v_\alpha) = \id_{\alpha,s(\alpha)} (v_\alpha) =
\Big\{ \matrix{ v_\alpha & {\rm if} \, \alpha_i =0,
 \cr
\id_{\alpha,\alpha+e_j}(v_{\alpha}) & {\rm if} \, \alpha_i=1, \alpha_j=0, 
 \cr 0 & {\rm if} \, \alpha_i=1, \alpha_j=1 .\cr}
\eqno{\rm (SA.16)}
$$
The components of $S_{ij} v$ are then given by
$$
(S_{ij}(v))_\beta = \pr_\beta(S_{ij} v) = \pr_\beta (\sum_\alpha 
S_{ij} v_\alpha) =
\Big\{ \matrix{ v_\beta & {\rm if} \, \beta_i =0,                \cr
\id_{\beta-e_j,\beta}(v_{\beta-e_j}) & {\rm if} \, \beta_i=1= \beta_j, \cr 
0 & {\rm if} \, \beta_i=1, \beta_j=0 .\cr}
\eqno{\rm (SA.17)}
$$
We say that $S_{ij}$ is a {\it positive shift} if $i>j$.
Note that $S \circ S=p_{1-e_i}p_{e_i^*}$ is the projection onto a maximal
subgeometry considered in Section MA.11 (Equation (MA.25)), and that
 $S_{ij}$ preserves the factors of the decomposition 
$E = \tilde V_{1-e_i} \oplus H_i$ (Section MA.11); it acts trivially on
the first factor and by a two-step nilpotent map $\eps_j$
on the second factor,
as may be summarized by the following matrix notation:
$$
S_{ij} = \pmatrix{ \1 &   \cr  & {0 \atop 1}{0 \atop 0} \cr} =
\pmatrix{ \1 & 0 \cr 0 & \eps_j \cr}.
$$
 The shift operators are not intrinsic
with respect to the whole general linear
group $\Gm^{1,k}(E)$ but only with respect to
a suitable subgroup.

\Proposition SA.13. 
The map $f^b$ commutes with the shift $S=S_{ij}$ if and only
if the following holds: for all $\alpha \in I_k$ such that
$\alpha_i =1$ and $\alpha_j=0$ and for all $\Lambda \in {\cal P}(\alpha)$
the following diagram commutes:
$$
\matrix{
V_{\Lambda^1} \times \ldots \times V_{\Lambda^\ell} &
{\buildrel b^\Lambda \over \to} & V_\alpha \cr
S \times \ldots \times S \downarrow \phantom{
\sigma \times \ldots \sigma} & & \phantom{\sigma} \downarrow S \cr
V_{\Lambda^1} \times \ldots \times V_{\Lambda^m \cup \{ i\}} \times \ldots
\times V_{\Lambda^\ell} &
{\buildrel b^{s \cdot \Lambda} \over \to} & V_{s \cdot \alpha} \cr}
\eqno({\rm SA.18})
$$
where $m$ is the (unique) index such that $i\in \Lambda^m$. In this
case, 
$$
s\Lambda = \{ s\Lambda^1,\ldots,s\Lambda^\ell \} = \{ \Lambda^1,
\ldots,\Lambda^m \cup \{ i \},\ldots,\Lambda^\ell \}
$$ 
is a partition of
$s \alpha =\alpha \cup \{i\}$, so that $b^{s \Lambda}$ is a matrix
coefficient of $f^b$. 

\Proof. On the one hand,
$$
\eqalign{
f_b(S_{ij}(v)) & = \sum_\alpha \sum_{\Lambda \in {\cal P}(\alpha)} 
b^\Lambda ((S_{ij} v)_{\Lambda^1},\ldots,(S_{ij} v)_{\Lambda^\ell} ) \cr
& = \sum_{\alpha \atop \alpha_i =0} b^\Lambda(v_{\Lambda^1},\ldots,
v_{\Lambda^\ell}) + 
\sum_{\alpha \atop \alpha_i =1=\alpha_j} 
b^\Lambda ((S_{ij} v)_{\Lambda^1},\ldots,(S_{ij} v)_{\Lambda^\ell} ) 
=: A + B, \cr}
$$
where $A$ corresponds to sum over all $\alpha$ with $\alpha_i =0$
(in this case we have also $(\Lambda^r)_i =0$ for all $i$, and hence
the shift $s_{ij}$ fixes these indices), $B$ corresponds to the sum
over all $\alpha$ with $\alpha_i = 1 = \alpha_j$, and the term
corresponding to the sum over all $\alpha$ with $\alpha_i=1$, $\alpha_j=0$
is zero since then there exists $r$ with $(\Lambda^r)_i=1$ and
$(\Lambda^r)_j=0$, and hence $(S_{ij} v)_{\Lambda^r} = 0$.
On the other hand, by definition of $S_{ij}$,
$$
\eqalign{
S_{ij}(f_b(v)) & = 
S_{ij}(\sum_\alpha \sum_{\Lambda \in {\cal P}(\alpha)}
b^\Lambda(v_\Lambda)) = 
\sum_\alpha \sum_{\Lambda \in {\cal P}(\alpha)}
S_{ij} (b^\Lambda(v_\Lambda)) \cr 
& = 
\sum_{\alpha \atop \alpha_i=0} \sum_{\Lambda \in {\cal P}(\alpha)}
b^\Lambda(v_\Lambda) + 
\sum_{\alpha \atop \alpha_i=1,\alpha_j=0} \sum_{\Lambda \in {\cal P}(\alpha)}
\id_{\alpha,\alpha+e_j}(b^\Lambda(v_\Lambda)) =: A + B', \cr}
$$
where the first term $A$ coincides with the first term of the preceding
calculation.
Comparing, we see that $f_b \circ S_{ij} = S_{ij} \circ f_b$ iff
$B=B'$. Let us have a look at $B$: if $\alpha_i=1=\alpha_j$, then
for a given $\Lambda \in {\cal P}(\alpha)$ two cases may arise:

\ssk
\item{(a)} either there exists $m=m(\Lambda)$ such that $i,j \in \Lambda^m$;
in this case 
$$
b^\Lambda ((S_{ij} v)_{\Lambda^1},\ldots,(S_{ij} v)_{\Lambda^\ell} ) =
b^\Lambda(v_{\Lambda^1},\ldots,\id_{\Lambda^m-e_j,\Lambda^m}
(v_{\Lambda^m-e_j}),\ldots,v_{\Lambda^\ell}),
$$ 
\item{(b)} there exists $m \not= n$ such that $i \in \Lambda^m$, $j \in
\Lambda^n$. In this case $(S_{ij} v)_{\Lambda^m} =0$, and the corresponding
term is zero. 

\ssk \nin Thus we are left with Case (a). In this case we let
$\gamma:=\alpha - e_j$ and define a partition
 $\Lambda':=\{ \Lambda^1,\ldots,\Lambda^m - e_j,
\ldots,\Lambda^\ell \}$; then $\Lambda' \in {\cal P}(\gamma)$ and
$\Lambda = s(\Lambda')$, $\alpha = s(\gamma)$, and the condition
$B=B'$ is seen to be equivalent to the condition from the claim.
\qed

\msk \nin
Let us define the following subgroups of $\Gm^{1,k}(E)$:
$$
\eqalign{
\Shi^{\ell,k}(E) &:=\{ f_b \in \Gm^{\ell,k}(E) | \, \forall i \not= j:
\, S_{ij} \circ f_b = f_b \circ S_{ij} \},\cr
\Shi^{\ell,k}_+(E) &:=\{ f_b \in \Gm^{\ell,k}(E) | \, \forall i > j:
\, S_{ij} \circ f_b = f_b \circ S_{ij} \}. \cr} 
\eqno{\rm (SA.19)}
$$
For $k=2$ the only  target space for non-trivial multilinear 
maps is $V_{11}$ and the only non-trivial partition is
$\Lambda = ((01),(10))$, which corresponds to Case (b) in the preceding 
proof, whence $\Gm^{1,2}(E) = \Shi^{1,2}(E)$ (cf.\ also Section BA.4).
We will see that 
the condition from Prop.\ SA.13 really has non-trivial consequences if
$k \geq 3$.

\Theorem SA.14.
Assume that $k>2$ and that the multilinear
connection $L$ is invariant under all 
shifts. Then $L$ is symmetrizable, and (if the integers are invertible in 
$\K$) $L_{\Sym}$ is again invariant under all shifts.
Equivalently, if $f_b \in \Shi^{1,k}(E)$, then $f^b$ is symmetrizable, 
and $f^{\Sym(b)}$ belongs to $\Shi^{1,k}(E) \cap \Vsm^{1,k}(E)$.

\Proof.
We prepare the proof by some remarks on shifts.
The shift maps $s_{ij}:I_k \to I_k$ are ``essentially injective" in the
sense that, if $\beta \not= 0$, then there is at most one $\alpha$
such that $s_{ij}(\alpha)=\beta$, and similarly for partitions.
Then we write $\alpha {\buildrel s_{ij} \over \to} \beta$, having in
mind that $\alpha$ is uniquely determined by $\beta$. Now let us
look at the following pattern of correspondences of (ordered) partitions 
of length $2$ for $k=3$:

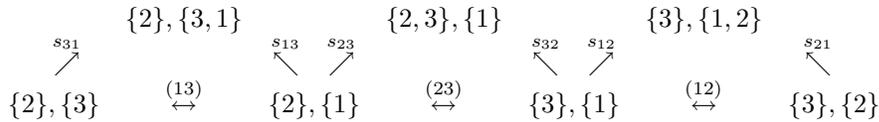
\msk
$$
\matrix{
 & \{2\}, \{ 3,1\} & & \{2,3\},\{1\}& & \{3\}, \{1,2\} & \cr
\phantom{ \nwarrow}  
{\buildrel s_{31} \over \nearrow} & & {\buildrel s_{13} \over \nwarrow} 
\quad {\buildrel s_{23} \over \nearrow} & & 
 {\buildrel s_{32} \over \nwarrow}\quad 
 {\buildrel s_{12} \over \nearrow} &  & 
{\buildrel s_{21} \over \nwarrow} \,
\phantom{ \nwarrow}  \cr 
\{2\},\{3\} & {\buildrel (13) \over \leftrightarrow} & \{2\},\{1\} &
 {\buildrel (23) \over \leftrightarrow} & \{3\}, \{1\} &
  {\buildrel (12) \over \leftrightarrow} & \{3\},\{2\} \cr}
$$

\msk
\nin If $f_b$ is invariant under all shifts, then the maps $b^\Lambda$ 
corresponding to the partitions of the preceding pattern are all conjugate
to each other under the respective shifts. But this implies that that
they are also conjugate to each other under the transpositions as indicated
in the bottom line, and by composing the three transpositions, we see that
also the effect of the transposition (32) on the partition $(\{2\},\{3\})$
is induced by the shift-invariance. 
(At this point, the reader may 
recognize the famous ``braid lemma" in disguise.)
Also, the partitions from the top line are conjugate to each other
under permutations (here the order is not important since their two 
elements are already distinguished by their cardinality). Summing up,
for $\ell(\Lambda)=2$,
all bilinear maps $b^{\sigma.\Lambda}$, $\sigma \in \Sigma_3$, are
conjugate to $b^\Lambda$ under $\sigma$, and according to SA.7,
this means that $f_b$ is symmetrizable. 
This argument can be generalized to arbitrary $k>2$ and $\Lambda$ with
$\ell(\Lambda)<k$ because for all such $\Lambda$ there exists $\Omega$ 
such that one of the relations $\Lambda {\buildrel s_{ij} \over \to}
\Omega$ or $\Omega {\buildrel s_{ij} \over \to} \Lambda$ holds
(the only partition for which no such relation holds is the longest one,
$\Lambda = \{ \{ 1\},\ldots,\{ k \} \}$). Then, using the above arguments,
one sees that $b^\Lambda$ and $b^{\tau \cdot \Lambda}$ for $\tau=(i,j)$
are conjugate, and so on, leading to invariance under a set of generating
transpositions of $\Sigma_k$. 

Finally, note that, if $L=L^b$, then $L_{\Sym}=L^{\Sym(b)}$, where
$\Sym(b)$ is obtained by symmetrising the highest coefficient of $b$.
By assumption, $f_b$ commutes with all shifts; then also all elements
of the $\Gm^{k-1,k}(E)$-coset commute with all shifts since elements
of $\Gm^{k-1,k}(E)$ commute with all shifts, by Prop.\ SA.13. 
Since $f_{\Sym(b)}$ belongs to this coset, it commutes with all shifts.
Summing up, $f_{\Sym(b)} \in \Shi^{1,k}(E) \cap \Vsm^{1,k}(E)$.
\qed

\msk \nin {\bf SA.15.} {\sl The intersection 
$\Vsm^{1,k}(E) \cap \Shi^{1,k}(E)$.}
Assume $f_b \in \Vsm^{1,k}(E) \cap \Shi^{1,k}(E)$.
Then all components $b^\Lambda$ of $f_b$ for with fixed length
 $\ell(\Lambda)=\ell$ are conjugate to each other under a combination of shifts
and permutations: by these operations, every
partition $\Lambda$ of length $\ell$
 can be transformed into the the ``standard partition''
$\{ \{ 1\},\ldots,\{ \ell \} \}$ of length $\ell$.
Thus $f_b$ is determined by a collection of $k$ symmetric multilinear maps
$d^j: V^j \to V$ where $V$ is the model space for all $V_\alpha$.
Let us denote by $\eps^\alpha:V \to V_\alpha$ the isomorphism with
the model space (so that $\id_{\alpha,\beta} = \eps^\alpha \circ
(\eps^\beta)^{-1}$). Then $f_b$ is given by
$$
f_b(v) = \sum_\alpha \eps^\alpha (\sum_{j=1}^{\vert \alpha \vert}
\sum_{\Lambda \in {\cal P}_j(\alpha)} d^j (v_\Lambda) ).
$$

\Theorem SA.16.
\item{(1)} 
A multilinear map $f_b \in \Gm(E)$ belongs to
$ \Shi^k(E)\cap \Vsm^k(E)$ if and only if there are symmetric 
multilinear maps $d^\ell:V^\ell \to V$ such that, if $\ell(\Lambda)=\ell$,
$$
b^\Lambda(\eps^{\Lambda_1} v_{\Lambda_1},\ldots,\eps^{\Lambda_\ell}
v_{\Lambda_\ell}) = \eps^{\uline \Lambda}
d^l(v_{\Lambda_1},\ldots, v_{\Lambda_\ell}) .
$$
\item{(2)} If $2,\ldots,k$ are invertible in $\K$, then
a multilinear map $f_b \in \Gm(E)$ belongs to
$ \Shi^k(E)\cap \Vsm^k(E)$ if and only if there 
is a polynomial map $F:V \to V$ of degree at most $k$ such that $f_b$
 is given by the following purely algebraic analog of {\rm Theorem 7.5}:
$$
f_b(v)=\sum_\alpha \sum_{\Lambda \in {\cal P}(\alpha)}
(d^{\ell (\Lambda)} F)(0)(v_{\Lambda}),
$$
where $(d^\ell F)(0)$ is the $\ell$-th differential of $F$ at the origin.
We write this  also in the form $f_b = (T^k F)_0$.

\Proof.
The condition from (1) is necessary since, as remarked above, 
by permutations and shifts every
partition $\Lambda$ can be transformed into the the standard partition
$\hat \Lambda:=\{ \{ 1\},\ldots,\{ \ell \} \}$ of length $\ell = \ell(\Lambda)$.
Then we let $d^\ell :=b^{\hat \Lambda}$.
The condition is also sufficient: with arbitrary choices of symmetric
multilinear maps $d^\ell$, the formula defines a multilinear map $f_b$;
this map satisfies the shift-invariance and the permutation-invariance
conditions (SA.13 and SA.4) and hence belongs to the intersection of
both groups.

(2) Given $f_b \in  \Shi^k(E)\cap \Vsm^k(E)$, we define $d^\ell$ as
in (1) and define a polynomial map $F:V \to V$ by
$$
F(x):=\sum_{\ell=1}^k {1 \over \ell !} d^\ell(x,\ldots,x).
$$
Then we have $(d^\ell F)(0)=d^\ell$, and hence $F$ satisfies the requirement from
(2). Conversely, given $F$, we let $d^\ell :=(d^\ell F)(0)$ and then define
$f_b$ as in Part (1).
\qed

\nin 
The preceding theorem shows that the group 
$\Vsm^{1,k}(E) \cap \Shi^{1,k}(E)$ 
is a purely algebraic version of the ``truncated polynomial group'',
denoted in [KMS93, Ch.\ 13.1] by $B_1^k$. Its Lie algebra (see next
chapter, or [KMS93]) is the truncated (modulo degre $> k$)
Lie algebra of polynomial vector fields on $V$.

\vfill \eject

\subheadline
{PG. Polynomial groups}

In the following, we will use the notion of {\it polynomial mappings between
modules over a ring or field} $\K$, which we define as finite sums 
$p=\sum_{k=0}^m p_k$ of homogeneous polynomial mappings $p_k:V \to W$,
where $p_k(x)=b_k(x,\dots,x)$ is given by evaluating  a
(not necessary symmetric) $\K$-multilinear map $b_k:V^k \to W$ on the
diagonal. However, in general it is not possible to recover the
homogeneous parts $p_k$ from the mapping $p:V \to W$.
Therefore one either has to work with a more formal concept of 
polynomials (as developed by N.\ Roby in [Rob63]), or one has to make
assumptions on the base ring. In order to simplify our presentation, we 
will choose the second option and
make rather strong assumptions on $\K$: in Sections PG.1 -- PG.4 we will
assume that $\K$ is an infinite field, and from Section PG.5 on we will assume
that $\K$ is a field of characteristic zero. 
Under these assumptions, the homogeneous parts can be recovered from $p$
by algebraic
differential calculus for polynomials -- this is classical, and it can also
be recovered in the general framework of differential calculus developed
in [BGN04] since polynomial maps over infinite fields
form a ``$C^0$-concept'' in the sense of [BGN04] (cf.\ also Appendix G.1).
Another choice of
assumptions could be to assume that $\K$ is a topological ring with dense unit
group, and to admit only continuous polynomials.

\msk \nin  
{\bf PG.1.} {\sl Definition of polynomial groups.}
A {\it polynomial group (over an infinite field $\K$) of degree at most $k$}
(where $k\in \N$) is a $\K$-vector space $M$
together with a group structure $(M,m,i,e)$ such that $e=0$, the product
map $m:M \times M \to M$ and the inversion $j:M \to M$ are
 polynomial, and all iterated product maps, for $i \in \N$,
$$
m^{(i)}: M^i \to M, \quad (x_1,\ldots,x_i) \mapsto x_1 \cdots x_i
\eqno {\rm (PG.1)}
$$
(which are polynomial maps since so is $m=m^{(2)}$) are polynomial maps
of degree bounded by $k$.

\msk \nin
{\bf PG.2.}
Our three main examples of polynomial groups are:

\ssk
\item{(1)} The general multilinear group $M=\Gm^{1,k}(E)$, with its
natural chart, is a polynomial group of degree at most $k$:
with respect to a fixed multilinear
connection $L=L^0$ on $E$, every $f \in
\Gm^{1,k}$ is of the form $f=f^b = \id + X^b$ with a multilinear family
$b$, and in this way $\Gm^{1,k}(E)$ is identified with the $\K$-module
$M={\rm Mult}_{>1}$ of multilinear families that are singular of order $1$.
The multiplication map $m:M \times M \to M$ is described by Prop.\ MA.13:
it is clearly polynomial in $(g,f)$, with degree bounded by $k$.
Applying twice the matrix multiplication rule MA.13, we get the iterated
product map
$$
(m^{(3)}(h,g,f))^\Lambda = \sum_{\Xi \preceq \Omega \preceq \Lambda}
h^\Xi \circ g^{\Xi\vert \Omega} \circ f^{\Omega\vert \Lambda}
\eqno{\rm (PG.2)}
$$
which is again polynomial of degree bounded by $k$:
applied to an element $v \in E$, (PG.2) is a sum of multilinear terms;
every occurence of $f,g$ or $h$ corresponds to a non-trivial contraction;
but in a multilinear space of degree $k$, one cannot contract more than
$k$ times, so the degree is bounded by $k$. Similar arguments apply to
$m^{(j)}$. Finally, the inversion procedure used in the proof of Theorem
MA.6 is also polynomial (details are left to the reader).

\item{(2)} For any Lie group $G$, the group $(T^k G)_e$, equipped with
the $\K$-module structure coming from left or right trivialization,
 is a polynomial group of degree at most $k$ -- this is clear from
Theorem 24.7 since the Lie algebra $(T^k \g)_e$ is nilpotent of order $k$.
(Examples (1) and (2) make sense in arbitrary characteristic.)

\item{(3)}
If $\K$ is of characteristic zero, then any nilpotent Lie algebra over $\K$,
equipped with the Campbell-Hausdorff multiplication,
is a polynomial group of degree bounded by the index of nilpotence of the
Lie algebra.

\msk \nin
{\bf PG.3.} {\sl The Lie bracket.} We write
$$
m(x,y)=\sum_{i=0}^k m_i(x,y) = \sum_{p+q \leq k} m_{p,q}(x,y)
\eqno {\rm (PG.3)}
$$
for the decomposition of $m:M \times M \to M$ into homogeneous polynomial maps
$m_i:M \times M \to M$, resp.\ into parts $m_{p,q}(x,y)$, homogeneous
in $x$ of degree $p$ and in $y$ of degree $q$. The constant term is
$m_0(x,y)=m_0(0,0)=0$, and the linear term is
$m_1(x,y)=m_1((x,0)+(0,y))=m_1(x,0)+m_1(0,y)=x+y$,
whence $m_{1,0}(x,y)=x$, $m_{0,1}(x,y)=y$.
For $j>1$, we have $m_{j,0}(x,y)=0 = m_{0,j}(x,y)$ since
$m_{j,0}(x,y)=m_{j,0}(x,0)$ is the homogeneous component of degree $j$ of
$x \mapsto m(x,0)=x$. Thus
$$
m(x,y)=x+y+m_{1,1}(x,y) +\sum_{p+q= 3,\ldots, k \atop p,q\geq 1} m_{p,q}(x,y)
$$
with a bilinear map $m_{1,1}:M \times M \to M$. Similar arguments 
(cf. [Bou72, Ch.III, Par.\ 5, Prop. 1]) show that
$$
\eqalign{
x^{-1} = i(x)  = -x + m_{1,1}(x,x) & \mod \deg 3 \cr
xyx^{-1} = y + m_{1,1}(x,y)-m_{1,1}(y,x) & \mod \deg 3 \cr
c_1(x,y):=xyx^{-1}y^{-1} =m_{1,1}(x,y)-m_{1,1}(y,x) & \mod \deg 3 ,\cr}
\eqno{\rm (PG.4)}
$$
(where $\mod \deg 3$ has the same meaning as in [Bou72, loc.cit.]),
and one proves that
$$
[x,y]:=m_{1,1}(x,y)-m_{1,1}(y,x) 
\eqno{\rm (PG.5)}
$$
is a Lie bracket on $M$. If one defines the iterated group commutator by
induction,
$$
c_{j+1}(x_1,\ldots,x_{j+2}):= c_1(x_1,c_j(x_2,\ldots,x_{j+2})),
$$
it is easily seen that 
$$
c_j(x_1,\ldots,x_{j+1}) = [x_1,[x_2,\ldots,[x_j,x_{j+1}] \ldots ]] \,
\mod \deg (j+2).
$$
The maps $c_j$ are polynomials, and their degree is uniformely bounded
by some integer $N \in \N$ (write $c_j = m_r \circ \iota_j$,
where $\iota_j:M^{j+1} \to M^r$ is a polynomial map whose degree
is equal to the degree of inversion $j:M \to M$ and which is a certain
 composition
of diagonal imbeddings and copies of $j$). Since 
the iterated Lie bracket is of degree
$j+1$, it follows that it vanishes as soon as $j+1 > N$, and hence the
Lie bracket is nilpotent.
Let us now describe this nilpotent Lie algebra
for the examples mentioned above (PG.2, (1) -- (3)):

\ssk
\item{(1)} Recall Formula (MA.19) (or\ Prop.\ MA.13)
 for the components of the product
$m(g,f)= g \circ f$ in the group $\Gm^{1,k}(E)$. It is linear in
$g$, but non-linear in $f$. The part which is linear in $f$ is obtained
by summing over all partitions $\Omega$ such that $\Lambda$ is a refinement
of $\Omega$ {\it of degree one}, i.e.\ there is exactly one element of 
$\Omega$ which is a union of at least two elements of $\Lambda$, whereas
all other elements of $\Omega$ are also elements of $\Lambda$. 
Thus, if we define  the {\it degree} of the refinement $(\Lambda,\Omega)$ 
to be the number of non-trivial contractions of $\Lambda$ induced by
$\Omega$:
$$
\deg (\Omega \vert \Lambda) := {\rm card} \{ \omega \in \Omega \vert \,
\omega \notin \Lambda \},
$$
then
$$
m_{1,1}(g,f)^\Lambda =
 \sum_{\Omega \prec \Lambda \atop \deg(\Omega | \Lambda)=1}
g^\Omega \circ f^{\Omega | \Lambda}.
$$
Viewing $f$ and $g$ as polynomial mappings $E \to E$ and with the usual
definition of differentials of polynomial maps, this can be written
$m_{1,1}(g,f)(x)=dg(x) \cdot f(x)$, and thus
$$
[g,f](x) = dg(x) \cdot f(x) - df(x) \cdot g(x)
$$
is the usual Lie bracket of the algebra of polynomial vector fields on $E$.
We denote the Lie algebra thus obtained by
 $\gm^{1,k}(E)$. It is a {\it graded Lie algebra}; the grading
depends on the linear structure $L^0$, whereas the {\it associated
filtration} is an intrinsic feature, given by
$$
\gm^{1,k}(E) \supset \gm^{2,k}(E) \supset \ldots \supset \gm^{k-1,k}(E),
$$
where $\gm^{j,k}(E)={\rm Mult}_{>j}(E)$ is the $\K$-module of multilinear
families that are singular of order $j$.

\item{(2)}  For $M=(T^k G)_e$, we recover the Lie algebra $(T^k \g)_0$ of
$(T^k G)_e$, as follows easily from the fundamental commutation rule
(24.1).

\item{(3)} Since the first terms of the Campbell-Hausdorff formula  are
$X+Y+{1 \over 2}[X,Y] + \ldots$, we get the Lie bracket we started with:
 ${1 \over 2} [X,Y] - {1 \over 2} [Y,X] = [X,Y]$.

\msk \nin {\bf PG.4.} {\sl The power maps.} We decompose the iterated
product maps into multihomogeneous parts,
$$
m^{(j)} = \sum_{ p_1,\ldots,p_j \geq 0 \atop
p_1 + \ldots +p_j \leq k} m^{(j)}_{p_1,\ldots,p_j},
$$
(where the $m^{(j)}_{p_1,\ldots,p_j}$ can be calculated in terms of
$m_{p,q}$, but we will not need the explicit formula)
and (following  the notation from [Bou72] and from [Se65, LG 4.19]) we let
$$
\psi_j: M \to M, \quad x \mapsto \sum_{p_1,\ldots, p_j >0}
 m^{(j)}_{p_1,\ldots, p_j}(x,\ldots,x).
\eqno {\rm (PG.6)}
$$
Since the degree of $m^{(j)}$ is bounded by $k$ and $\sum_i p_i$ is
the degree of a homogeneous component of $m^{(j)}$, it follows that
$\psi_j = 0$ for $j >k$.
By convention, $\psi_0(x)=0$; then $\psi_1(x)=x$ and
$$
\eqalign{
x^2 & =m(x,x) =x+x+\sum_{i>1} m_i(x,x)=2 \psi_1(x) + \psi_2(x), \cr
x^3 & =m^{(3)}(x,x,x) = \sum_{p_1,p_2,p_3 \geq 0}
 m^{(3)}_{p_1,p_2,p_3} (x,x,x) \cr
& = (\sum_{p_1,p_2,p_3 \not=0} +
\sum_{p_1=0,p_2\not=0,p_3\not=0} +
\sum_{p_2=0,p_1\not=0,p_3\not=0} +
\sum_{p_3=0,p_2\not=0,p_1\not=0} +
\cr & 
\sum_{p_1=0,p_2=0,p_3\not=0} +
\sum_{p_1=0,p_3=0,p_2\not=0} +
\sum_{p_2=0,p_3=0,p_1\not=0} +
\sum_{p_1=p_2=p_3=0})(
 m^{(3)}_{p_1,p_2,p_3} (x,x,x)) \cr
&=\psi_3(x) + 3 \psi_2(x) + 3 \psi_1(x) + \psi_0(x) \cr
x^n & =  \sum_{j=0}^k \pmatrix{n \cr j\cr} \psi_j(x). \cr}
$$
See [Bou72, III, Par.\ 5, Prop.\ 2] for the details of the proof
(mind that $\psi_j = 0$ for $j>k$).
The last formula holds for all $n\in \N$ and in arbitrary characteristic.

\msk \nin {\bf PG.5.} {\sl One-parameter subgroups.}
From now on we assume that $\K$ is a field of characteristic zero.  Then,
for fixed $x \in M$, we let
$$
f:=f_x: \K \to M, \quad t \mapsto x^t:=
 \sum_{j=0}^k \pmatrix{t \cr j\cr} \psi_j(x).
\eqno {\rm (PG.7)}
$$
Clearly, $f$ is a polynomial map (of degree bounded by $k$), and we have
just seen that $f(t+s)=f(t) f(s)$ for all $t,s \in \N$.
For fixed  $s \in \N$,
both sides of this equality depend polynomially on $t \in \K$;
they agree for $t \in \N$ and hence for all $t \in \K$
 (since $\K$ is a field of characteristic zero).
Since $f(n)=x^n$ for $n \in \N$, we may say that $f_x$ is a
{\it one-parameter subgroup through $x$}.
In particular, it follows that
$$
x^{-1} =   \sum_{j=0}^k \pmatrix{-1 \cr j\cr} \psi_j(x) =
 \sum_{j=0}^k (-1)^j  \psi_j(x),
\eqno {\rm (PG.8)}
$$
and hence the inversion map is a polynomial map of degree bounded by $k$.
Note that $f_x$ is not a one-parameter subgroup ``in direction of $x$'' since
the derivative $f_x'(0)$ is not equal to $x$ but is given by the linear term in
$$
\eqalign{
f_x(t) & = tx + {t (t-1) \over 2} \psi_2(x) + 
{t (t-1) (t-2) \over 3!} \psi_3(x) + \ldots \cr
&  =
t(x-{1 \over 2} \psi_2(x) + {1 \over 3} \psi_3(x) \pm \ldots ) 
\, \mod \,(t^2), \cr}
$$
and thus $f_x$ is a one-parameter subgroup in direction
$f_x'(0) = \sum_{p=1}^k {(-1)^{p-1} \over p} \psi_{p}(x).$
Thus, in order to define ``the one-parameter subgroup in direction $y$'',
all we need is to invert the polynomial 
$\sum_{p=1}^k {(-1)^{p-1} \over p} \psi_{p}(x)$, if possible.
 This is essentially the
content of the following theorem:

\Theorem PG.6. Assume $M$ is a polynomial group over
a field $\K$ of characteristic zero.
We denote by $\m$ a formal copy of the $\K$-vector space $M$ (where we may
forget about the group structure).
Then there exists a unique polynomial
$\exp: \m \to M$ such that
\item{(1)}
for all $X \in \m$ and $n \in \Z$, $\exp(nX)=(\exp X)^n$,
\item{(2)}
the linear term of $\exp$ is the identity map  $\m = M$.
\ssk \nin
The polynomial $\exp$ is bijective, and its inverse map $\log:M \to \m$
is again polynomial. We have the explicit formulae
$$
\eqalign{
\exp(x) & 
= \sum_{p=1}^k {1 \over p!} \psi_{p,p}(x) \cr
\log(x) & = \sum_{p=1}^k {(-1)^{p-1} \over p} \psi_{p}(x) \cr}
$$
where $\psi_{p,m}$ is the homogeneous component of degree $m$ of $\psi_p$.
                                              
\Proof. We proceed in several steps.

\ssk
Step 1. Uniqueness of the exponential map.
We claim that two polynomial homomorphisms $f,g:\K \to M$
with $f'(0)=g'(0)$ are equal
 (here we use the standard definition of the derivative of a polynomial). 
Indeed, this follows from general arguments
given in [Se65, p.\ LG 5.33]: $f,g$ being homomorphisms, they satisfy
the same formal differential equation with same initial condition and
hence must agree. 
It follows that, for every $v \in M$, there exists at most one polynomial
 homomorphism  $\gamma_v:\K \to M$ such that $\gamma_v'(0)=v$.
On the other hand,
Property (1) of the exponential map implies (by the same ``integer density'' 
argument as above) that
$\exp(t X)=(\exp X)^t$ for all $t \in \K$ and hence
$\exp((t+s)X)=\exp(tX) \cdot \exp(sX)$ for all $t,s \in \K$.
Invoking  Property (2) of the exponential map, it follows that
 the exponential map, if it exists, has to be given by
$\exp(v)=\gamma_v(1)$. 

\ssk
Step 2. We define the logarithm by the polynomial expression given in the
claim and 
establish the functional equation of the logarithm.
As remarked before stating the theorem, $\log(x) = f_x'(0)$.
This implies the functional equation
$$
\forall n \in \Z: \quad \log(x^n)=n \log(x).
$$
Indeed, from $f_x(n)=x^n=f_{x^n}(1)$ we get $f_x(nt)=f_{x^n}(t)$ for
all $t \in \K$ and hence
$$
\log(x^n)={d \over dt}|_{t=0} f_{x^n}(t)=
{d \over dt}|_{t=0} f_{x}(nt)= n \log(x).
$$
This implies $\log(x^t)=t \log(x)$ for all $t \in \K$. 

\ssk
Step 3.  We define $\exp$ by the polynomial expression given in the claim and
 show that $\log \circ \exp = \id_\m$.
In fact, here the proof from [Bou72, III. Par.\ 5 no.\ 4] carries over word
by word.
We briefly repeat the main arguments : for $t \in \K^\times$,
$$
\log(t x)=t \log((tx)^{t^{-1}}) = t \log (\sum_{r \leq m} t^{m-r} \phi_{r,m}
(x)),
\eqno{\rm (PG.9)}
$$
where the polynomials $\phi_{r,m}$ are defined by the expansion
$$
x^t = \sum_{r,m} t^r \phi_{r,m}(x)
$$
into homogeneous terms of degree $r$ in $t$ and degree $m$ in $x$.
From the definition of $x^t$ and of
 the components $\psi_{p,m}$ one gets in particular
$$
\phi_{1,m}(x)=\sum_{p \leq m} {(-1)^{p-1} \over p} \psi_{p,m}(x), \quad \quad
\phi_{m,m}(x)={1 \over m!} \psi_{m,m}(x).
$$
One notes also that $\psi_{p,m} =0$ if $m<p$. 
Taking ${d \over dt}|_{t=0}$ in (PG.9), we get
$$
x = \log (\sum_m \phi_{m,m}(x)) = 
\log (\sum_m {1 \over m!} \psi_{m,m}(x)) = \log(\exp (x)).
$$

\ssk
Step 4. We show that $\exp \circ \log = \id_M$. We  let
$$
\gamma_v(t) := f_{\exp(v)}(t).
$$
Then $\gamma_v'(0)=
f_{\exp(v)}'(0)=\log (\exp(v)) = v$ by Step 3, and hence the preceding 
definition is consistent with the definition of $\gamma_v$ in Step 1.
We have
$$
\gamma_v(1) = f_{\exp(v)}(1)= \exp(v).
$$
The uniqueness statement on one-parameter subgroups from Step 1 shows
that $f_x = \gamma_{\log(x)}$ since both have derivative $\log(x)$ at $t=0$.
Thus $x=f_x(1)=\gamma_{\log(x)}(1) = \exp(\log(x))$.

\ssk
Step 5. We prove that $\exp$ satisfies (1) and (2). In fact,  Property
(2) follows from the fact that $\Psi_1(x)=\Psi_{1,1}(x)=x$. Since
$\exp=(\log)^{-1}$,
the functional equation (1) follows from the functional equation of the
logarithm established in Step 2. (One may also adapt the usual ``analytic''
argument:
from the uniqueness statement on one-parameter subgroups from Step 1
we get $\gamma_{tv}(1)=\gamma_v(t)$, which yields $\exp(t v)=(\exp v)^t$
and hence the functional equation (1).)
\qed

\nin
Let us add some comments on the theorem and on its proof.
We have chosen notation $\m$ and $M$ such that on the side of $\m$ only the
$\K$-vector space structure is involved  (see Theorem PG.8 below
for relation with the Lie algebra structure) and on $M$ only the group
structure and its ``regularity''. 
In the theory of general formal
 groups, the existence and uniqueness
of formal power series $\exp$
and $\log$ with the desired properties is proved by quite different
arguments: in  [Bou72, Ch.III,
Par.\ 4, Thm.\ 4 and Def.\ 1] the existence of an exponential map
is established by using the theory of  universal envelopping algebras
of Lie algebras, and the argument from
 [Se65, LG 5.35. -- Cor.\ 2] relies heaviliy on the existence of
a Campbell-Hausdorff multiplication.

\msk \nin
{\bf PG.7.} {\sl Examples.}
Let us give a more explicit description of the exponential maps
for our examples PG.2 (1) -- (3).

\msk
\item{(1)} Let $X=X^b \in M = \gm^{1,k}(E)$, where $b$ is
 a multilinear family $b$, singular of order $1$. Then, for $\ell(\Lambda)>1$,
$$
\exp(X^b)^\Lambda =
 b^\Lambda + \sum_{j=2}^k {1\over j!}
  \sum_{\Omega^{(1)} \prec \ldots \prec \Omega^{(j)} =\Lambda  \atop
\deg(\Omega^{(i)} \vert \Omega^{(i+1)} )=1 }
b^{\Omega^{(1)}} \circ b^{\Omega^{(1)} \vert \Omega^{(2)} } \circ \ldots
\circ b^{\Omega^{(j-1)} \vert \Lambda}.
\eqno {\rm (PG.10)}
$$
Moreover, if $E$ is of tangent or bundle type
(Chapter SA), then the very symmetric and the (weakly) symmetric
subgroups of $\Gm^{1,k}(E)$ are again polynomial groups, and the
polynomial $\exp$ commutes with the action of the
symmetric group (this follows from functoriality of $\exp$).
Thus the exponential map of $\Sm^{1,k}(E)$ and of $\Vsm^{1,k}(E)$ is
given by restriction of the one of $\Gm^{1,k}(E)$.
In particular, we get a combinatorial formula for the exponential
map of the group 
$\Vsm^{1,k}(E) \cap \Gm^{1,k}(E)$ which, as seen at the end of Chapter SA,
is a version of ``jet group'' from
 [KMS93, Section 13.1]. 

\item{(2)} See Theorem 25.2. 

\item{(3)} If $M$ is a Lie algebra with Campbell-Hausdorff multiplication,
then $\exp$ is the identity map
(here $\psi_j = 0$ for $j \geq 2$, and this property characterizes the
``canonical chart'', cf.\ [Bou72, III. Par.\ 5, Prop.\ 4]).

\msk \nin
Concerning Example (1), we add the remark that,
even if the integers are not invertible in $\K$, it is possible to express
the group law of the general multilinear group
$\Gm^{1,k}(E)$ in terms of the Lie bracket, by a formula
similar to the one given for higher order tangent groups in Theorem 24.7.
In any case, it can be shown that the filtration  of the Lie algebra
$\gm^{1,k}(E)$ has an analog on the group level:
 the chain of subgroups from part (3) of
Theorem MA.6 is a {\it filtration} of $\Gm^{1,k}(E)$
by normal subgroups 
in the sense that the group commutator satisfies
$[\Gm^{i,k}(E),\Gm^{j,k}(E)] \subset \Gm^{i+j,k}(E)$.

\Theorem PG.8. Under the assumptions of the preceding theorem,
the multiplication on a polynomial group is the Campbell-Hausdorff
multiplication with respect to the polynomial $\exp$ from the preceding
theorem. In other words,
 $x * y :=\log(\exp (x) \exp(y))$ is given by the Campbell-Hausdorff
formula, with the Lie bracket defined by {\rm (PG.4)}.
The Campbell-Hausdorff multiplication is again polynomial. 

\Proof. The first claim follows from  a general result on formal groups
characterizing the exponential map as the unique formal homomorphism from
the Campbell-Hausdorff group chunk to the group having the identity as
linear term
(cf.\ [Bou72, Ch.III, Par. 4, Th.\ 4 (v)], see also [Se65, LG 5.35]).
In our case, the Campbell-Hausdorff multiplication is polynomial since
since so are $\exp$, $\log$ and the ``original'' multiplication.
\qed

\msk \nin {\bf PG.9} {\sl The projective limit and formal groups.}
We could define now a {\it formal (power series) group} as the
projective limit $M_\infty$ of a projective system
$\pi_{\ell,k}:M_k \to M_\ell$, $k \geq \ell$, of polynomial groups 
$M_k$, $k \in \N$,
where $M_k$ is of degree (at most) $k$. 
A formal power series group in the usual sense (see, e.g., [Haz78], [Se65])
gives rise to a sequence of truncated polynomial groups of which it is
the projective limit.
The exponential maps and logarithms $\exp_k$, resp.\ $\log_k$,
then define, via their projektive limits,
 formal power series $\exp_\infty$ and $\log_\infty$
for $M_\infty$. 
This way of proving the existence of the exponential and logarithm series
has the advantage of using only ``elementary'' methods which are
close to the usual notions of analysis on Lie groups.

\vfill \eject

\def\title{ } 
\def\rightheadline{\hfil{\eightrm \title}\hfil\tenrm\folio}

\def\entries{

\[ALV91 Alekseevski, D.V., Lychagin, V.V and A.M. Vinagradov,
``Geometry I: Basic Ideas and Concepts of Differential Geometry",
vol.\ {\bf 28} of the Encyclopedia of Mathematical Sciences
(Ed.\ R.V. Gamkrelidze), Springer-Verlag, Berlin 1991

\[B57 Berger, M., {\it Les espaces sym\'etriques non-compacts},
Ann. Sci. Ec. Norm. Sup. (3) {\bf 74} (1957), 85 -- 177

\[Bes78 Besse, A.L. ``Manifolds all of whose Geodesics are Closed'',
Springer EMG {\bf 93}, New York 1978

\[BBG98 Bertelson, M., Bielavsky, P. et S. Gutt, {\it Parametrizing
equivalence classes of invariant star products}, Lett. Math. Phys.
{\bf 46} (4) (1998), 339 -- 345

\[Be96 Bertram, W., {\it Un th\'eor\`eme de Liouville pour les alg\`ebres
de Jordan}, Bull.\ Soc.\ Math.\ Fran\c caise {\bf 124} (1996), 299 --
327

\[Be00 Bertram, W., ``The Geometry of Jordan and Lie Structures'',
Sprin\-ger LNM {\bf 1754}, Berlin 2000

\[Be02 Bertram, W. {\it Generalized projective geometries: general
theory and equivalence with Jordan structures}, 
Advances in Geometry {\bf 2} (2002), 329--369

\[Be04 Bertram, W. {\it From linear algebra via affine algebra to
projective algebra},
Linear Algebra and its Applications {\bf 378} (2004), 109 -- 134

\[Be07 Bertram, W., ``Calcul diff\'erentiel topologique \'el\'ementaire'',
book in preparation

\[BGN04 Bertram, W., H. Gl\"ockner and K.-H. Neeb, {\it Differential
Calculus over General Base Fields and Rings},
Expo. Math. {\bf 22} (2004), 213--282. 
arXiv: math.GM/0303300

\[BeNe04 Bertram, W. and K.-H. Neeb, {\it Projective completions of
Jordan pairs. Part I: The generalized projective geometry of a
Lie algebra}. J. of Algebra {\bf 227} (2004), 474 -- 519.
 arXiv: math.RA/0306272

\[BeNe05 Bertram, W. and K.-H. Neeb, {\it Projective completions of
Jordan pairs. Part II: Manifold structures and symmetric spaces}.
Geom.\ Ded.\ {\bf 112} (2) (2005), 73 -- 113.
   arXiv: math.GR/0401236

\[BGV92 Berline, N., Getzler E. and M. Vergne,
``Heat Kernels and Dirac Operators'', Springer-Verlag, 
New York 1992

\[Bor02 Bordemann, M., {\it Sur l'existence d'une prescription d'ordre
naturelle projectivement invariante}, preprint 2002, 
arXiv: math.DG/0208171

\[Bou67  Bourbaki, N. ``Vari\'et\'es diff\'erentielles et analytiques --
Fascicule de r\'e\-sul\-tats'', Hermann, Paris 1967 -- 1971

\[Bou72 Bourbaki, N. ``Groupes et alg\`ebres de Lie. Chapitres 1 -- 3.'',
 Hermann, Paris 1972

\[Bue95 Buekenhout, F. (editor), `` Handbook of Incidence Geometry'',
Elsevier, Amsterdam 1995

\[Ca71 Cartan, H., ``Calcul Diff\'erentiel'', Hermann, Paris 1971

\[CFT01 Cattaneo, A., Felder, G. and L. Tomassini, {\it From local to global
deformation quantization of Poisson manifolds}, preprint, arXiv
math. QA/ 0012228v3

\[Chap03 Chaperon, M., ``Calcul diff\'erentiel et calcul int\'egral'',
Dunod, Paris 2003

\[Ch01 Chapoton, F., {\it Alg\`ebres de pr\'e-Lie et alg\`ebres de Hopf
li\'ees \`a la renormalisation}, C.R.\ Acad.\ Sci.\ Paris, t. {\bf 332} 
(2001), S\'erie I, p.\ 1 -- 4

\[Did06 Didry, M., {\it
Construction of groups associated to Lie- and Leibniz-algebras},
preprint, Nancy 2006 (Journal of Lie Theory, to appear)

\[D73 Dieudonn\'e, J., ``Introduction to the Theory of Formal Groups",
Marcel Dekker, New York 1973


\[GD71 Grothendieck, A. and J. Dieudonn\'e, ``El\'ements de G\'eom\'etrie
Alg\'e\-bri\-que. I.'', Springer, Berlin 1971


\[Gi03 Giordano, P., {\it Differential Geometry in spaces of mappings
using nilpotent infinitesimals}, preprint, arXiv:
math.DG/0308119v1

\[Gl03a Gl\"ockner, H., {\it Implicit Functions from Topological Vector
Spaces to Banach Spaces}, Israel J.\ Math.\ 
{\bf 155} (2006), 205 -- 252; see also
 Preprint Nr. 2271 TU Darmstadt, 2003.
 arXiv: math. GM/0303320

\[Gl03b  Gl\"ockner, H., {\it Bundles of locally convex spaces, group
actions, and hypo\-continuous bilinear mappings}, submitted

\[Gl03c Gl\"ockner, H., {\it Every smooth $p$-adic Lie group admits a
compatible  analytic structure}, to appear in Forum Math.\
{\bf 18} (2006), 45 -- 84.
 (also TU Darmstadt Preprint 2307, December 2003; arXiv:math. GR/0312113)

\[Gl04a  Gl\"ockner, H., {\it Examples of differentiable mappings into
non-locally convex spa\-ces},
TU Darmstadt Preprint 2272, March 2003, 6 pp, also arXiv: math. FA/0307040; 
 Topol.\ Proceedings.\ {\bf 28} (2004), No. 2, 47 -- 486

\[Gl04b  Gl\"ockner, H., {\it Tensor products in the category of topological
vector spaces are not associative}, 
Comment.\ Math.\ Univ.\ Carol.\ {\bf 45} (2004), no. 4, 607 -- 614

\[Gl04c  Gl\"ockner, H., {\it Fundamentals of direct limit Lie theory},
 Compositio Math.\ {\bf 141} (2005), 1551-1577.
 (cf. TU Darmstadt Preprint 2324, March 2004; arXiv: math. GR/0403093)

\[Gl04d  Gl\"ockner, H., {\it Smooth Lie groups over local fields of positive
 characteristic need not be analytic},
 J.\  Algebra {\bf 285} (2005), 356 -- 371 (cf. arXiv: math. GR/0408239)

\[Gl04e  Gl\"ockner, H., {\it Lie groups over non-discrete topological
fields}, TU Darmstadt Preprint 2356, July 2004; also arXiv: math. GR/0408008 

\[Haz78 Hazewinkel, M., ``Formal Groups and Applications", Academic Press,
New York 1978

\[Hel78 Helgason, S., ``Differential Geometry and Symmetric Spaces'',
Academic Press 1978

\[Hel84 Helgason, S., ``Groups and Geometric Analysis'',
Academic Press 1984

\[Ho01 Holdgr\"un, H.S., ``Analysis. Band 2: Differential- und Integralrechnung
mehrerer Variablen'', Leins-Verlag, G\"ottingen 2001

\[KM87 Kainz, G. and P.\ Michor, {\it Natural transformations in
differential geometry}, Czech.\ Math.\ Journal {\bf 37} (112) (1987),
584 -- 607


\[Ko81 Kock, A., ``Synthetic Differential Geometry'', London Math. Soc.
Lect. Note Series {\bf 51}, Cambridge Univ. Press, 1981

\[Koe69 Koecher, M., {\it Gruppen und Lie-Algebren von rationalen
Funktionen}, Math. Z. {\bf 109}, 349--392 (1969)

\[KMS93  Kol\'a\v r, I., Michor P. W. and J. Slov\'ak, 
``Natural Operations in Differential Geometry", Springer-Verlag,
Berlin 1993

\[KoNo63 Kobayashi, S. and S.  Nomizu, ``Foundations of Differential
Geometry, vol. I'', Wiley, New York 1963

\[La99 Lang, S., ``Fundamentals of Differential Geometry'', Graduate
Texts in Math. {\bf 191}, Springer-Verlag, Berlin 1999

\[Lau77 Laugwitz, D., ``Differentialgeometrie'', Teubner,
Stuttgart 1977

\[Lau86  Laugwitz, D., ``Zahlen und Kontinuum'', BI Wissenschafts\-verlag,
Mann\-heim 1986

\[Lav87 Lavendhomme, R., ``Le\c cons de G\'eom\'etrie Diff\'erentielle
Synth\'etique Na\-\"\i \-ve'', Mo\-nographies de math\'ematique {\bf 3},
Louvain-la-Neuve 1987 

\[Lo67 Loos, O., {\it Spiegelungsr\"aume und homogene symmetrische
R\"aume}, Math. Z. {\bf 99} (1967), 141 -- 170 

\[Lo69 Loos, O., ``Symmetric Spaces I'', Benjamin, New York 1969



\[Mac79 Macdonald, I.\ G., ``Symmetric Functions and Hall Polynomials'',
Clarendon Press, Oxford 1979

\[MR91 Moerdijk, I. and G.E. Reyes, ``Models for Smooth Infinitesimal
Analysis'', Springer-Verlag, New York 1991

\[Ne02a
 Neeb, K.-H., {\it A Cartan-Hadamard Theorem for Banach-Finsler Manifolds},
Ge\-om.  Dedicata {\bf 95} (2002), 115 -- 156

\[Ne02b Neeb, K.-H., ``Nancy lectures on infinite dimensional Lie groups'',
Lecture Notes, Nancy 2002 

\[Ne06 Neeb, K.-H., {\it Towards a Lie theory of locally convex groups},
Jap.\ J.\ of Math.\ {\bf 1} (2006), 291--468

\[N88 Nel, L.D., {\it Categorial Differential Calculus for infinite
dimensional Spa\-ces}, Cahi\-ers de Top. et G\'eom. Diff. Cat\'egoriques,
{\bf 29}, t. 4 (1988), 257--286

\[Nes03 Nestruev, J., ``Smooth Manifolds and Observables'', Springer-Verlag,
New York 2003

\[P62 Pohl, W.F., {\it Differential Geometry of Higher Order},
Topology {\bf 1} (1962), 169--211

\[Pos01 Postnikov, M.M., ``Riemannian Geometry'', vol. {\bf 91} of
the Encyclopedia of Mathematical Sciences, Springer-Verlag,
New York 2001

\[Ro66 Robinson, A, ``Non-Standard Analysis'', North Holland,
Amsterdam 1966

\[Rob63 Roby, N., {\it Lois polyn\^omes et lois formelles en th\'eorie de
modules}, Ann. Sci. Ec. Norm. Sup. 3e s\'erie t. {\bf 80} (1963),
213 --  348

\[Sch84 Schikhof, W.H., ``Ultrametric Calculus -- An Introduction to
$p$-adic Analysis'', Cambridge University Press, Cambridge 1984

\[Se65 Serre, J.-P., ``Lie Algebras and Lie Groups'', Benjamin,
New York 1965

\[Sh97  Sharpe, R.W.,  ``Differential Geometry  Cartan's Generalization
of Klein's Erlangen Program", Springer, New York 1997

\[Sp79 Spivak, M., ``A comprehensive Introduction to Differential
Geometry'', vol. 5 (Scd. Ed.), Publish and Perish, Houston 1979

\[W53 Weil, A., {\it Th\'eorie des points proches sur les vari\'et\'es
diff\'erentiables}, Colloque de G\'eom. Diff., Strasbourg 1953, pp. 111--117
(in: A. Weil, \OE uvres  scientifiques, Vol. 2, Springer,
   New York  1980, pp. 103--109)

\[T04 Tu, L.W., {\it Une courte d\'emonstration de la formule de
Campbell-Haus\-dorff}, J.\ of Lie Theory {\bf 14}, no.\ 2 (2004), 501--508

\[Tr67 Tr\`eves, F., ``Topological Vector Spaces, Distributions and
Kernels'', Academic Press, New York 1967

}

\references

\vfill\eject

\end